\documentclass[a4paper,10pt]{article}
\usepackage{authblk}

\usepackage[ngerman,english]{babel}
\usepackage[latin1]{inputenc}
\usepackage{csquotes}
\usepackage{amsfonts,amsmath,amsthm}
\usepackage{empheq}
\usepackage{cancel}
\usepackage[titletoc,title]{appendix}
\usepackage[backend=bibtex8,doi=false,eprint=false,firstinits=true,isbn=false,style=numeric-comp,maxnames=99]{biblatex}
\makeatletter
\def\blx@maxline{77}
\makeatother

\bibliography{Bibliography.bib}

\usepackage{cases}
\usepackage{mathabx}
\usepackage{bbm}
\usepackage{xfrac}
\usepackage{fancyhdr}
\usepackage{color}
\usepackage[colorinlistoftodos,textsize=small]{todonotes}
\usepackage[colorlinks]{hyperref}
\definecolor{blue75}{rgb}{0,0,.75}
\definecolor{green75}{rgb}{0,.75,0}
\hypersetup{colorlinks=true, urlcolor=blue75,linkcolor=blue75,citecolor=green75,pdfstartview=FitB,bookmarksopen=true,bookmarksopenlevel=1}
\usepackage{graphicx,subcaption}
\usepackage[a4paper, left=2.5cm, right=2.0cm, top=2.5cm,bottom=2cm]{geometry}
\usepackage{constants}
\newcommand{\parenthezises}[1]{\arabic{#1}}
\newconstantfamily{C}{
	symbol=C,
	format=\parenthezises,
	reset={section}
}
\newconstantfamily{M}{
	symbol=M,
	format=\parenthezises,
	reset={section}
}
\newconstantfamily{B}{
	symbol=B,
	format=\parenthezises,
	reset={section}
}
\newconstantfamily{xi}{
	symbol=\xi,
	format=\parenthezises,
	reset={section}
}
\usepackage{enumerate}
\usepackage{graphicx}
\graphicspath{ {images/} }
\usepackage{wrapfig}
\usepackage{figbib}
\allowdisplaybreaks
\usepackage[capitalise]{cleveref}

\crefdefaultlabelformat{{\it #2#1#3}}

\crefname{equation}{}{}

\crefname{enumi}{}{}
\creflabelformat{enumi}{{(#2#1#3)}}

\crefname{section}{{\it Section}}{{\it Sections}}
\crefname{subsection}{{\it Subsection}}{{\it Subsections}}
\crefname{subsubsection}{{\it Paragraph}}{{\it Paragraphs}}

\newtheorem{Theorem}{Theorem}[section]
\crefname{Theorem}{{\it Theorem}}{{\it Theorems}}

\crefname{Definition}{{\it Definition}}{{\it Definitions}}
\newtheorem{Lemma}[Theorem]{Lemma}
\crefname{Lemma}{{\it Lemma}}{{\it Lemmas}}

\crefname{Proposition}{{\it Proposition}}{{\it Propositions}}

\crefname{Assumption}{{\it Assumption}}{{\it Assumptions}}

\crefname{Corollary}{{\it Corollary}}{{\it Corollaries}}

\theoremstyle{definition}
\newtheorem{Remark}[Theorem]{Remark}
\crefname{Remark}{{\it Remark}}{{\it Remarks}}

\crefname{Notation}{{\it Notation}}{{\it Notations}}

\crefname{Example}{{\it Example}}{{\it Examples}}

\newcommand{\be}{\begin{equation} \label}
	\newcommand{\ee}{\end{equation}}
\newcommand{\bea}{\begin{eqnarray}\label}
	\newcommand{\eea}{\end{eqnarray}}
\newcommand{\bas}{\begin{eqnarray*}}
	\newcommand{\eas}{\end{eqnarray*}}
\newcommand{\bit}{\begin{itemize}}
	\newcommand{\eit}{\end{itemize}}
\newcommand{\R}{\mathbb{R}}

\begin{document}
	\enlargethispage{10mm}

	\title{Mathematical modeling of Buruli ulcer spread}
	\author[1]{Shimi Chettiparambil Mohanan} 
		\author[1]{Christina Surulescu}
		\author[2]{Neslihan Nesliye Pelen}
		\affil[1]{{\small RPTU Kaiserslautern-Landau,  Postfach 3049, 67653 Kaiserslautern, Germany}} 
	\affil[2]{{\small Missouri University of Science and Technolgy, Rolla, MO65401, USA}}
	\date{\today}
	\maketitle
	
	\begin{abstract} We propose two approaches to a model deduction for Buruli ulcer spread in a tissue, prove global existence of solutions to the obtained macroscopic cross-diffusion PDE-ODE systems, and perform numerical simulations to illustrate the behavior of solutions under various scenarios and compare the outcome of the two approaches.
	\end{abstract}
	
	\section{Introduction}\label{sec1}
	
	\noindent
	Buruli ulcer (BU) is caused by Mycobacterium ulcerans, a member of the family of bacteria causing tuberculosis and leprosy. It is a neglected tropical or subtropical disease \cite{neglectedburuli} characterised by destructive skin lesions and tissue damage. Delay of treatment results in prolonged periods in the hospital \cite{Asiedu} and long-term disability \cite{ab}. It often affects geographical regions with limited healthcare resources, posing significant challenges for diagnosis, treatment, and control, while demonstrating higher incidence rates in children compared to adults \cite{bar}. \\[-2ex]
	
	\noindent
	Advanced stages of BU require surgery, however, the true extent of the lesions cannot be assessed by commonly available medical techniques and this often leads to post-surgery recurrences \cite{Muelder}. Therefore, mathematical models are called upon to help make predictions about the spread of bacteria and necrotic matter and hence guide the resection. To get a mathematical model, it is important to understand the structure and the pathogenesis of the illness in the organism. Observations have been made by investigating the pathogenesis of BU and models should be set up by  making use of these observations.\\[-2ex]
	
	\noindent 
	Mathematical modeling of the spread of M. ulcerans bacteria within tissue is less developed; to our knowledge there hasn't been proposed any reasonable space-time mathematical model; the setting in \cite{nyarko} does not make much sense. Ordinary differential equation (ODE) disease spread models within a population have been studied in \cite{Nyabadza2015, zhao2021mathematical}. However, the spread of bacteria within the wound is crucial for the treatment and has to be studied. We present two novel RDTE models specifically designed to address this issue. In particular, we employ two approaches to characterize the dynamics of bacteria spread: multiscale modeling using kinetic transport equations within the kinetic theory of active particles framework \cite{bellomo2017quest} and a random jump description (as e.g. in \cite{doi:10.1098/rsif.2008.0014, randomjumpKPH, randomjumpSAO}). \\[-2ex]
	
	\noindent
	The pathogenesis of Buruli ulcer involves the infiltration of Mycobacterium ulcerans through insect bites into subcutaneous tissue.  Although the exact mode of transmission is unclear \cite{JSSPMPHHA}, there is strong evidence that BU does not transfer from person to person. Instead, it might be due to bites by some water-living insects \cite{ABT02,Walsh}. The bacteria then proliferate, producing a toxin called mycolactone \cite{Portaels}. The toxin spreads from the initial site of infection, penetrating the surrounding areas and diffusing further, resulting in the necrosis of affected tissue \cite{Wansbrough}. This creates an environment conducive to the continued proliferation of the pathogen \cite{Walsh}. The bacteria invade regions with high concentrations of normal \cite{Oliviera} and necrotic tissue, where they further proliferate and produce more toxin, facilitating the sustained spread of lesions. \\[-2ex]
	
	\noindent
	The rest of this note is structured as follows: Section \ref{sec:models} contains the two approaches to deduction of macroscale dynamics of bacteria, normal and necrotic tissue, and mycolactone. Section \ref{sec:analysis} is concerned with the analysis of the obtained macroscopic formulations, in terms of existence of solutions. In Section \ref{sec:simulations} we perform numerical simulations to illustrate the solution behavior for the macroscale systems and compare the two approaches. 
		
	\section{Modeling}\label{sec:models}
	
	\subsection{Approach via kinetic transport equations (velocity jumps)}
	
	\noindent
We start on the mesoscopic scale on which groups of cells which share the same velocity regime are considered. We are interested in the dynamics of cell densities $ u(t,x,\vartheta) $, where $t>0$ and $x\in \mathbb{R}^d$ represent time and space variables, while $\vartheta\in \Theta=s \mathbb{S}^{d-1}$ denotes velocity. We assume, e.g., as in \cite{conte2021mathematical,conte2023mathematical,corbin2021modeling,engwer2015glioma,engwer2016effective,kumar2021multiscale} for different, but mathematically related problems, that the (average) speed $s>0$ is constant and only the direction of movement matters; $\mathbb{S}^{d-1}$ represents the unit sphere in $\mathbb{R}^d$.\\[-2ex]

\noindent
In our modeling attempt we intend to pay attention to the interactions between bacterial cells and their surrounding environment, including necrotic matter, normal tissue, and toxin. These interactions influence various aspects of bacterial behaviour, such as movement and growth rates, as well as potential consequences on necrotisation of normal tissue. In the sequel $U,V, N$ denote the macroscopic densities of bacteria, normal, and necrotic matter, respectively, while $M$ is the concentration of mycolactone. All these quantities are space-time dependent. \\[-2ex]

\noindent
 We consider the following  kinetic transport equation (KTE) describing velocity jumps of bacteria: 
\begin{equation}\label{-_-}
	\partial_t u+\nabla_x \cdot (\vartheta u)= \mathcal{L}_u[\eta] u+\mu_u(U,V,N)u,
\end{equation}
where  the turning operator, $\mathcal{L}_u[\eta] u$ characterizes bacteria reorientation. Since the cells are attracted by normal tissue and necrotic matter, we propose the following turning operator: 
$$\mathcal{L}_u[\eta]u=-\eta(x,\vartheta,V,N)u(t,x,\vartheta)+\int_{\Theta} \frac{1}{|\Theta|}\eta(x,\vartheta',V,N)u(\vartheta')d\vartheta'.$$ 
Here we considered for simplicity a uniform density function over the unit sphere in order to describe the turning kernel. In \cite{conte2021mathematical,conte2023mathematical,corbin2021modeling,engwer2015glioma,engwer2016effective,kumar2021multiscale} it took into account the orientation distribution of the underlying fibrous structure of the tissue. Our simplification is twofold motivated: on the one hand, bacteria are so small, that they can spread largely independently on how the tissue is oriented; on the other hand, we also aim at describing tissue evolution, but cannot do that on the mesocopic level, since we would eventually obtain a two-scale cross-diffusion system coupling PDEs, ODEs, and an integro-differential ODE featuring tissue orientation as a further variable. A similar system arising from such modeling process was obtained in \cite{corbin2021modeling} and raised several  challenges as far as numerics was concerned, while its analysis is still not addressed.\\[-2ex]

\noindent
The turning rate $\eta $ is is chosen in the form (see \cite{conte2021mathematical},\cite{sym12111870} for a similar approach in the context of glioma spread): 
$$\eta(x,\vartheta,V,N)=\eta_0(x) e^{-a(V,N)(\frac{\gamma_1}{K_N} D_t N+\frac{\gamma_2}{K_V} D_t V)},$$ 
where $\eta_0(x)$ is the turning rate in the absence of external cues, while $K_V, K_N >0$ are the carrying capacities of $V$ and $N$, respectively. The coefficient functions $\frac{\gamma_1}{K_N}a(V, N)$ and $\frac{\gamma_2}{K_V} a(V, N)$ are related to interactions of  bacteria with normal tissue and with necrotic matter, respectively. The constants $\gamma_1,\gamma_2 > 0$ account for the change in the turning rate per unit of change in $dy^*/dN$ and $dy^*/dV$, respectively. Thereby, $y^*$ denotes the steady-state of bacteria receptor binding dynamics to $N$ and $V$ (from mass action kinetics on the microscale) and is given as follows:
\begin{align*}
	y^* =\frac{\frac{N}{K_N}+ \frac{V}{K_V}}{\frac{N}{K_N}+ \frac{V}{K_V}+1}.
\end{align*}
Then we have
$$\frac{dy^*}{dN}= \frac{\frac{1}{K_N}}{(1 +\frac{N}{K_N }+\frac{V}{K_V})^2}$$
$$\frac{dy^*}{dV} = \frac{\frac{1}{K_V}}{(1+\frac{N}{K_N }+\frac{V}{K_V})^2}.$$
The function $a(V,N)$ is chosen as 
$$a(V,N):= \frac{1 }{(1+\frac{N}{K_N }+\frac{V}{K_V})^2}.$$ 
\noindent
 The turning rate $\eta$ also depends on the pathwise gradients of necrotic and normal cells: 
$$\begin{array}{lll}D_t N &=&\partial_t N +\vartheta \cdot \nabla N
	\\
	D_t V&=&\partial_t V +\vartheta \cdot \nabla V.\end{array},$$  
	which is inspired by \cite{conte2021mathematical,othmer2002diffusion,kumar2021multiscale}.  
The last term in (\ref{-_-}) represents proliferation of bacteria.\\[-2ex]

\noindent
Due to the high dimensionality of equation (\ref{-_-}) we aim to deduce an effective macroscopic equation for bacteria density. This is done by rescaling  the time and space variables by $t\rightarrow \epsilon^2 t$ and $x\rightarrow \epsilon x$ where $\epsilon$ is a small dimensionless parameter. This is motivated by the fact that the motion of bacteria is diffusion dominated; we already mentioned that the orientation of tissue fibers does not play a relevant role, and in the present framework this would have been the only source of drift dominance. Since proliferation is much slower than migration we rescale the source term also by $\epsilon^2$. Therewith, (\ref{-_-}) becomes
\begin{equation}
	\label{1}
	\epsilon^2 \partial_t u+\epsilon \nabla_x \cdot (\vartheta u)=\mathcal{L}_u[\eta]u+\epsilon^2\mu_u(U,V,N) u.
\end{equation} 
\noindent  $\eta$ after rescaling is given as follows:
$$\eta(x,\vartheta,V,N)=\eta_0(x) exp\left(-a(V,N)(\epsilon^2\frac{\gamma_1}{K_N}\partial_t N+\epsilon\frac{ \gamma_1}{K_N} \vartheta \cdot \nabla N+\epsilon^2 \frac{\gamma_2}{K_V}\partial_t V+\epsilon \frac{\gamma_2}{K_V} \vartheta \cdot \nabla V)\right).$$ 

\noindent  As mentioned previously, $U(t,x)$ is the macroscopic density of bacteria, i.e. the zeroth order moment of $u$ with respect to $\vartheta$:
$$U(t,x)=\int_{\Theta} u(t,x,\vartheta) d\vartheta.$$ 
We will subsequently use the Hilbert expansions:
$$u(t,x,\vartheta)=\sum_{k=0}\epsilon^k u_k, \quad U(t,x)=\sum_{k=0}\epsilon^k U_k$$
and Taylor expansions for functions in \eqref{1}:
\begin{equation}
	\begin{aligned}
		\mu_u(U,V,N)&=\mu_u(U_0,V,N)+\partial_U \mu_u (U_0,V,N)(U-U_0)+O|U-U_0|^2\\
		\eta(x,\vartheta,V,N)&=\eta_0(x)\bigg[1-\epsilon a(V,N) (\frac{ \gamma_1}{K_N} \vartheta \cdot \nabla N +\frac{\gamma_2}{K_V} \vartheta \cdot \nabla V) \\ 
		& \quad + \epsilon^2\bigg(-a(V,N)(\frac{ \gamma_1}{K_N}\partial_t N+\frac{\gamma_2}{K_V}\partial_t V)+\frac{a(V,N)^2}{2}\big((\frac{ \gamma_1}{K_N} \vartheta \cdot \nabla N)^2
		\\
		&\quad +(\frac{\gamma_2}{K_V} \vartheta \cdot \nabla V)^2+ 2\frac{ \gamma_1}{K_N}  \frac{\gamma_2}{K_V}(\vartheta \cdot \nabla N ) (\vartheta \cdot \nabla V)\big) \bigg)+O(\epsilon ^3)\bigg].
	\end{aligned}
\end{equation}
By equating  powers of $\epsilon$ in \eqref{1}, we obtain \\

\noindent $\epsilon^0:$

\begin{align}
	0=-\eta_0u_0+\frac{1}{|\Theta|} \int_{\Theta} \eta_0(x)u_0 (\vartheta')d\vartheta' = \eta_0(\frac{1}{|\Theta|}U_0-u_0), \label{2}
\end{align} 
\noindent where $\frac{1}{|\Theta|}=\frac{s^{1-d}}{|{\mathbb{S}}^{d-1}|}$.\\ \\

\noindent $\epsilon^1:$
\begin{align}
	\nabla_x \cdot (\vartheta u_0)&= \eta_0(x)\bigg[\frac{1}{|\Theta|} U_1(t,x)-u_1(t,x,\vartheta)+a(V,N)\big(\frac{ \gamma_1}{K_N} \vartheta \cdot \nabla N+\frac{\gamma_2}{K_V} \vartheta \cdot \nabla V\big) u_0\nonumber
	\\
	&\quad - \frac{1}{|\Theta|} \ a(V,N) \int_{\Theta} u_0(\vartheta')\vartheta'd\vartheta'\cdot(\frac{ \gamma_1}{K_N}\nabla N+\frac{\gamma_2}{K_V}\nabla V \big)\bigg].\label{3}
\end{align}

\noindent $\epsilon^2:$
\begin{equation}
	\label{4}
	\begin{aligned}
		\partial_t u_0+\nabla_x \cdot (\vartheta u_1)=& \eta_0(x) \bigg[ \bigg(a(V,N)(\frac{ \gamma_1}{K_N}\partial_t N+\frac{\gamma_2}{K_V}\partial_t V) -\frac{ a(V,N)^2}{2}\big((\frac{ \gamma_1}{K_N} \vartheta \cdot \nabla N)^2\\[5pt]
		&  +(\frac{\gamma_2}{K_V} \vartheta \cdot \nabla V)^2 + 2\frac{ \gamma_1}{K_N} \frac{\gamma_2}{K_V}(\vartheta \cdot \nabla N ) (\vartheta \cdot \nabla V)\big) \bigg) u_0
		\\[5pt]
		&+ a(V,N)(\frac{ \gamma_1}{K_N} \vartheta \cdot \nabla N +\frac{\gamma_2}{K_V} \vartheta \cdot \nabla V )u_1- u_2 \bigg] 
		\\[5pt]
		&+\frac{1}{|\Theta|}\eta_0(x)\bigg[ U_2-a(V,N )\int_{\Theta}\vartheta'u_1(\vartheta')d\vartheta'\cdot (\frac{ \gamma_1}{K_N}\nabla N+\frac{\gamma_2}{K_V} \nabla V) 
		\\[5pt]
		& -a(V,N)(\frac{ \gamma_1}{K_N}\partial_t N+\frac{\gamma_2}{K_V}\partial_t V)U_0\ +\frac{a(V,N)^2}{2}\int_{\Theta}  \bigg( (\frac{ \gamma_1}{K_N} \vartheta' \cdot \nabla N )^2 
		\\[5pt]
		&+(\frac{\gamma_2}{K_V}\vartheta' \cdot \nabla V )^2 +  2\frac{ \gamma_1}{K_N}  \frac{\gamma_2}{K_V} (\vartheta' \cdot \nabla N ) (\vartheta' \cdot \nabla V)\bigg)  u_0(\vartheta')d\vartheta' \bigg] \\[5pt]
		&+\mu_u(U_0,V,N)u_0.
	\end{aligned}
\end{equation}
\noindent By using \eqref{2}, we have
\begin{equation}\label{u0}
	u_0=\frac{1}{|\Theta|}U_0,
\end{equation}
which means $u_0$ depends on the constant speed $s$.
\\ Now let us consider equation \eqref{3}. Using \eqref{u0} we get
\begin{align}\label{epsilon1-eqn-simplified}
	\nabla_x \cdot (\vartheta u_0)&= \eta_0(x)\bigg[\frac{1}{|\Theta|} U_1(t,x)-u_1(t,x,\vartheta)+a(V,N)\big(\frac{ \gamma_1}{K_N}  \vartheta \cdot \nabla N+\frac{\gamma_2}{K_V}  \vartheta \cdot \nabla V\big) u_0\bigg]
\end{align}
and therefore, as e.g. in \cite{engwer2015glioma,painter2013mathematical} we have the compact Hilbert-Schmidt operator
\begin{equation*}
	\begin{aligned}
		\mathcal{L}_u[\eta_0(x)]u_1&=-\eta_0(x)u_1+\frac{1}{|\Theta|}\eta_0(x)U_1\\
		&=\nabla_x \cdot (\vartheta u_0)-a(V,N)\big(\frac{ \gamma_1}{K_N}  \vartheta \cdot \nabla N+\frac{\gamma_2}{K_V}\vartheta \cdot \nabla V\big) u_0\
	\end{aligned}
\end{equation*}
defined on the $L^2$ space. We
can therefore calculate its pseudo-inverse to get an expression for $u_1$ :
\begin{equation}\label{5buruli2} u_1=-\frac{1}{\eta_0(x)}\nabla_x\cdot (\vartheta u_0)+a(V,N)(\frac{ \gamma_1}{K_N}  \vartheta \cdot \nabla N+\frac{\gamma_2}{K_V}\vartheta \cdot\nabla V)u_0\end{equation} 

\noindent Again by using \eqref{5buruli2} in \eqref{epsilon1-eqn-simplified}, we get
\begin{equation}\label{u1}
	U_1=0.
\end{equation}
Integrating \eqref{4} with respect to $\vartheta$, we get the macroscopic equation for bacteria as:
	
	\begin{equation}\label{6}\partial_t U_0+\int_{\Theta}\nabla_x\cdot(\vartheta u_1)d\vartheta=\mu_u(U_0,V,N) U_0.\end{equation}
	To compute the integral on left-hand side, by using \eqref{5buruli2} we have 
	\begin{equation}
\begin{aligned}
	\int_{\Theta}\nabla_x\cdot(\vartheta u_1)d\vartheta&=\int_{\Theta}\nabla_x \cdot  \left[ \vartheta \bigg(-\frac{1}{\eta_0(x)}\nabla_x \cdot (\vartheta u_0)+a(V,N )(\frac{ \gamma_1}{K_N} \vartheta \cdot \nabla N+\frac{\gamma_2}{K_V}\vartheta \cdot \nabla V )u_0\bigg)\right]d\vartheta\\[5pt]
	&=\nabla_x \cdot \big[ \int_{\Theta}-\frac{1}{\eta_0(x)}\vartheta \otimes \vartheta\nabla_x (\frac{1}{|\Theta|}U_0)d\vartheta\big]\\[5pt]
	&\quad + \nabla_x \cdot\big[ \frac{a(V,N) }{|\Theta|} \int_{\Theta} \vartheta\otimes \vartheta d\vartheta (\frac{ \gamma_1}{K_N}\nabla N+\frac{\gamma_2}{K_V} \nabla V ) U_0 \big]
	\\[5pt]
	&=-\nabla \cdot (\mathbb{D}_{u}\nabla U_0)+\nabla \cdot \big(\chi_a(V,N)(\frac{ \gamma_1}{K_N}\nabla N+\frac{\gamma_2}{K_V} \nabla V)U_0\big)
\end{aligned}
\end{equation}
where $\mathbb{D}_{u}:=\frac{s^2}{d\eta_0(x)}\mathbb{I}_d$, $\chi_a(V,N):=\frac{s^2 a(V,N)}{d}\mathbb{I}_d=a(V,N)\eta_0(x)\mathbb{D}_{u}.$ 
\\ \\
Therefore, the macroscopic equation for bacteria is obtained as 
\begin{equation}\label{13}
\partial_t U_0-\nabla\cdot (\mathbb{D}_{u}(x) \nabla U_0)+\nabla\cdot(\chi_a(V,N)(\frac{ \gamma_1}{K_N}\nabla N+\frac{\gamma_2}{K_V} \nabla V)U_0)
= \mu_u(U_0,V,N)U_0.
\end{equation} \normalsize \\ Here, we take the growth rate as $$\mu_u(U_0,V,N)=\alpha_u\frac{N}{K_N+N}\left(1-\frac{U_0}{K_U}-\frac{V}{K_V}-\frac{N}{K_N}\right),$$
where $\alpha_u>0$ is the proliferation rate of bacteria and $K_U >0$ is the carrying capacity for the bacteria. We consider a logistic growth slowed down by the competition for space between bacteria, normal tissue, and necrotic matter. There is also a growth enhancement factor that depends on necrotic matter for bacteria dynamics, as mentioned in \cite{palomino1998effect}.\\[-2ex]

\noindent
The above macroscopic equation for bacteria density is coupled with macroscopic equations describing the dynamics of mycolactone, normal tissue, and necrotic matter. Hence we get the following system:
\begin{equation}
	\label{14}
	\begin{aligned}
		\partial_t U&=\nabla \cdot (\mathbb{D}_{u}(x) \nabla U)-\nabla \cdot \big(\chi_a(V,N)(\frac{ \gamma_1}{K_N}\nabla N+\frac{\gamma_2}{K_V}\nabla V)U\big)  
		\\
		&\quad \quad  +\alpha_u\frac{N}{K_N+N}U\left(1-\frac{U}{K_U}-\frac{V}{K_V}-\frac{N}{K_N}\right) , \quad \quad   t>0, \ x \in  \Omega,
		\\
		\partial_t M &= D_m \Delta M+{\delta} \frac{U}{K_U+U}-\lambda M , \quad \quad   t>0, \ x \in  \Omega,
		\\
		\partial_t V &= -\beta_1 V \frac{M}{K_M} , \quad \quad   t>0, \ x \in  \Omega,
		\\
		\partial_t N &= \beta_2 \frac{V}{K_V} \frac{M}{K_M}-\gamma N,  \quad \quad   t>0, \ x \in  \Omega,
	\end{aligned}
\end{equation}
\\
where $D_m, \ \delta, \ K_M, \ \lambda, \  \beta_1 , \ \beta_2 , \ \gamma$ are positive constants, and $\Omega\subset \R^d$ is a bounded, sufficiently smooth domain.\\[-2ex]

\noindent
The evolution of mycolactone is governed by the second equation in the system (\ref{14}), where diffusion serves as the primary mechanism driving its spread. Bacteria release mycolactone, which diffuses into wounds and surrounding tissues, causing necrosis. Over time, mycolactone degrades, altering its concentration profile.
\\[-2ex]

\noindent
The third equation in (\ref{14}) describes the dynamics of normal tissue, capturing its necrotisation by the action of mycolactone on it. The fourth equation focuses on necrotic tissue, featuring its natural decay and the ongoing necrotisation of adjacent healthy tissue.
\\[-2ex]

\noindent The above parabolic upscaling was actually performed with $x\in \R^d$. Going down to $\Omega $ and obtaining no-flux boundary conditions can be done (at least informally) similarly to \cite{plaza}, yielding
\begin{equation}\label{*_}
	\left(\mathbb{D}_{u}(x) \nabla U+\chi_a(V,N)(\frac{ \gamma_1}{K_N}\nabla N+\frac{\gamma_2}{K_V} \nabla V)U\right)\cdot \nu=0, \quad \nabla M\cdot \nu=0.
\end{equation} 
For the initial conditions we take
\begin{equation}\label{*-*}
	U(0,x)=u_0, \ V(0,x)=v_0, \ N(0,x)=n_0, \ M(0,x)=m_0, \ x \in  \Omega.
\end{equation}
\noindent
The obtained PDE for bacteria dynamics features double haptotaxis and Fickian diffusion, with a space-dependent diffusion tensor $\mathbb D_u$, which is also present in the tactic sensitivity function $\chi_a$. It depends on the average cell speed (squared) and is inversely proportional to the space dimension and the turning rate. In \cite{conte2021mathematical,conte2023mathematical,corbin2021modeling,engwer2015glioma,engwer2016effective,kumar2021multiscale,mohanan2024mathematicalmodeltissueregeneration} the parabolic upscaling from KTEs led to myopic diffusion, with the corresponding tensor depending on the orientation distribution of tissue fibers. Here the simpler choice of the turning kernel leads to omission of the supplementary drift with velocity $\nabla \cdot \mathbb D_u$.

	\subsection{Approach via random position jumps}\label{subsec:PJA}
	
	\noindent
	In this subsection we deduce a continuous, macroscopic PDE for bacteria dynamics by an alternative approach, namely assuming that cells are particles jumping with certain probabilities between the sites of a lattice, as developed, e.g., in \cite{randomjumpKPH, randomjumpSAO}. We will use the notations $u,m,v,n$ for the (macroscopic) volume fractions of bacteria, mycolactone, normal and necrotic tissue, respectively. As the informal approach is rather standard, we only sketch the deduction in the 1D case, for convenience. \\[-2ex]
	
	\noindent
	Imagine a 1D lattice with nodes $x_i$ and a particle (bacterium)  jumping from one node to its neighbours or resting at its current position (node).  We are interested in $u_i(t)$, the probability of finding a particle at $x_i$ at time $t$. Let $h$ and $\Delta t$ denote the lengths of space and time steps, respectively. Let $T_i^{\pm}$ denote the probabilities of jumping to the right (+) or to the left (-), respectively. Then we have
	\\
	\begin{equation}\label{u_i}
		u_i(t+\Delta t) = \bigg(T_{i-1}^+ u_{i-1}(t) +T_{i+1}^- u_{i+1}(t)\bigg) \Delta t + \bigg(1-\Delta t (T_i^+ +T_i^-)\bigg) u_i (t).
	\end{equation}
	\noindent
	When deciding to move, a cell 'tests' its environment by comparing - via cell receptors- the amount of haptoattractants (normal tissue and necrotic matter) at its current and at neighboring locations. Therefore, the jump probability per unit time typically depends on the amount of necrotic matter and normal tissue, as well as on the cell density. Thus, we define $ T_i^{\pm}$ as:
	\begin{equation*}
		T_i^{\pm}(u,v,n):=\lambda \bigg({a}(u_{i},v_{i})+\kappa(u_{i},n_{i})(\tau(n_{i\pm1})-\tau(n_i)) \bigg),
	\end{equation*}
	\noindent
	where $v_i, \  n_i$ represent the probabilities of normal tissue and necrotic matter at a node $x_i$ at time $t,$ respectively. Define $\kappa(u_i,n_i)$ as 
	\begin{equation}\label{kappa}
		\kappa (u_{i},n_{i})= \frac{K_u K_n }{K_u K_n + u_{i}n_{i}}.
	\end{equation}
		\noindent
	The choice of $\kappa$ to be inversely proportional to cell-necrotic matter complexes is motivated by the fact that the increase of such bindings reduces the tendency of bacteria $u$ to migrate towards regions with higher necrotic matter density $n$.\\[-2ex] 
	
	\noindent
	Now divide \eqref{u_i} by $\Delta t$ and informally pass to the limit for $\Delta t \rightarrow 0$
	\begin{equation*}
		\frac{\partial u_i}{\partial t}=T_{i-1}^+u_{i-1}+T_{i+1}^-u_{i+1}-(T_i^++T_i^-)u_i. 
	\end{equation*}
	\noindent
	With the Taylor expansions (for the simplification of notation we set ${T}_i:=T(u_i,v_i,n_i)$):
	\begin{equation*}
		\begin{aligned}
			{T}_{i+1}^ -u_{i+1}&={T}_i^-u_i+h\frac{\partial}{\partial x}({T}_i^-u_i)+\frac{h^2}{2}\frac{\partial^2}{\partial x^2}({T}_i^-u_i)+\frac{h^3}{3!}\frac{\partial^3}{\partial x^3}({T}_i^-u_i)+O(h^4)
			\\
			{T}_{i-1}^ + u_{i-1}&={T}_i^+u_i-h\frac{\partial}{\partial x}({T}_i^+u_i)+\frac{h^2}{2}\frac{\partial^2}{\partial x^2}({T}_i^+u_i)-\frac{h^3}{3!}\frac{\partial^3}{\partial x^3}({T}_i^+u_i)+O(h^4).
		\end{aligned}
	\end{equation*}
	\noindent
	Hence
	\begin{equation*}
		\begin{aligned}
			\frac{\partial u_i}{\partial t}=&{T}_i^+u_i-h\frac{\partial}{\partial x}({T}_i^+u_i)+\frac{h^2}{2}\frac{\partial^2}{\partial x^2}({T}_i^+u_i)-\frac{h^3}{3!}\frac{\partial^3}{\partial x^3}({T}_i^+u_i)
			\\
			&+{T}_i^-u_i+h\frac{\partial}{\partial x}({T}_i^-u_i)+\frac{h^2}{2}\frac{\partial^2}{\partial x^2}({T}_i^-u_i)+\frac{h^3}{3!}\frac{\partial^3}{\partial x^3}({T}_i^-u_i)
			\\
			&-(T_i^++T_i^-)u_i+O(h^4),
		\end{aligned}
	\end{equation*}
	from which follows
	\begin{equation*}
		\begin{aligned}
			\frac{\partial u_i}{\partial t}
			&=   h \frac{\partial}{\partial x} \big((T_i^- - T_i^+)u_i\big) + \frac{h^2}{2} \frac{\partial}{\partial {x^2}}\big((T_i^+ + T_i^-)u_i\big)+ O(h^3)\\
			&= \lambda \bigg(-2 h^2 \frac{\partial}{\partial x} \big( \kappa (u_i,n_i) \frac{(\tau(n_{i+1})- \tau(n_{i-1}))}{2h}u_i\big) +  \frac{h^2}{2} \frac{\partial}{\partial {x^2}}\big( 2 a (u_i, v_i)u_i \big)\bigg)+ O(h^3).
		\end{aligned}
	\end{equation*}
		
	\noindent
	Passing to the limit (again, informally) for $h\rightarrow 0$ with $2{\lambda}h^2\rightarrow D,$ 
	we get
		$$\frac{\partial u}{\partial t}=\frac{D}{2}({a}(u,v)u)_{xx}-D(\kappa(u,n) u \tau'(n) n_x)_x,$$
where $u,v,n$ represent the continuous limits (which are supposed to exist) of $u_i,v_i,n_i$.\\[-2ex]

	\noindent The higher dimensional version of this equation reads 
	\begin{equation*}
		\partial_t u=D\nabla.\left(\frac{1}{2}\left({a}(u,v)+\frac{\partial{a}}{\partial u}u\right)\nabla u+\frac{1}{2}\frac{\partial{a}}{\partial v}u\nabla v-\kappa u\frac{d\tau}{dn}\nabla n\right).
	\end{equation*}  
	\noindent
 With the choice ${a}(u,v)=\frac{K_u K_v}{K_uK_v+uv}$, 
 we have 
	\begin{equation*}
		\frac{D}{2}\left({a}(u,v)+\frac{\partial{a}}{\partial u}u\right)\nabla u+\frac{D}{2}\frac{\partial{a}}{\partial v}u\nabla v=\frac{DK_u^2 K_v^2}{2(K_uK_v+uv)^2}\nabla u-\frac{D K_u K_vu^2}{2(K_uK_v+uv)^2}\nabla v.
	\end{equation*}
	\noindent
	To find the chemotactic sensitivity function $\tau(n),$ we follow the idea in \cite{randomjumpSAO} and use the receptor binding kinetics of bacteria to necrotic matter.  Thus, we take $R$ as the receptors on bacteria which are free, $R_0$ as the receptors occupied by $n$, and let $R_T$ be the total amount of receptors on bacteria, which is conserved. Let $R_T=R+R_0$ be a constant.
	\\[-2ex]
	
	\noindent
	Using simple mass action kinetics for the receptor binding 
	\begin{equation*}
		R+\frac{n}{K_n}\underset{k{-}}{\overset{k^{+}}{\rightleftharpoons }} R_0,
	\end{equation*}
 we obtain an ODE for the dynamics of occupied receptors :
	\begin{equation} \label{ode-reaction}
		\frac{dR_0}{dt}=k^+\frac{n}{K_n}R-k^-R_0,
	\end{equation}
	\noindent
	Since attachment/detachment of receptors is very fast, we assume it to quickly reach its equilibrium. The steady state of (\ref{ode-reaction}) is
	\begin{equation*}
		R_0^*=\frac{{K}R_T\frac{n}{K_n}}{1+{K}\frac{n}{K_n}},
	\end{equation*}where ${K}=\frac{k^+}{k^-}.$
		\noindent 
	If we consider $\tau(n)$ to be a chemical mechanism to measure the concentration of chemoattractant $n$, we could assume $\tau(n)$ to be proportional to the amount of receptors occupied by $n$, that is
	\begin{equation*}
		\tau(n)={b} R_0^*,
	\end{equation*}where ${b}$ is a positive constant. Hence
	\begin{equation*}
		D\kappa(u,n) u\frac{d\tau}{dn}\nabla n=\frac{D\kappa(u,n){b}\frac{K}{K_n}R_Tu}{(1+{K}\frac{n}{K_n})^2}\nabla n,
	\end{equation*}
	is obtained. 
	\\
	
	\noindent
	Assuming the attachment rate $k^+$ and the detachment rate $k^-$ of bacteria to the necrotic matter to be almost the same, we get ${K} \approx 1$ and let $K_1:=D{b} {K}R_T,$ and using  \eqref{kappa} we have
	\begin{equation*}
		D\kappa(u,n) u\frac{d\tau}{dn}\nabla n=\frac{K_1 K_u  u}{(1+\frac{n}{K_n})^2 (K_u K_n +un)}\nabla n. 
	\end{equation*}
	
	\noindent
	The system describing the dynamics of $u,n,v,m$ then reads:
	\begin{equation}
		\begin{aligned}\label{burulinonlinear}
			\partial_tu&=\nabla \cdot (D_u(u,v)\nabla u)-\nabla \cdot (\chi_1(u,v)u\nabla v)-\nabla \cdot          (\chi_2(u,n)u\nabla n)  \\& \qquad \qquad +\alpha \frac{ n}{K_n+n}u(1-\frac{u}{K_u}-\frac{v}{K_v}-\frac{n}{K_n}),  \quad \quad   t>0, \ x \in  \Omega,\\
			\partial _tm&=D_m\Delta m+\delta \frac{u}{K_u+u}-\lambda m,  \quad \quad   t>0, \ x \in  \Omega,\\
			\partial_tv&=-\beta _1 v \frac{m}{K_m} , \quad \quad   t>0, \ x \in  \Omega,\\
			\partial_tn&=\beta _2\frac{v}{K_v}\frac{m}{K_m} - \gamma n ,  \quad \quad   t>0, \ x \in  \Omega,
		\end{aligned}
	\end{equation}
	
	\noindent where $$D_u(u,v)=\frac{DK_u^2 K_v^2 }{2(K_uK_v+uv)^2}, \ \chi_1(u,v)=\frac{D K_u K_v u}{2(K_uK_v+uv)^2}, \ \chi_2(u,n)=\frac{K_1 K_u  }{(1+\frac{n}{K_n})^2 (K_u K_n +un)}.$$
	As before, $\Omega \subset \R^d$ is a bounded, smooth enough domain, and we consider no-flux boundary conditions
	\begin{equation}\label{5}
		(D_u(u,v)\nabla u-\chi_1(u,v) \nabla v-\chi_2(n)u \nabla n)\cdot \nu =0, \ \nabla m \cdot \nu=0 \quad   \quad \quad   t>0, \ x \in  \partial\Omega,
	\end{equation} where $\nu$ represents the outer unit normal vector to $\partial \Omega$. 
	\\[-2ex]
	
	\noindent
	For the initial conditions we take
	\begin{equation}\label{6}u(0,x)=u_0,\  v(0,x)=v_0, \ n(0,x)= n_0, \ m(0,x)=m_0 \quad \quad x \in \Omega .\end{equation}
	
\section{Existence of solutions}\label{sec:analysis}

\noindent
In this section we prove the global existence of solutions to the obtained macroscopic systems. We will begin with system \eqref{14}, \eqref{*_}, \eqref{*-*}; the proof for system \eqref{burulinonlinear}, \eqref{5}, \eqref{6} is very similar. \\[-2ex]

\noindent
Before performing the analysis, we nondimensionalise system \eqref{14}. Therefore, we take
\begin{equation}\label{dim-var}
	\tilde{u}=\frac{U}{K_U}, \ \tilde{v}=\frac{V}{K_V}, \ \tilde{n}=\frac{N}{K_N}, \ \tilde{m}=\frac{M}{K_M}, \ \tilde{x}={x}\sqrt{\frac{\alpha_u}{D_m}}, \ \tilde{t}=t \alpha_u,
\end{equation}
where $K_M$ is the maximum possible level of mycolactone. By using the transformations in \eqref{dim-var}, system \eqref{14} becomes
\begin{equation}\begin{aligned}
		\label{nondimsysburuli2}
		\tilde{u}_{\tilde{t}}&=\nabla \cdot \left({\tilde{\mathbb{D}_u}}\nabla \tilde{u}\right)-\nabla \cdot \left(\tilde{a}(\tilde{v}, \tilde{n}) \tilde{\mathbb{D}_u} \tilde{\gamma_1} \tilde{u}\nabla \tilde{n} \right)-\nabla \cdot \left(\tilde{a}(\tilde{v}, \tilde{n}) \tilde{\mathbb{D}_u} \tilde{\gamma_2}\tilde{u} \nabla \tilde{v}\right)\\
		&\quad \ + \frac{\tilde{n}}{1+\tilde{n}} \tilde{u} \left(1-\tilde{u}-\tilde{v}-\tilde{n}\right)
		\\
		\tilde{m}_{\tilde{t}}&=\Delta \tilde{m}+\frac{\tilde{\delta} \tilde{u}}{1+\tilde{u}}-\tilde{\lambda} \tilde{m}
		\\
		\tilde{v}_{\tilde{t}}&=-\tilde{\beta}_1\tilde{m}\tilde{v}
		\\
		\tilde{n}_{\tilde{t}}&=\tilde{\beta}_2\tilde{m}\tilde{v}-\tilde{\gamma}\tilde{n},
	\end{aligned}
\end{equation}
where 
\begin{align*} 
	&\widetilde{\mathbb{D}_u}=\frac{\mathbb{D}_{u}}{D_m},\qquad \tilde{\gamma_1}=\gamma_1 \eta_0 , \qquad \tilde{\gamma_2}=\gamma_2 \eta_0 ,  \qquad \tilde{\delta}=\frac{\delta}{ K_M \alpha_u}, \qquad \tilde{\beta}_1=\frac{{\beta}_1}{\alpha_u},\\ 
	&  \qquad \tilde{\beta_2}=\frac{\beta_2 }{K_N \alpha_u},  \qquad \tilde{\lambda}=\frac{\lambda  }{\alpha_u} \qquad \tilde{\gamma} = \frac{\gamma}{\alpha_u}, \qquad \tilde{a}(\tilde{v}, \tilde{n})=\frac{1}{(1+\tilde{n} +\tilde{v})^2}.
\end{align*}

\noindent To simplify the notation we omit all tildes in system \eqref{nondimsysburuli2} and continue with this system along with the initial conditions (\ref{*-*}) and boundary conditions (\ref{*_}). 
\\

\noindent \textbf{Assumption: (A1)}  Let $ \mathbb{D}_{u}\in \big(C^{2,\rho}(\Omega)\cap C(\bar{\Omega})\big)^{d \times d},$ where $\rho\in(0,1)$ and observe that there exists $\Lambda>0$ such that for any $\xi\in \mathbb{R}^d$ and $x\in \Omega;$ $$\xi^{t}\mathbb{D}_{u}(x)\xi\geq \Lambda|\xi|^2$$ is satisfied.

\subsection{Local existence of a solution to \eqref{14}}
The following result establishes the local existence of solutions to (\ref{nondimsysburuli2}). The proof is similar to that in \cite{TAOchemohapto}.
\begin{Lemma}\label{lemma1buruli2}
	Let $\Omega\subset \mathbb{R}^d$ $(d\geq 1)$ be a bounded domain with a smooth boundary. Suppose that $u_0\in W^{1,\infty}(\Omega), \ m_0, \ n_0, \ v_0 \in W^{2,\infty}(\Omega)$ satisfy the compatibility conditions on boundary and are non-negative. Let \textbf{(A1)} be satisfied. Then there exists a quadruple of non-negative functions $(u,m,v,n)$ $\in C^0([0,T_{\max}) \times\bar{\Omega} ) \cap C^{2,1}((0,T_{\max})\times{\Omega})$ 
	which	solve system \eqref{nondimsysburuli2} classically in $(0,T_{\max})\times\Omega$. Moreover, we have the following dichotomy:
	\begin{equation}
		\begin{aligned}\label{dichoburuli2} 
			&\text{either } T_{max}=\infty \  \text{or } \\
			&  \limsup_{t\nearrow T_{\max}}(||u(t)||_{L^{\infty}(\Omega)}+||m(t)||_{W^{1,\infty}(\Omega)}+||v(t)||_{W^{1,\infty}(\Omega)}+||n(t)||_{W^{1,\infty}(\Omega)})=\infty.
		\end{aligned}
	\end{equation}
\end{Lemma}
\begin{proof} 
	We define the following closed, bounded, and convex set: 
	\begin{equation*}
		Q:=\{\bar{u}\in L^{\infty}((0,T) \times {\Omega}):0\leq \bar{u}\leq ||u_0||_{L^{\infty}(\Omega)}+1 \text{  a.e in } \Omega \times (0,T) \}
	\end{equation*}
	with $T>0$ small, to be fixed below. We consider a fixed point problem with the operator $F$ defined on $Q$ such that $F({\bar{u}})=u$, where $u$ is the first component of the solution $(u,m,v,n)$ to the following system for $ x \in \Omega , \ t \in (0,T)$:
	
	\begin{equation}\label{15}
		\begin{aligned}
			{u}_{{t}}&=\nabla \cdot \left({{\mathbb{D}_u}}\nabla {u}\right)-\nabla \cdot \left({a}({v}, {n}) {\mathbb{D}_u} {\gamma_1} {u}\nabla {n} \right)-\nabla \cdot \left({a}({v}, {n}) {\mathbb{D}_u} {\gamma_2}{u} \nabla {v}\right)\\
			&\quad \ + \frac{{n}}{1+{n}} {u} \left(1-{\bar u}-{v}-{n}\right)
			\\[5pt]
			\partial_t m &= \Delta m+\delta \frac{\bar u}{1+\bar u}-\lambda m
			\\[5pt]
			\partial_t v &= -\beta_1 v m
			\\[5pt]
			\partial_t n &= \beta_2 v m-\gamma n \\[5pt]
			\big(\mathbb{D}_{u} &\nabla u-a(v,n)\mathbb{D}_{u}u( \gamma_1 \nabla n+\gamma_2 \nabla v) \big)\cdot \vartheta=0, \quad \nabla m\cdot \vartheta=0 \quad x \in \partial \Omega , \ t> 0 \\[5pt]
			u(0,x)&=u_0 (x), \ m(0,x)=m_0 (x),\ v(0,x)=v_0 (x), \ n(0,x)=n_0 (x), \quad x \in \Omega.
		\end{aligned}
	\end{equation} 
	For any $q> \frac{d+2}{2}$, by using Theorem A.1 in \cite{Giesselmann}, there exists a solution $m(t,x)$ such that 
	\begin{equation}
		\label{m-lq}||\nabla m||_{L^{q}((0,T) \times \Omega)}\leq \tilde{c_1},
	\end{equation} 
	where $\tilde{c_1}$ is a positive constant which depends on $||u_0||_{L^{\infty}(\Omega)},||m_0||_{W^{1,\infty}(\Omega)},||n_0||_{W^{1,\infty}(\Omega)} $, and  $||v_0||_{W^{1,\infty}(\Omega)}$. The following constants $\tilde{c_2}, c_3, c_4, c_5 ,..$ also depend on the above values. \\
	
	\noindent	By solving the ODE for $v$ we get
	\begin{equation}
		\label{17buruli2}
		v(t,x)=v_0(x)\exp\left(-\beta_1\int_0^t m(s,x) ds \right),  \quad \quad t \in (0,T), \ x \in \Omega.  
	\end{equation} 
	
	\noindent From \eqref{m-lq} and \eqref{17buruli2}, we get
	\begin{equation}
		\label{v-lq}||\nabla v||_{L^{q}((0,T) \times \Omega)}\leq \tilde{c_2}.
	\end{equation} 
	\noindent Using the variation of constants formula for the fourth equation of system \eqref{nondimsysburuli2}, we obtain the following representation:
	\begin{equation}\label{18buruli2}
		n(t,x)=e^{-\gamma t} n_0(x)+\beta_2\int_0^t m(s,x)v(s,x) e^{-\gamma(t-s)}ds, \quad \quad t \in (0,T), \ x \in \Omega.
	\end{equation} 
	\noindent Using the previous estimates \eqref{m-lq} and \eqref{v-lq}, we have
	\begin{equation}\label{n-lq}
		||\nabla n||_{L^{q}((0,T) \times \Omega)}\leq c_3.
	\end{equation}
	
	\noindent Hence, we can rewrite the first equation in system  \eqref{nondimsysburuli2} as
	\begin{equation*}
		\partial_t u=\nabla \cdot (\mathbb{D}_{u}(x) \nabla u - g(t,x)u) +\frac{n}{n+1}u\left(1-\bar u-v-n\right)  
	\end{equation*}
	where $ ||g||_{L^{q}((0,T) \times \Omega)}\leq c_4$. Using the parabolic theory(Thm V.2.1 and V.1.1 from \cite{LAD}), 
	$||u||_{L^{\infty}((0,T) \times \Omega)}\leq c_5$ and we get $c_6$ such that $||u||_{C^{\beta, \frac{\beta}{2}}((0,T) \times \Omega)}\leq c_6$  for some $\beta >0$.  Hence $||u||_{L^{\infty}((0,T) \times \Omega)}\leq ||u_0||_{L^{\infty}(\Omega)} + c_7 T^{\frac{\beta}{2}}$ for some $T<1$ such that $ c_7 T^{\frac{\beta}{2}} <1$. 
	\\
	
	\noindent
	Now, consider the fixed point problem $F(\bar{u})=u$. By using the above, we have $F(Q)\subset Q.$ Since $F$ is continuous and $\overline{F(Q)}$ is a compact subset in $L^{\infty}((0,T) \times \Omega),$ by applying Schauders's fixed point theorem we have that there exists at least one fixed point $u\in Q.$ Hence by parabolic regularity theory $(u,v,n,m)$ is obtained as the classical solution of system \eqref{nondimsysburuli2} in $C^0([0,T) \times \bar{\Omega})\cap C^{2,1}((0,T)\times {\Omega}).$ By applying the maximum principle in \cite{protter2012maximum}, we get the non-negativity of $u$. Using Theorem A.3 in \cite{lenz} we get the non-negativity of $m$. By using the non-negativity of $m$, the assumption on $v_0,$ and \eqref{17buruli2}, the non-negativity of $v$ is obtained. By using the positivity of $\gamma$, the assumption on $n_0$, \eqref{18buruli2}, and by the non-negativity of $m \text{ and } v$ we get the nonnegativity of $n$.
	\\
	
	\noindent Define\begin{equation}\label{Aburuli2}
		A :=||u(0)||_{L^{\infty}(\Omega)}+||m(0)||_{W^{1,\infty}(\Omega)}+||v(0)||_{W^{1,\infty}(\Omega)}+||n(0)||_{W^{1,\infty}(\Omega)}.
	\end{equation}
	and let 
	$$T_{max}= \sup{\{T \in (0,\infty ] \ : \ \text{there exists a solution in } C^0([0,T_{max}) \times \bar{\Omega}))\cap C^{2,1}((0,T_{\max})\times {\Omega})\}}$$ be the maximum existence time. Let $T_{max} < \infty$ and
	$$\limsup_{t\nearrow T_{\max}}(||u(t)||_{L^{\infty}(\Omega)}+||m(t)||_{W^{1,\infty}(\Omega)}+||v(t)||_{W^{1,\infty}(\Omega)}+||n(t)||_{W^{1,\infty}(\Omega)})<\infty.$$
	Hence there exists $A'>0$ such that 
	$$||u(t)||_{L^{\infty}(\Omega)}+||m(t)||_{W^{1,\infty}(\Omega)}+||v(t)||_{W^{1,\infty}(\Omega)}+||n(t)||_{W^{1,\infty}(\Omega)} < A' \quad \text{for all } t < T_{max}.$$ 
	Next we choose $T_1$ sufficiently small such that the above procedure works in $[0,T_1]$ after replacing $A$ (defined in (\ref{Aburuli2}))
	by $A'$. Define $T_2 = T_{max} -\frac{ T_1}{2}$. As $T_1$ depends only on $A'$, we may choose $u(T_2,.), m(T_2,.),v(T_2,.) ,n(T_2,.)$ as new initial values. Hence by the above procedure, we obtain a solution in $[T_2, T_{max} + \frac{T_1}{2}]$. This contradicts the maximality of $T_{max}$. Hence (\ref{dichoburuli2}) holds.
	
\end{proof}	

\subsection{Boundedness of the solution components for \eqref{14}, global existence}
\noindent To prove the global existence of solutions we will begin by proving the boundedness of functions $u, \ m, \ v, \text{ and } n$ satisfying (\ref{nondimsysburuli2}), (\ref{*-*}), and (\ref{*_}). We follow the proof in \cite{kumar2021multiscale}.
\begin{Remark}
	The boundedness of the solution components is independent of the local existence of the solution.   However, for convenience, the same maximal time $T_{max}$ from the local existence proof is often used in the boundedness proof. This allows for consistency in the analysis while treating the existence and boundedness as separate properties.
\end{Remark}

\begin{Lemma}\label{lem1buruli2}
	There exists $C_m>0$ such that
	
	\begin{equation*}
		\begin{aligned}
			||m||_{L^{\infty}([0,T_{\max})\times\Omega)}&\leq ||m_0||_{L^{\infty}(\Omega)}+\frac{\delta}{\lambda}\\ 
			||\nabla m||_{L^{\infty}([0,T_{\max})\times\Omega)}&\leq C_m(||\nabla m_0||_{L^{\infty}(\Omega)}+1),
		\end{aligned}
	\end{equation*}

\end{Lemma}
\begin{proof}
	\noindent 
	Let $pm^{p-1}$ be a test function for the $m$-equation of system \eqref{nondimsysburuli2} for $p>1$. Integrating over $\Omega,$  we have for any $\epsilon\in(0,1)$
	
	\begin{equation*}
		\begin{aligned}
			\frac{d}{dt}\int_{\Omega} m^p&=- p(p-1)  \int_{\Omega} m^{p-2}|\nabla m|^2+p\delta\int_{\Omega}m^{p-1}\frac{u }{1+u}-\lambda p \int_{\Omega} m^p
			\\ &= -\frac{4(p-1)}{p} \int_{\Omega}|\nabla m^{p/2}|^2+p\delta\int_{\Omega}m^{p-1}\frac{u}{1+u}-\lambda p \int_{\Omega} m^p
			\\ &\leq -\frac{4(p-1)}{p}  \int_{\Omega}|\nabla m^{p/2}|^2+p\delta\int_{\Omega}|m|^{p-1}-\epsilon\lambda p \int_{\Omega} m^p.
		\end{aligned}
	\end{equation*}
	
	\noindent For some $A>0$ and $\psi\in(0,1)$, we have
	
	\begin{equation*}
		\begin{aligned}
			p\delta\int_{\Omega}|m|^{p-1}&=p\delta\int_{ \{x\in\Omega \ : \ |m(x)|\leq A\}}|m|^{p-1}+p\delta\int_{\{x\in\Omega \ : \ |m(x)|> A\}}|m|^{p-1}
			\\
			&\leq p\delta A^{p-1}|\Omega|+p\delta\psi\int_{\Omega}|m|^{p-1}.
		\end{aligned}
	\end{equation*}
	
	\noindent	If we choose $A:=\frac{\delta(1-\psi)^{\frac{1}{p-1}}}{\lambda(1-\epsilon)},$ then we get
	\begin{equation}\label{16}
		\frac{d}{dt}\int_{\Omega} m^p+\epsilon\lambda p\int_{\Omega} m^p\leq \frac{p\delta^p|\Omega|}{(1-\epsilon)^{p-1}\lambda^{p-1}}.
	\end{equation} 
	Multiplying both sides of \eqref{16} with $e^{\epsilon\lambda p t}$ and integrating both sides from $0$ to $t$ for any $t\in[0,T_{\max})$, we obtain
	$$\int_{\Omega}m^p\leq \int_{\Omega}m_0^p+\frac{\delta^p|\Omega|}{\epsilon\lambda^{p}(1-\epsilon)^{p-1}}.$$
	Additionally,
	\begin{equation*}
		\begin{aligned}
			||m(t)||_{L^{\infty}(\Omega)}&=\lim_{p\rightarrow \infty}\left(\int_{\Omega}m^p(t)\right)^{\frac{1}{ p}}
			\\
			&\leq \lim_{p\rightarrow \infty} \left(\int_{\Omega}m_0^p+\frac{\delta^p|\Omega|}{\epsilon\lambda^{p}(1-\epsilon)^{p-1}}\right)^{\frac{1}{ p}}\\
			&\leq||m_0||_{L^{\infty}(\Omega)}+\lim_{p\rightarrow \infty}\frac{\delta |\Omega|^{\frac{1}{p}}}{\lambda \epsilon^{\frac{1}{ p}}(1-\epsilon)^{\frac{p-1}{ p}} }.
			\\
			&\leq||m_0||_{L^{\infty}(\Omega)}+\frac{\delta}{\lambda}.
		\end{aligned}
	\end{equation*} 
	
	\noindent	Because of the arbitrariness of $\epsilon\in(0,1),$ we get the desired result.
	By using the second equation of system \eqref{nondimsysburuli2}, we have for all $t\in(0,T_{\max})$
	$$m(t)= e^{ t \Delta} m_0+\int_0^t e^{ (t-s)\Delta}\big(\delta \frac{u}{1+u} - \lambda m\big) ds,$$ where $e^{ t \Delta}$ denotes the Neumann heat semigroup. By using the above estimate along with Lemma 1.3 (ii) and 1.3 (iii) from \cite{winkler}, we obtain for all $t \in (0,T_{max})$
	\begin{equation*}
		\begin{aligned}
			||\nabla m(t)||_{L^{\infty}(\Omega)}&\leq ||\nabla e^{ t \Delta} m_0||_{L^{\infty}(\Omega)}+\int_0^t ||\nabla e^{ (t-s)\Delta}\big(\frac{\delta u}{1+u} - \lambda m\big)||_{L^{\infty}(\Omega)}  ds
			\\
			&\leq k_1 e^{- \lambda_1 t } ||\nabla m_0||_{L^{\infty}(\Omega)}+k_2\int_0^t  e^{- \lambda_1(t-s)}(1+(t-s)^{\frac{-1}{2}}) || \frac{ \delta u}{1+u} - \lambda m||_{L^{\infty}(\Omega)} ds
			\\
			&\leq k_1 e^{- \lambda_1 t } ||\nabla m_0||_{L^{\infty}(\Omega)}+k_2\int_0^t  e^{- \lambda_1(t-s)}(1+(t-s)^{\frac{-1}{2}}) ( \delta + \lambda ||m||_{L^{\infty}(\Omega)}) ds
			\\
			&\leq C_m (||\nabla m_0||_{L^{\infty}(\Omega)}+1),
		\end{aligned}
	\end{equation*} 
	\normalsize where $ \lambda_1>0$ denotes the first non-zero eigenvalue of $-\Delta$ in $\Omega\subset \mathbb{R}^d$ under Neumann boundary conditions. 
	\\
\end{proof}
\begin{Remark}\label{rem1buruli2}
	By using \eqref{17buruli2}, non negativity of $m$, the assumption on $v_0$ and Lemma \ref{lemma1buruli2}, we can conclude that there exist constants $C_1, c_1> 0$ such that $||v||_{L^{\infty}([0,T_{\max})\times\Omega)}\leq C_1.$ By using \eqref{17buruli2} and Lemma \ref{lem1buruli2}, we also get 
	$$||\nabla v||_{L^{\infty}([0,T_{\max})\times\Omega)}\leq c_1.$$
\end{Remark}
\begin{Remark}
	\label{rem2buruli2}
	\noindent	By using  \eqref{18buruli2}, the positivity of $\gamma$, the assumption on $n_0$, Lemma \ref{lem1buruli2}, and Remark \ref{rem1buruli2}, there exist constants $C_2, c_2> 0$ such that $||n||_{L^{\infty}([0,T_{\max})\times\Omega)}\leq C_2$. Using \eqref{18buruli2}, Lemma \ref{lem1buruli2}, and Remark \ref{rem1buruli2}, we also get $$||\nabla n||_{L^{\infty}([0,T_{\max})\times\Omega)}\leq c_2.$$
\end{Remark}

\begin{Lemma}\label{lem2buruli2}
	Let \textbf{(A1)} be satisfied. Then for any $p>1,$  there exists $C>0$ such that for $t\in(0,T_{\max}),$ we have
	$$||u(t)||_{L^{p}(\Omega)}\leq C.$$ 
\end{Lemma}
\begin{proof}
	Let $pu^{p-1}$ be a test function for the first equation of system \eqref{nondimsysburuli2}. Then integrating over $\Omega$, using the no flux boundary conditions of the system and Young's inequality, and denoting $||\mathbb{D}_{u}||_{L^{\infty}(\Omega)}:=D_0$, we have
	\begin{equation*}
		\begin{aligned}
			\frac{d}{dt}\int_{\Omega}u^p=&-\frac{4(p-1)}{p}\int_{\Omega}(\nabla u^{\frac{p}{2}})^\top\mathbb{D}_{u}(x)\nabla u^{\frac{p}{2}}
			+(p-1)\int_{\Omega}(\nabla u^{{p}})^{\top} a(v,n)\mathbb{D}_{u}(x)(\gamma_1\nabla n+\gamma_2\nabla v)
			\\
			&+\int_{\Omega} p\frac{n}{n+1}u^p\left(1-u-v-n\right)
			\\ 
			\leq& -\frac{4(p-1)\Lambda}{p}\int_{\Omega} |\nabla u^{\frac{p}{2}}|^2+ \frac{4(p-1)\Lambda}{p}\int_{\Omega} |\nabla u^{\frac{p}{2}}|^2+\tilde{C}\int_{\Omega}  u^p
			+ p \int_{\Omega}  u^p,
			\\ 
			\leq& (\tilde{C}+ p)\int_{\Omega}u^p,
		\end{aligned}
	\end{equation*} where 
	\begin{equation*}
		\begin{aligned}
			\tilde{C}&=\frac{(p-1)p}{4\Lambda} D_0^2 (\gamma_1 ||\nabla n||_{L^{\infty}}+ \gamma_2 ||\nabla v||_{L^{\infty}})^2,
		\end{aligned}
	\end{equation*}
	
	\noindent By Gronwall's inequality, we obtain
	$$\int_{\Omega}u^p\leq \exp(t( \tilde{C}+p)) \int_{\Omega} u_0^p .$$
	Thus, we get the desired result
	\begin{equation}
		||u(t)||_{L^p(\Omega)} \leq C(p,T_{max}, ||u_0||_{L^p(\Omega)}).
	\end{equation}
\end{proof}
\begin{Remark}\label{rem3}
	By the usual Moser iteration process as in \cite{MARINO2019154}, we can prove that there exists $\tilde{C_1}(T_{max})>0$ such that
	$$||u(t)||_{L^\infty (\Omega)} \leq \tilde{C_1} \quad \text{for all } t \in (0,T_{max}). $$
\end{Remark}

\begin{Theorem}\label{thm2}
	Under the assumptions of Lemma \ref{lemma1buruli2}, system \eqref{nondimsysburuli2} has a non-negative classical solution $(u,v,n,m)$ which is global in time.
\end{Theorem}
\begin{proof}
	This is now an immediate consequence of the previous results.
\end{proof}

\subsection{Global existence of solutions to \eqref{burulinonlinear}}

\noindent
The global existence of solutions to \eqref{burulinonlinear} follows the same lines. The fixed-point mapping $\phi $ is defined on the same set $Q$ as above: let $ A= ||u_0||_{L^{\infty}({\Omega})}$. For $\bar{u} \in Q$, let $\phi(\bar{u}) \equiv u$ denote the first component of the solution of the following system for $ x \in \Omega, \ {t}  \in (0,T)$:

\begin{equation}\label{15}
	\begin{aligned}
		\partial_t u &=\nabla \cdot ( {D_u} (\bar{u},v)\nabla u)-\nabla \cdot \left( {\chi_1 }(\bar{u},v))u\nabla v\right)-\nabla \cdot \left( {\chi_2}(\bar{u},n)u \nabla {n}\right) \\
		& \quad \quad \quad +  u{f}(\bar{u},v,n)
		\\
		\partial_t m &= D_m \Delta m+\delta \frac{\bar{u}}{K_u + \bar{u}}-\lambda m
		\\
		\partial_t v &= -\beta_1 v m
		\\
		\partial_t n &= \beta_2 v m-\gamma n \\
	\end{aligned}
\end{equation} 
\noindent
with boundary conditions
\begin{equation*}
	\bigg( {D_u} (\bar{u},v)\nabla u-{\chi_1 }(\bar{u},v))\nabla v- {\chi_2}(\bar{u},n)\nabla {n}  \bigg)\cdot \nu =0, \quad \nabla m\cdot \nu=0 \quad \ t> 0 ,\  x \in \partial \Omega , 
\end{equation*}
\noindent
and initial conditions
\begin{equation*}
	u(0,x)=u_0 (x), \ m(0,x)=m_0 (x),\ v(0,x)=v_0 (x), \ n(0,x)=n_0 (x), \quad x \in \Omega,
\end{equation*}
where $u_0, \ m_0$ are non negative and $v_0 ,\ n_0 >0.$ Then Schauder's fixed point theorem and standard parabolic regularity show  the local  existence of a solution in $C^0([0, T) \times \bar{\Omega})\cap C^{2,1}((0,T)\times {\Omega}).$ Analogous estimates as above ensure the boundedness of such solution, while a dichotomy analogous to that in Lemma \ref{lemma1buruli2} leads to the global existence. The solution-dependent motility coefficients $D_u,\ \chi_1,\ \chi _2$ are no particular challenge, since they are well-behaved in terms of boundedness and non-degeneracy. For details of this proof we refer to \cite{diss-shimi}.

\section{Numerical simulations}\label{sec:simulations}

\noindent
In this section we perform numerical simulations of the macroscopic systems for Buruli ulcer spread obtained in Section \ref{sec:models}. We employ the finite difference method (FDM) to discretize the equations. The diffusion terms are treated with the standard central difference scheme, while for the taxis terms we use a first-order upwind discretisation. For the time derivatives an implicit-explicit (IMEX) scheme is used, thereby treating the diffusion parts implicitly and discretizing the taxis and source terms with an explicit Euler method.\\[-2ex]

\noindent
In order to conduct simulations of the models, a set of realistic parameters is required. To establish an estimated range for these parameters, we consulted several relevant studies in the literature.  
The simulations are performed obtaining the non dimensionalized parameters using the following dimensional parameter set:
	\begin{center}
	\begin{table}[hbt]
		\begin{tabular}{ |c|c|c| } 
		\hline
		$D_u$ & $10^{-4} \ mm^2/h$ & \cite{Licata2016}\\ 
		$D_m$ & $0.086 \ mm^2/h$ & \cite{alghamdi2024combinedexperimentalmathematicalstudy} \\ 
		$\delta $ & $1/h$ &  \cite{alghamdi2024combinedexperimentalmathematicalstudy} \\ 
		$\lambda$ & $0.1/h$ & \cite{alghamdi2024combinedexperimentalmathematicalstudy}  \\ 
		$\beta_1$ & $0.3/h$ &  \cite{Shi2016,King2003,McGann2009}, this work \\
		$\beta_2$ & $0.3/h$ & \cite{Shi2016, King2003,McGann2009}, this work\\
		$\gamma $ &$3*10^{-4}/h$ & this work\\
		$\alpha_u$ & $.005/h$ & \cite{Roberts2019}, this work \\
		$\gamma_1 $ &$10^{-5}- 10^{-3}h$ & this work\\
		$\gamma_2 $ &$10^{-5}- 10^{-3}h$ & this work\\
		$\gamma_3 $ &$10^{-4}h$ & this work\\
		$\eta_0 $ &$10/h$ & this work\\
		$K_U $ &$10^{4}  \ cells \ mm^{-2}$ & \cite{Roberts2019}, this work \\
		$K_M $ &$10^{4}\ mol/L$ & this work \\
		$K_V $ &$10^{4} \ cells \ mm^{-2}$ & this work\\
		$K_N $ &$10^{4}  \ cells \ mm^{-2}$ & this work \\
		\hline
				\end{tabular}
				\caption{Parameters used in the simulations}
				\end{table}
					\end{center}
\subsection{Simulations of system \eqref{14}}

\noindent
We consider a variety of scenarios, where $\gamma_1$ and $\gamma_2$ take different values. Furthermore, we simulate a setting in which bacteria perform chemotaxis towards mycolactone, as well as another one, in which initially there is less normal tissue.\\[-1ex]

\noindent
\textbf{Scenario 1:}  $\gamma_1 \ = \ \gamma_2 \ = \ 10^{-5}h$\\

\noindent
\textbf{Scenario 2:} $\gamma_1 \ = \ 10^{-3}h \ > \ \gamma_2 \ = \ 10^{-5}h$\\

\noindent
\textbf{Scenario 3:} $\gamma_1 \ =  \ 10^{-5}h \ < \ \gamma_2 \  = \ 10^{-3}h$\\

\noindent
\textbf{Scenario 4:} Bacteria additionally perform a linear chemotaxis towards mycolactone. To this aim we consider for this scenario the macroscopic equation of bacteria in the form
\begin{equation*}
	\partial_t u =\nabla \cdot (\mathbb{D}_{u}(x) \nabla u)-\nabla \cdot(a(v,n)u(\gamma_1\nabla n+\gamma_2 \nabla v)) -\nabla \cdot({\gamma_3} u \nabla m)+\frac{n}{n+1}u\left(1-{u}-v-n\right),    
\end{equation*}
where $ \gamma_3 = 10^{-4}mm^2/h$.
\\[-2ex]

\begin{figure}[htbp!]
	\centering
	\begin{subfigure}{0.24\textwidth}
		\centering
		\includegraphics[width=\textwidth]{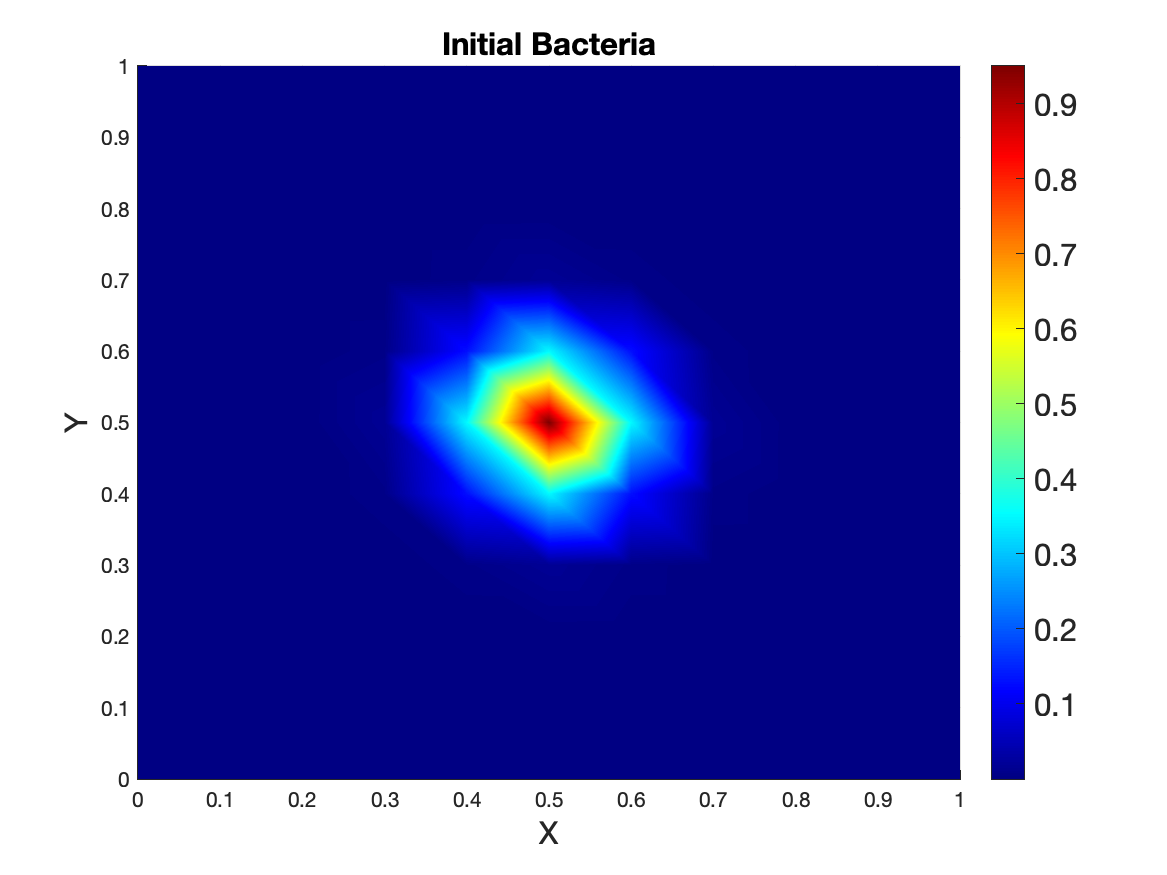}
		\caption{Bacteria }
		\label{} 
	\end{subfigure} 
	\begin{subfigure}{0.24\textwidth}
		\centering
		\includegraphics[width=\textwidth]{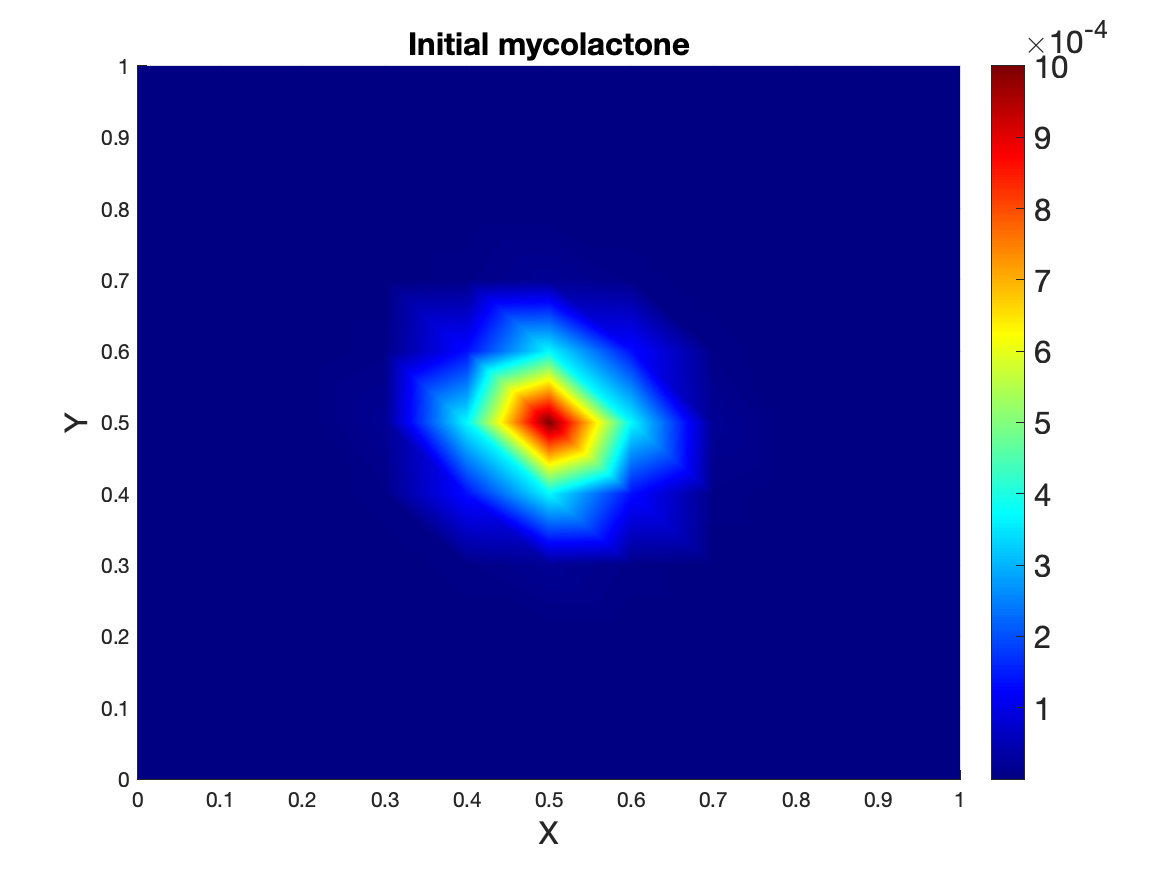}
		\caption{Mycolactone}
	\end{subfigure}    
	\begin{subfigure}{0.24\textwidth}
		\centering
		\includegraphics[width=\textwidth]{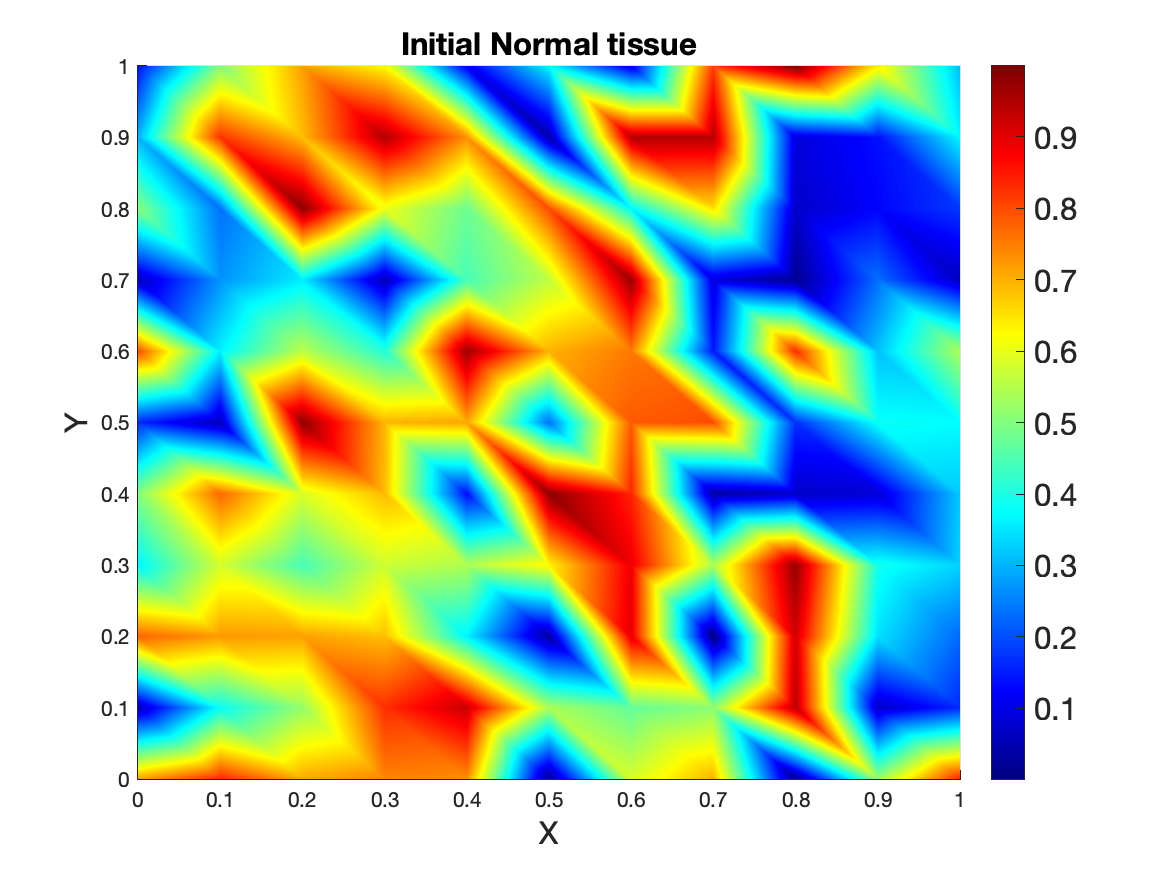}
		\caption{Normal tissue}
	\end{subfigure} 
	\begin{subfigure}{0.24\textwidth}
		\centering
		\includegraphics[width=\textwidth]{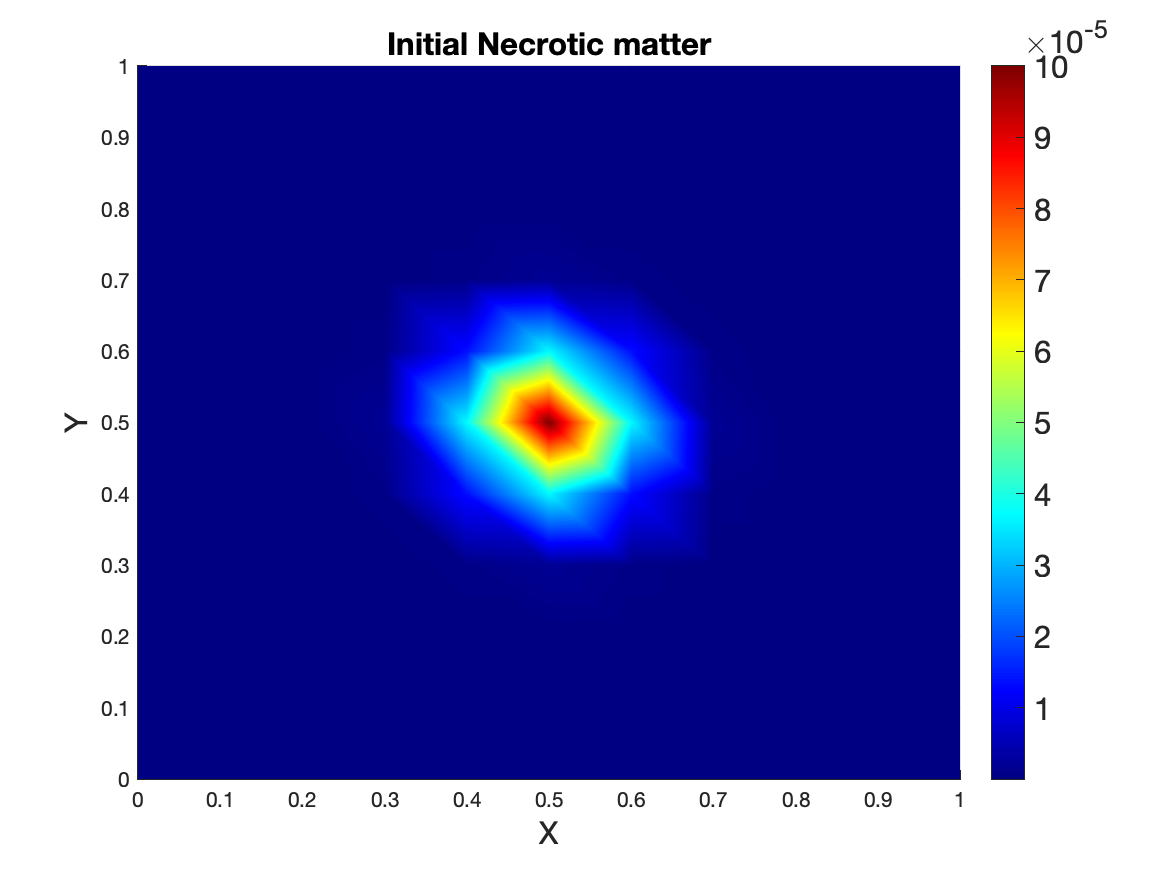}
		\caption{Necrotic matter}
	\end{subfigure} 
	\caption{Initial conditions for bacteria, mycolactone, normal tissue, and necrotic matter.}
	\label{IC-big}
\end{figure}

\noindent
We consider the non-dimensionalised initial conditions 
\begin{equation}
	\begin{aligned}
		\label{ic}
		u(0,x,y)  &= 0.95 \exp(-\frac{(x-0.5)^2+(y-0.5)^2}{.01}), \qquad \ x,y\in[0,1],
		\\
		m(0,x,y) &= 0.001 \exp(-\frac{(x-0.5)^2+(y-0.5)^2}{.01}), \qquad   x,y \in[0,1],
		\\
		v(0,x,y) & = \mathcal{U}, 
		\\
		n(0,x,y) &=0.0001\exp(-\frac{(x-0.5)^2+(y-0.5)^2}{.01}), \qquad x,y\in[0,1],
	\end{aligned}
\end{equation}
where $\mathcal{U}$ is the uniform distribution in $(0,1).$ They are illustrated in Figure \ref{IC-big}.

\noindent
We apply the no-flux boundary conditions in \eqref{*_} to the boundary of $\Omega=(0,1)\times(0,1)$.
According to the initial conditions, the space is predominantly occupied with normal tissue and bacteria are provided at the centre of the domain (assumed site of insect bite), from where they spread into the tissue. As a consequence, toxin and necrotic matter are initially rather concentrated, as well.\\[-2ex]

\begin{figure}[htbp!]
	\centering
	
	\begin{subfigure}{0.24\textwidth} 		
		\centering
		\includegraphics[width=\textwidth]{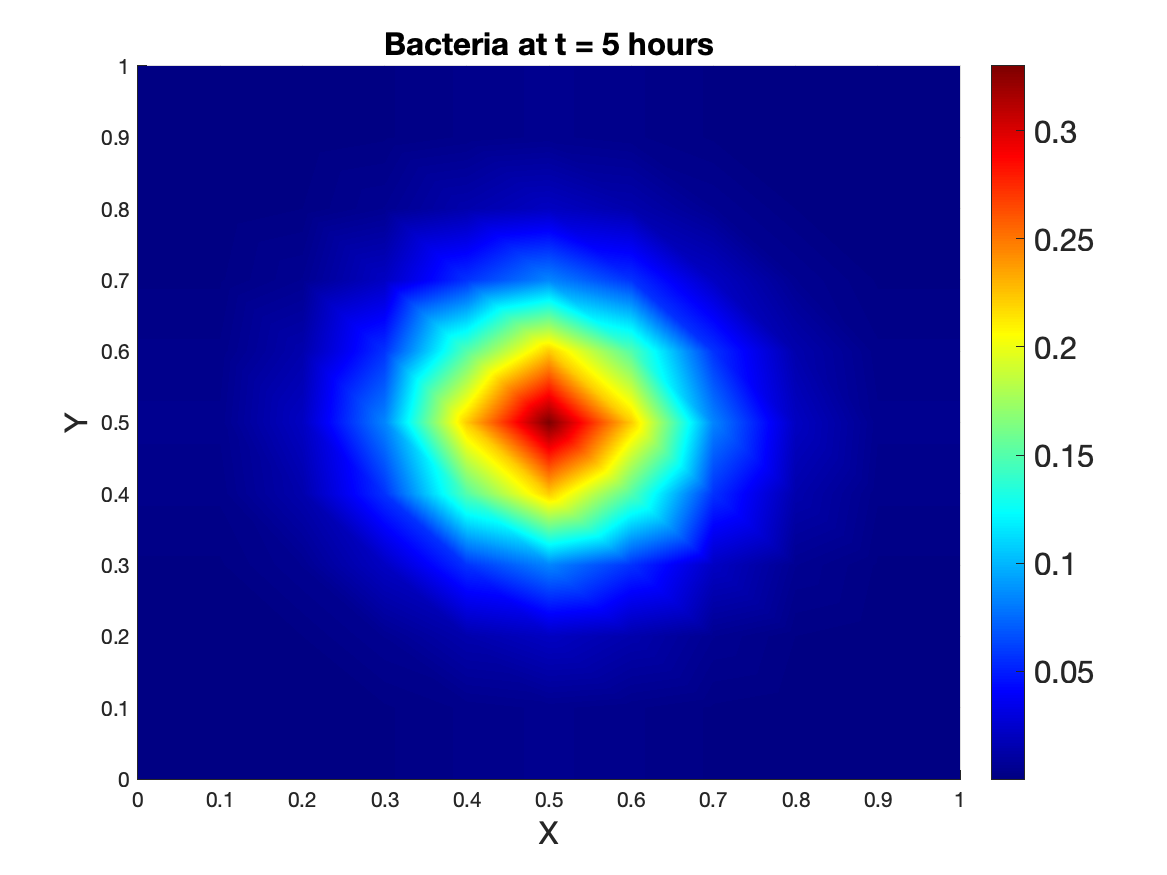}
		\caption{Bacteria at $t=5$}
		\label{}
	\end{subfigure} 
	\begin{subfigure}{0.24\textwidth} 		
		\centering
		\includegraphics[width=\textwidth]{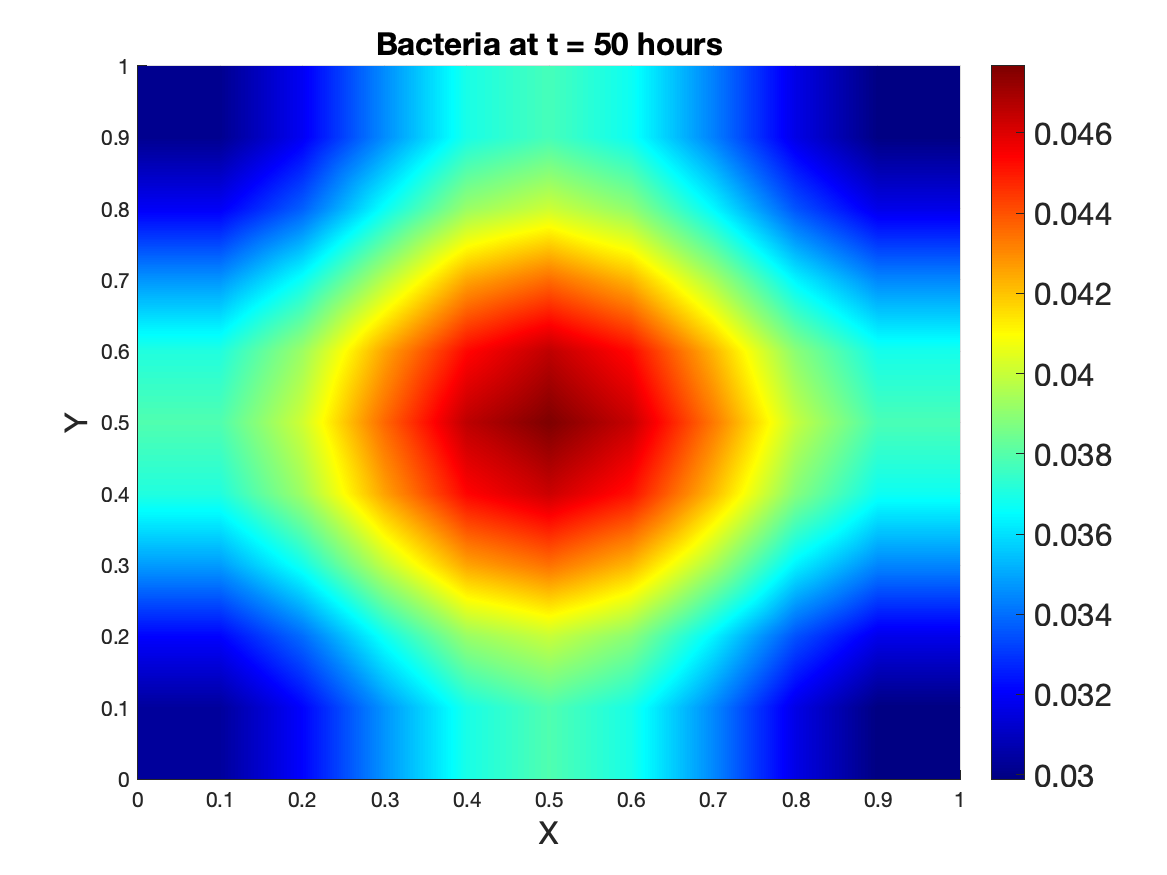}
		\caption{Bacteria at $t=50$}
	\end{subfigure}
	\begin{subfigure}{0.24\textwidth} 		
		\centering
		\includegraphics[width=\textwidth]{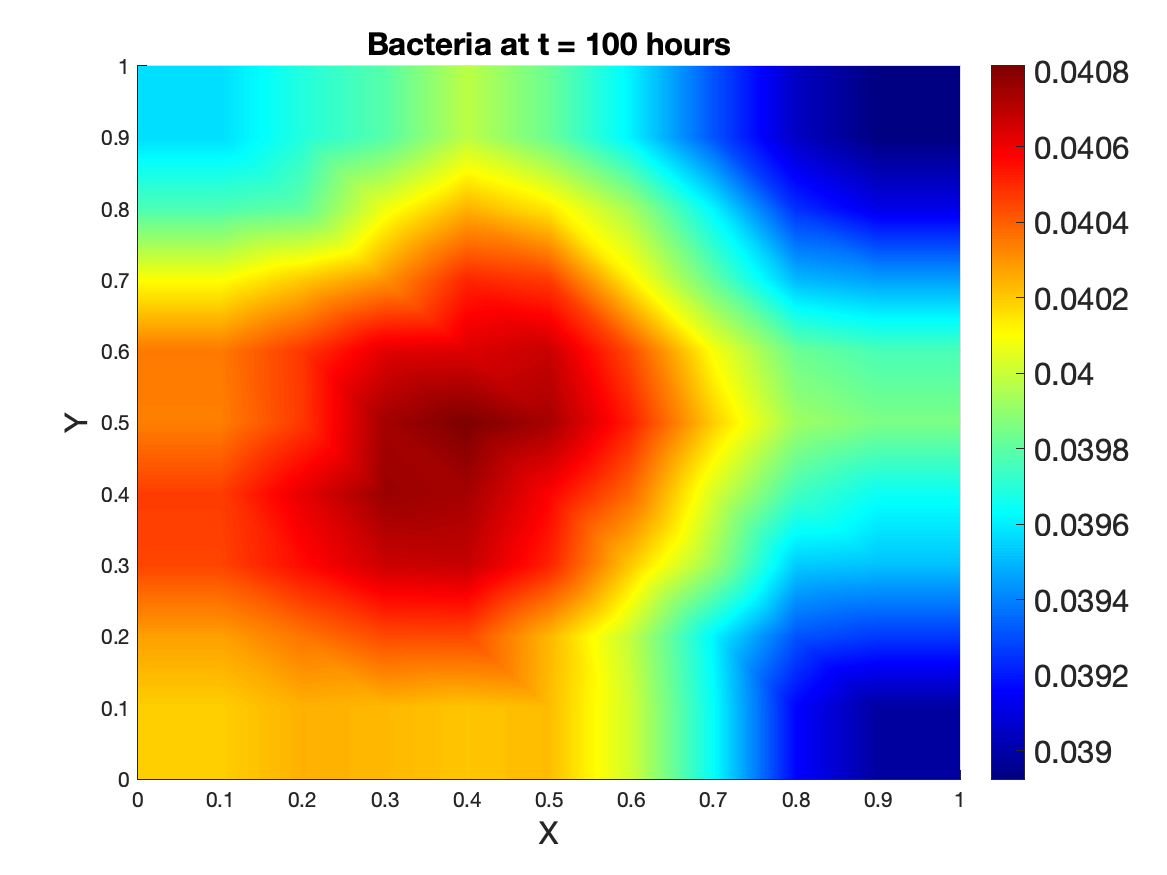}
		\caption{Bacteria at $t=100$}
	\end{subfigure}
	\begin{subfigure}{0.24\textwidth} 		\centering
		\includegraphics[width=\textwidth]{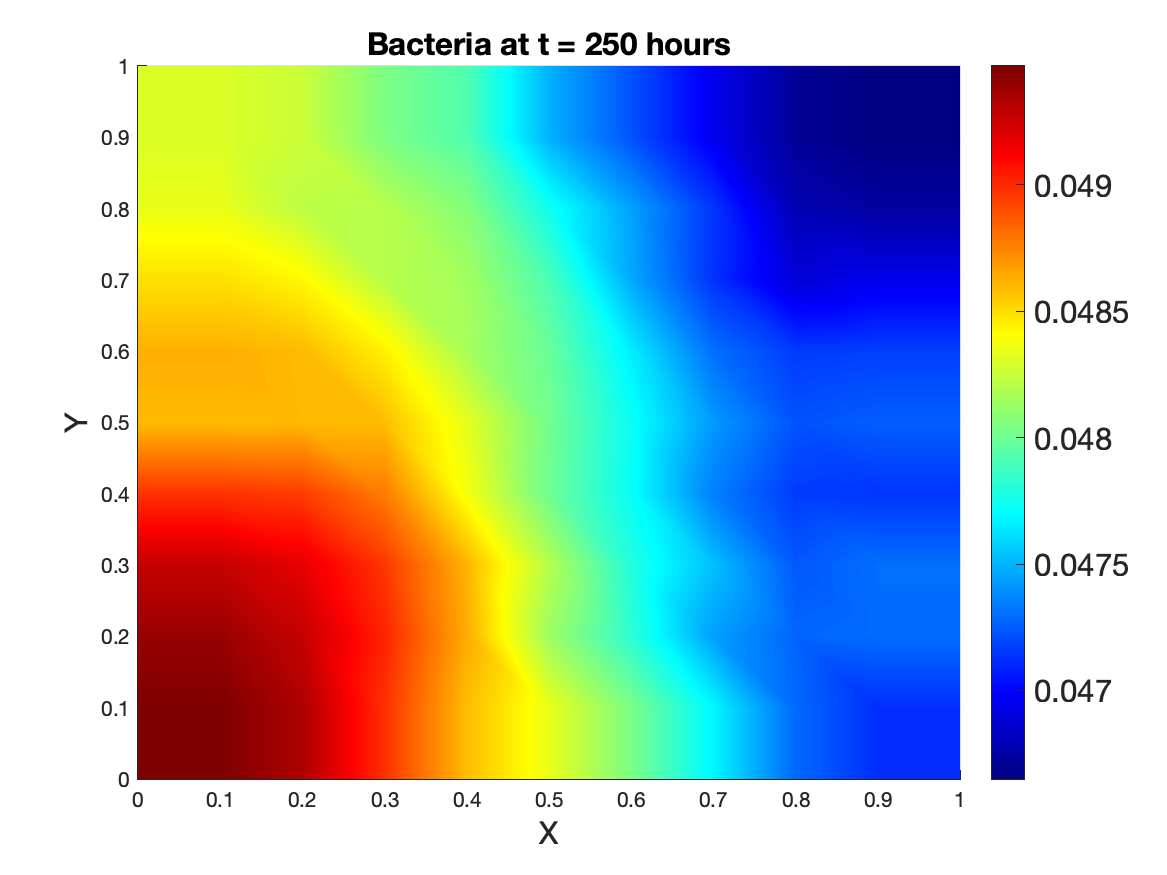}
		\caption{Bacteria at $t=250$}
	\end{subfigure}     \vspace{0.5cm}     \\
	
	\begin{subfigure}{0.24\textwidth} 		
		\centering
		\includegraphics[width=\textwidth]{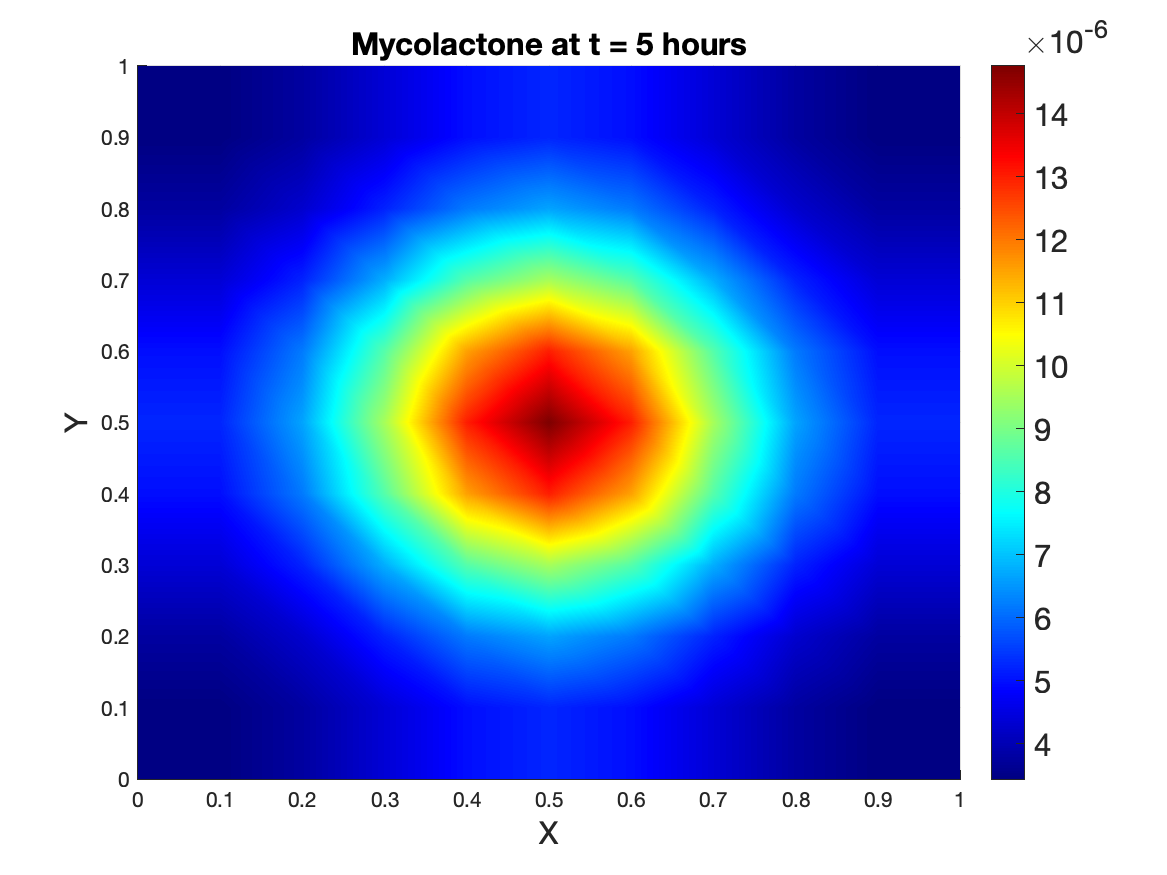}
		\caption{Mycolactone at \\ $t=5$}
		\label{}
	\end{subfigure}
	\begin{subfigure}{0.24\textwidth} 		\centering 
		\includegraphics[width=\textwidth]{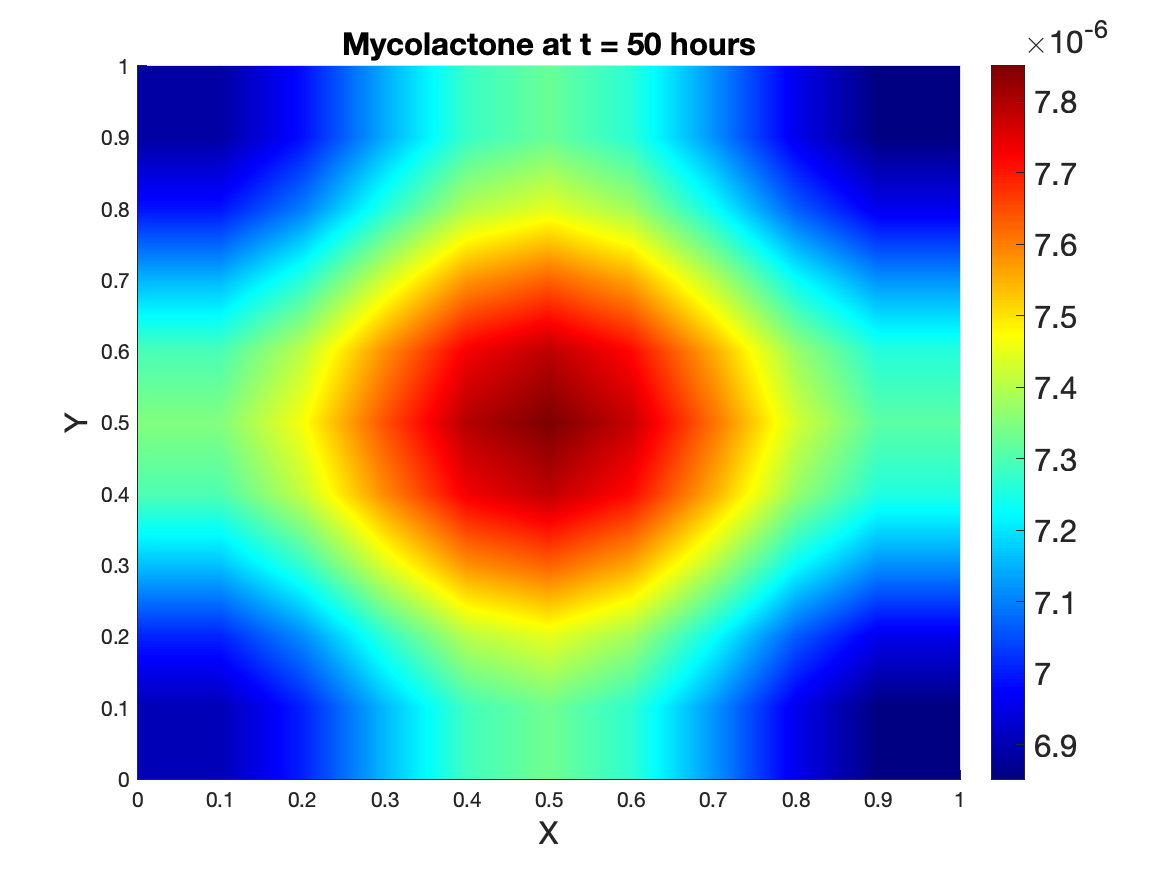}
		\caption{Mycolactone at \\ $t=50$}
	\end{subfigure}
	\begin{subfigure}{0.24\textwidth} 		\centering 
		\includegraphics[width=\textwidth]{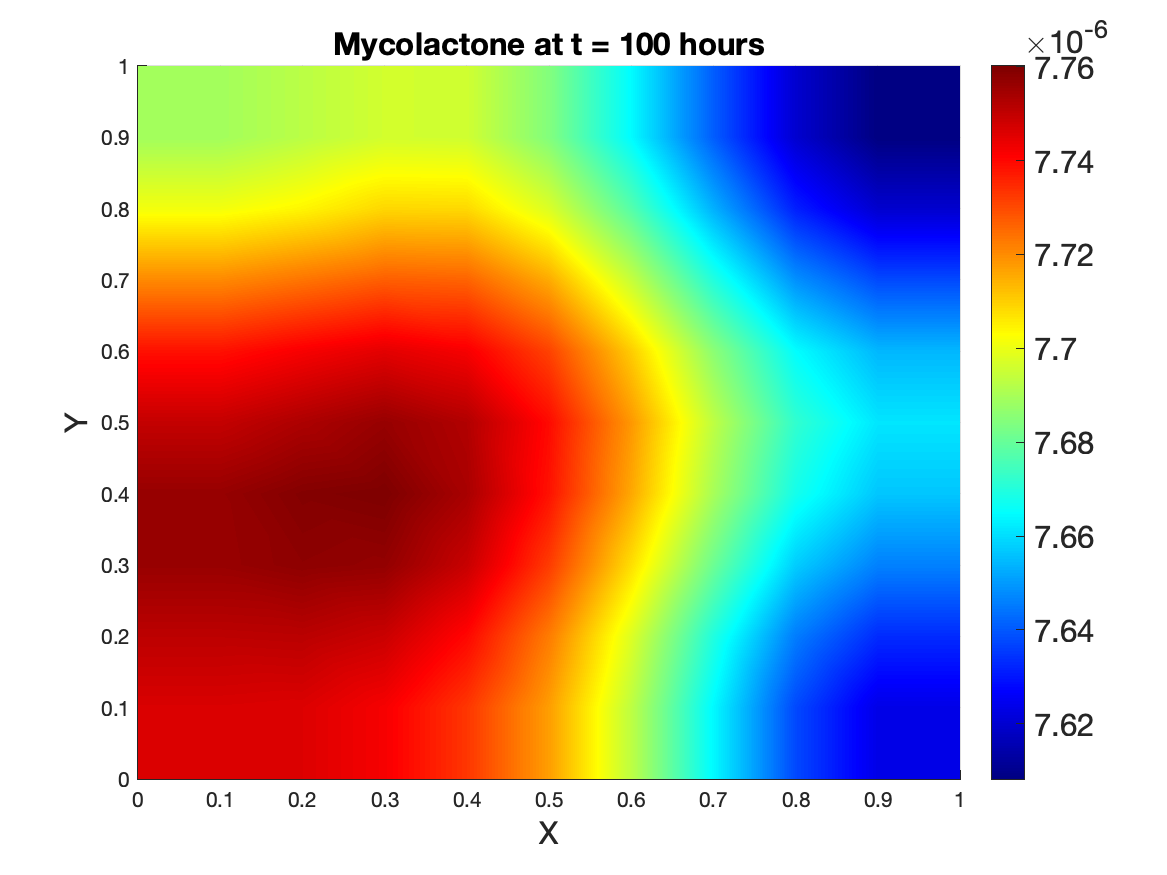}
		\caption{Mycolactone at \\ $t=100$}
	\end{subfigure}
	\begin{subfigure}{0.24\textwidth} 		
		\centering
		\includegraphics[width=\textwidth]{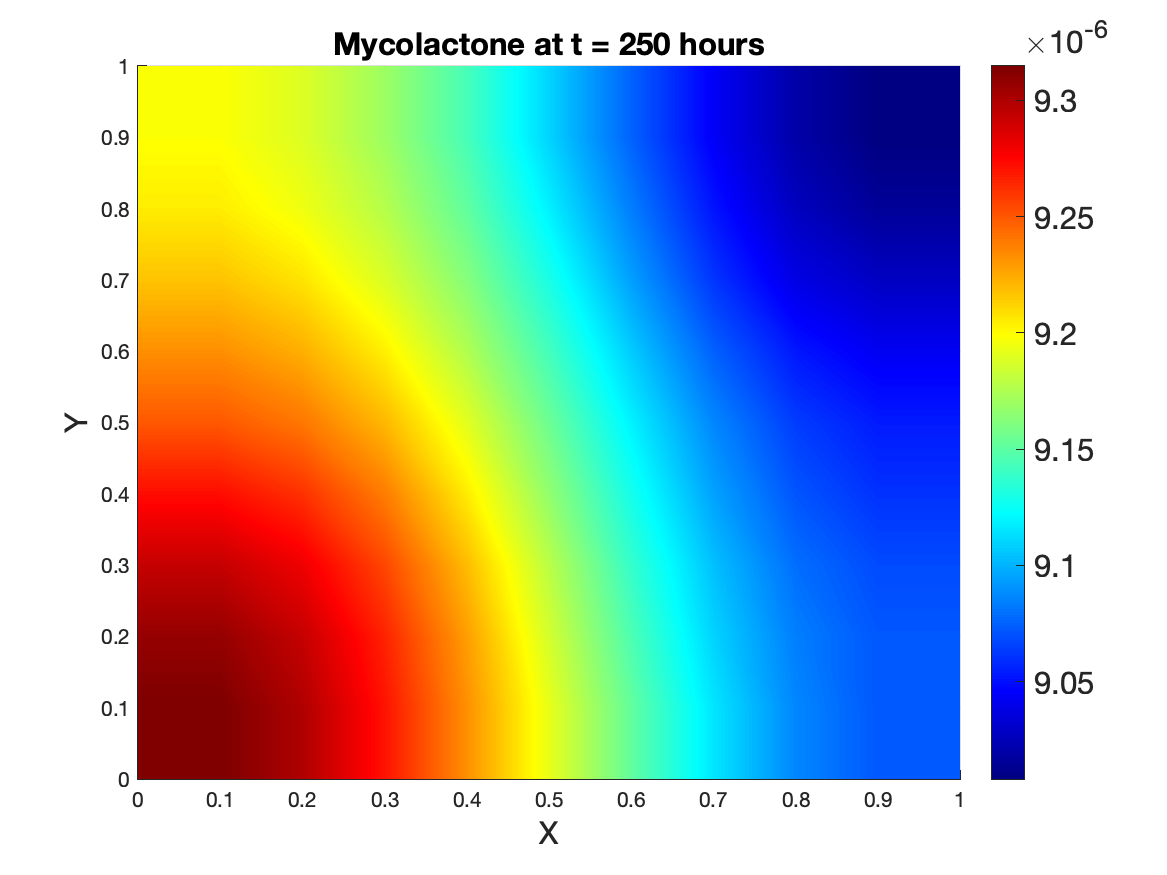}
		\caption{Mycolactone at \\ $t=250$}
	\end{subfigure}     \vspace{0.5cm}     \\
	
	\begin{subfigure}{0.24\textwidth} 		\centering
		\includegraphics[width=\textwidth]{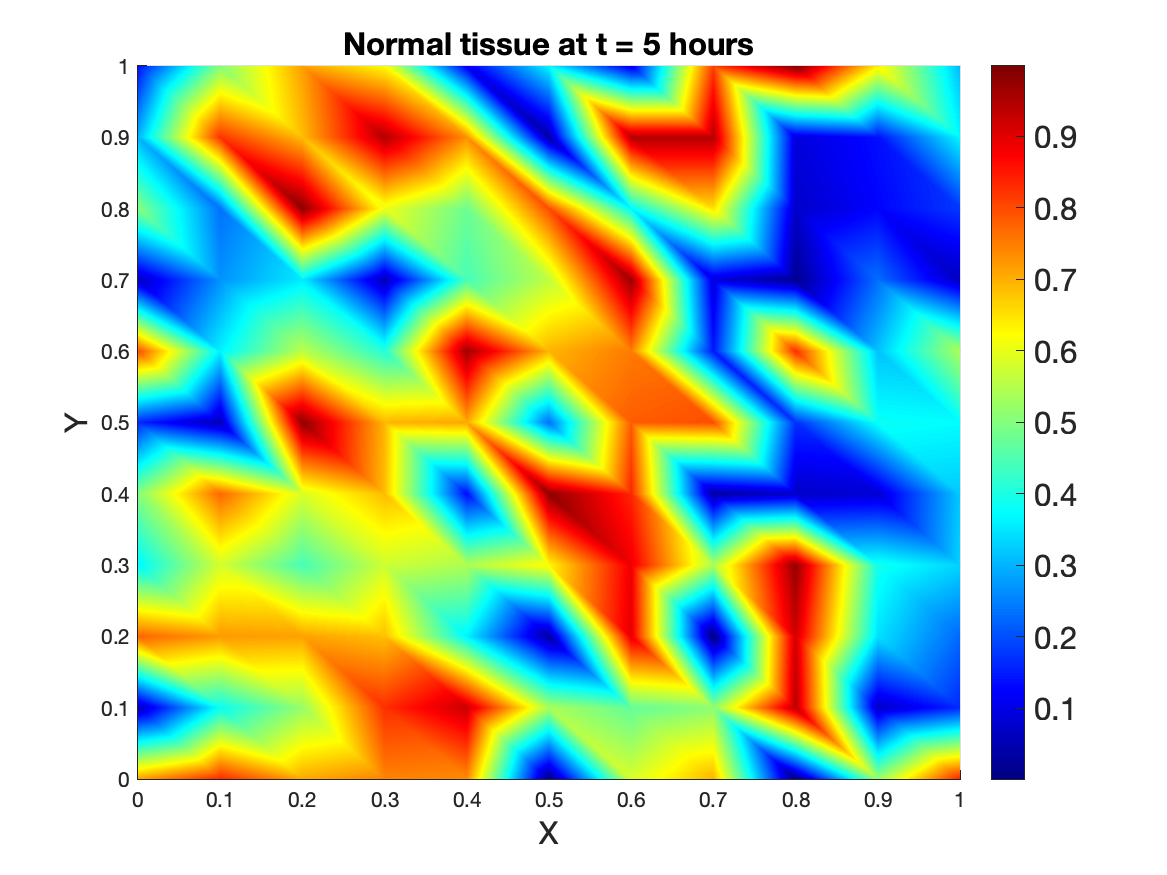}
		\caption{Normal tissue at \\ $t=5$}
		\label{}
	\end{subfigure}	
	\begin{subfigure}{0.24\textwidth} 		
		\centering
		\includegraphics[width=\textwidth]{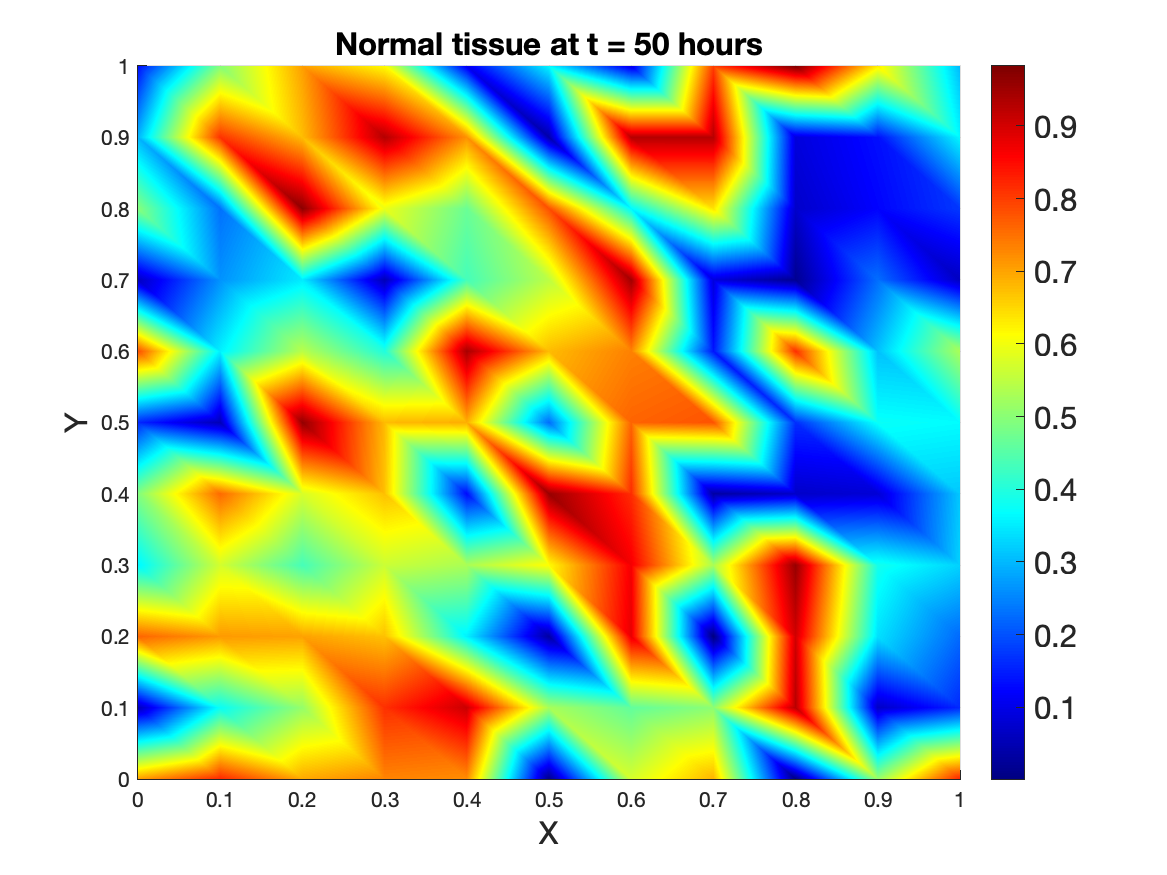}
		\caption{Normal tissue at \\ $t=50$}
	\end{subfigure}	
	\begin{subfigure}{0.24\textwidth} 		
		\centering
		\includegraphics[width=\textwidth]{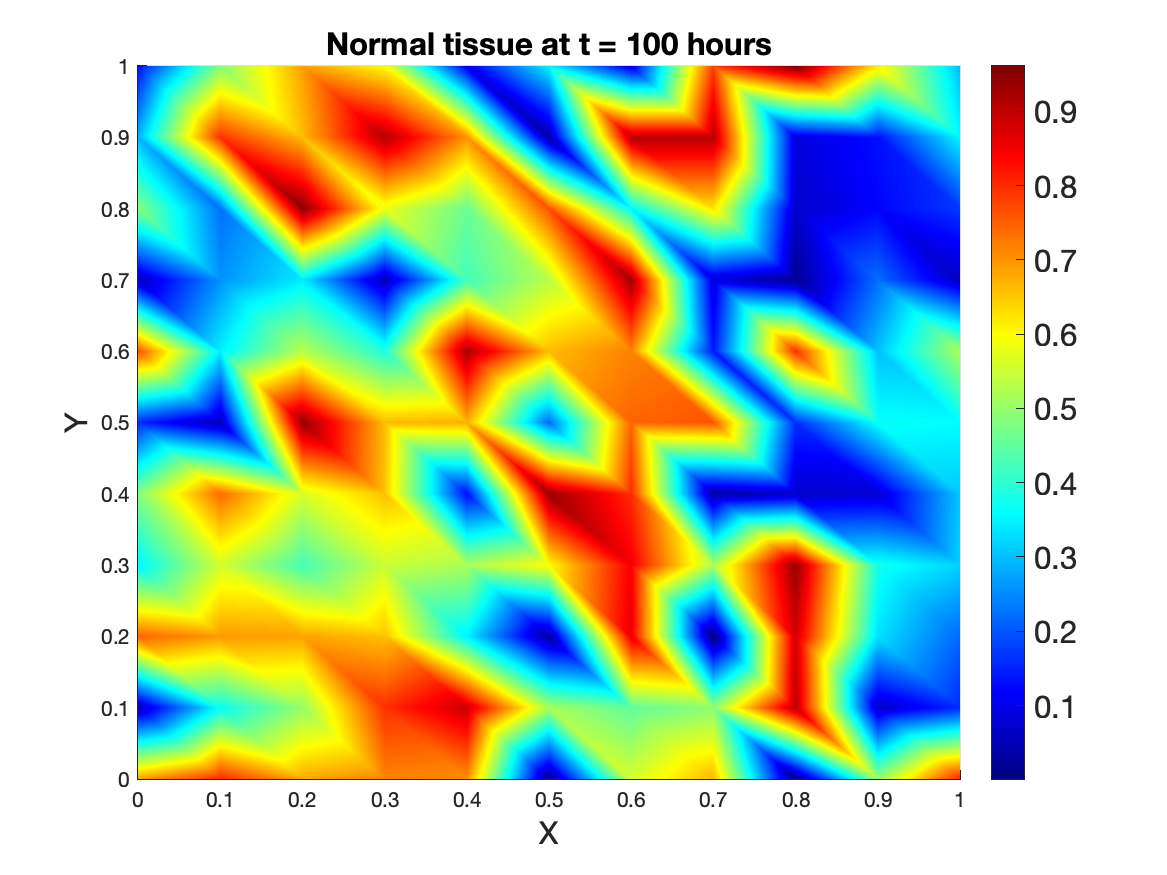}
		\caption{Normal tissue at \\ $t=100$}
	\end{subfigure}	
	\begin{subfigure}{0.24\textwidth} 	
		\centering
		\includegraphics[width=\textwidth]{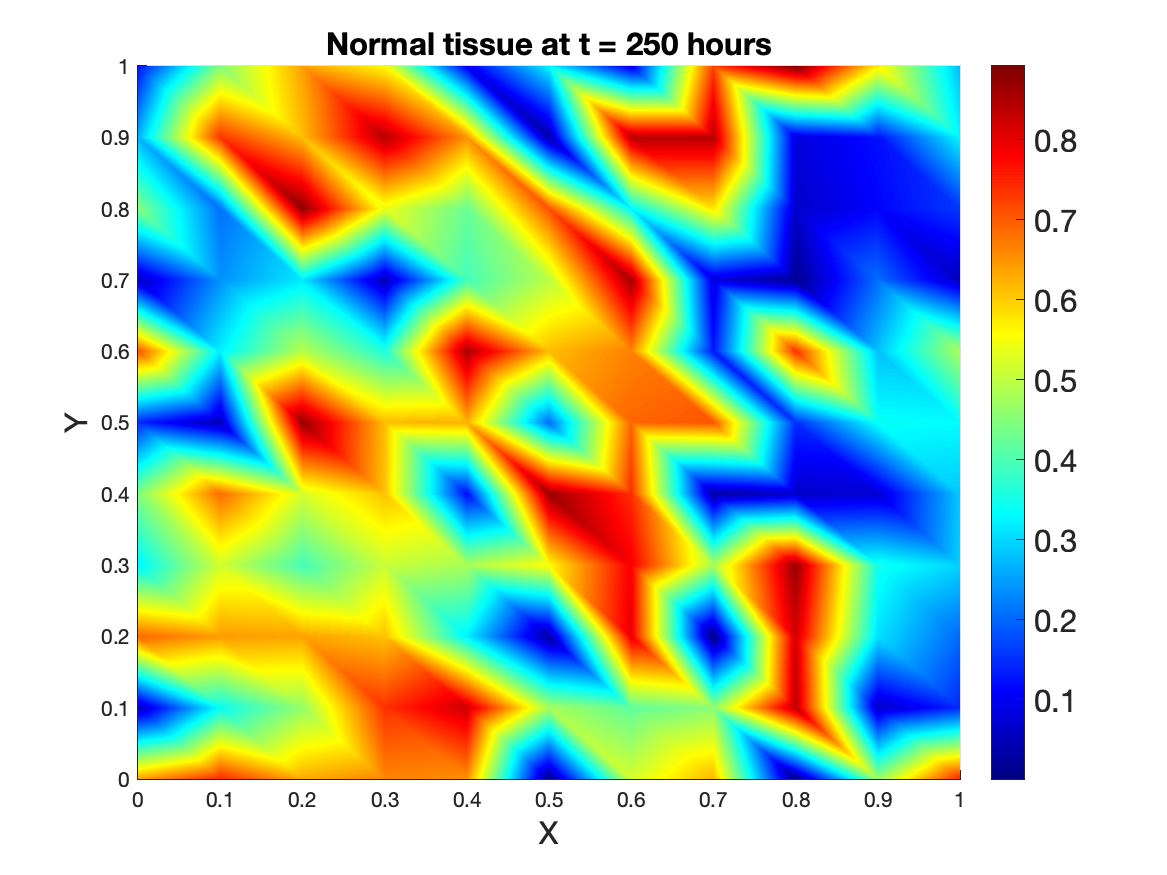}
		\caption{Normal tissue at \\ $t=250$}
	\end{subfigure}
	\vspace{0.5cm} \\
	\begin{subfigure}{0.24\textwidth} 		
		\centering
		\includegraphics[width=\textwidth]{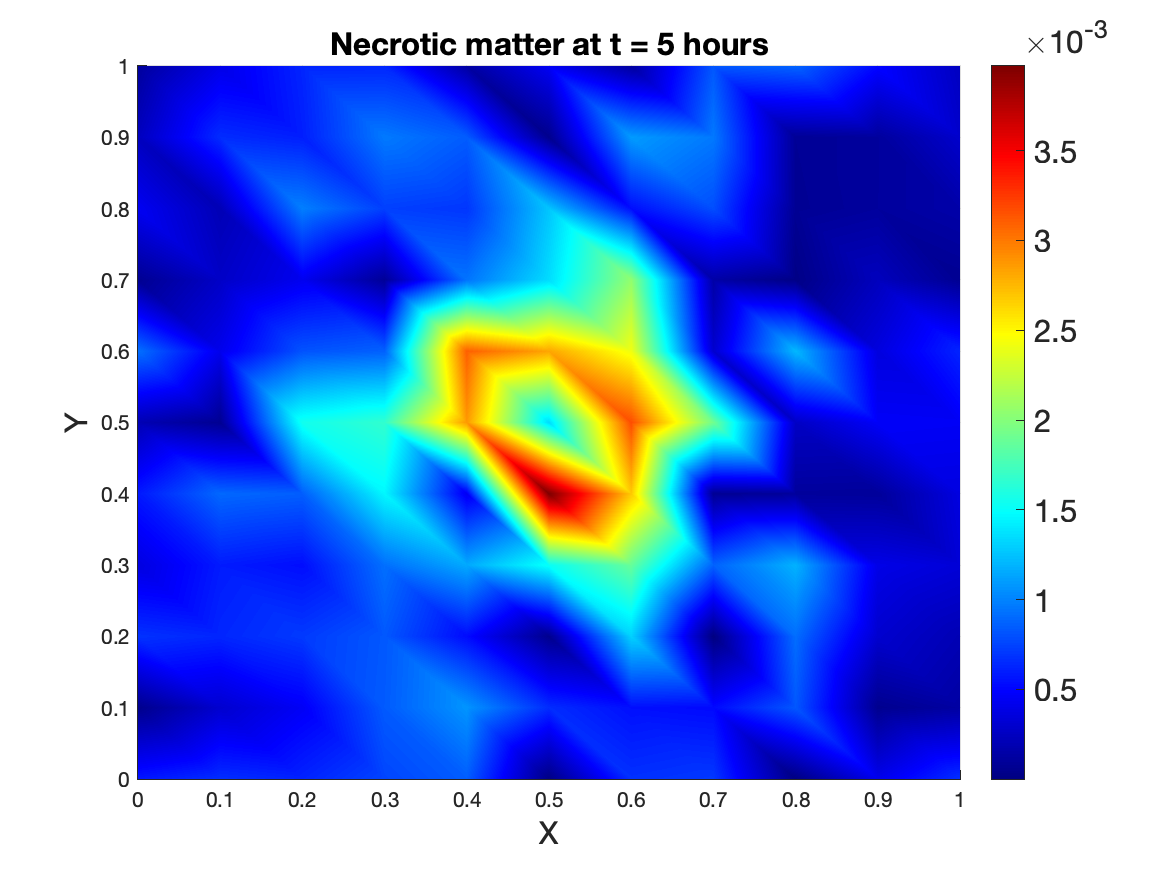}
		\caption{Necrotic matter \\ at $t=5$}
		\label{}
	\end{subfigure}	 	
	\begin{subfigure}{0.24\textwidth} 		
		\centering
		\includegraphics[width=\textwidth]{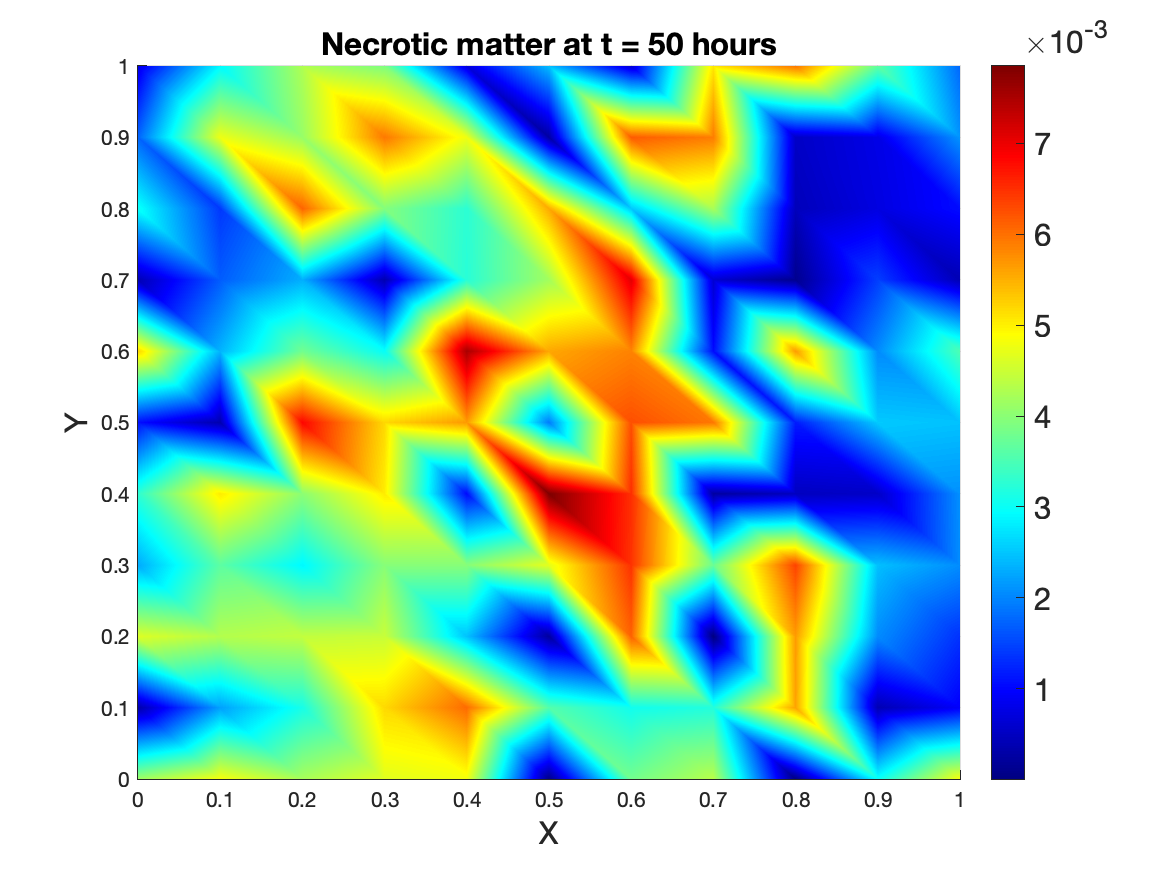}
		\caption{Necrotic matter at \\ $t=50$}
	\end{subfigure}	
	\begin{subfigure}{0.24\textwidth} 		
		\centering
		\includegraphics[width=\textwidth]{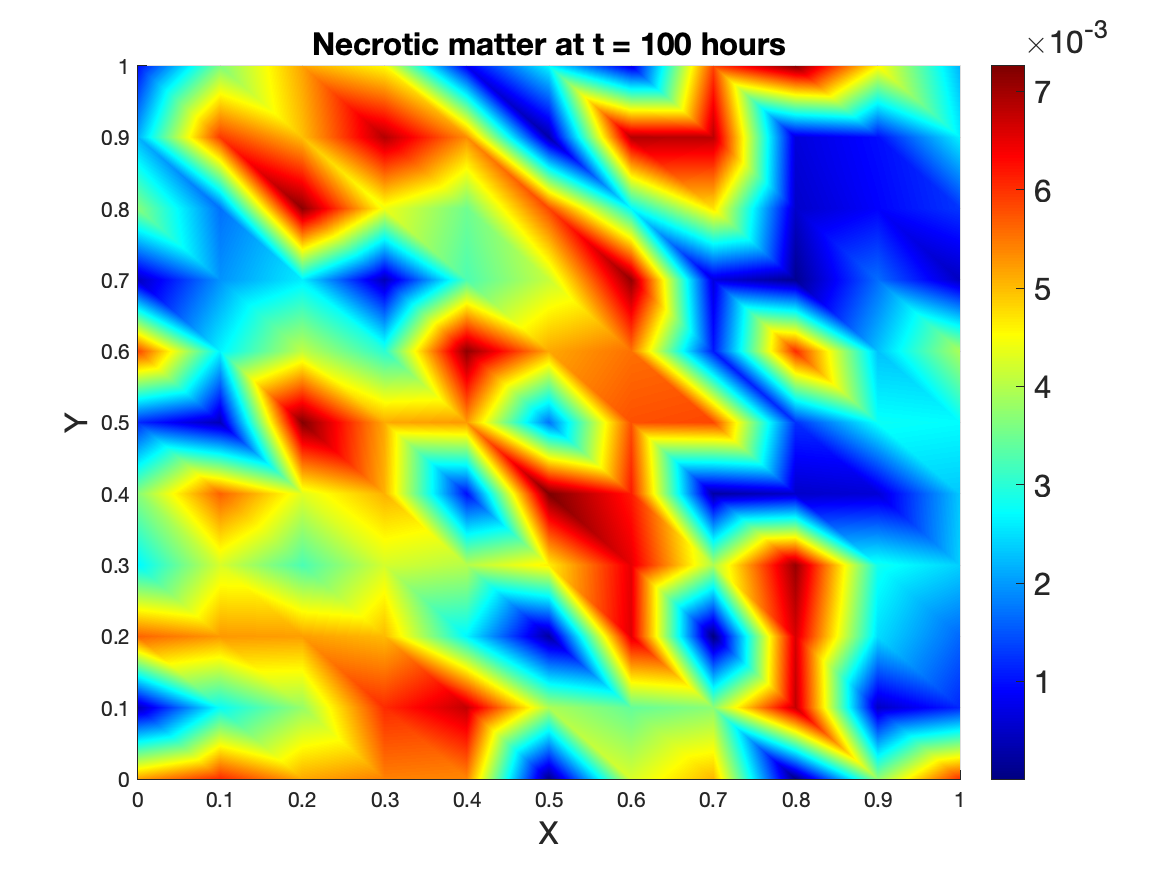}
		\caption{Necrotic matter at \\ $t=100$}
	\end{subfigure}	
	\begin{subfigure}{0.24\textwidth} 		\centering
		\includegraphics[width=\textwidth]{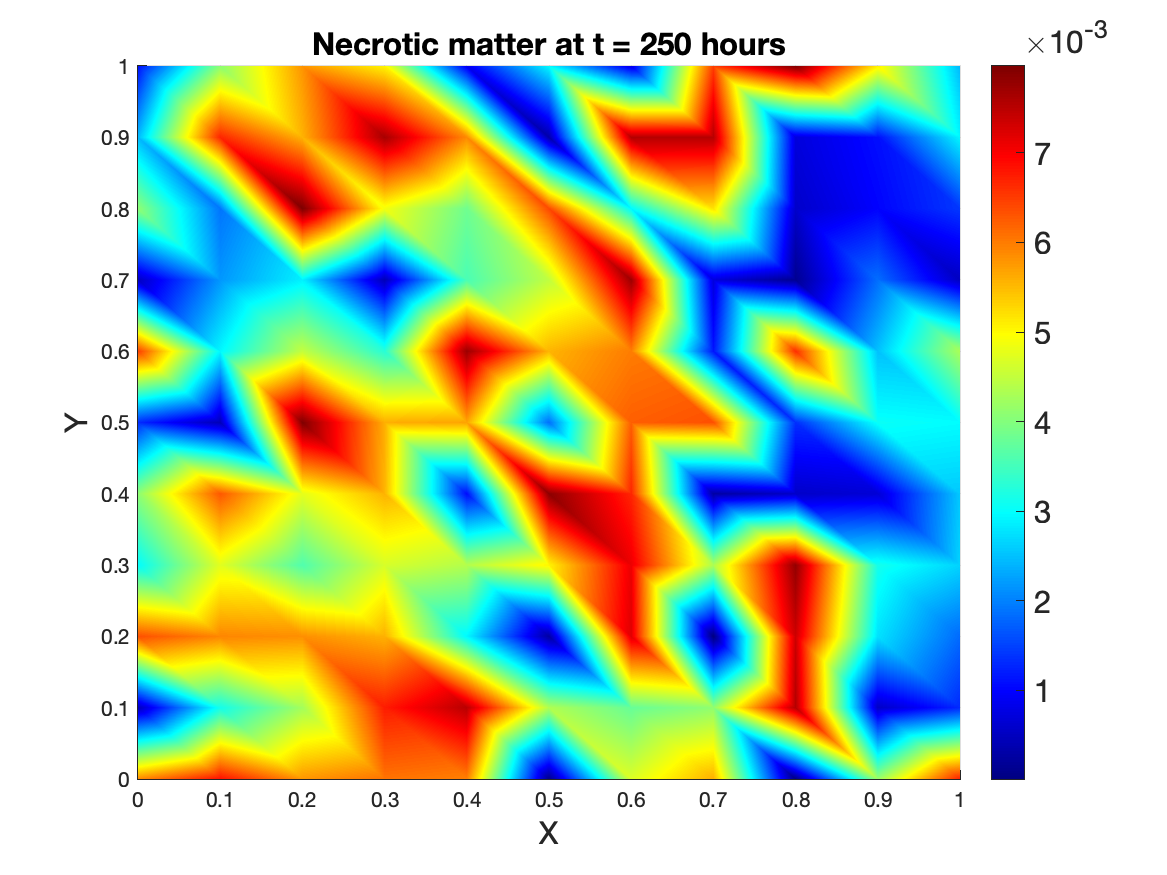}
		\caption{Necrotic matter at \\ $t=250$}
	\end{subfigure}     \vspace{0.5cm}     \\
	\caption{Scenario 1: Bacteria, mycolactone, normal tissue, and necrotic matter at different times for $\gamma_1 = \gamma_2.$}
	\label{fig:gamma1=gamma2}
\end{figure}

\noindent In Figure \ref{fig:gamma1=gamma2} we observe the dynamics of bacteria, mycolactone, normal tissue, and necrotic matter for Scenario 1 when $\gamma_1 = \gamma_2$. While the progression patterns are similar across all cases, they exhibit varying densities. The simulation shows a biphasic progression of the disease. Initially, over approximately 100 hours (about 4 days), there is a decrease in bacteria, mycolactone, and necrotic tissue. This early decline likely stems from a combination of factors: the initial lack of necrotic tissue to support bacterial growth, the host's immune response, and localized depletion of nutrients. Crucially, the absence of sufficient necrotic matter in early stages limits bacterial proliferation, as bacteria thrives in necrotic environments. However, after this period, all three components increase, marking a turning point in infection progression. This second phase probably represents the establishment of a favorable environment for bacteria: accumulated necrotic tissue enhances bacterial growth, increased mycolactone levels suppress the local immune response, and bacteria potentially spread to unaffected areas with fresh resources. A positive feedback loop ensues, where bacterial growth leads to more mycolactone production, causing more tissue necrosis, which further promotes bacterial proliferation. This pattern underscores the critical role of necrotic tissue in disease progression - initially as a growth-limiting factor, later as a growth-promoter. It also suggests a crucial early window for therapeutic intervention, before the infection becomes self-sustaining. The model's insights highlight potential treatment strategies, including preventing tissue necrosis, inhibiting mycolactone production, and emphasizing early-stage treatments to disrupt the establishment of an expanding infection.\\[-2ex]

\begin{figure}[htbp!]
	
	\begin{subfigure}{0.24\textwidth} 		\centering
		\includegraphics[width=\textwidth]{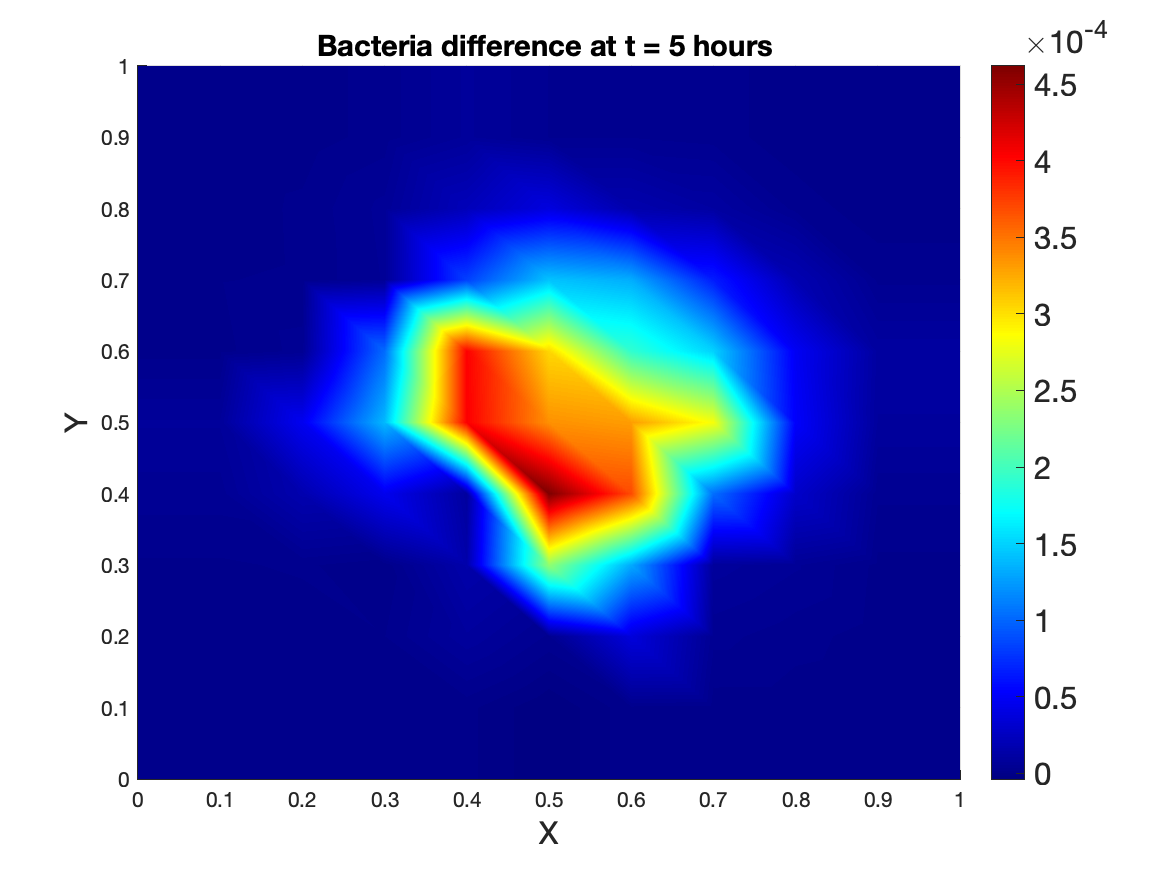}
		\caption{Bacteria at $t=5$}
		\label{}
		
	\end{subfigure} 
	\begin{subfigure}{0.24\textwidth} 		\centering
		\includegraphics[width=\textwidth]{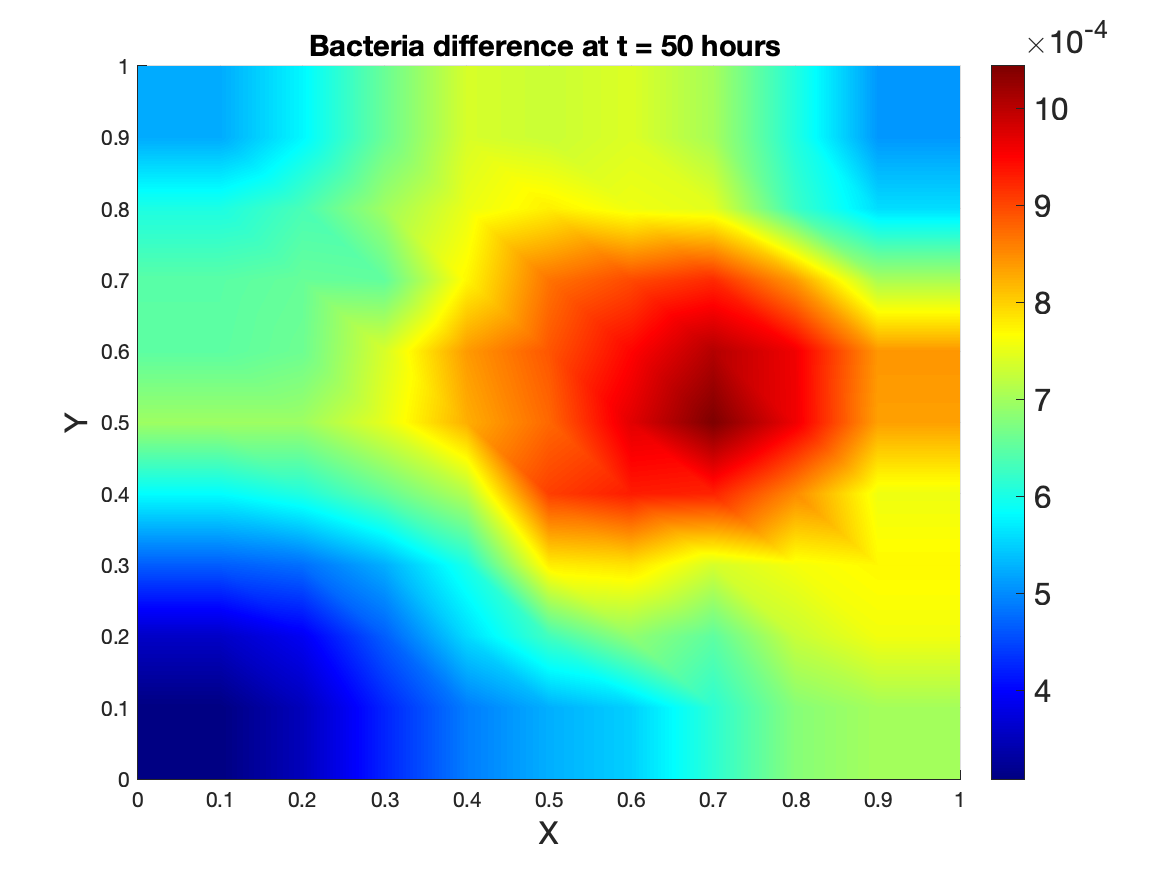}
		\caption{Bacteria at $t=50$}
		
	\end{subfigure} 
	\begin{subfigure}{0.24\textwidth} 		\centering
		\includegraphics[width=\textwidth]{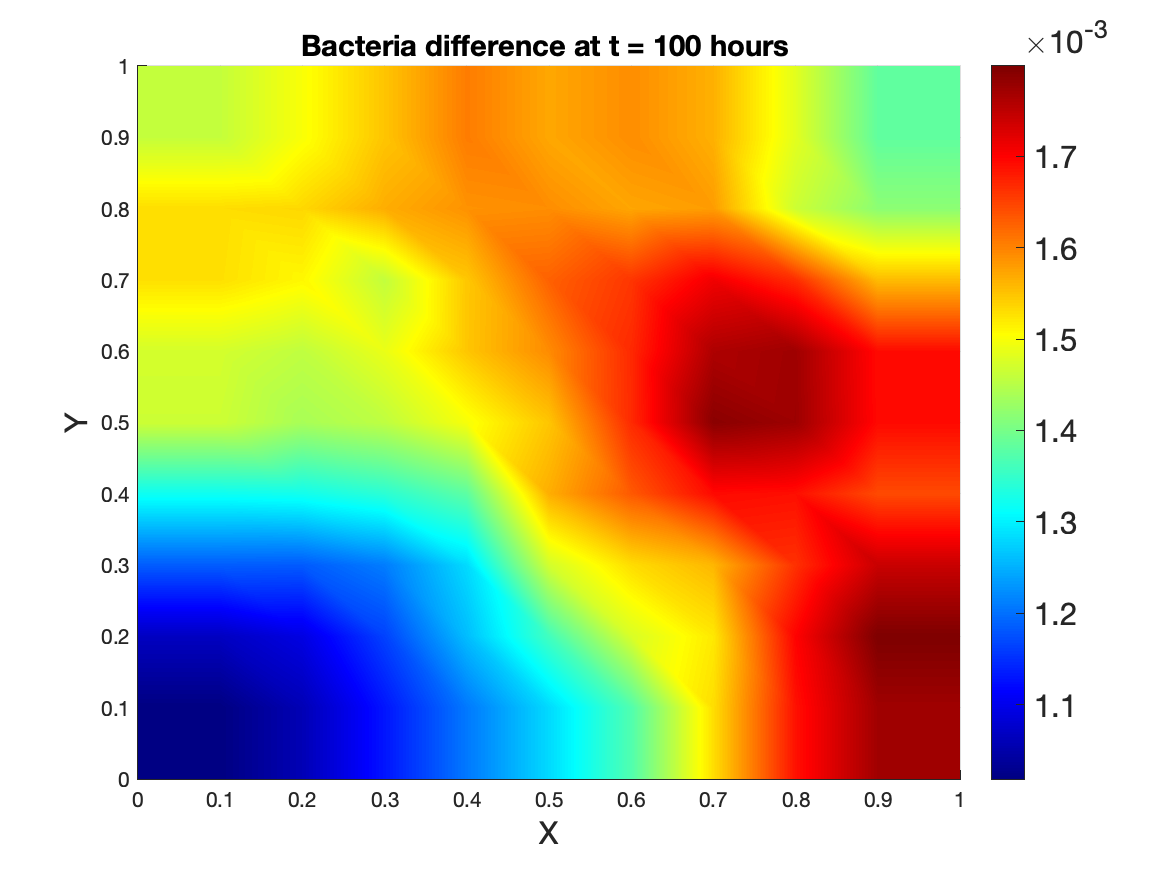}
		\caption{Bacteria at $t=100$}
		
	\end{subfigure} 
	\begin{subfigure}{0.24\textwidth} 		\centering
		\includegraphics[width=\textwidth]{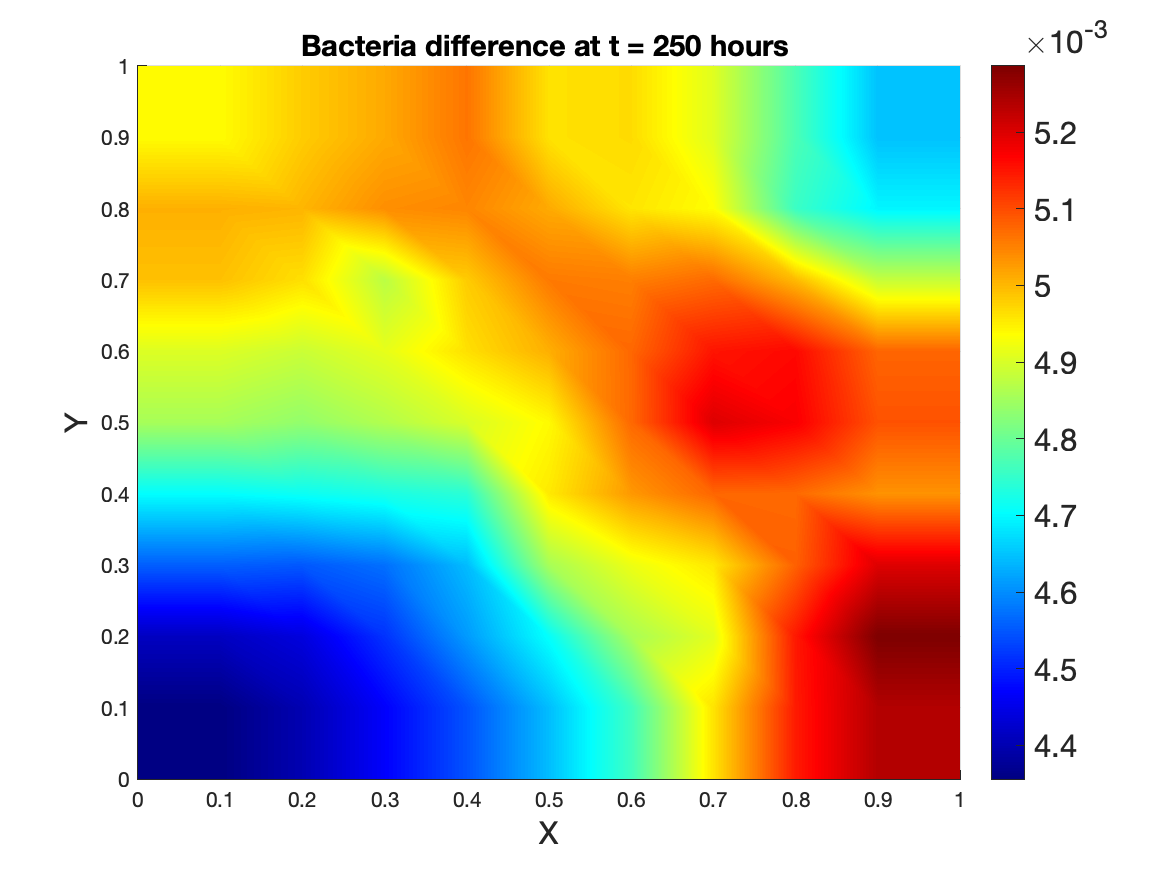}
		\caption{Bacteria at $t=250$}
		
	\end{subfigure}     \vspace{0.5cm}     \\
	
	\begin{subfigure}{0.24\textwidth} 		\centering
		\includegraphics[width=\textwidth]{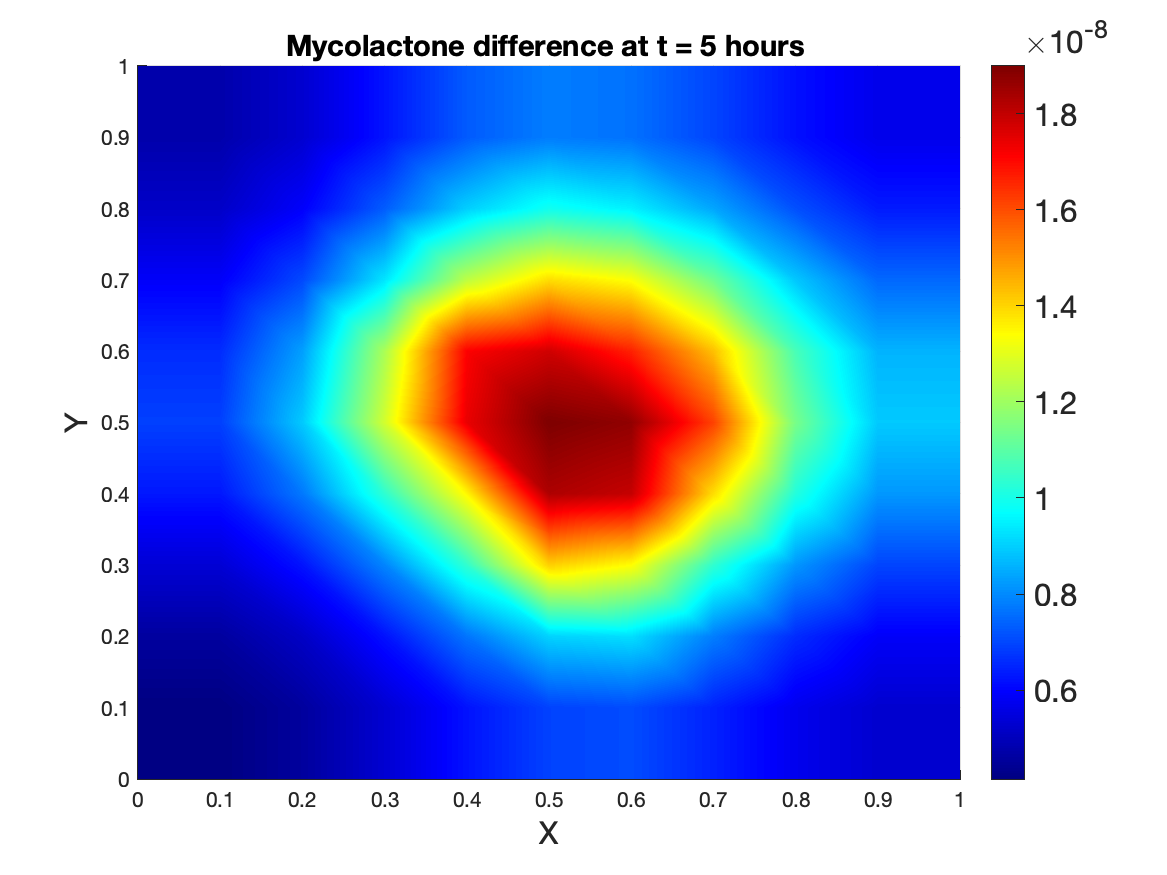}
		\caption{Mycolactone at \\ $t=5$}
		\label{}
	\end{subfigure} 
	\begin{subfigure}{0.24\textwidth} 		\centering 
		\includegraphics[width=\textwidth]{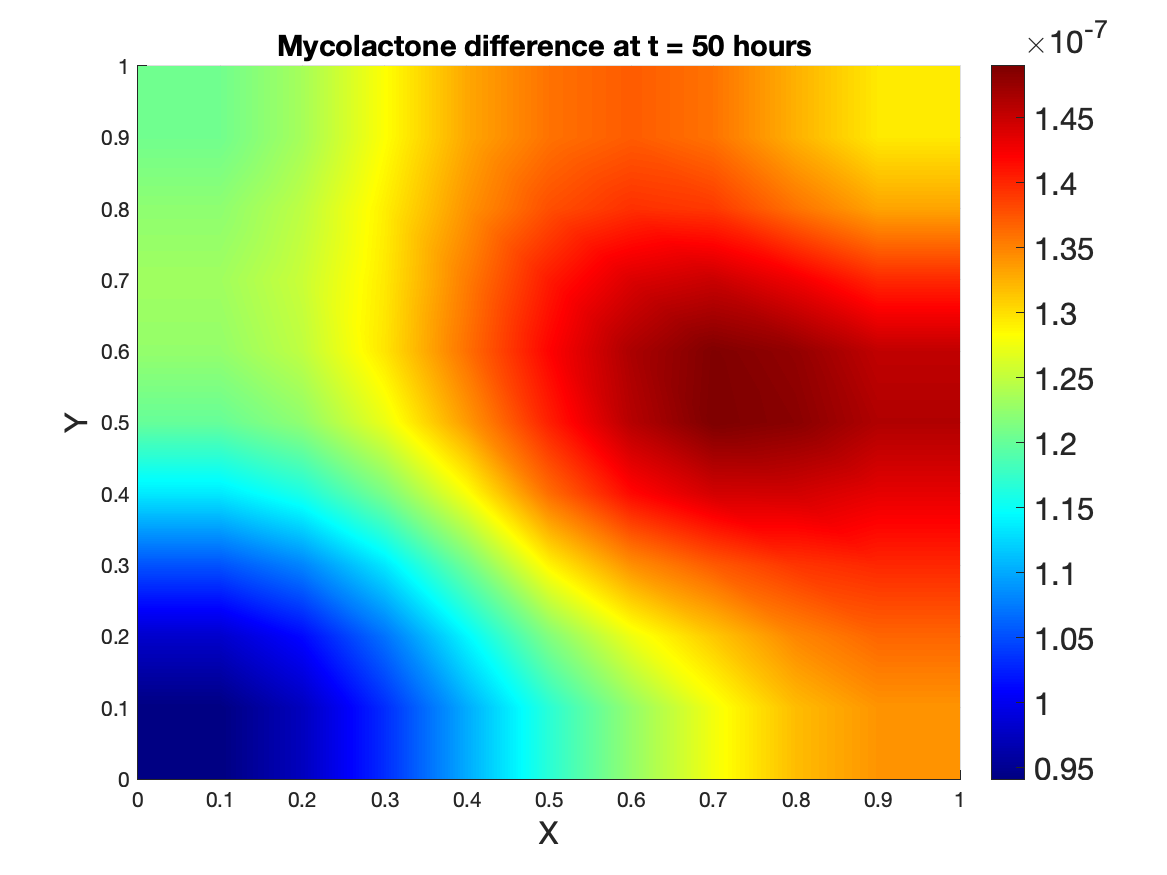}
		\caption{Mycolactone at \\ $t=50$}
		
	\end{subfigure} 
	\begin{subfigure}{0.24\textwidth} 		\centering 
		\includegraphics[width=\textwidth]{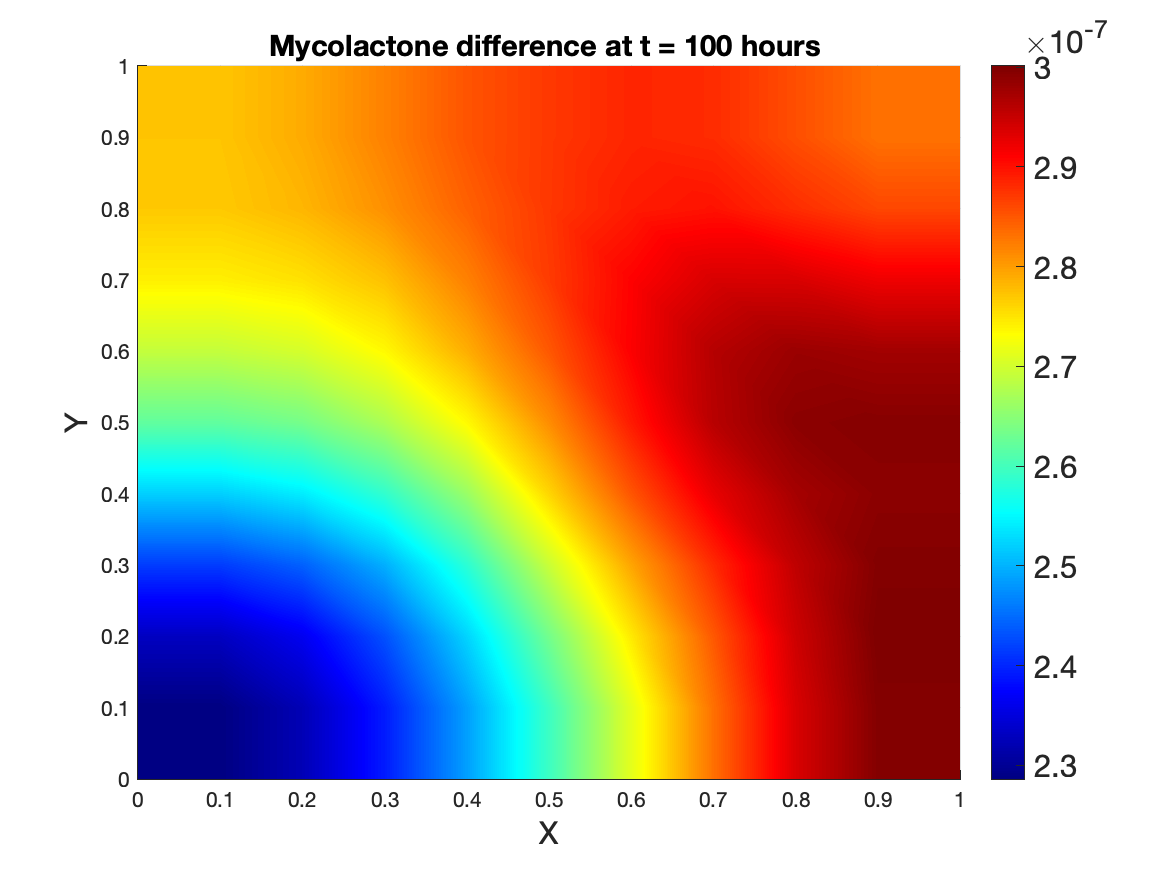}
		\caption{Mycolactone at \\ $t=100$}
		
	\end{subfigure}
	\begin{subfigure}{0.24\textwidth} 		\centering
		\includegraphics[width=\textwidth]{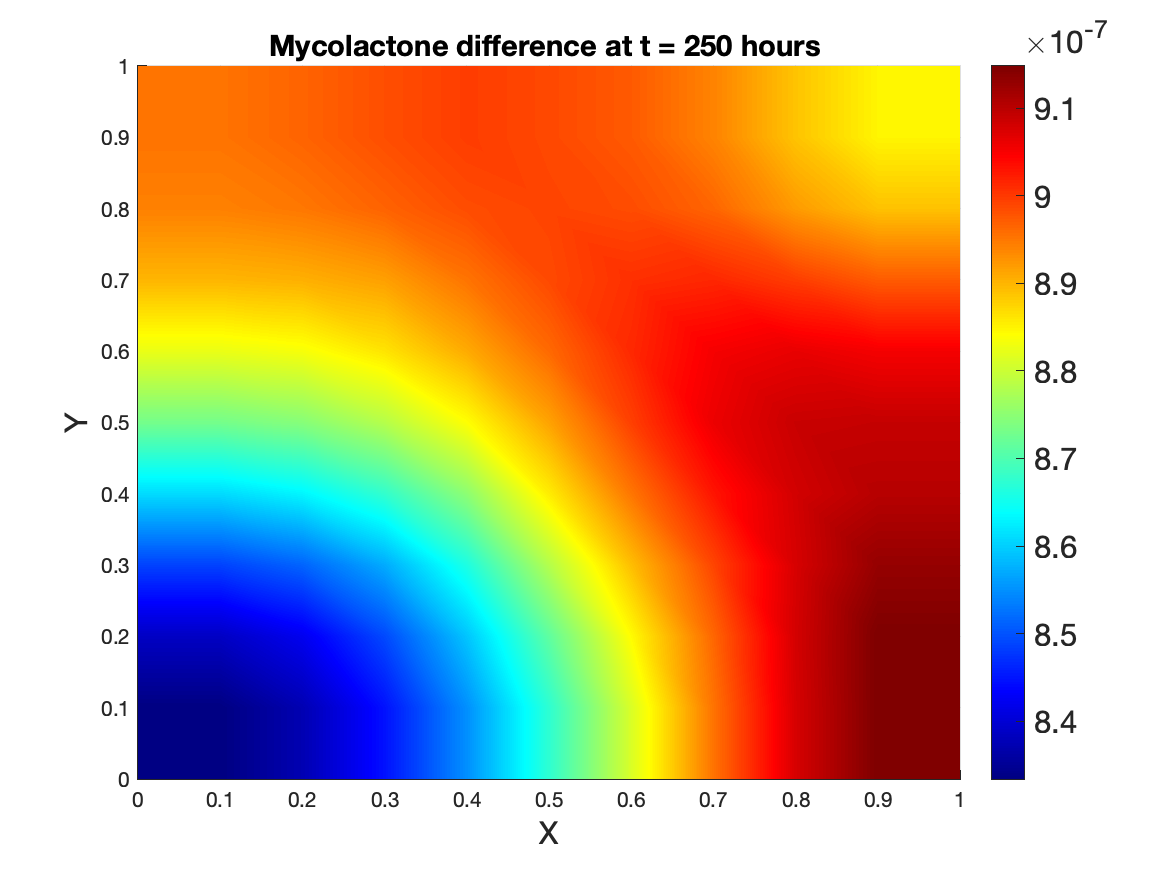}
		\caption{Mycolactone at \\ $t=250$}
	\end{subfigure}     \vspace{0.5cm}     \\
	
	\begin{subfigure}{0.24\textwidth} 		\centering
		\includegraphics[width=\textwidth]{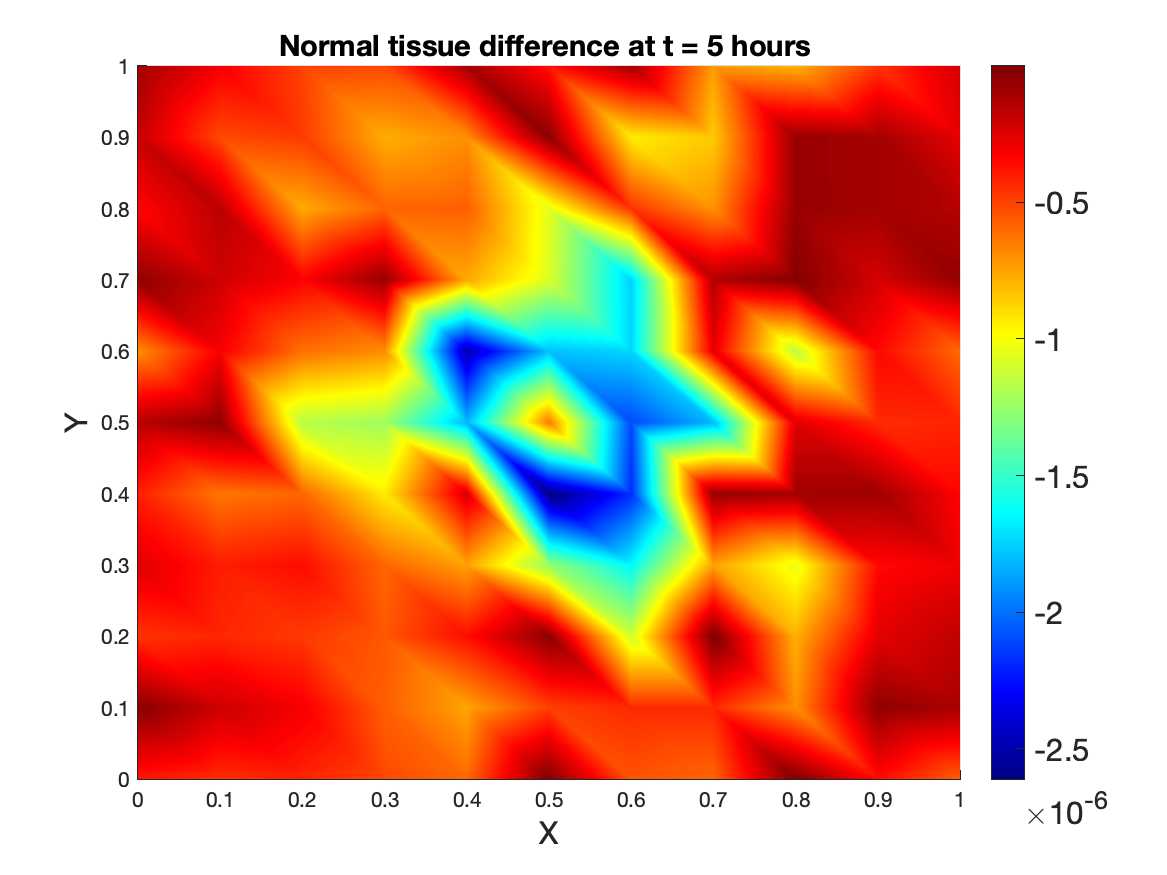}
		\caption{Normal tissue at \\ $t=5$}
		\label{}
	\end{subfigure} 
	\begin{subfigure}{0.24\textwidth} 		\centering
		\includegraphics[width=\textwidth]{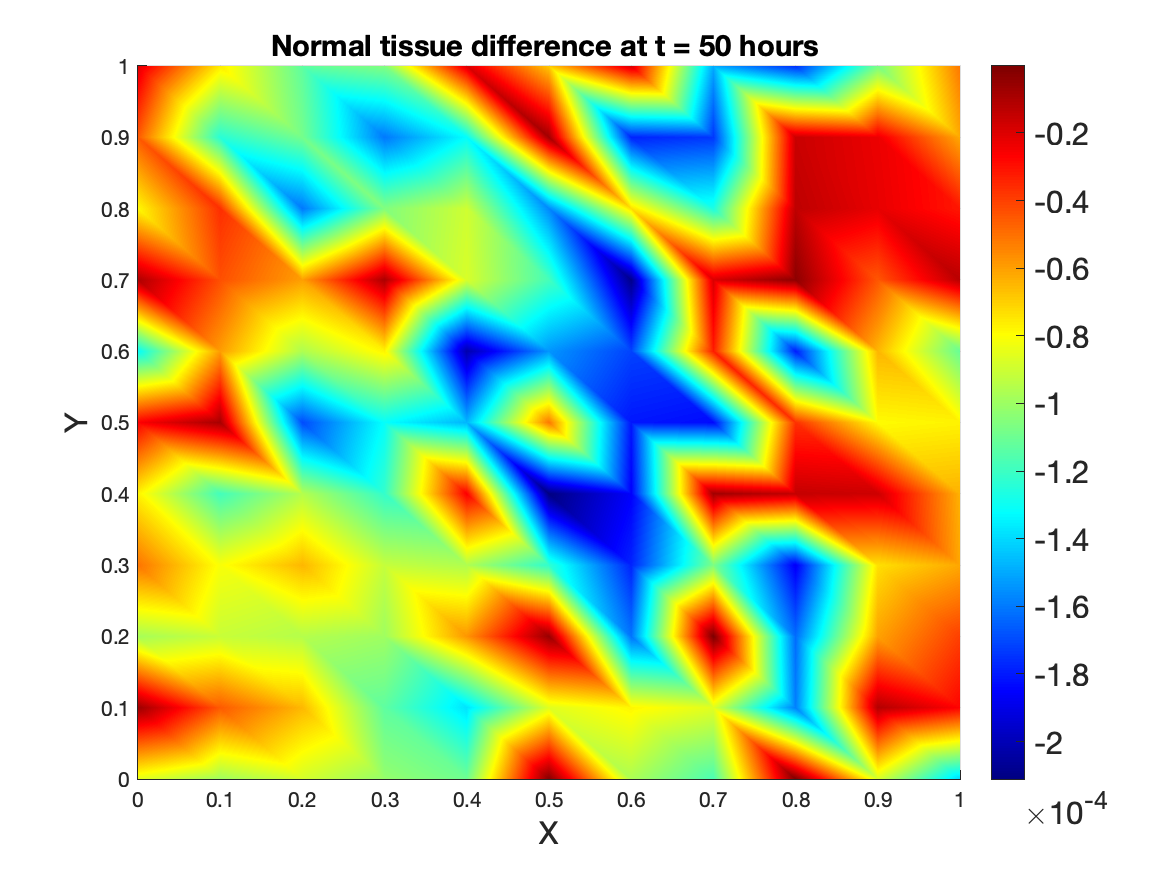}
		\caption{Normal tissue at \\ $t=50$}
	\end{subfigure}
	\begin{subfigure}{0.24\textwidth} 		\centering
		\includegraphics[width=\textwidth]{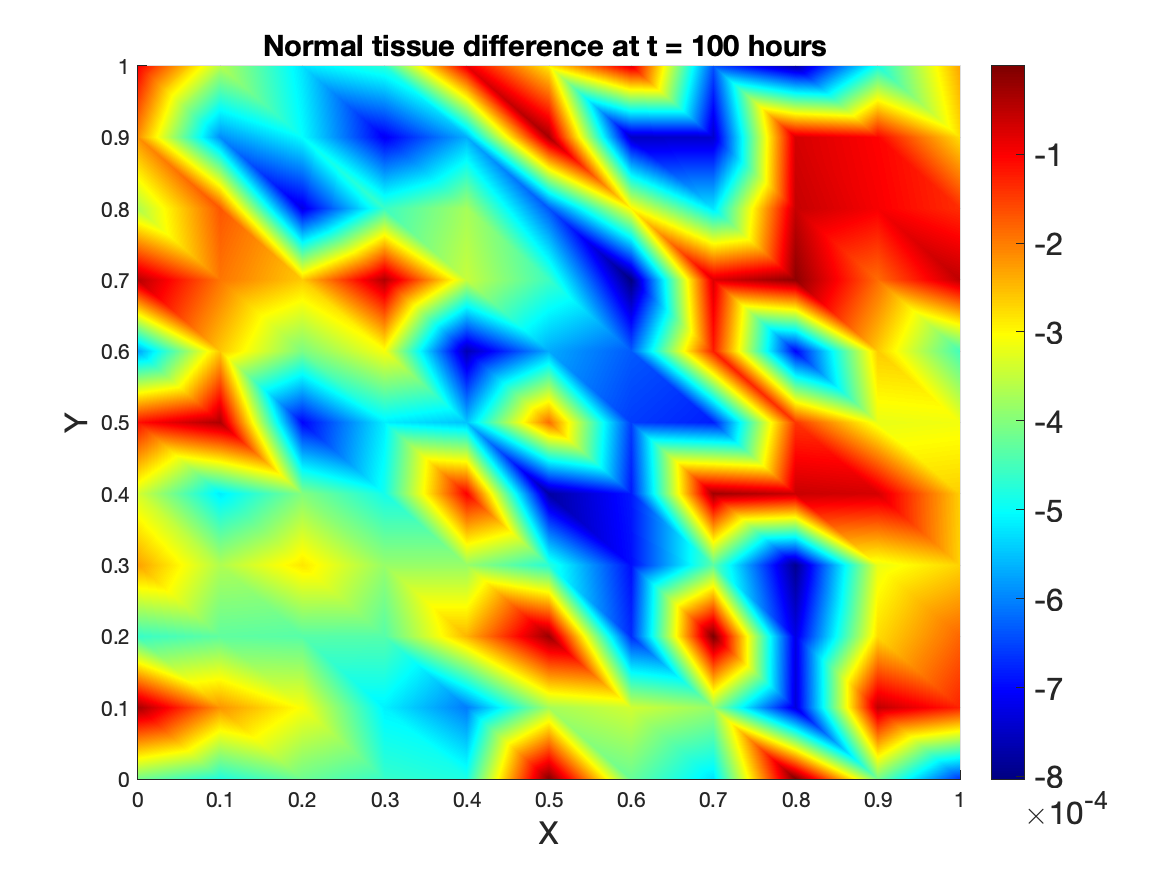}
		\caption{Normal tissue at \\ $t=100$}
	\end{subfigure}
	\begin{subfigure}{0.24\textwidth} 		\centering
		\includegraphics[width=\textwidth]{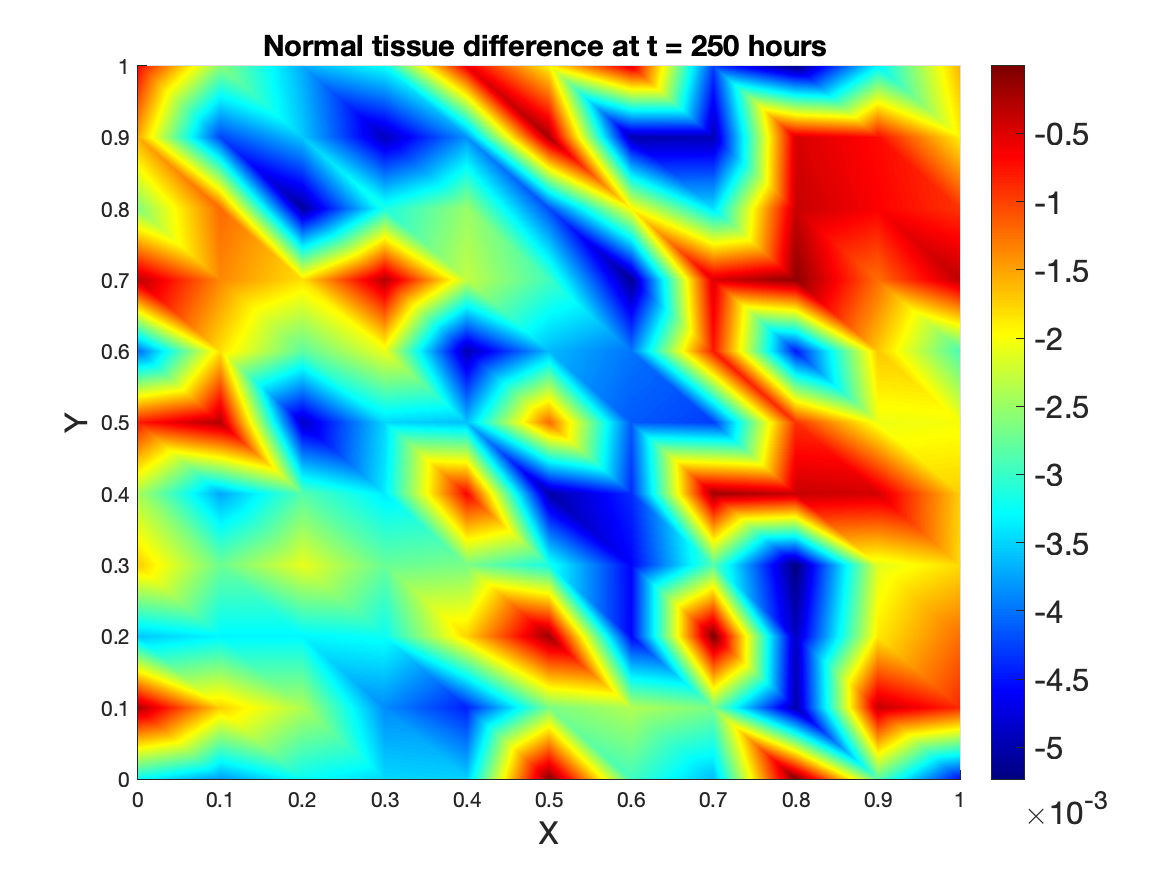}
		\caption{Normal tissue at \\ $t=250$}
	\end{subfigure}\\
	
	\begin{subfigure}{0.24\textwidth} 		\centering
		\includegraphics[width=\textwidth]{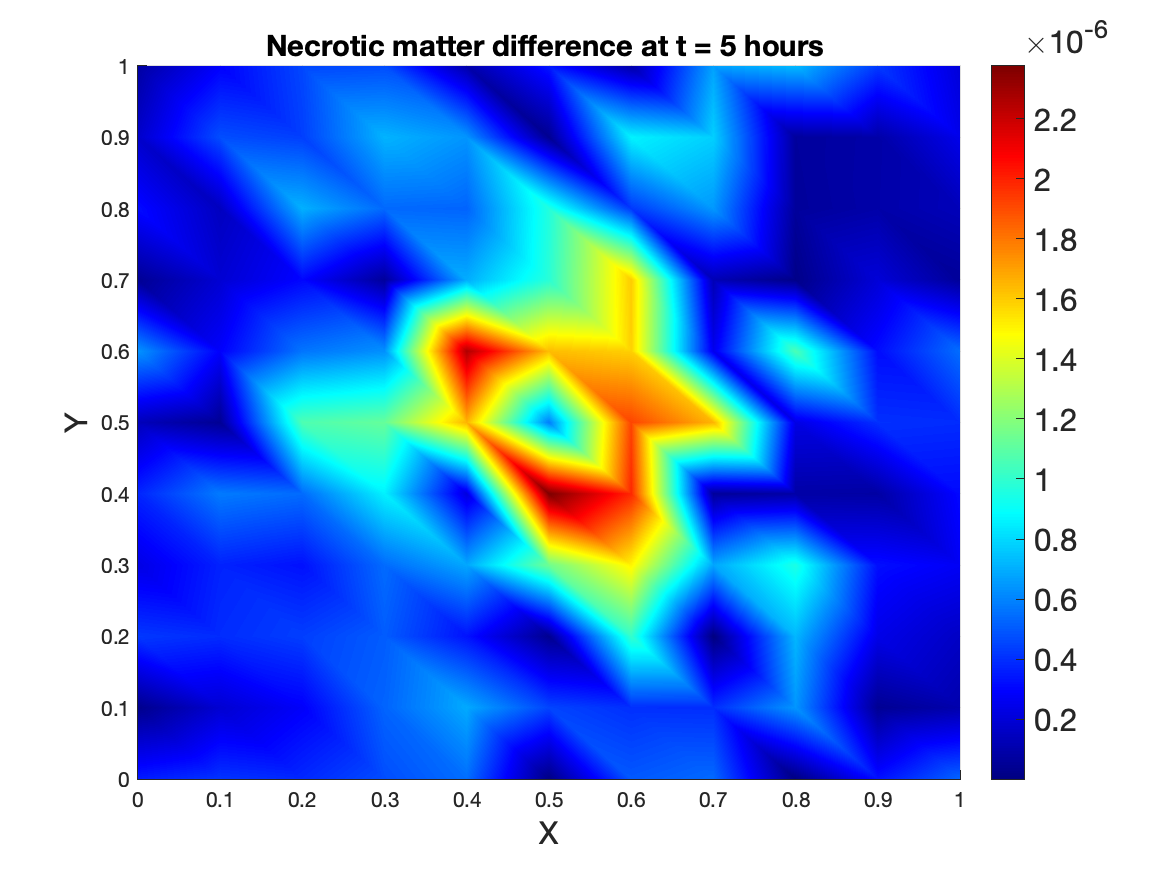}
		\caption{Necrotic matter \\ at $t=5$}
		\label{}
	\end{subfigure}  	
	\begin{subfigure}{0.24\textwidth} 		\centering
		\includegraphics[width=\textwidth]{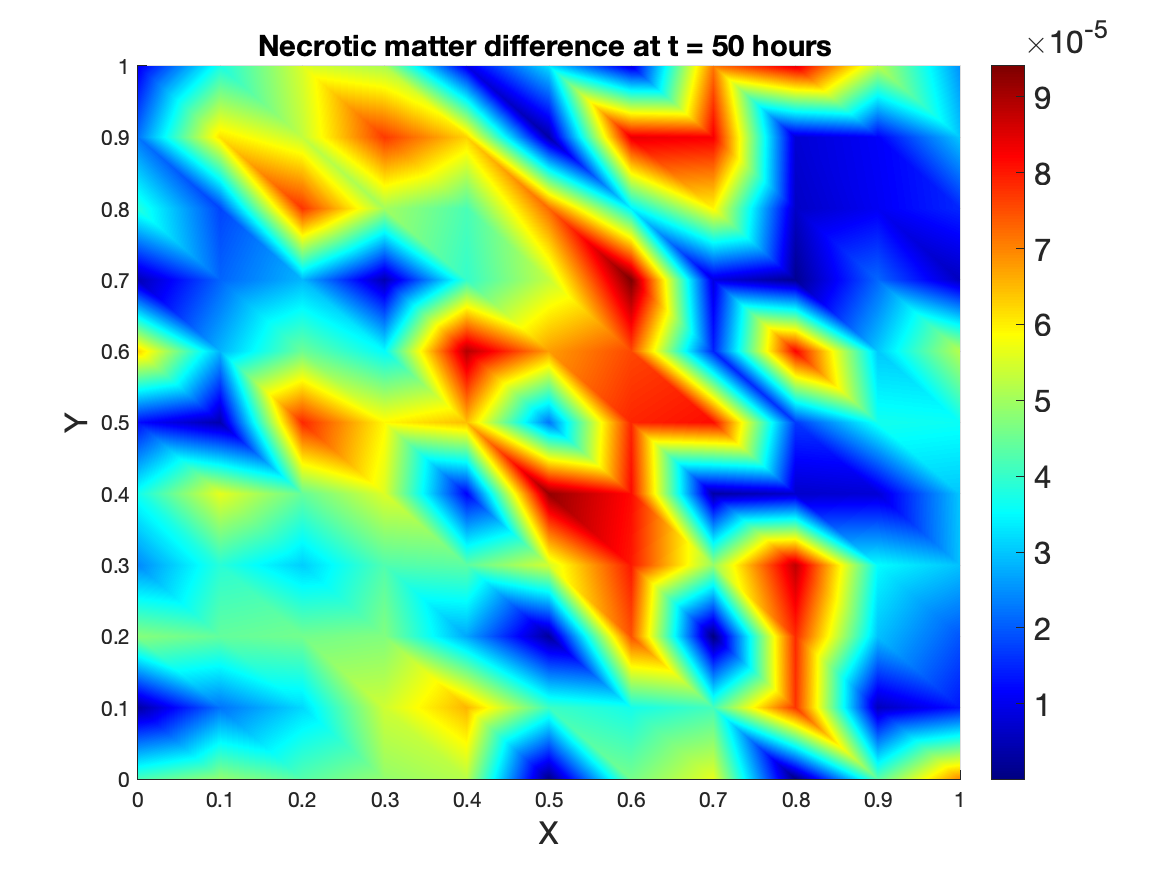}
		\caption{Necrotic matter at \\ $t=50$}
		
	\end{subfigure} 
	\begin{subfigure}{0.24\textwidth} 		\centering
		\includegraphics[width=\textwidth]{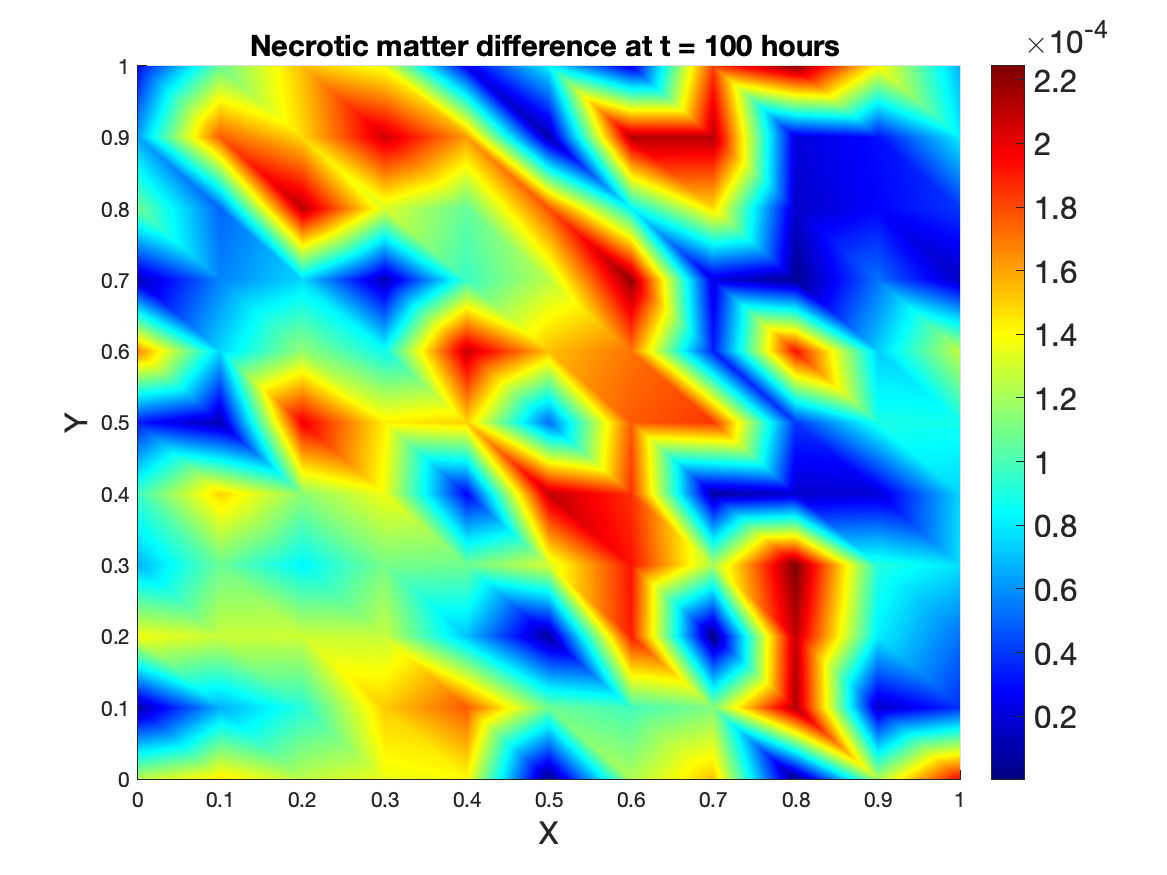}
		\caption{Necrotic matter at \\ $t=100$}
		
	\end{subfigure} 
	\begin{subfigure}{0.24\textwidth} 		\centering
		\includegraphics[width=\textwidth]{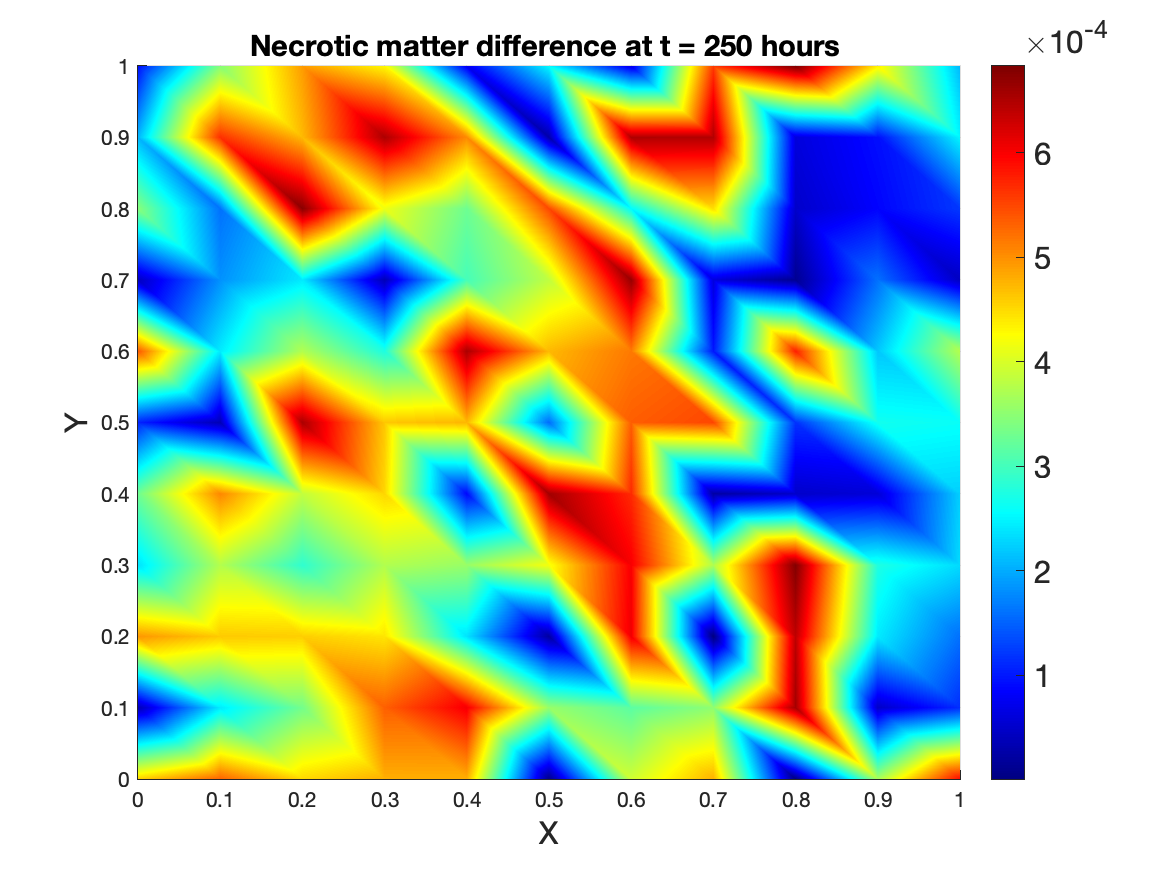}
		\caption{Necrotic matter at \\ $t=250$}
	\end{subfigure}     \vspace{0.5cm}     \\
	\caption{ Difference between densities of bacteria, mycolactone, normal tissue, and necrotic matter at different times for Scenario 2 and Scenario 3 (Scenario 3 - Scenario 2).}
	\label{fig:comparison1}
\end{figure}

\begin{figure}[htbp!]
	
	\begin{subfigure}{0.24\textwidth} 		\centering
		\includegraphics[width=\textwidth]{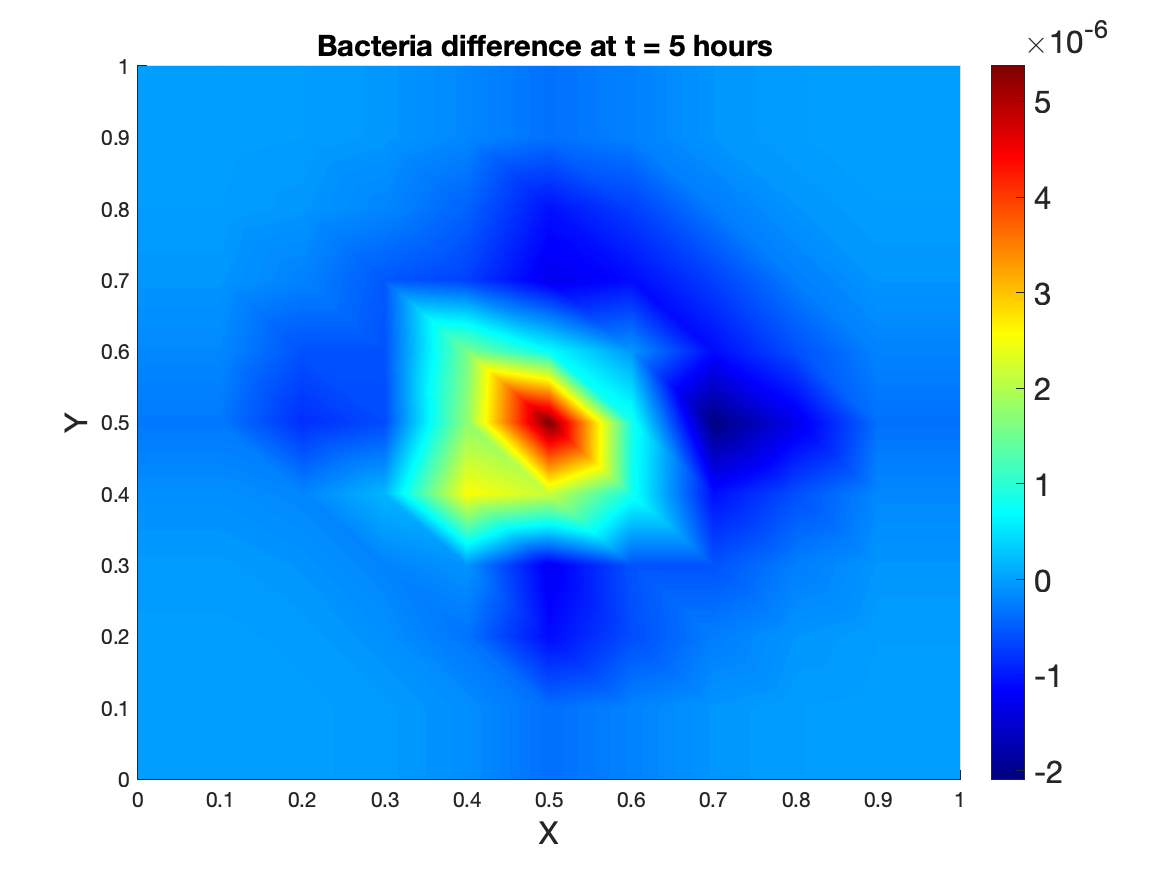}
		\caption{Bacteria at $t=5$}
		\label{}
		
	\end{subfigure} 
	\begin{subfigure}{0.24\textwidth} 		\centering
		\includegraphics[width=\textwidth]{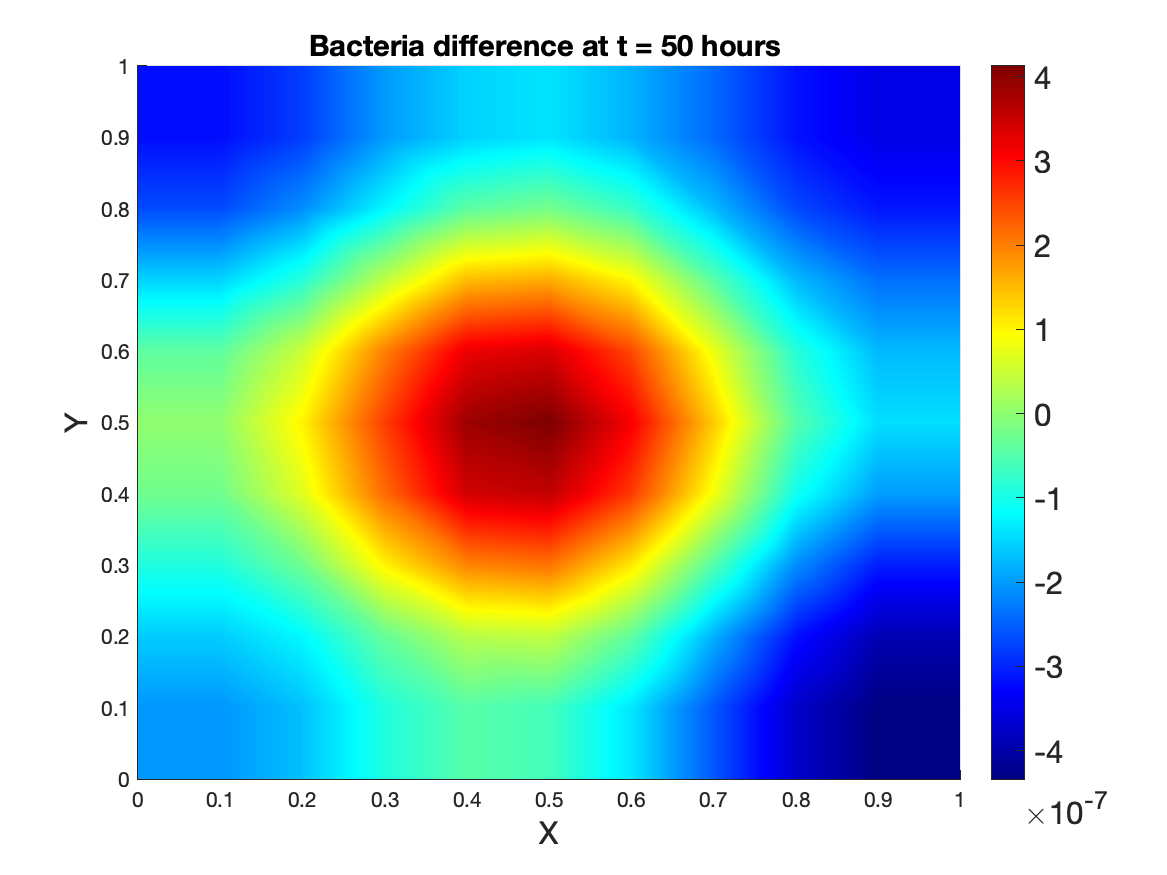}
		\caption{Bacteria at $t=50$}
		
	\end{subfigure} 
	\begin{subfigure}{0.24\textwidth} 		\centering
		\includegraphics[width=\textwidth]{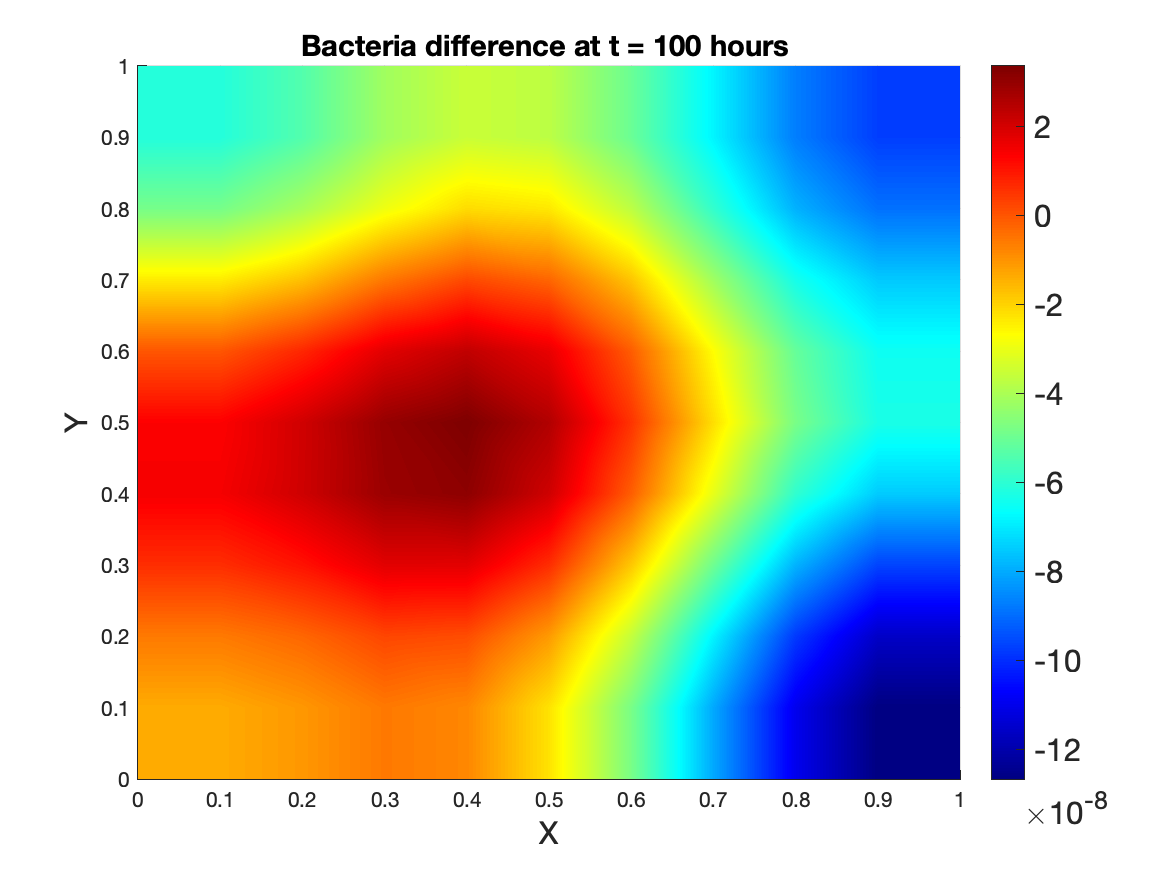}
		\caption{Bacteria at $t=100$}
		
	\end{subfigure} 
	\begin{subfigure}{0.24\textwidth} 		\centering
		\includegraphics[width=\textwidth]{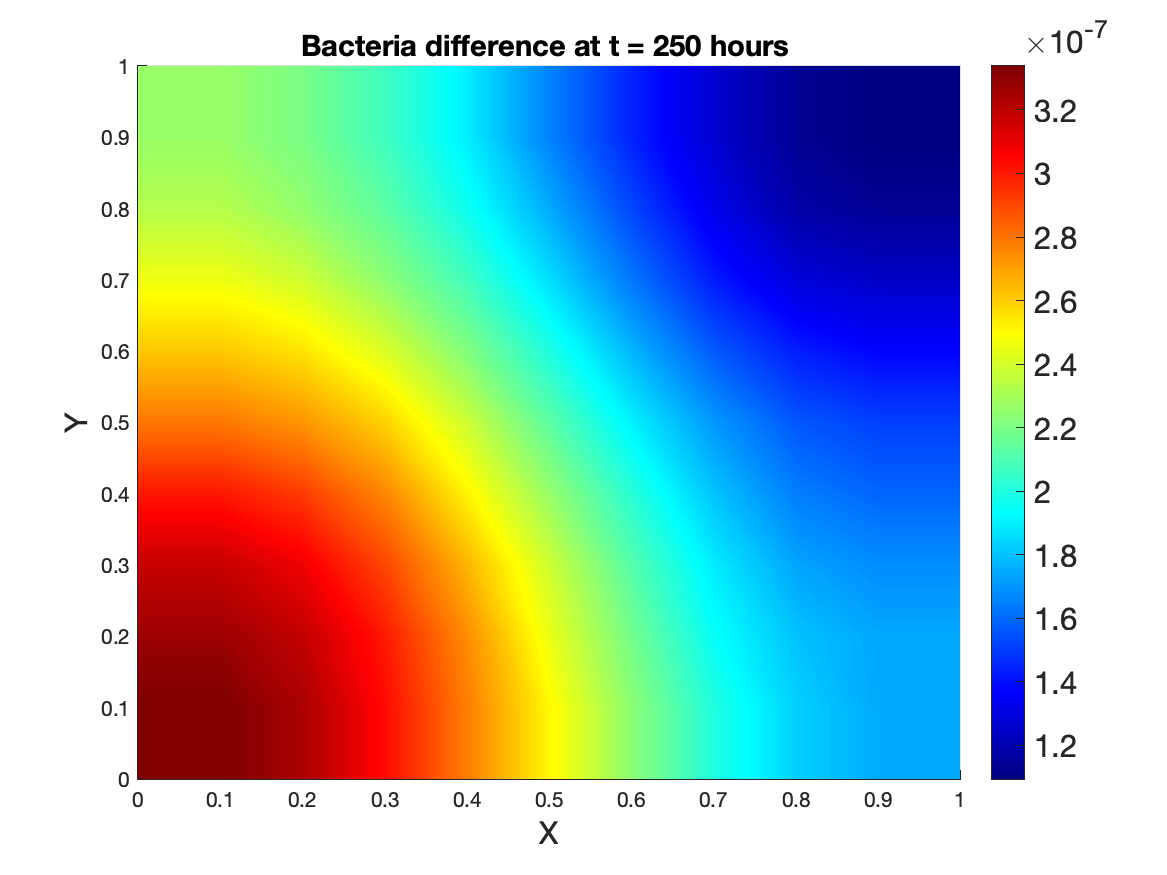}
		\caption{Bacteria at $t=250$}
		
	\end{subfigure}     \vspace{0.5cm}     \\
	
	\begin{subfigure}{0.24\textwidth} 		\centering
		\includegraphics[width=\textwidth]{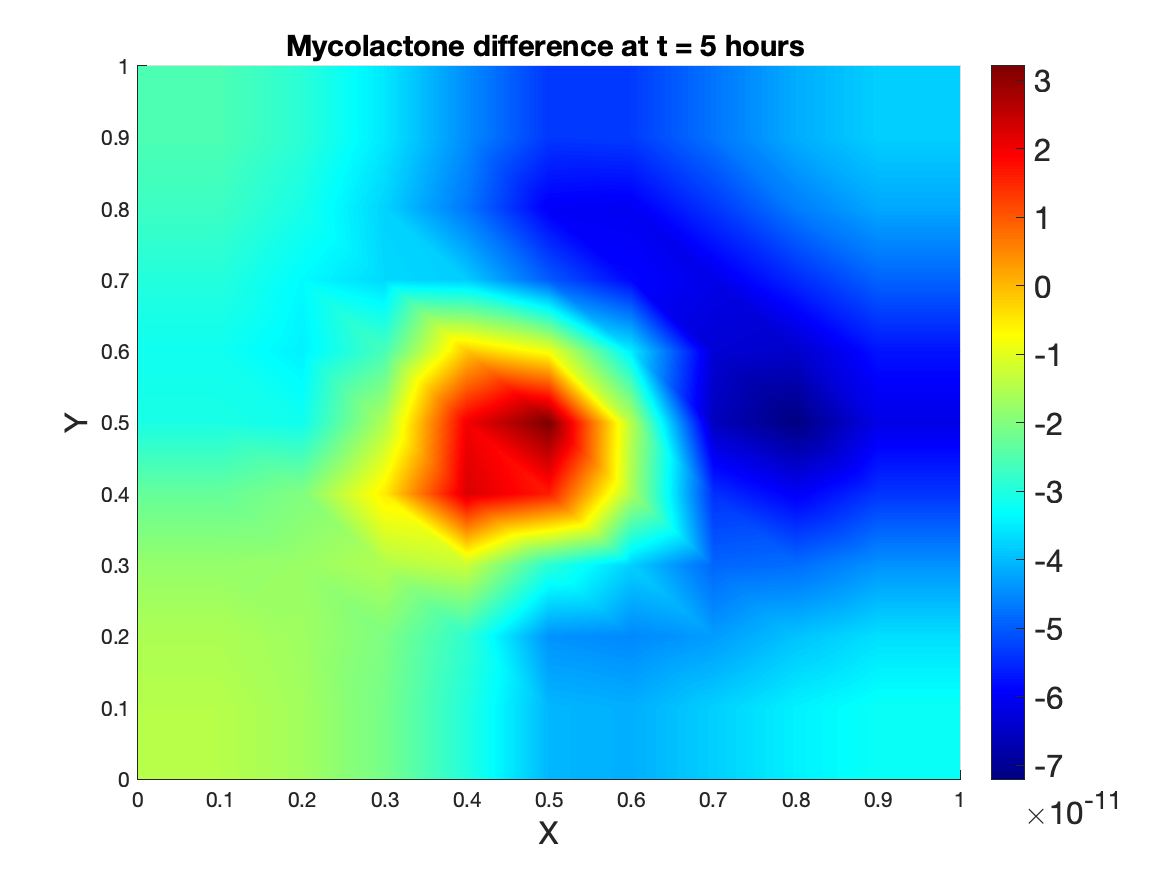}
		\caption{Mycolactone at \\ $t=5$}
		\label{}
	\end{subfigure} 
	\begin{subfigure}{0.24\textwidth} 		\centering 
		\includegraphics[width=\textwidth]{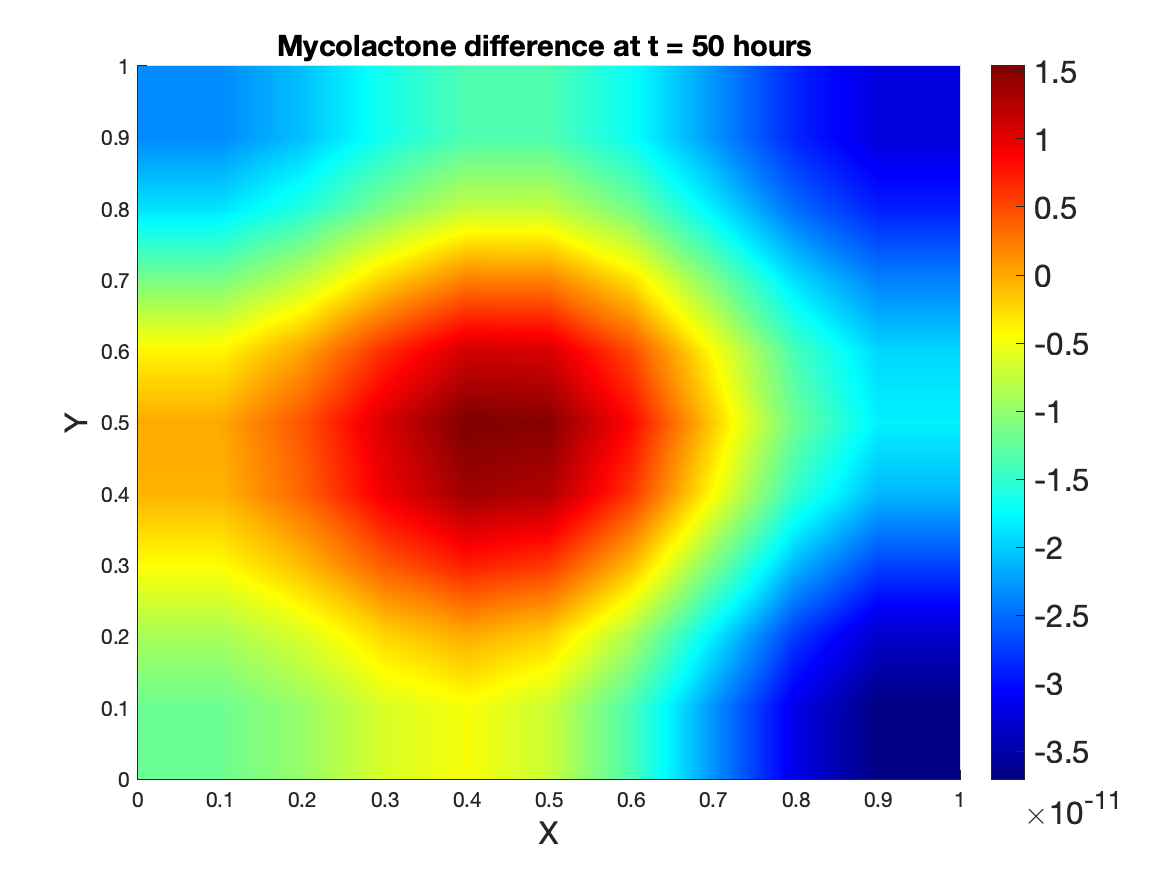}
		\caption{Mycolactone at \\ $t=50$}
		
	\end{subfigure} 
	\begin{subfigure}{0.24\textwidth} 		\centering 
		\includegraphics[width=\textwidth]{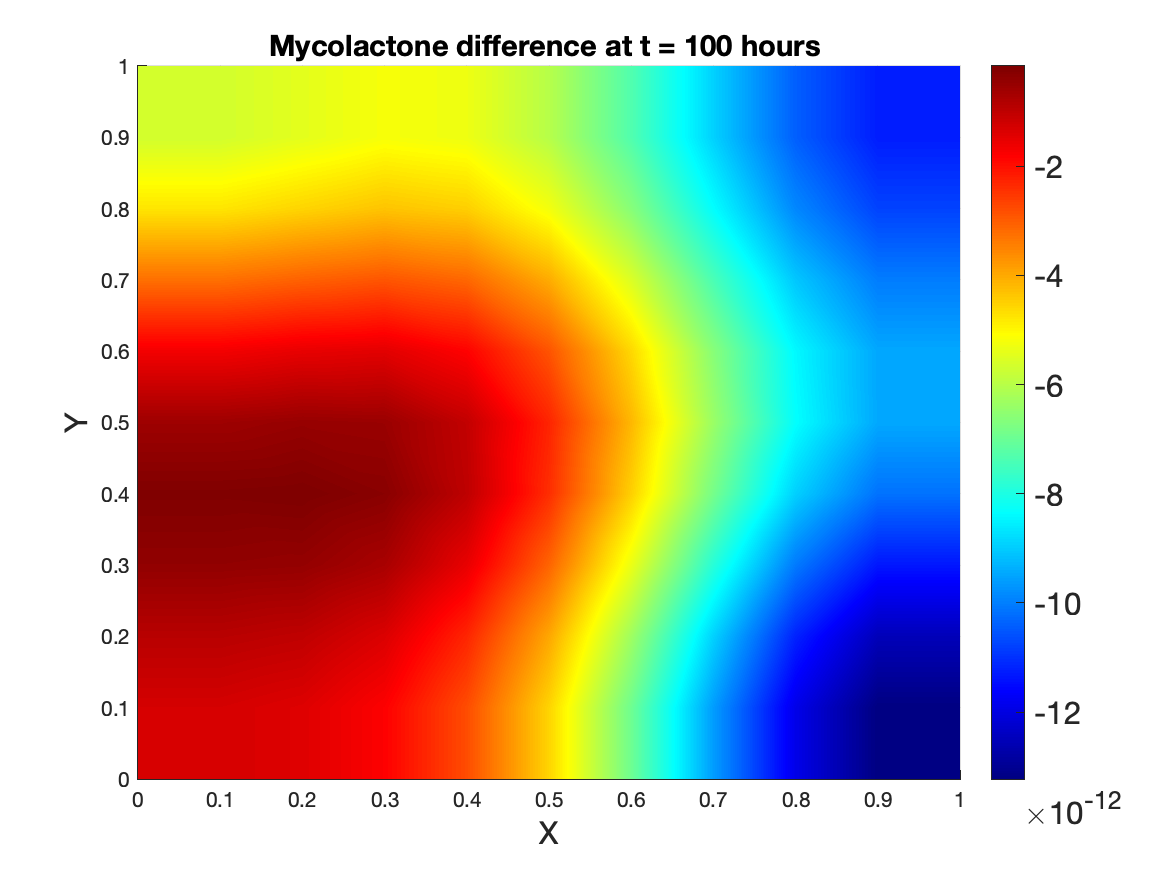}
		\caption{Mycolactone at \\ $t=100$}
		
	\end{subfigure}
	\begin{subfigure}{0.24\textwidth} 		\centering
		\includegraphics[width=\textwidth]{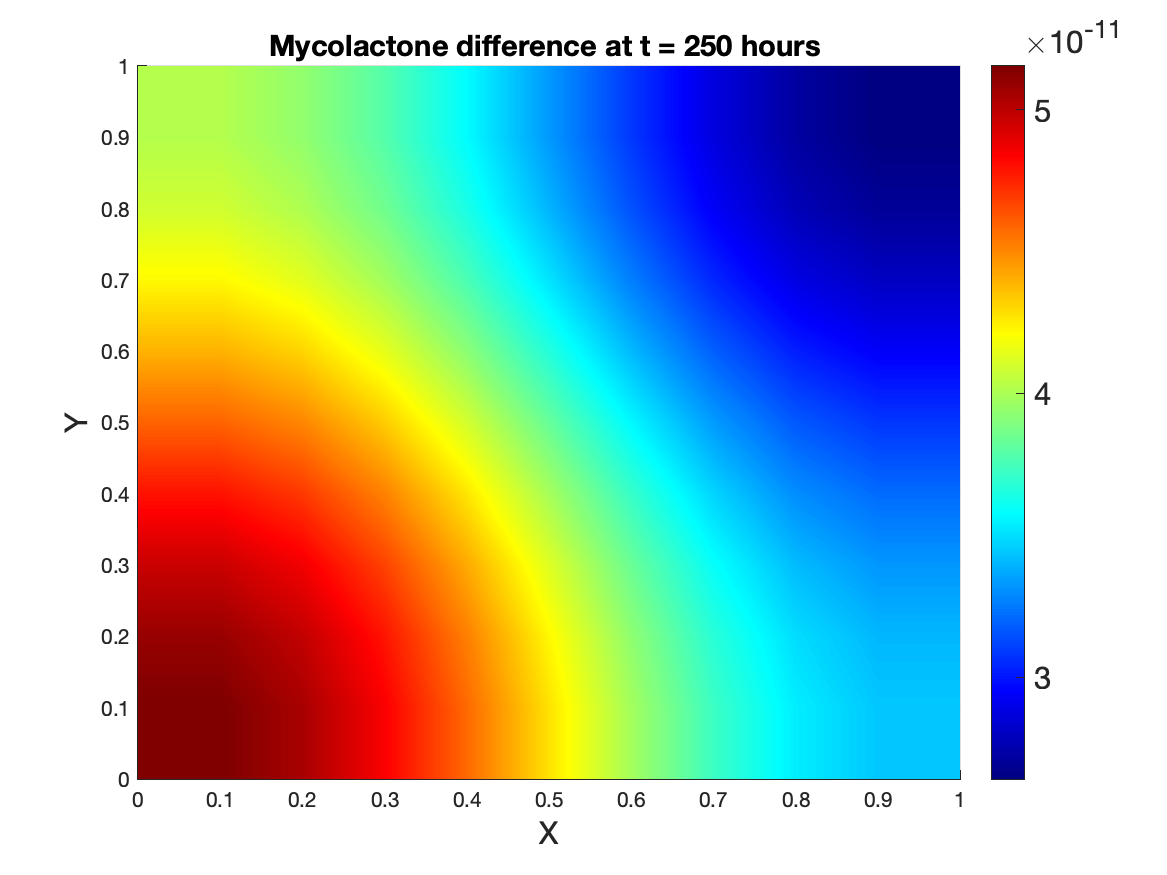}
		\caption{Mycolactone at \\ $t=250$}
	\end{subfigure}     \vspace{0.5cm}     \\
	
	\begin{subfigure}{0.24\textwidth} 		\centering
		\includegraphics[width=\textwidth]{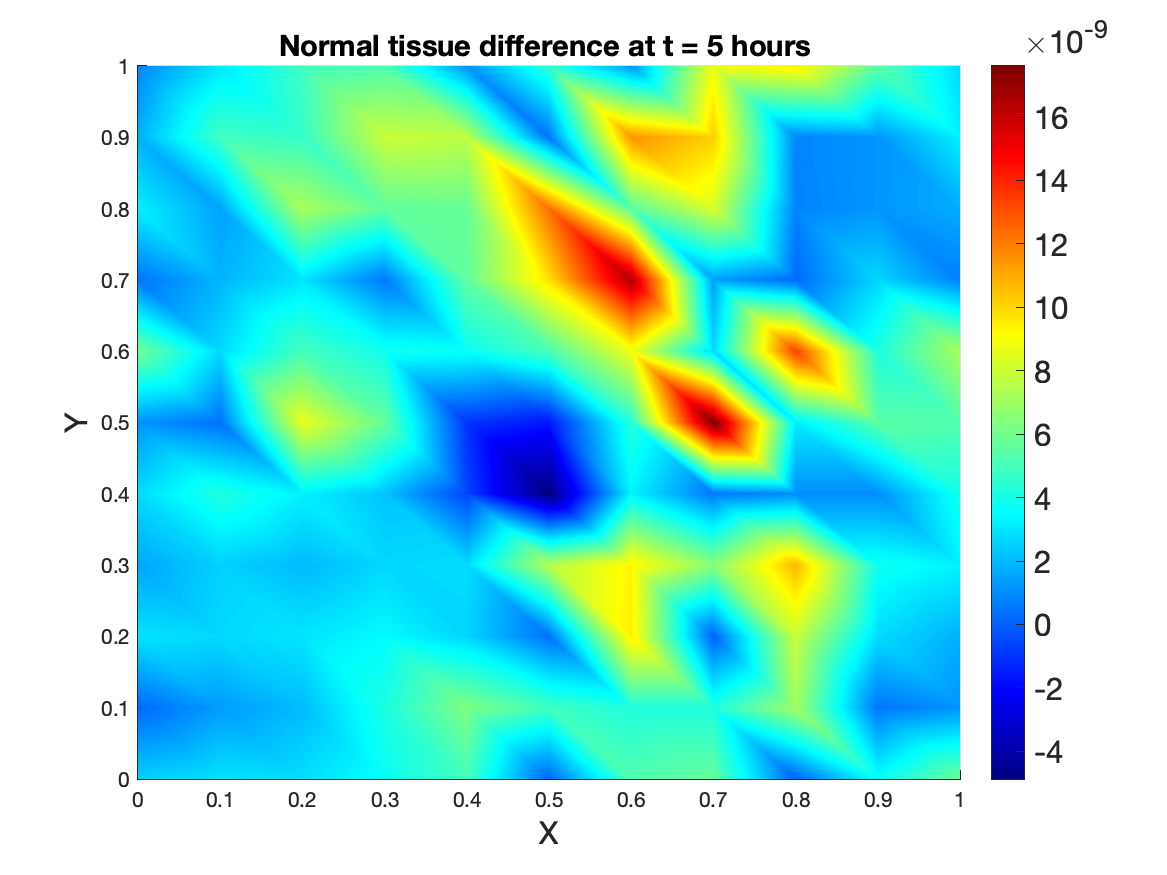}
		\caption{Normal tissue at \\ $t=5$}
		\label{}
	\end{subfigure} 
	\begin{subfigure}{0.24\textwidth} 		\centering
		\includegraphics[width=\textwidth]{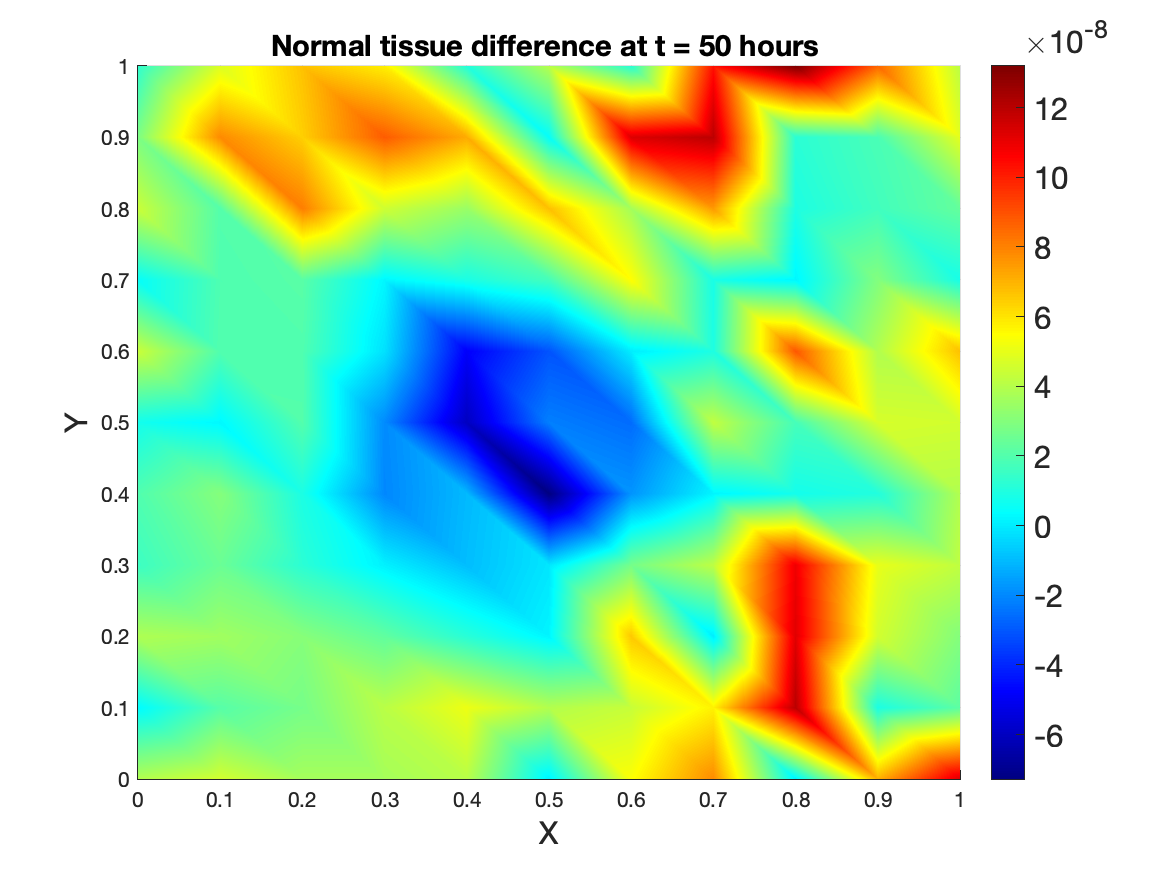}
		\caption{Normal tissue at \\ $t=50$}
	\end{subfigure}
	\begin{subfigure}{0.24\textwidth} 		\centering
		\includegraphics[width=\textwidth]{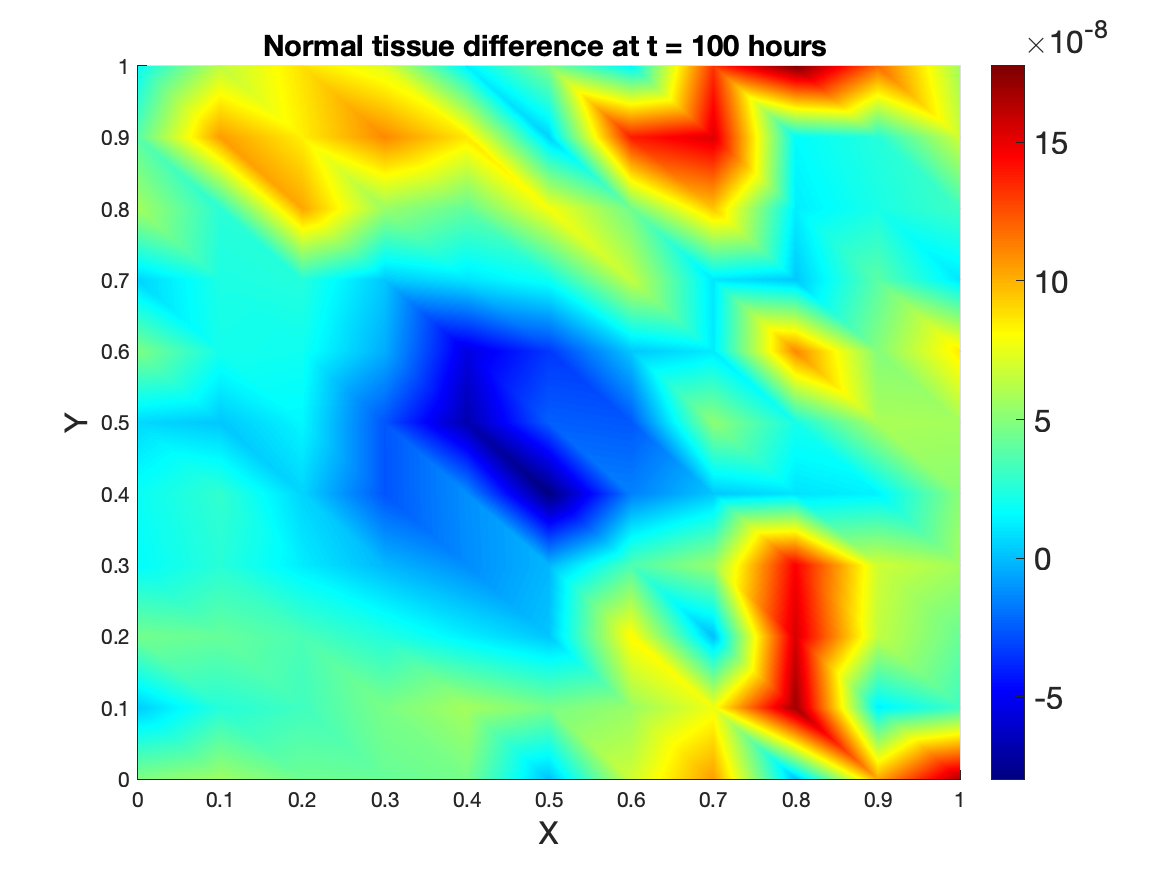}
		\caption{Normal tissue at \\ $t=100$}
	\end{subfigure}
	\begin{subfigure}{0.24\textwidth} 		\centering
		\includegraphics[width=\textwidth]{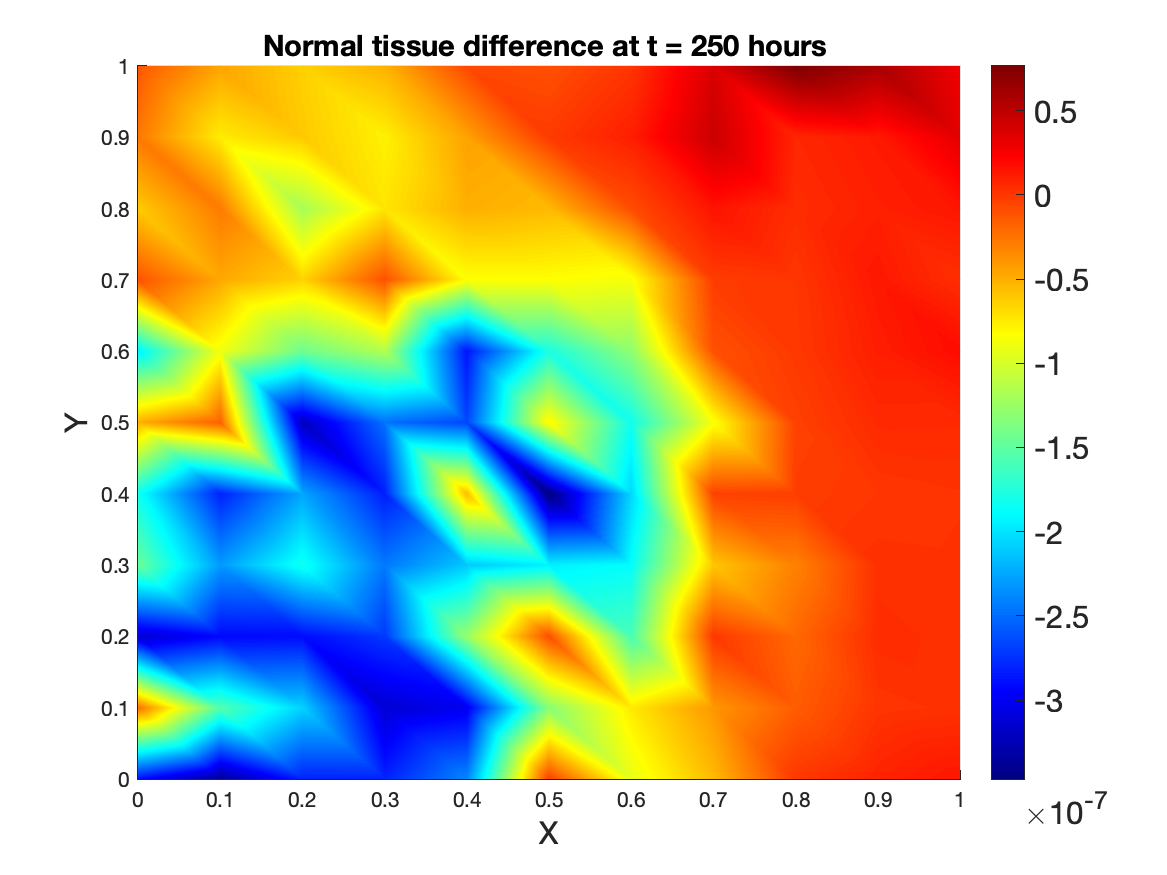}
		\caption{Normal tissue at \\ $t=250$}
	\end{subfigure}\\
	
	\begin{subfigure}{0.24\textwidth} 		\centering
		\includegraphics[width=\textwidth]{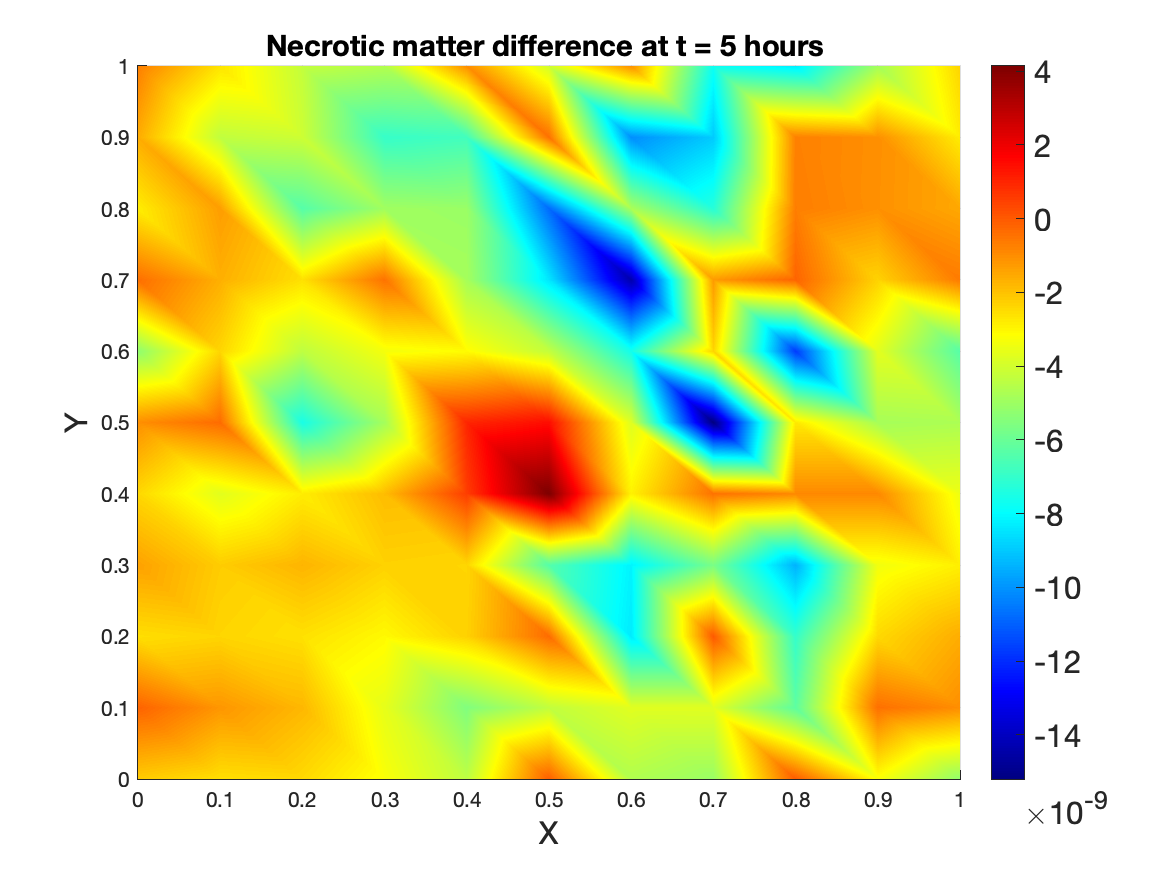}
		\caption{Necrotic matter \\ at $t=5$}
		\label{}
	\end{subfigure}  	
	\begin{subfigure}{0.24\textwidth} 		\centering
		\includegraphics[width=\textwidth]{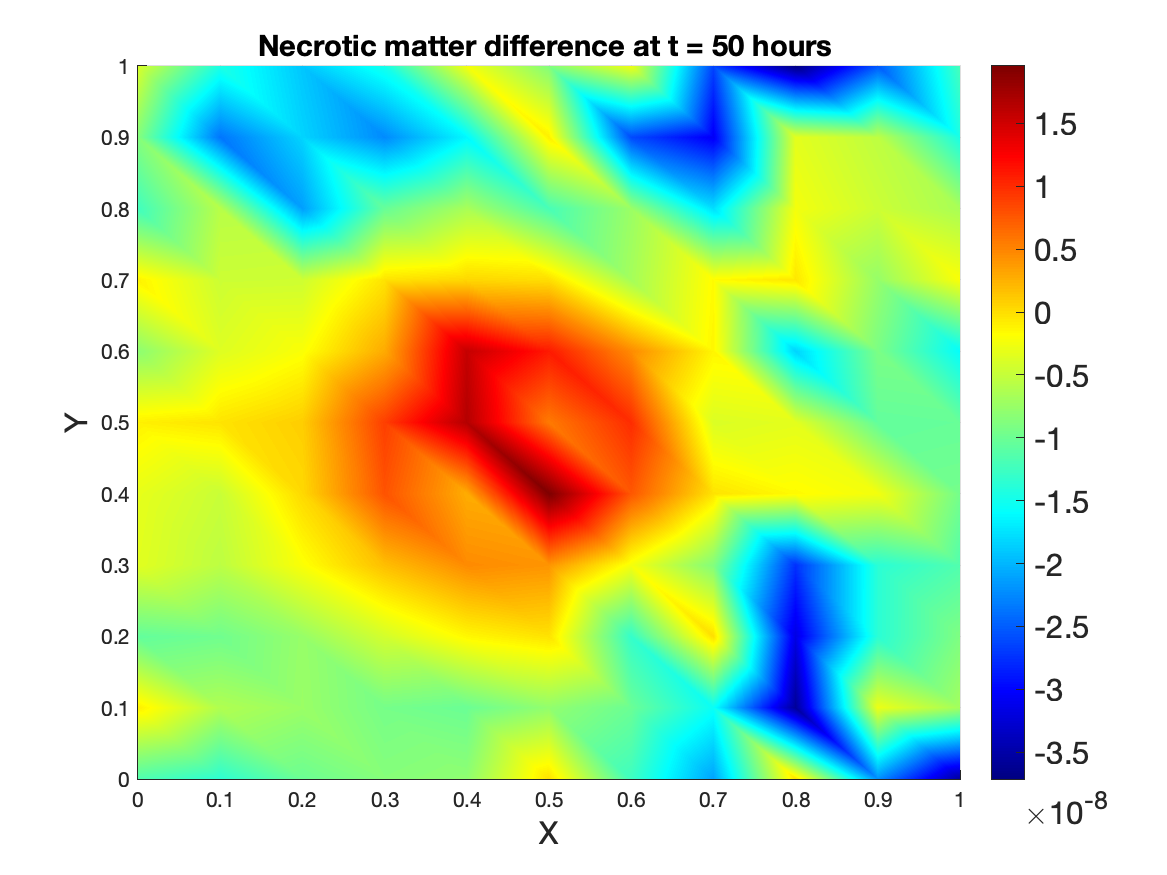}
		\caption{Necrotic matter at \\ $t=50$}
		
	\end{subfigure} 
	\begin{subfigure}{0.24\textwidth} 		\centering
		\includegraphics[width=\textwidth]{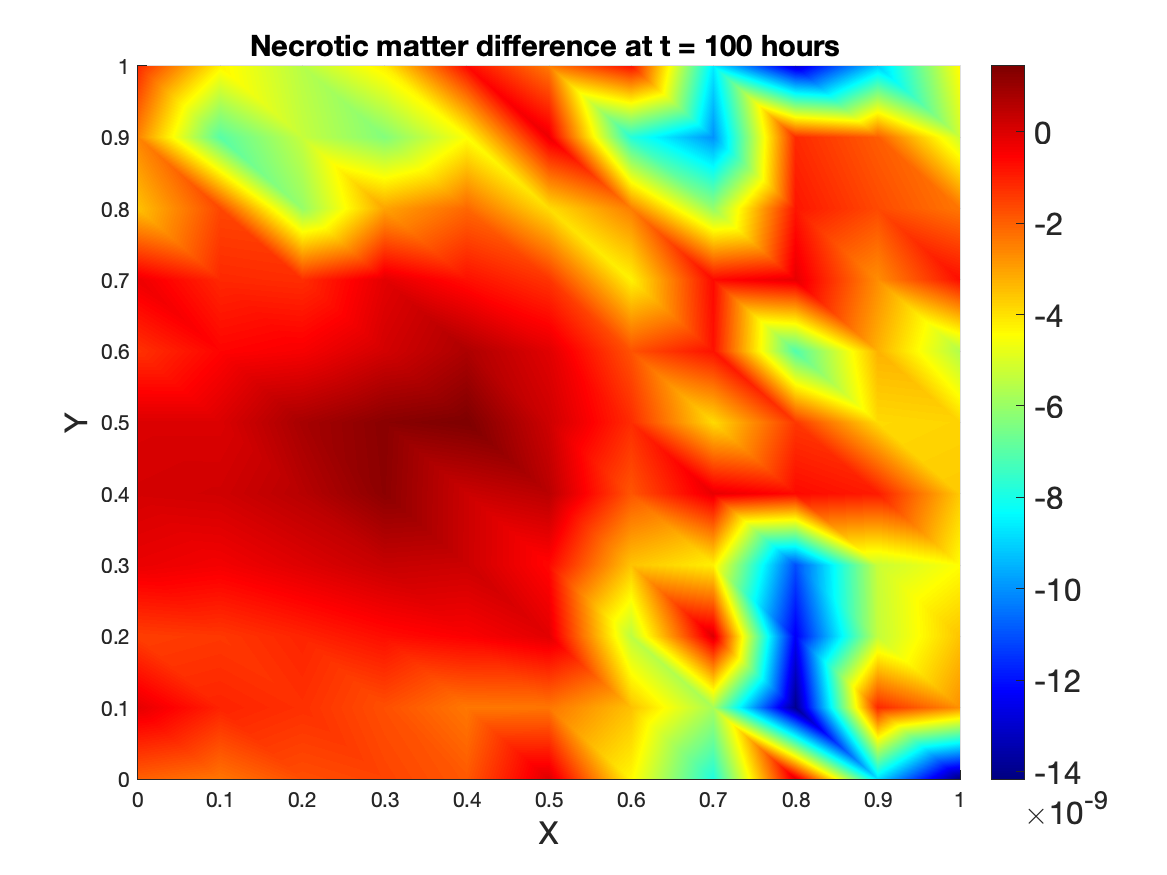}
		\caption{Necrotic matter at \\ $t=100$}
		
	\end{subfigure} 
	\begin{subfigure}{0.24\textwidth} 		\centering
		\includegraphics[width=\textwidth]{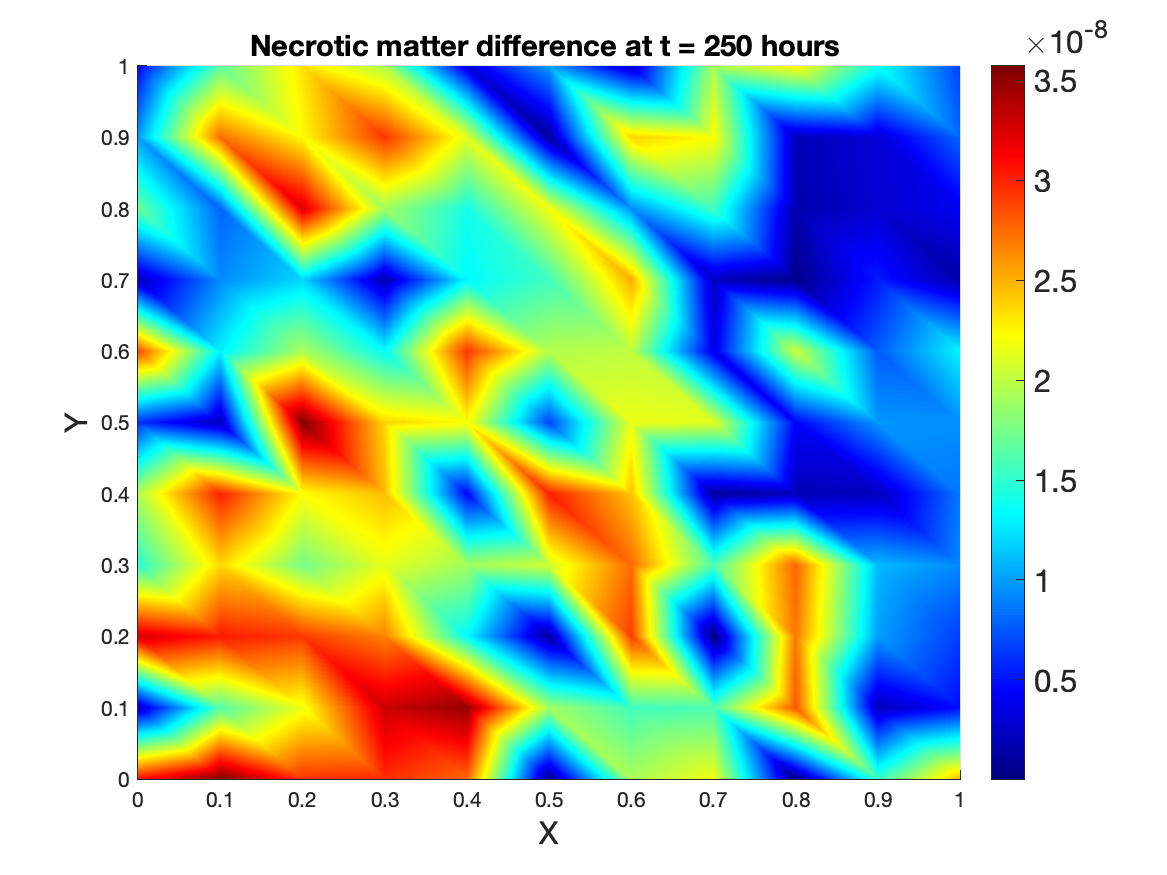}
		\caption{Necrotic matter at \\ $t=250$}
	\end{subfigure}     \vspace{0.5cm}     \\
	\caption{ Difference between densities of bacteria, mycolactone, normal tissue, and necrotic matter at different times for Scenario 4 and Scenario 1 (Scenario 4 - Scenario 1).}
	\label{fig:comparison2}
\end{figure}

\noindent In order to identify possibly significant differences in the dynamics of each solution component, we show the comparison of Scenario 2 to Scenario 3 in Figure \ref{fig:comparison1}, and of Scenario 1 to Scenario 4 in Figure \ref{fig:comparison2}.\\[-2ex]

\noindent
When comparing Scenario 2 ($\gamma_1 > \gamma_2)$ and Scenario 3 ($\gamma_2 > \gamma_1)$ (i.e., Scenario 3- Scenario 2), the bacteria density difference is rather small at early times. The absolute difference gets larger with advancing time. This is the case for all solution components, however for some of them (normal tissue, necrotic matter and bacteria) being more substantial within the simulated time span. Consequently, necrotzsation is stronger and faster when $\gamma_2 > \gamma_1$: the cells follow existing tissue rather than necrotic matter, which is slowly built up.
\\[-2ex]

\noindent
When comparing Scenario 1 (without chemotaxis) to Scenario 4 (with chemotaxis towards mycolactone) in fact we perform Scenario 4 - Scenario 1, 
several key differences emerge in the dynamics of bacterial proliferation, mycolactone concentration, normal tissue, and necrotic matter. The bacterial density difference increases from \( t = 5 \) to \( t = 50 \), decreases slightly between \( t = 50 \) and \( t = 100 \), and then rises again by \( t = 250 \), suggesting a fluctuating but overall enhanced bacterial proliferation in Scenario 4. This is driven by chemotaxis, where bacteria actively migrate towards regions with higher mycolactone concentrations, thus amplifying their spread. Mycolactone and necrotic matter differences exhibit a similar behaviour. Tissue degradation is more localised around areas with high bacterial concentrations. Notably, there are negative differences in normal tissue density where bacterial concentrations are higher, indicating that in Scenario 4, necrotisation is more severe in regions of bacterial accumulation, or at least there are some higher cell aggregates - probably due to the cells being attracted to the sites with higher mycolactone concentration. Bacteria accumulate in areas with high toxin concentration and express even more mycolactone, which in turn leads to enhanced tissue degradation. Overall, however, the differences are so small, that one could give up chemotaxis towards mycolactone, in order not to complicate unncecessarily the model.\\[-2ex]

\begin{figure}[htbp!]
	\centering
	\begin{subfigure}{0.24\textwidth}
		\centering
		\includegraphics[width=\textwidth]{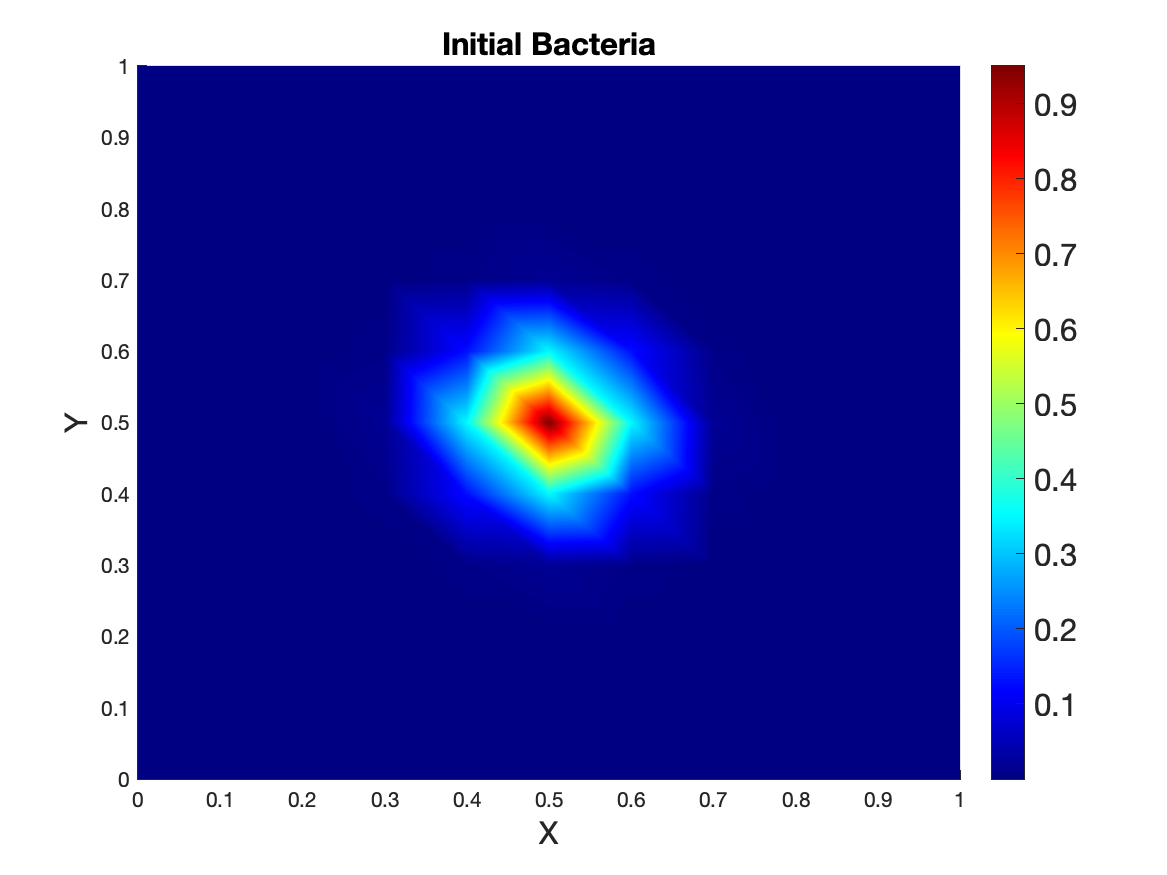}
		\caption{Bacteria }
		\label{} 
	\end{subfigure} 
	\begin{subfigure}{0.24\textwidth}
		\centering
		\includegraphics[width=\textwidth]{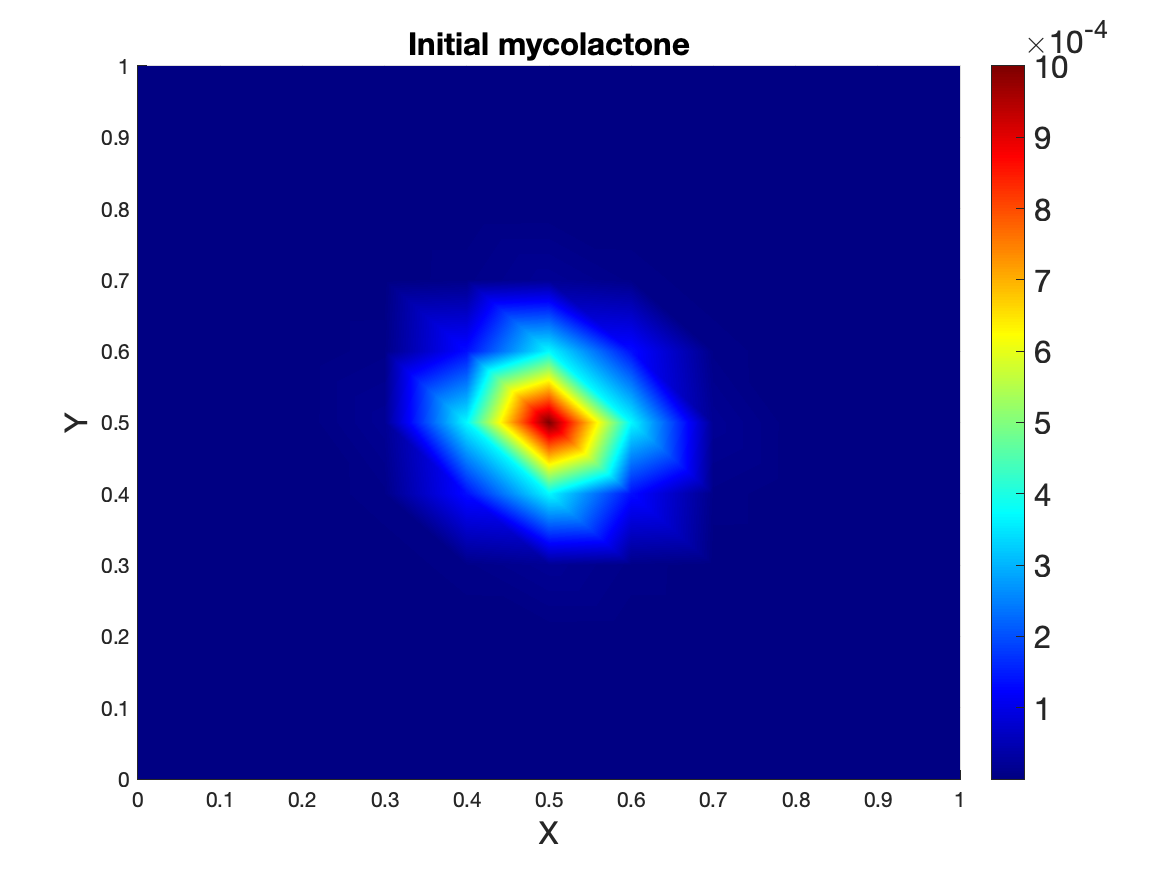}
		\caption{Mycolactone}
		
	\end{subfigure} 
	\begin{subfigure}{0.24\textwidth}
		\centering
		\includegraphics[width=\textwidth]{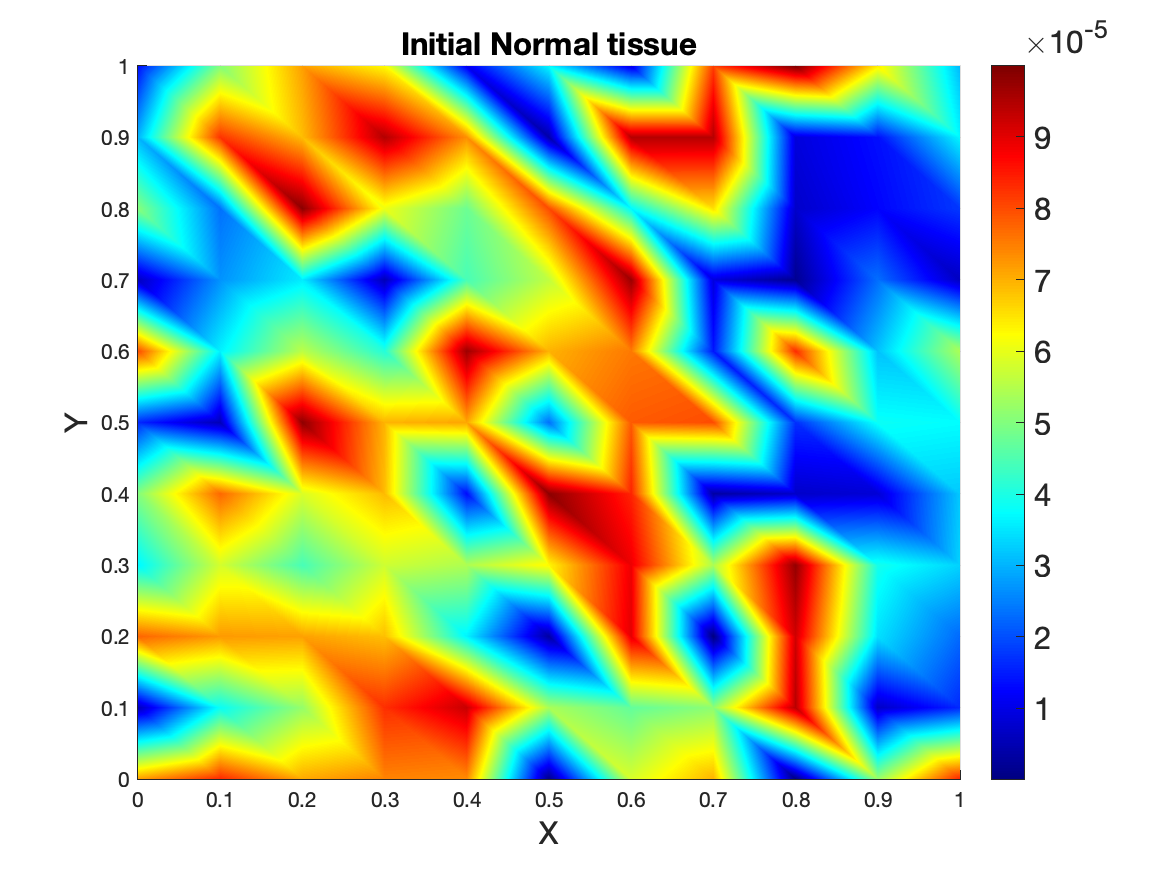}
		\caption{Normal tissue}
		
	\end{subfigure} 
	\begin{subfigure}{0.24\textwidth}
		\centering
		\includegraphics[width=\textwidth]{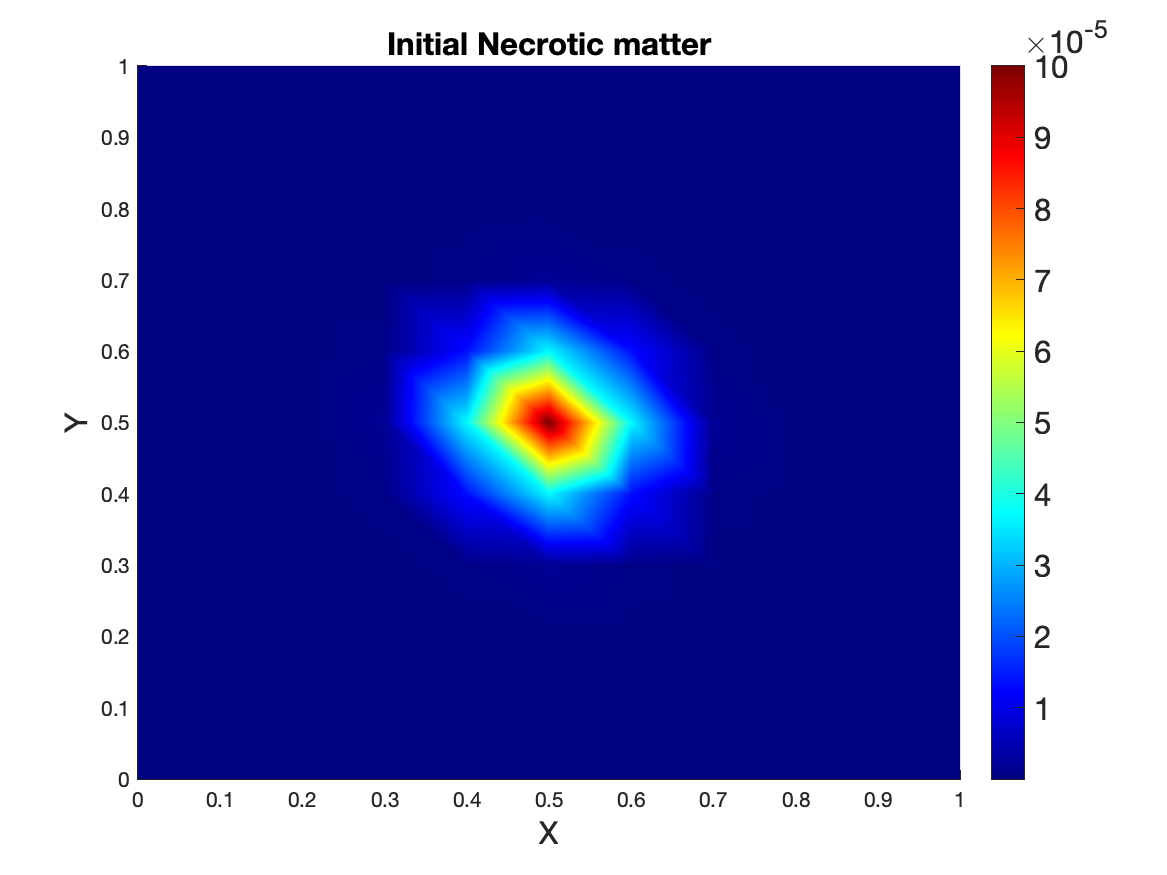}
		\caption{Necrotic matter}
		
	\end{subfigure} 
	\caption{Initial conditions for bacteria, mycolactone, normal tissue, and necrotic matter in  Scenario 5}
	\label{IC-small}
\end{figure}

\noindent
We also consider a \textbf{Scenario 5}, with much less initial density of normal tissue. Here we take $\gamma_1$ and $\gamma_2$ equal. The initial conditions in this scenario are given as follows:
\begin{equation}
	\begin{aligned}
		\label{ic-smallnorma}
		u(0,x,y)  &= 0.95 \exp(-\frac{(x-0.5)^2+(y-0.5)^2}{.01}), \qquad \ x,y\in[0,1],
		\\
		m(0,x,y) &= 0.001 \exp(-\frac{(x-0.5)^2+(y-0.5)^2}{.01}), \qquad   x,y \in[0,1],
		\\
		v(0,x,y) & = 0.0001*\mathcal{U}, 
		\\
		n(0,x,y) &=0.0001\exp(-\frac{(x-0.5)^2+(y-0.5)^2}{.01}), \quad x,y\in[0,1],
	\end{aligned}
\end{equation}
where $\mathcal{U}$ represents a uniform distribution in $(0,1).$ They are illustrated in Figure \ref{IC-small}.

\begin{figure}[htbp!]
	
	\begin{subfigure}{0.24\textwidth} 		\centering
		\includegraphics[width=\textwidth]{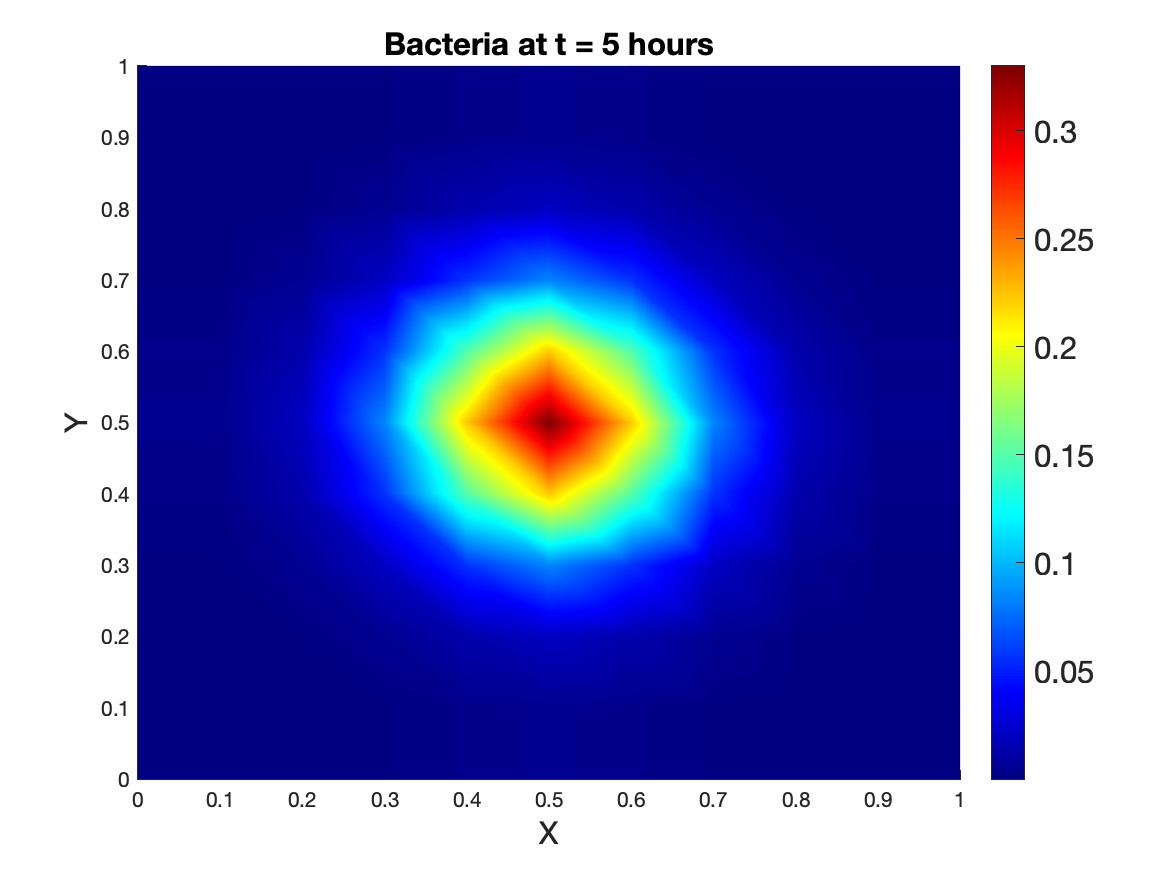}
		\caption{Bacteria at $t=5$}
		\label{}
		
	\end{subfigure} 
	\begin{subfigure}{0.24\textwidth} 		\centering
		\includegraphics[width=\textwidth]{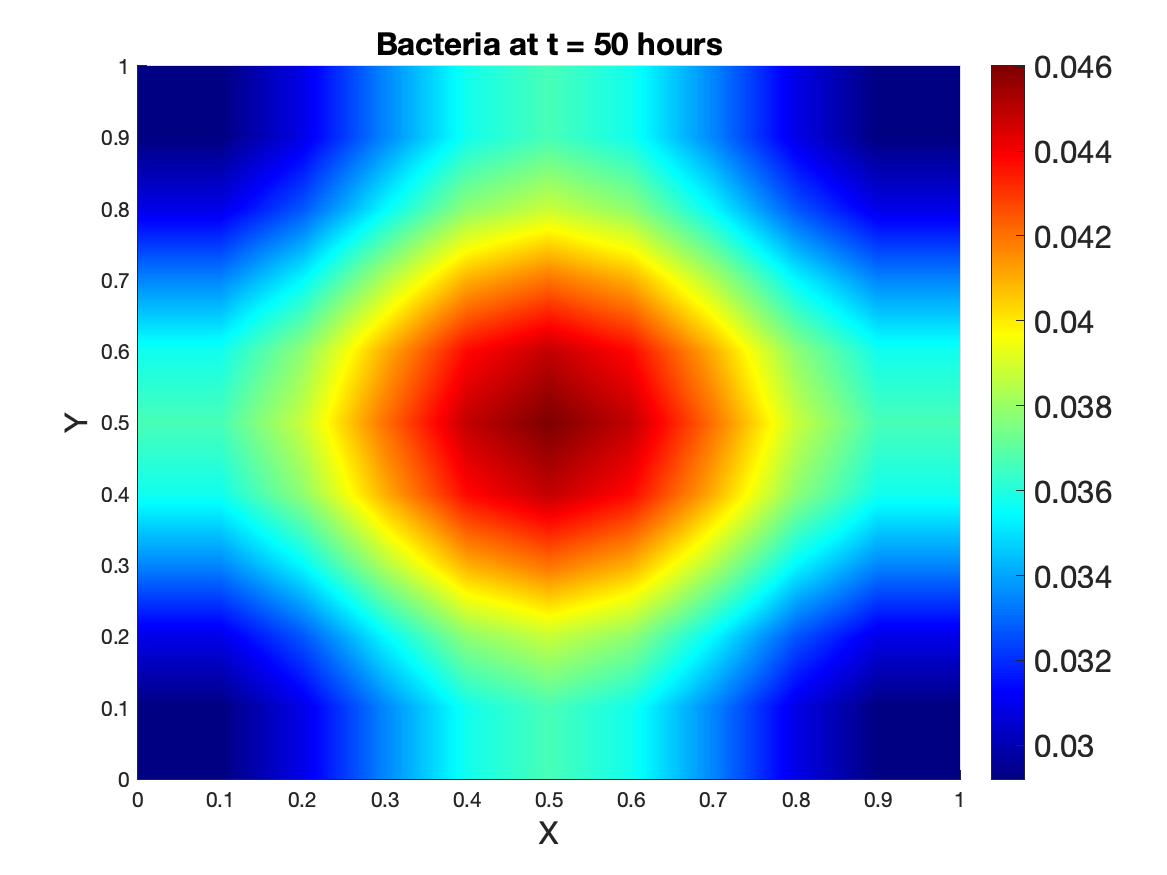}
		\caption{Bacteria at $t=50$}
		
	\end{subfigure} 
	\begin{subfigure}{0.24\textwidth} 		\centering
		\includegraphics[width=\textwidth]{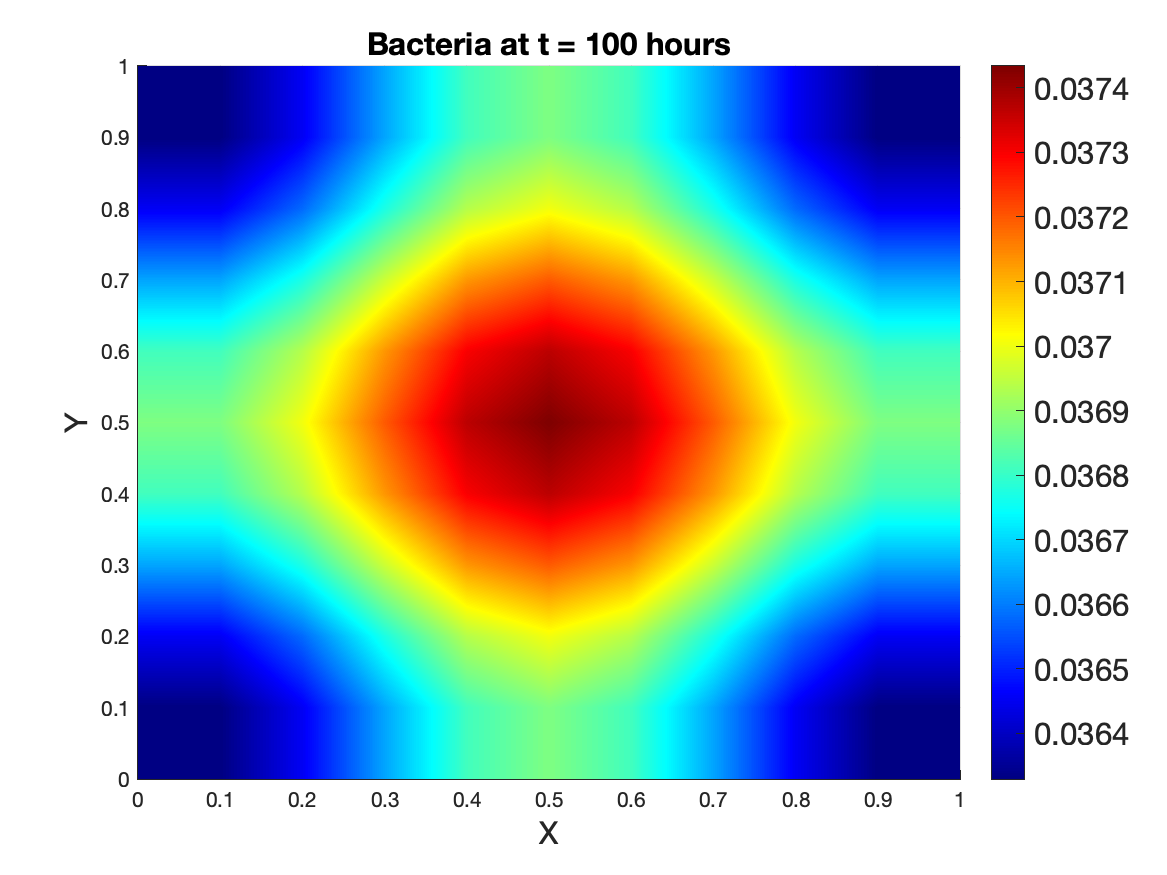}
		\caption{Bacteria at $t=100$}
		
	\end{subfigure} 
	\begin{subfigure}{0.24\textwidth} 		\centering
		\includegraphics[width=\textwidth]{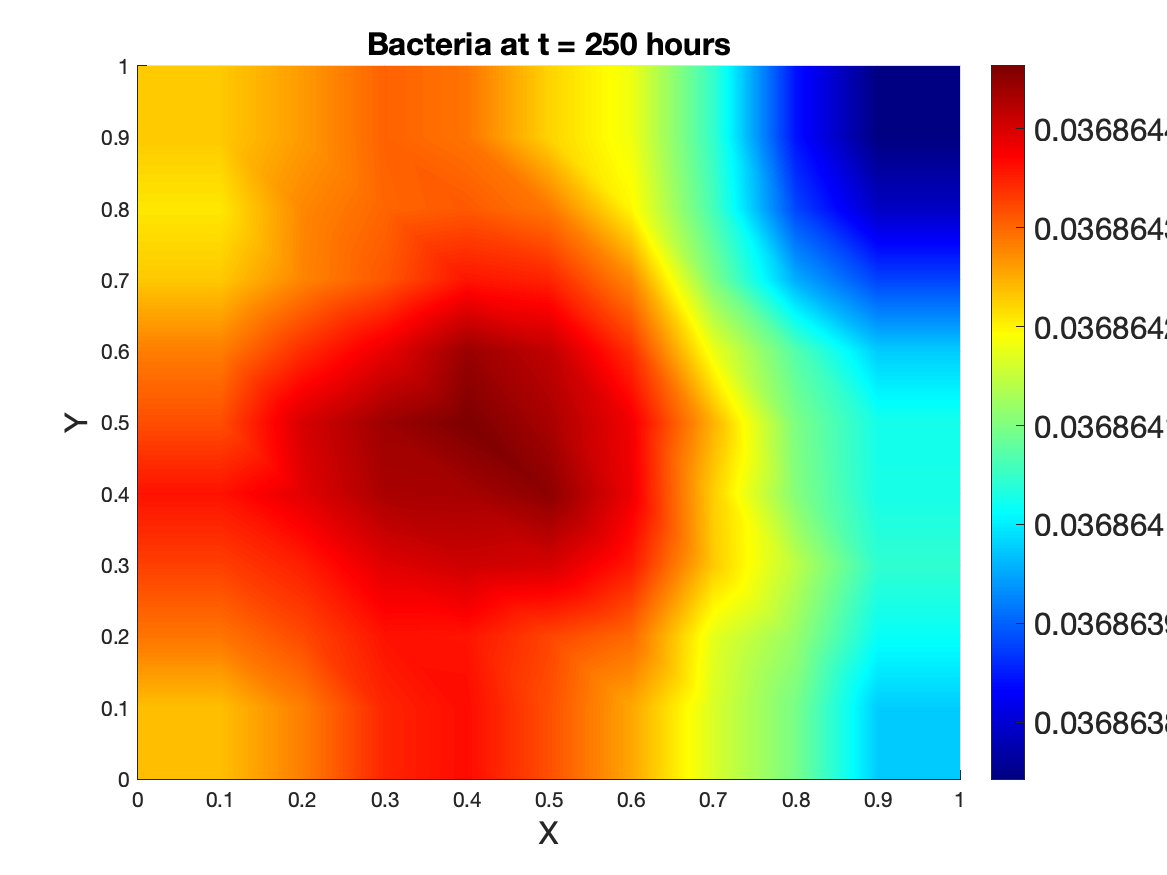}
		\caption{Bacteria at $t=250$}
		
	\end{subfigure}     \vspace{0.5cm}     \\
	
	\begin{subfigure}{0.24\textwidth} 		\centering
		\includegraphics[width=\textwidth]{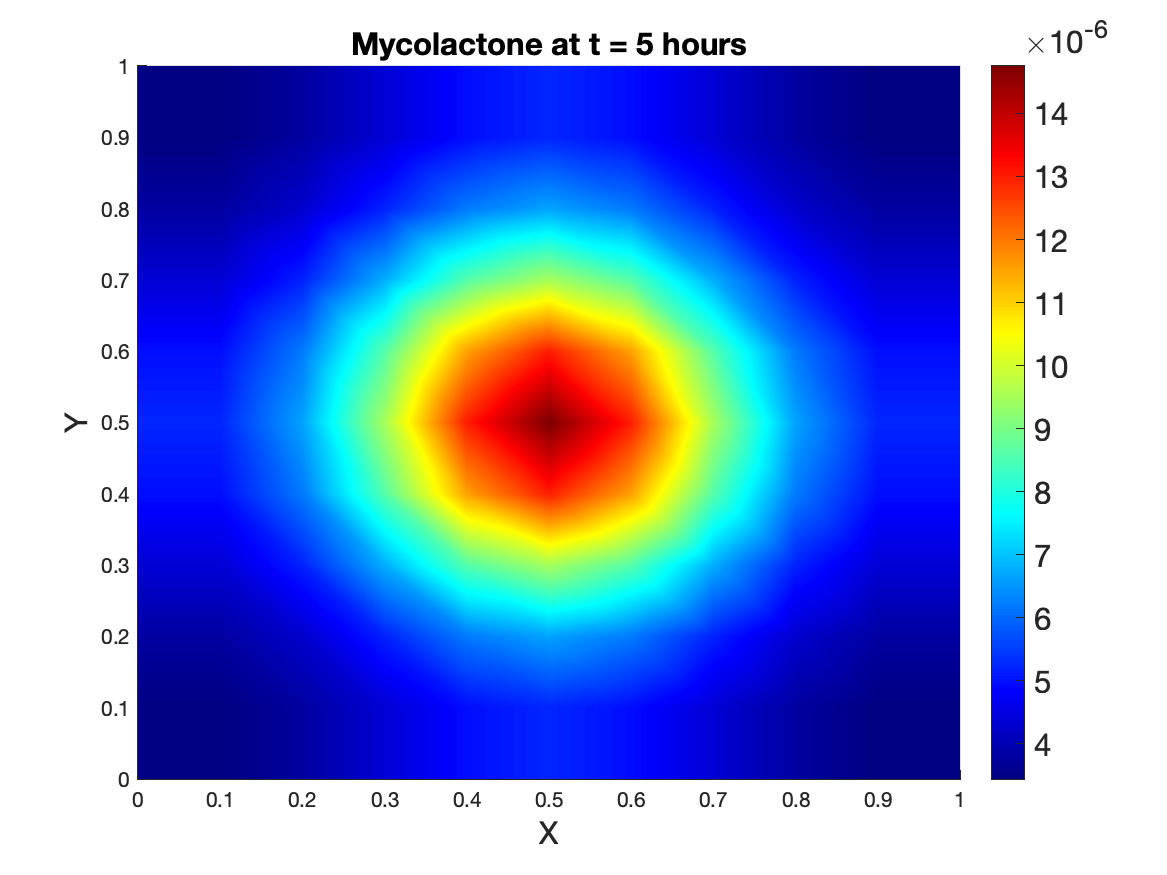}
		\caption{Mycolactone at \\ $t=5$}
		\label{}
	\end{subfigure} 
	\begin{subfigure}{0.24\textwidth} 		\centering 
		\includegraphics[width=\textwidth]{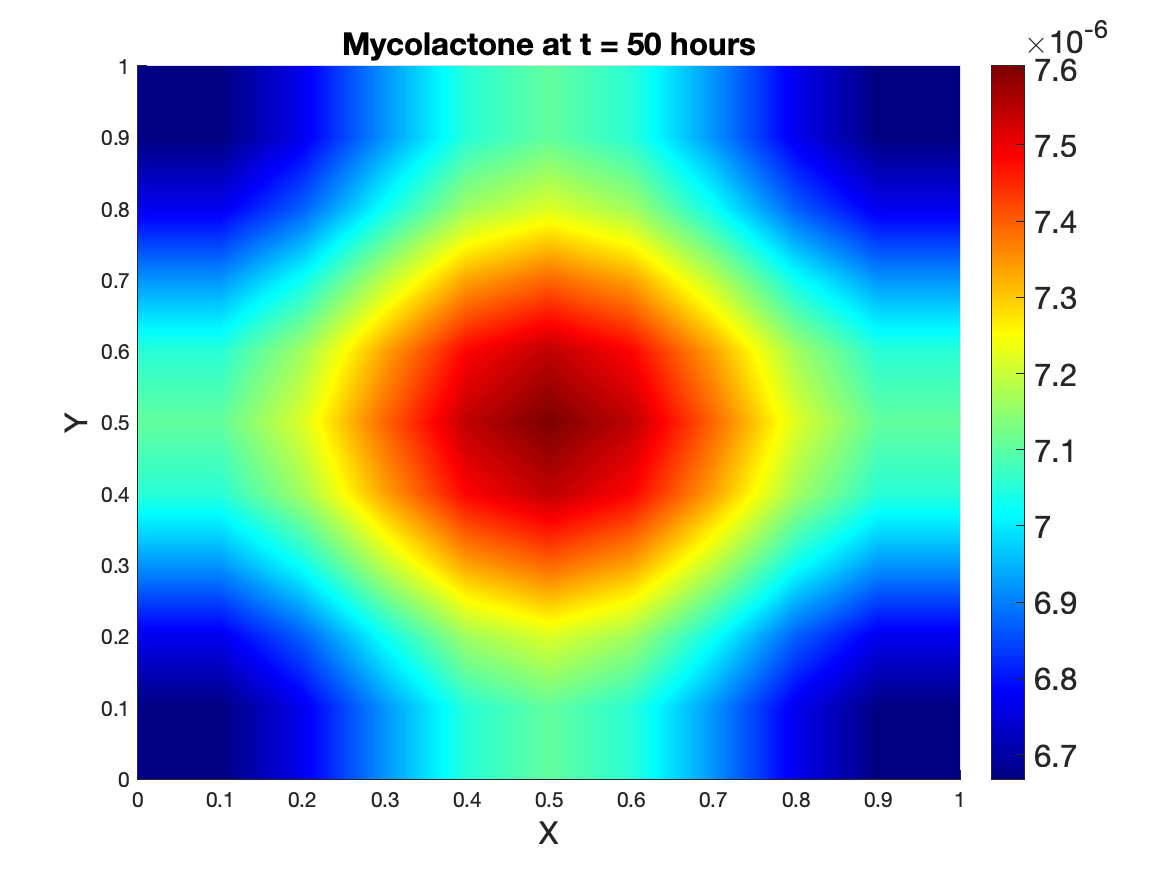}
		\caption{Mycolactone at \\ $t=50$}
		
	\end{subfigure} 
	\begin{subfigure}{0.24\textwidth} 		\centering 
		\includegraphics[width=\textwidth]{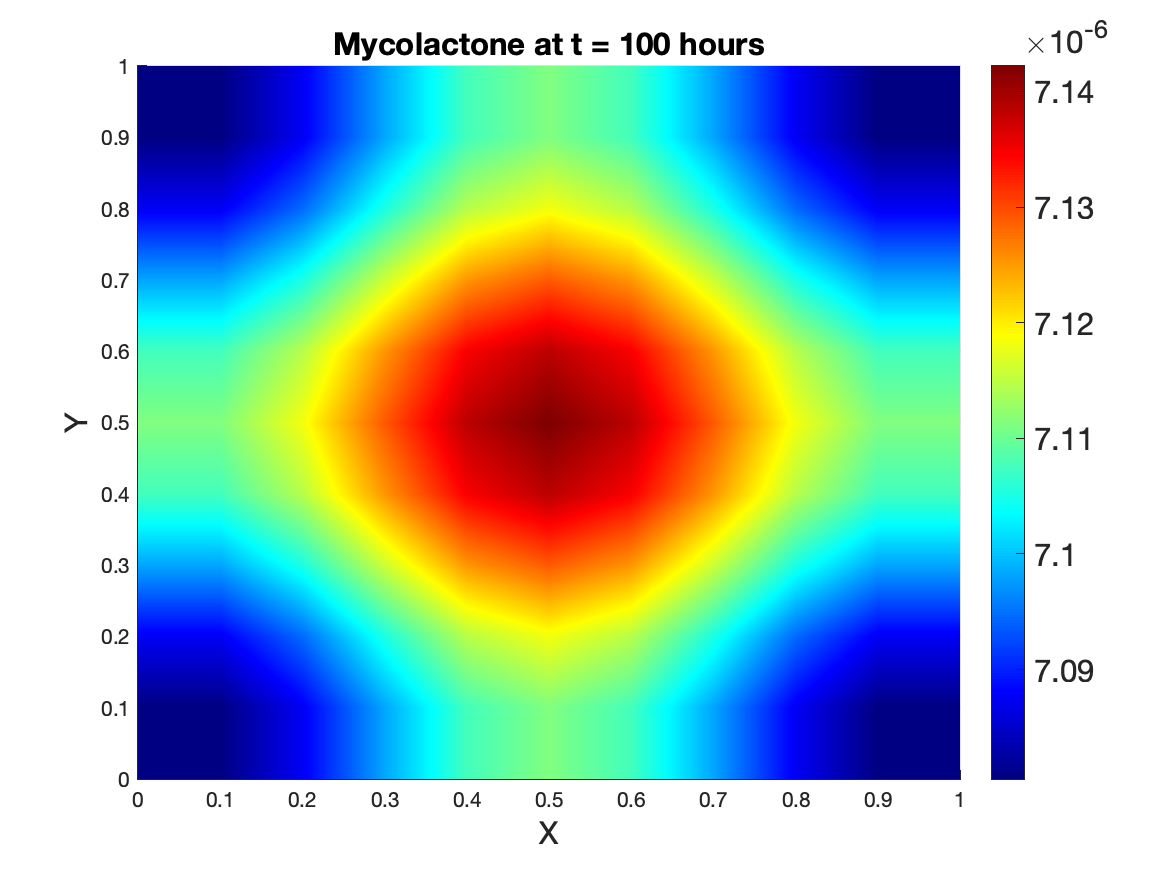}
		\caption{Mycolactone at \\ $t=100$}
		
	\end{subfigure}
	\begin{subfigure}{0.24\textwidth} 		\centering
		\includegraphics[width=\textwidth]{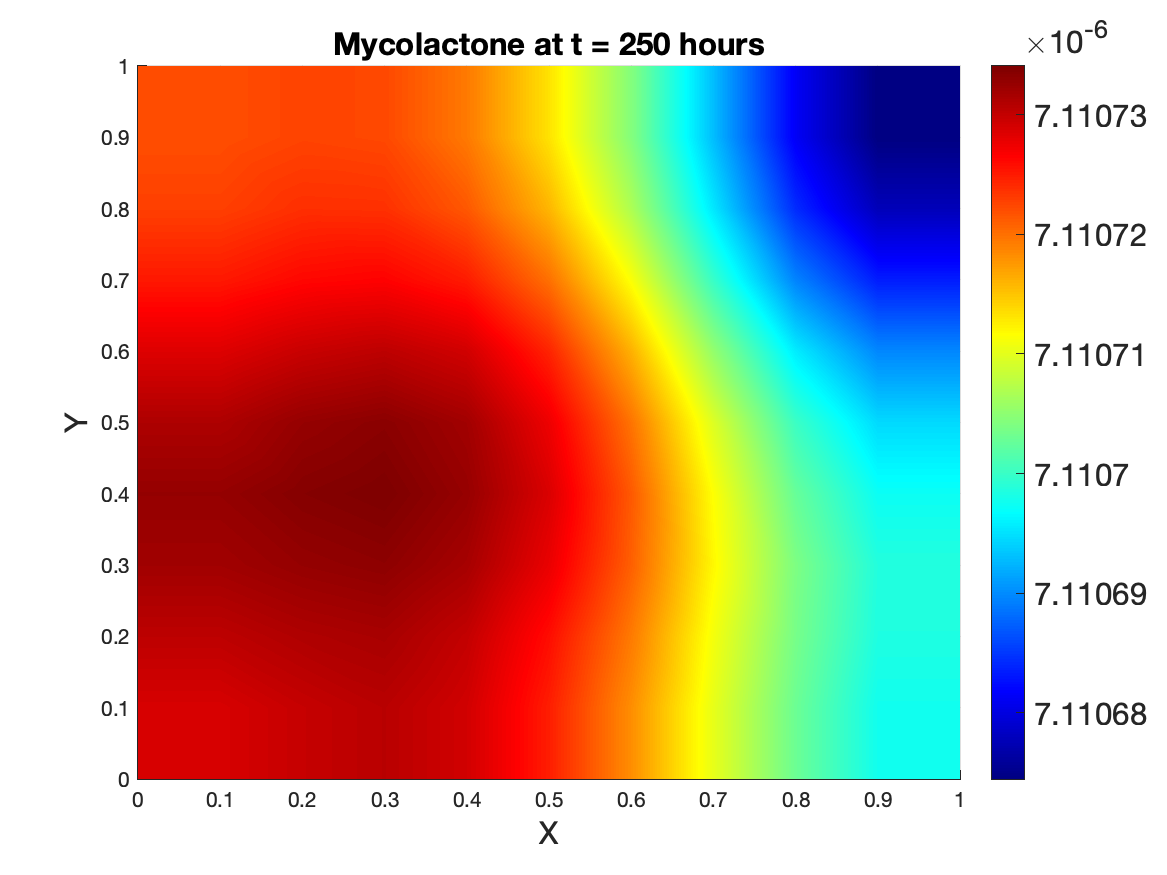}
		\caption{Mycolactone at \\ $t=250$}
	\end{subfigure}     \vspace{0.5cm}     \\
	
	\begin{subfigure}{0.24\textwidth} 		\centering
		\includegraphics[width=\textwidth]{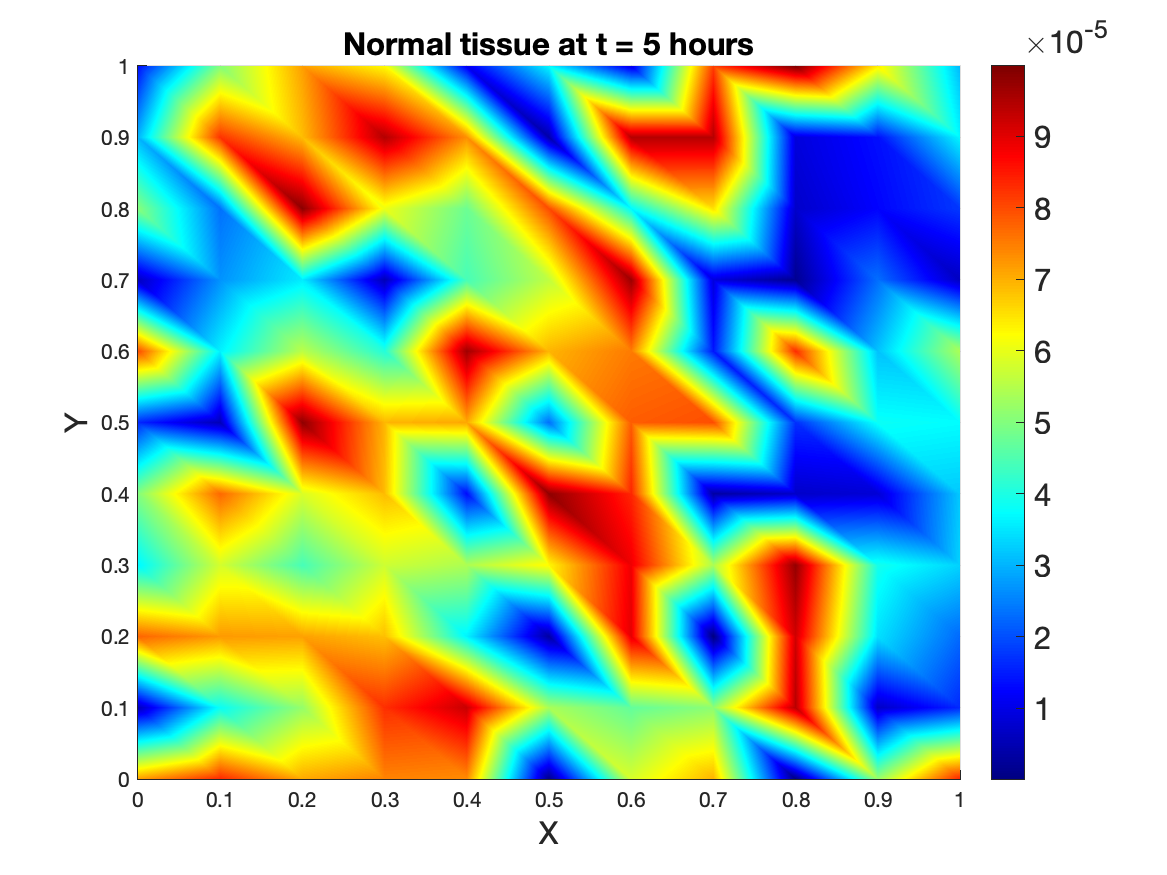}
		\caption{Normal tissue at \\ $t=5$}
		\label{}
	\end{subfigure} 
	\begin{subfigure}{0.24\textwidth} 		\centering
		\includegraphics[width=\textwidth]{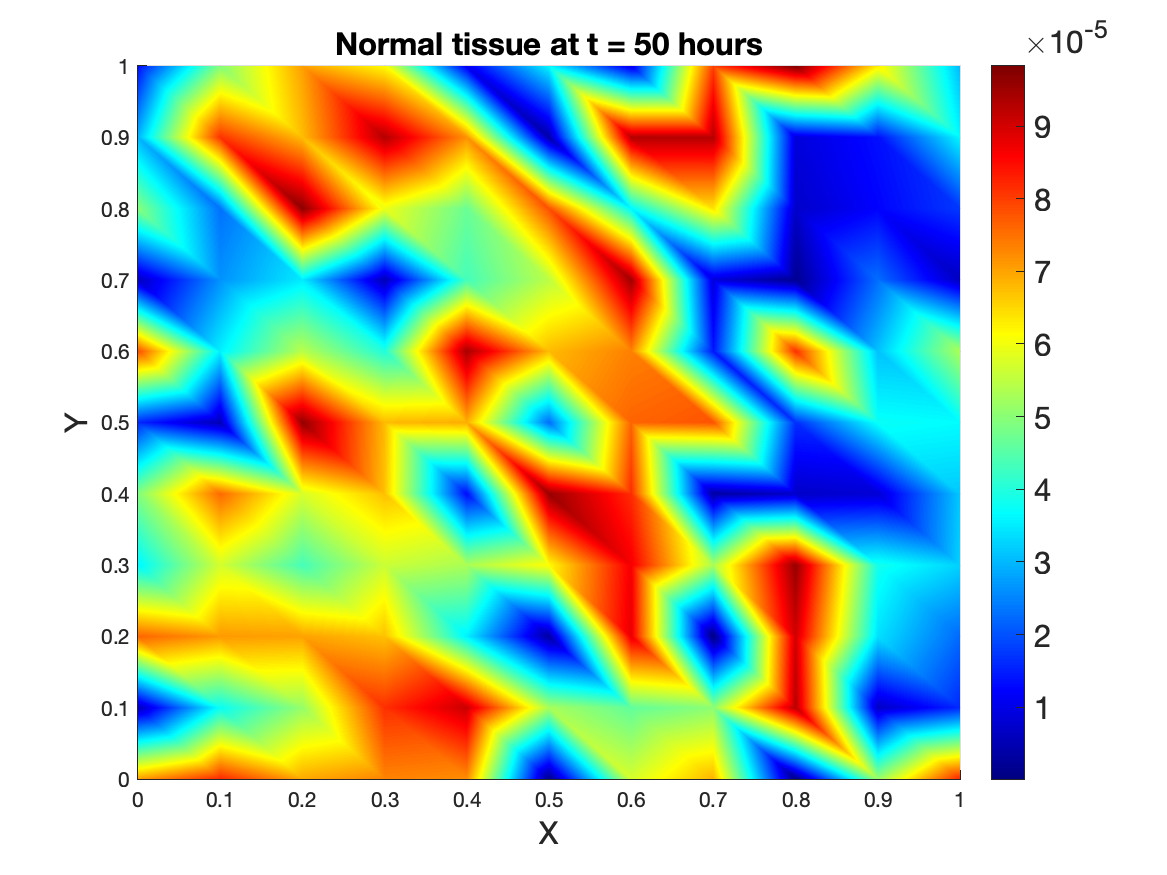}
		\caption{Normal tissue at \\ $t=50$}
	\end{subfigure}
	\begin{subfigure}{0.24\textwidth} 		\centering
		\includegraphics[width=\textwidth]{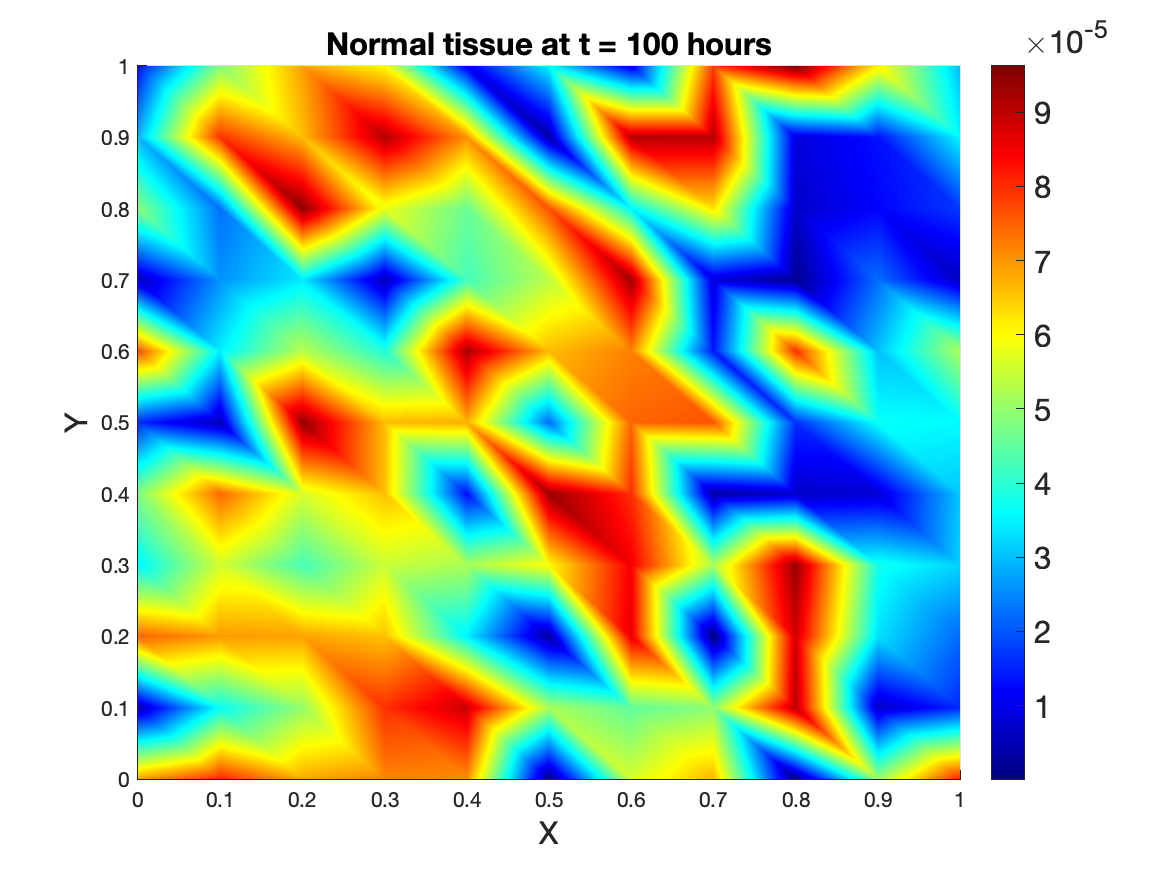}
		\caption{Normal tissue at \\ $t=100$}
	\end{subfigure}
	\begin{subfigure}{0.24\textwidth} 		\centering
		\includegraphics[width=\textwidth]{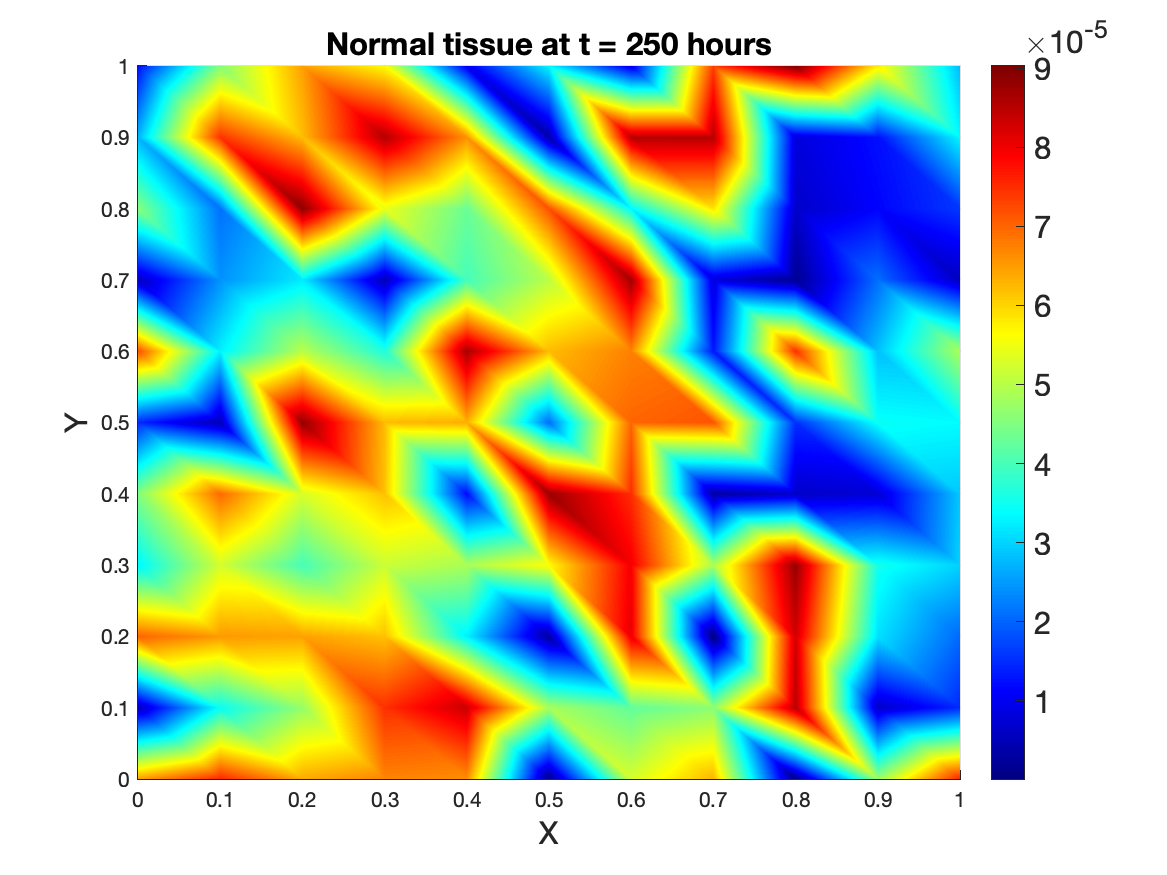}
		\caption{Normal tissue at \\ $t=250$}
	\end{subfigure}\\
	
	\begin{subfigure}{0.24\textwidth} 		\centering
		\includegraphics[width=\textwidth]{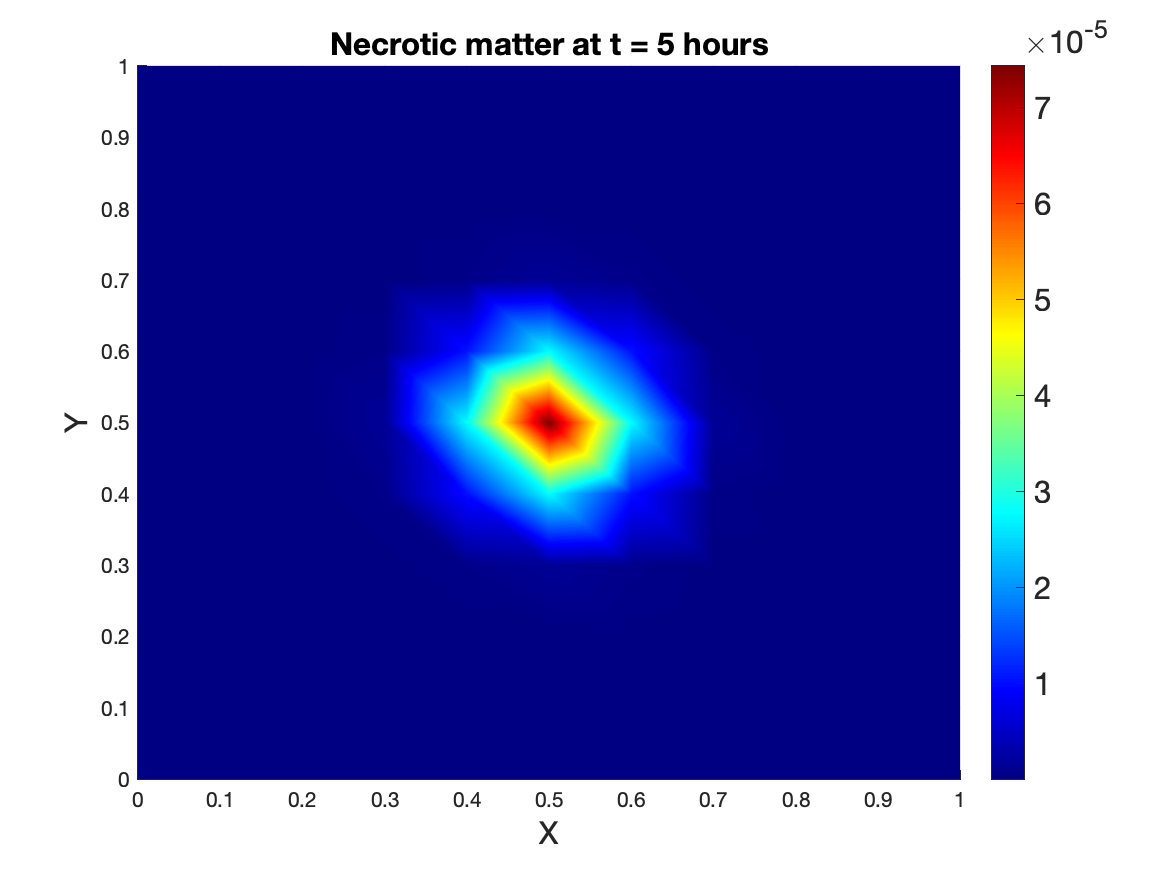}
		\caption{Necrotic matter \\ at $t=5$}
		\label{}
	\end{subfigure}  	
	\begin{subfigure}{0.24\textwidth} 		\centering
		\includegraphics[width=\textwidth]{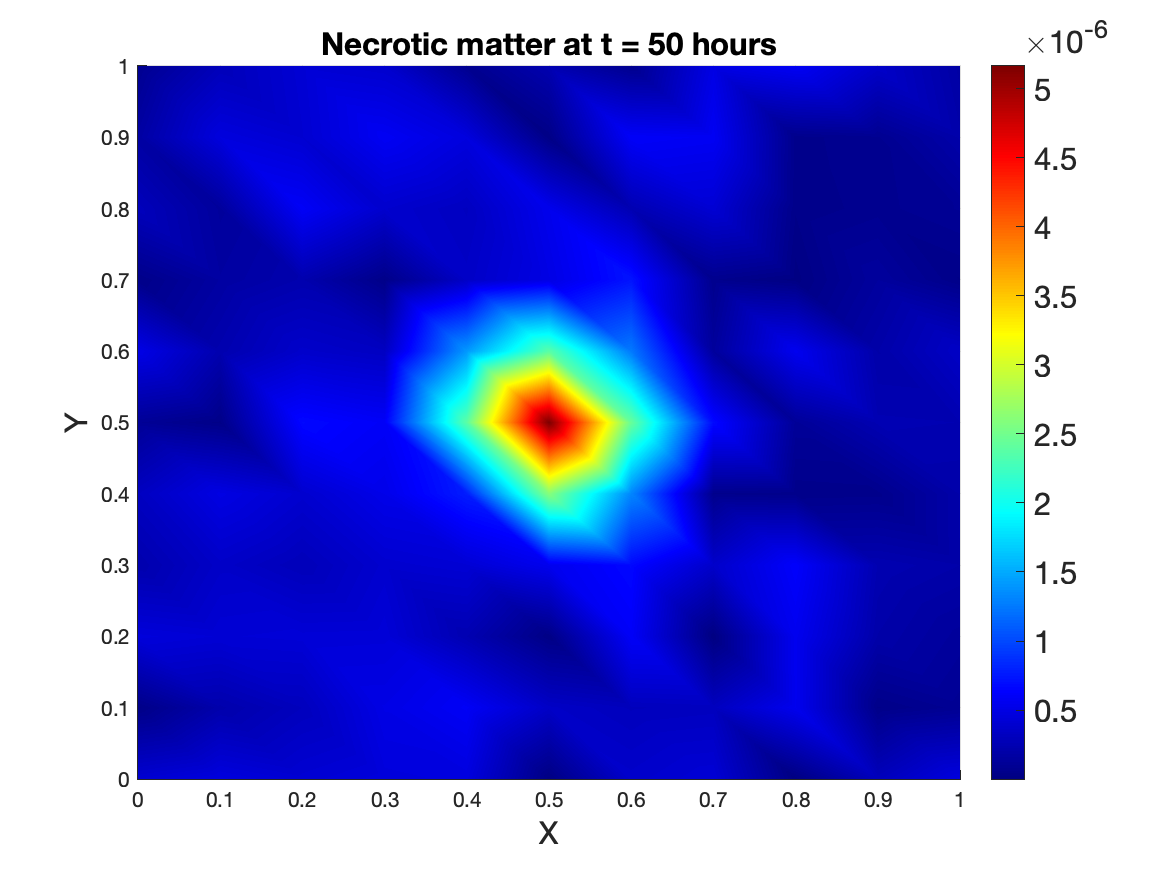}
		\caption{Necrotic matter at \\ $t=50$}
		
	\end{subfigure} 
	\begin{subfigure}{0.24\textwidth} 		\centering
		\includegraphics[width=\textwidth]{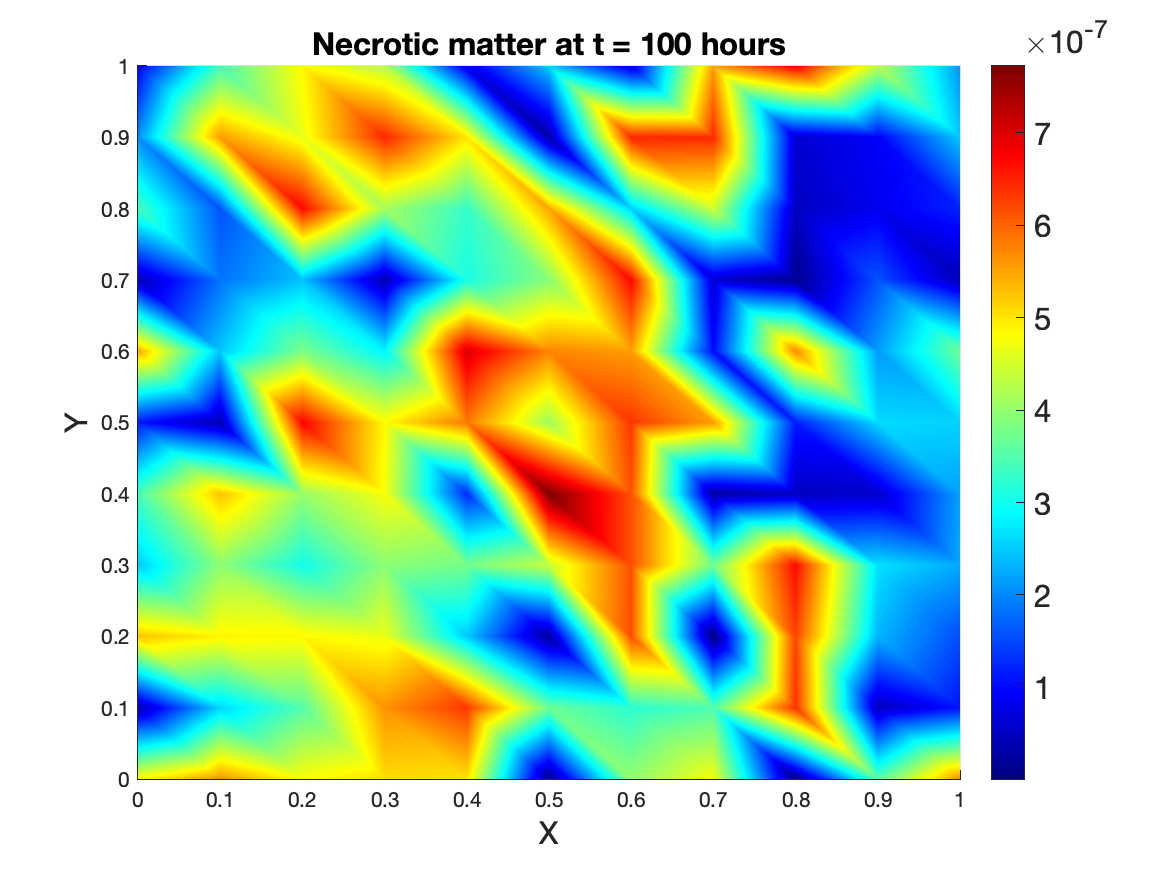}
		\caption{Necrotic matter at \\ $t=100$}
		
	\end{subfigure} 
	\begin{subfigure}{0.24\textwidth} 		\centering
		\includegraphics[width=\textwidth]{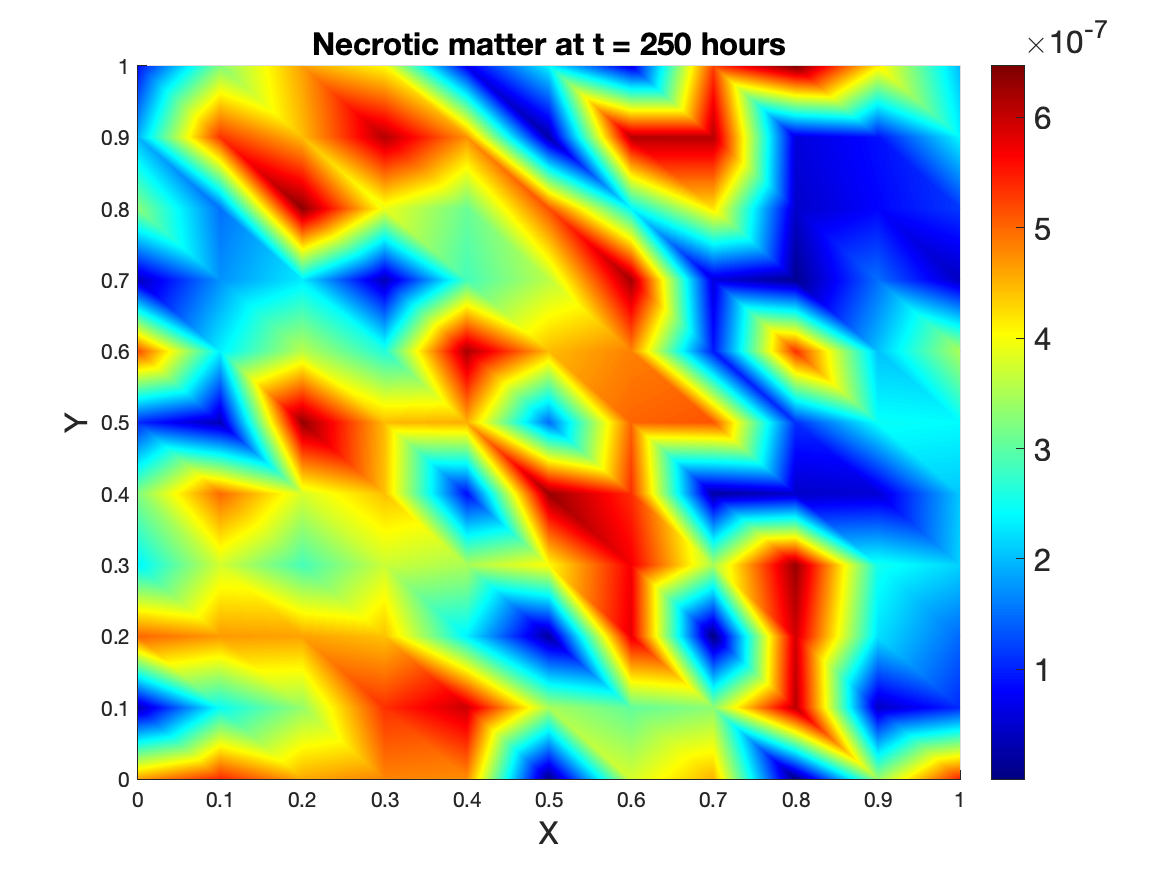}
		\caption{Necrotic matter at \\ $t=250$}
	\end{subfigure}     \vspace{0.5cm}     \\
	\caption{ Scenario 5: Bacteria, mycolactone, normal tissue, and necrotic matter at different times for less initial density of normal tissue.}
\end{figure}

\noindent In contrast to Scenarios 1-4, where abundant normal tissue fosters bacterial proliferation, an intriguing pattern emerges when initial normal tissue levels are insufficient. Instead of growth, bacteria undergo degradation over time, accompanied by increased diffusion. This decline in bacteria shows the importance of sufficiently available normal tissue at the beginning of the process. In Scenario 5,  where normal tissue is lacking, bacteria proliferation is hindered, leading to a gradual reduction in cell density over time. While proliferation is subdued for bacteria, diffusion becomes pronounced due to the absence of adequate normal tissue or necrotic matter to stimulate taxis.\\[-2ex]

\noindent
The diminished population of bacteria results in a corresponding decrease in mycolactone levels. Nevertheless, the ongoing necrotization of normal tissue continues, however at a rather slower pace compared to scenarios with sufficient initial tissue levels. This observation highlights the close interrelation between dynamics of bacteria, normal tissue, and mycolactone production. \\[-2ex]

\begin{figure}[htbp!]
	\centering
	
	\includegraphics[width=45mm,scale=0.5]{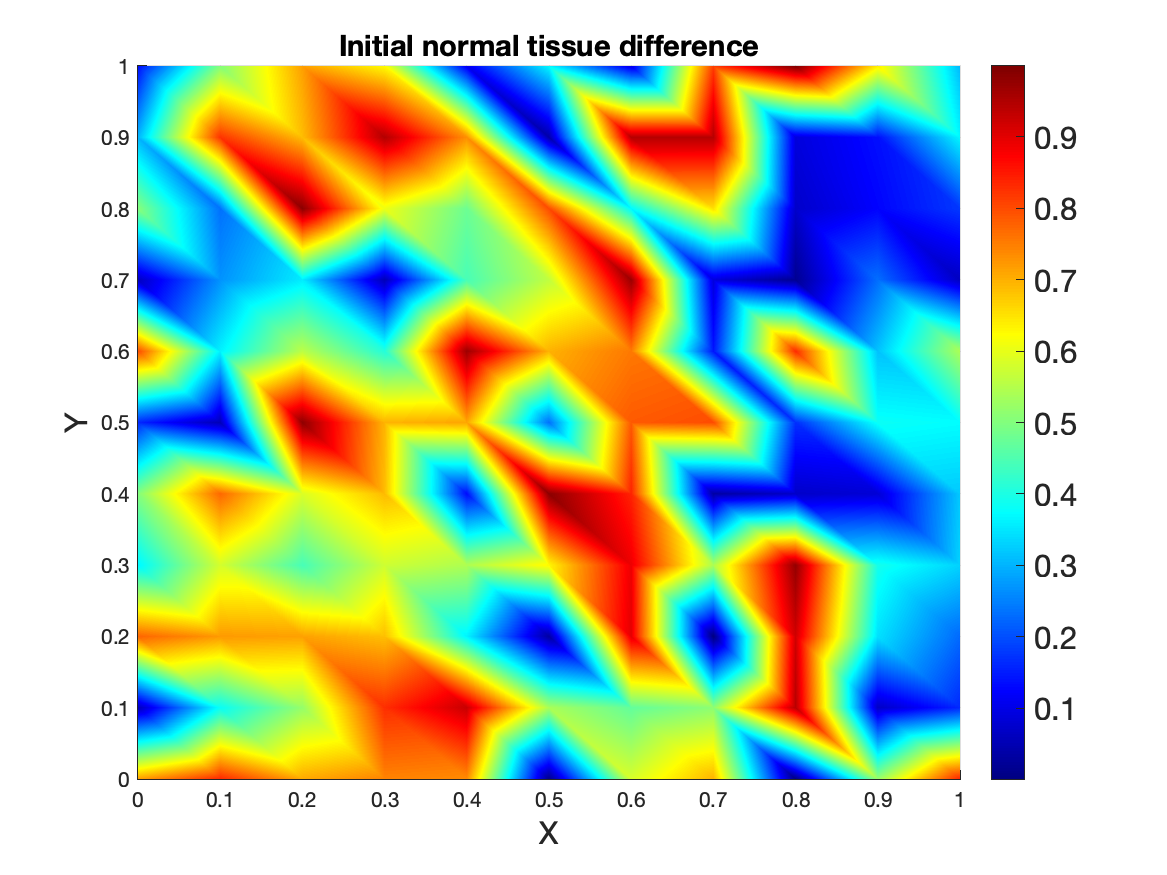}
	
	\caption{ Difference between the initial conditions for normal tissue, in Scenarios 1 and 5.}
	\label{fig:IC-comparison5and1}
\end{figure}

\begin{figure}[htbp!]
	
	\begin{subfigure}{0.24\textwidth} 		\centering
		\includegraphics[width=\textwidth]{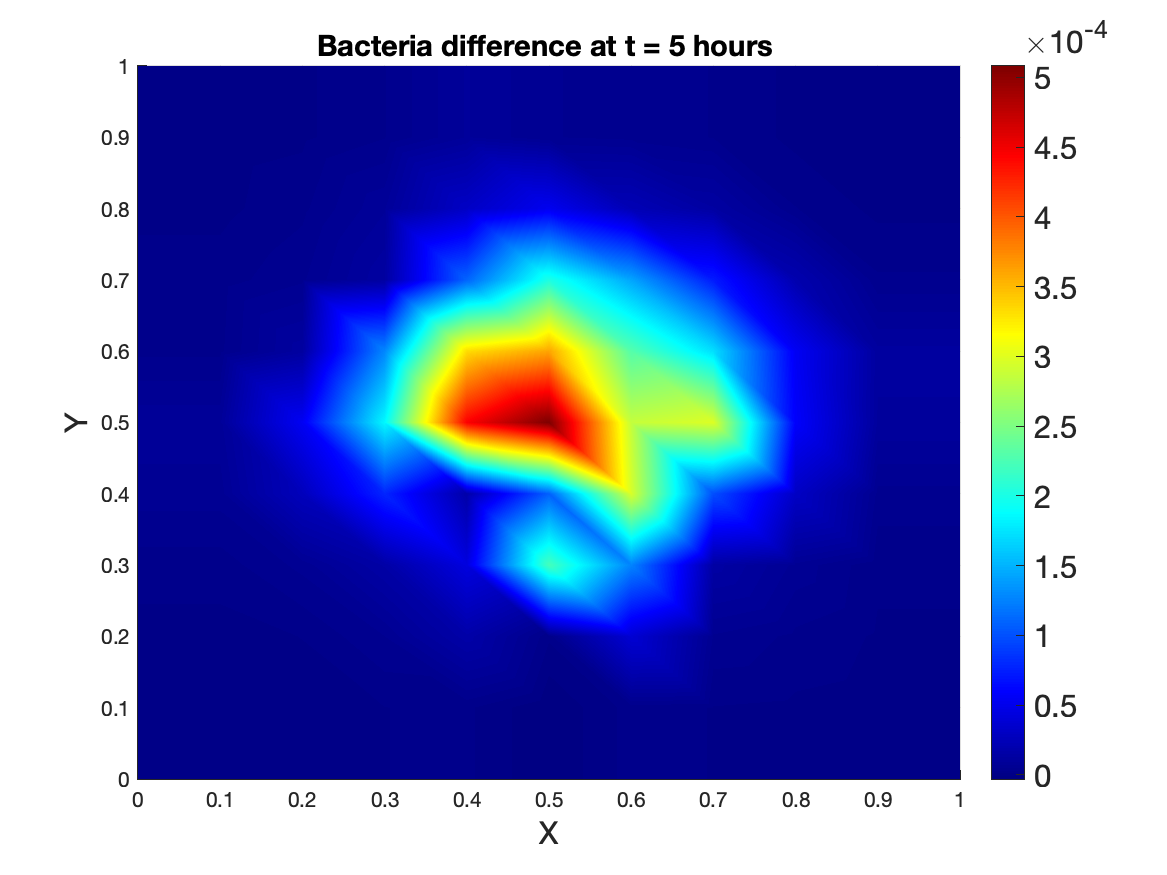}
		\caption{Bacteria at $t=5$}
		\label{}
		
	\end{subfigure} 
	\begin{subfigure}{0.24\textwidth} 		\centering
		\includegraphics[width=\textwidth]{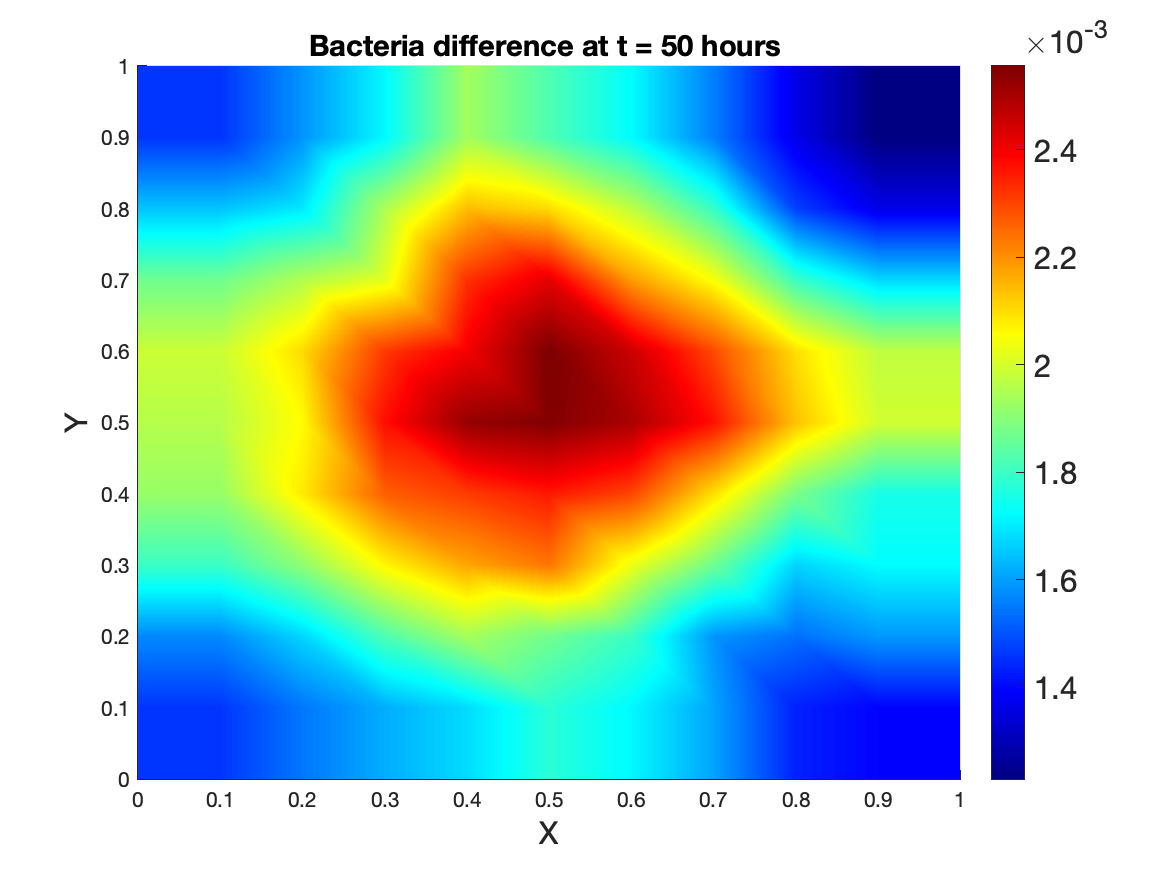}
		\caption{Bacteria at $t=50$}
		
	\end{subfigure} 
	\begin{subfigure}{0.24\textwidth} 		\centering
		\includegraphics[width=\textwidth]{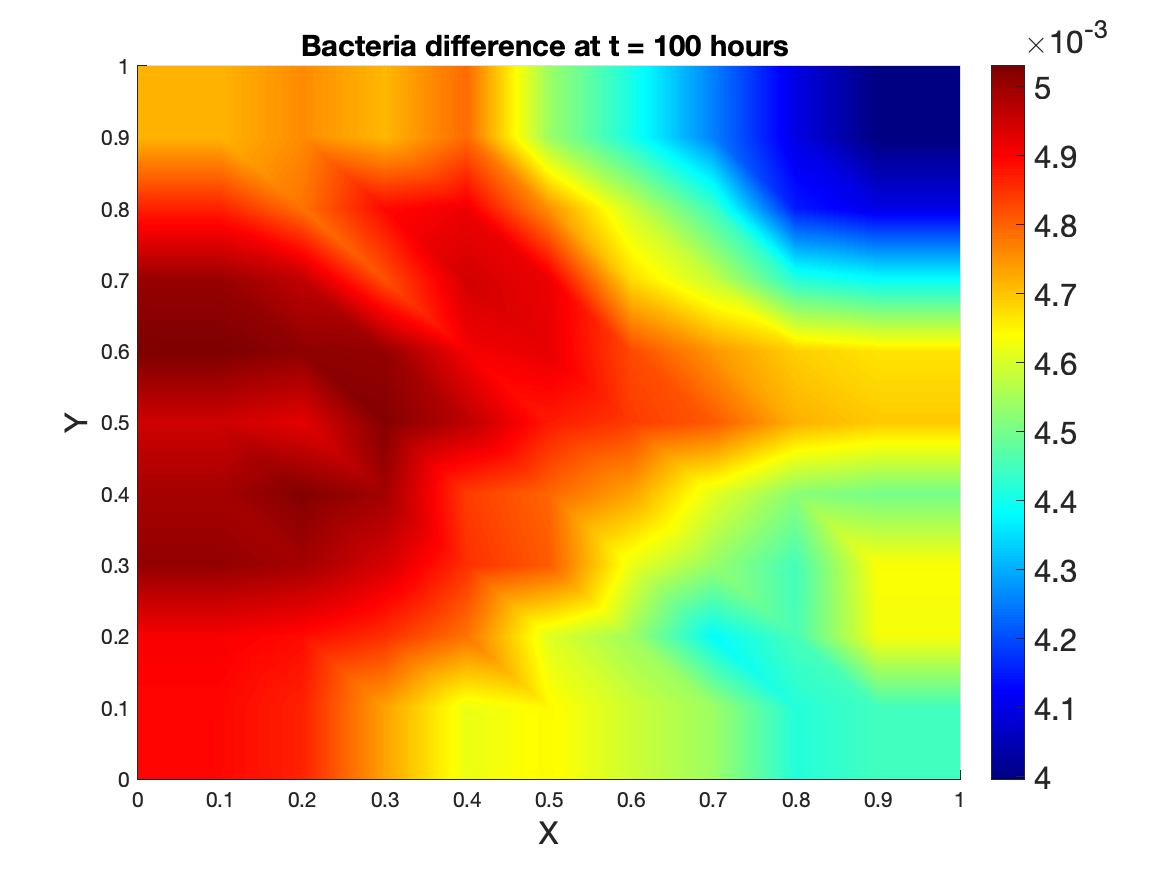}
		\caption{Bacteria at $t=100$}
		
	\end{subfigure} 
	\begin{subfigure}{0.24\textwidth} 		\centering
		\includegraphics[width=\textwidth]{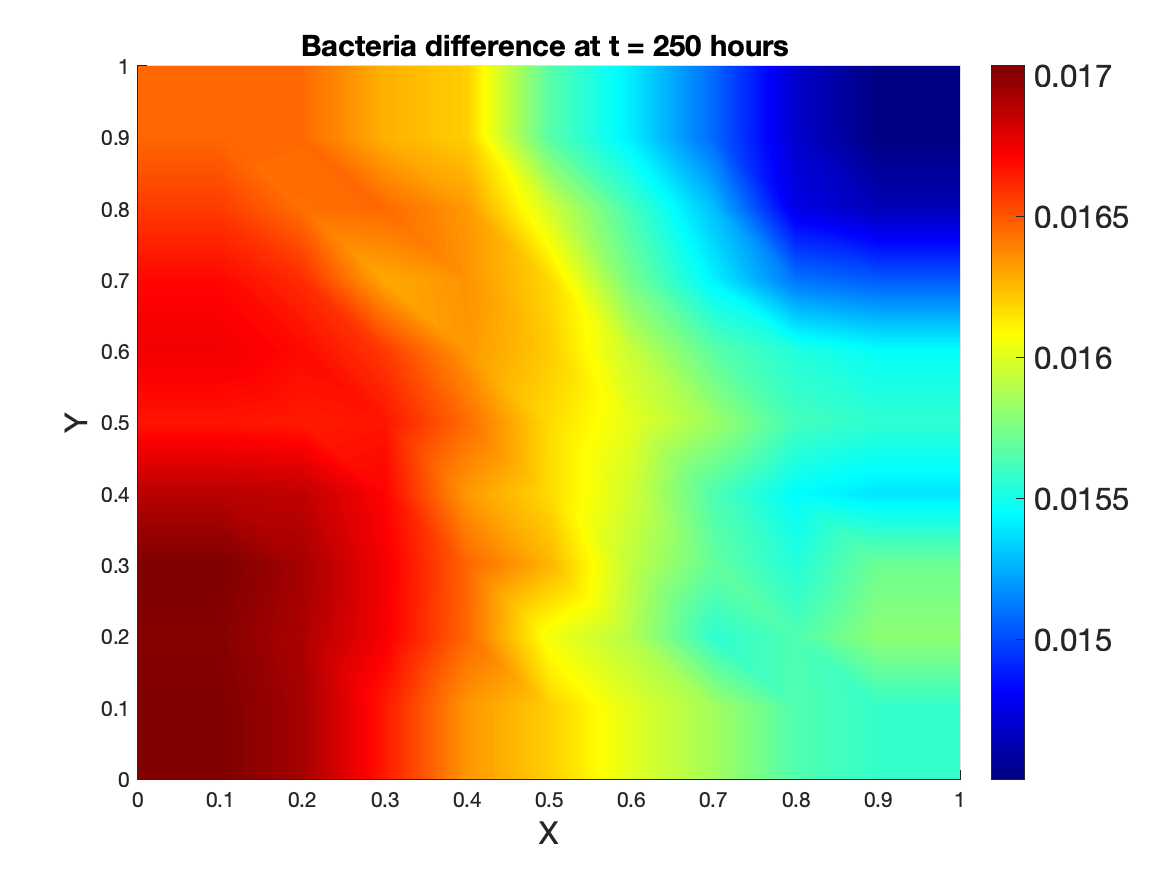}
		\caption{Bacteria at $t=250$}
		
	\end{subfigure}     \vspace{0.5cm}     \\
	
	\begin{subfigure}{0.24\textwidth} 		\centering
		\includegraphics[width=\textwidth]{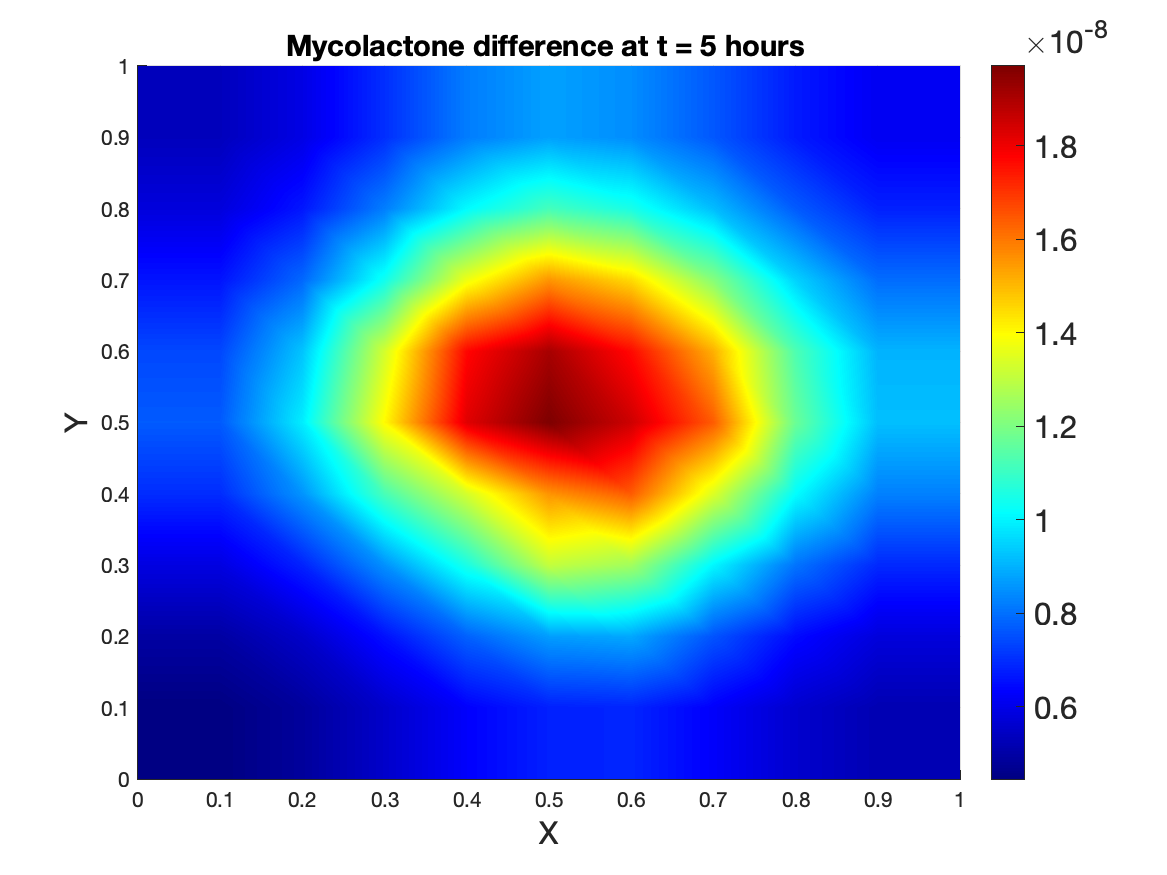}
		\caption{Mycolactone at \\ $t=5$}
		\label{}
	\end{subfigure} 
	\begin{subfigure}{0.24\textwidth} 		\centering 
		\includegraphics[width=\textwidth]{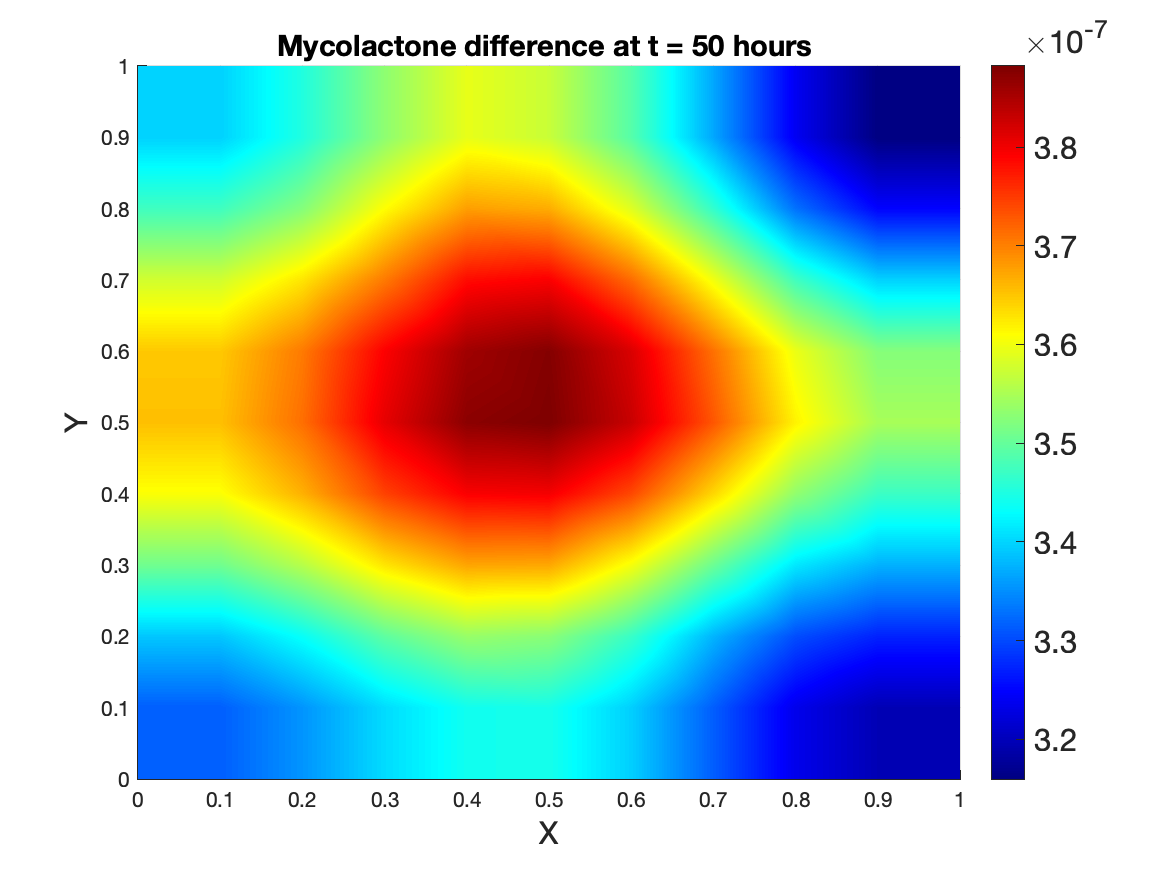}
		\caption{Mycolactone at \\ $t=50$}
		
	\end{subfigure} 
	\begin{subfigure}{0.24\textwidth} 		\centering 
		\includegraphics[width=\textwidth]{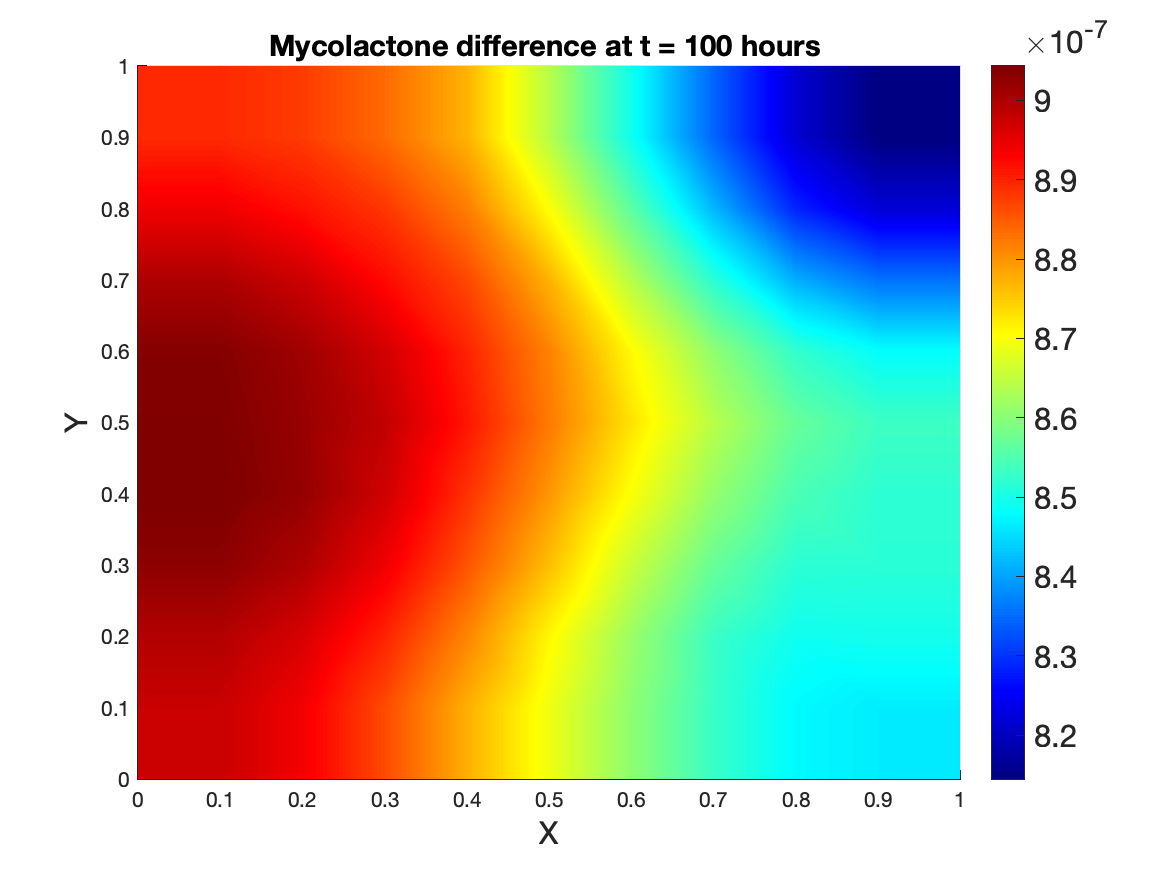}
		\caption{Mycolactone at \\ $t=100$}
		
	\end{subfigure}
	\begin{subfigure}{0.24\textwidth} 		\centering
		\includegraphics[width=\textwidth]{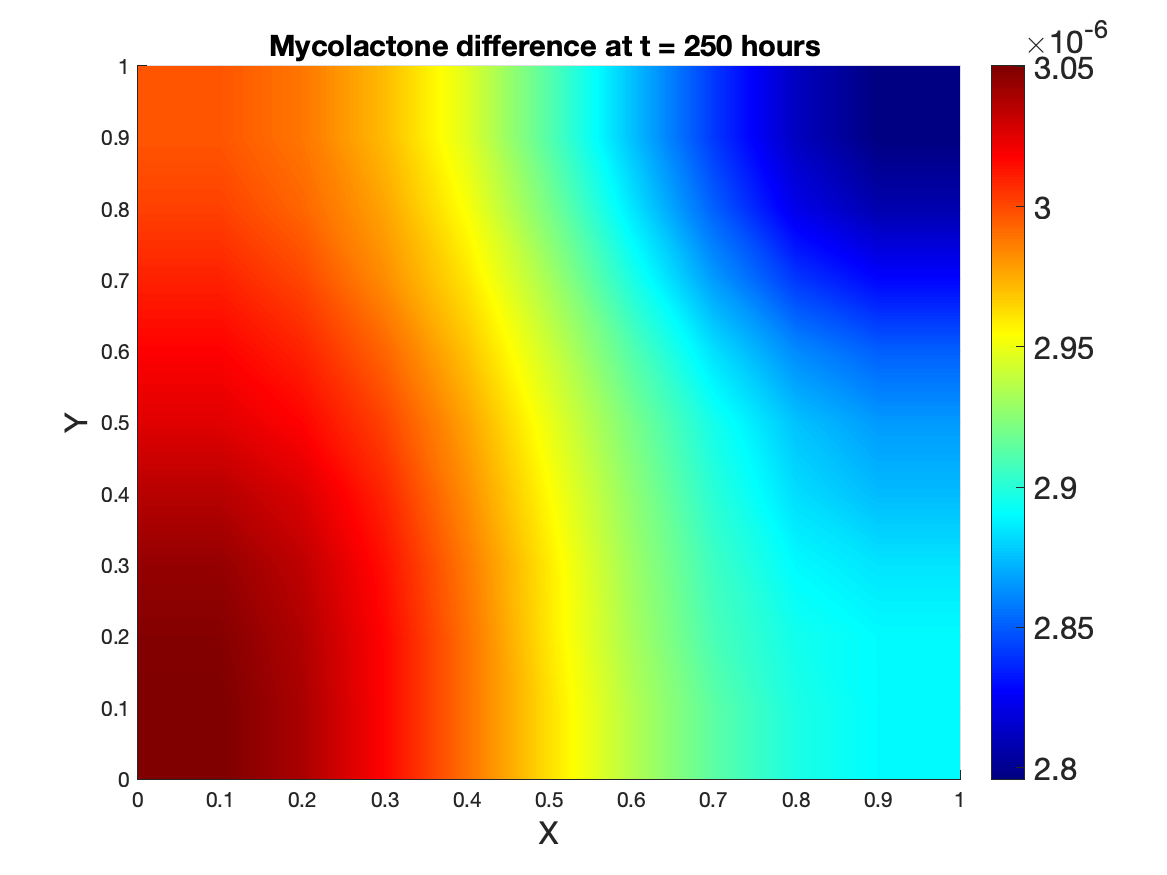}
		\caption{Mycolactone at \\ $t=250$}
	\end{subfigure}     \vspace{0.5cm}     \\
	
	\begin{subfigure}{0.24\textwidth} 		\centering
		\includegraphics[width=\textwidth]{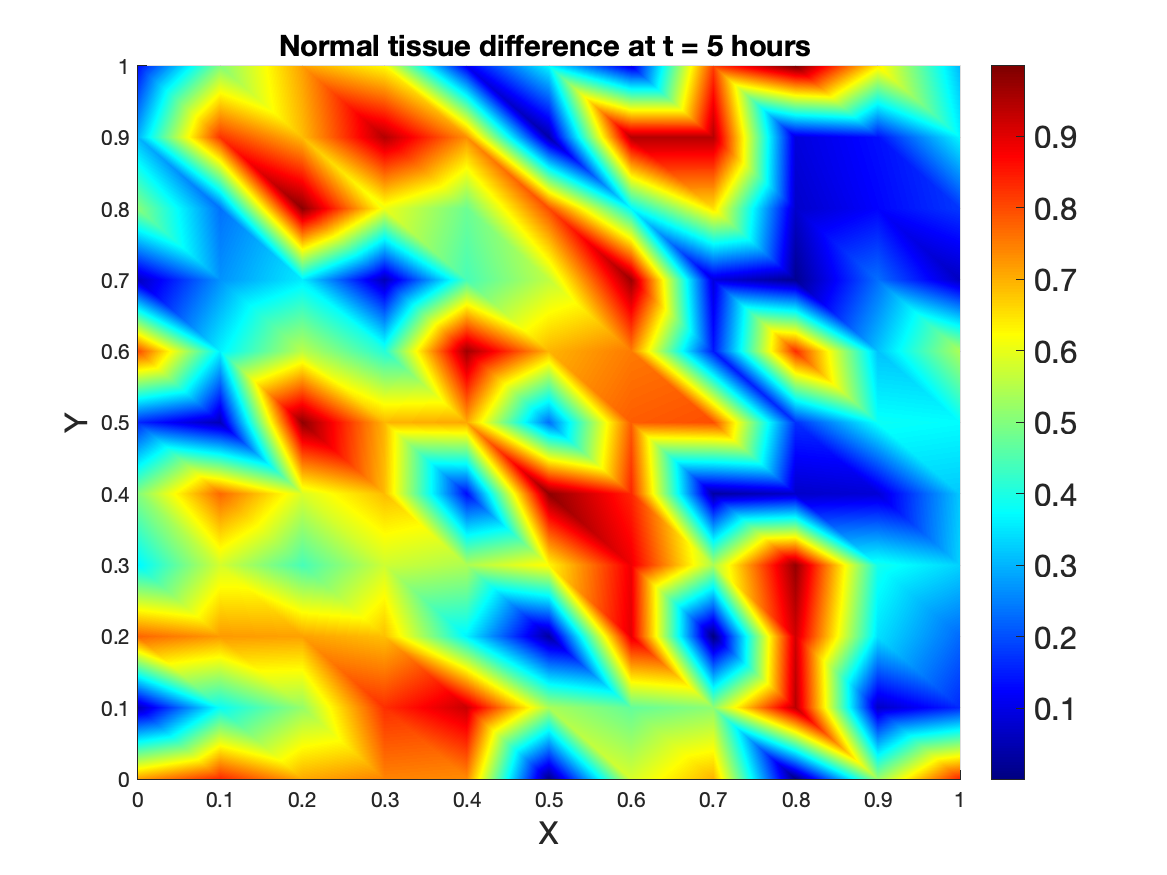}
		\caption{Normal tissue at \\ $t=5$}
		\label{}
	\end{subfigure} 
	\begin{subfigure}{0.24\textwidth} 		\centering
		\includegraphics[width=\textwidth]{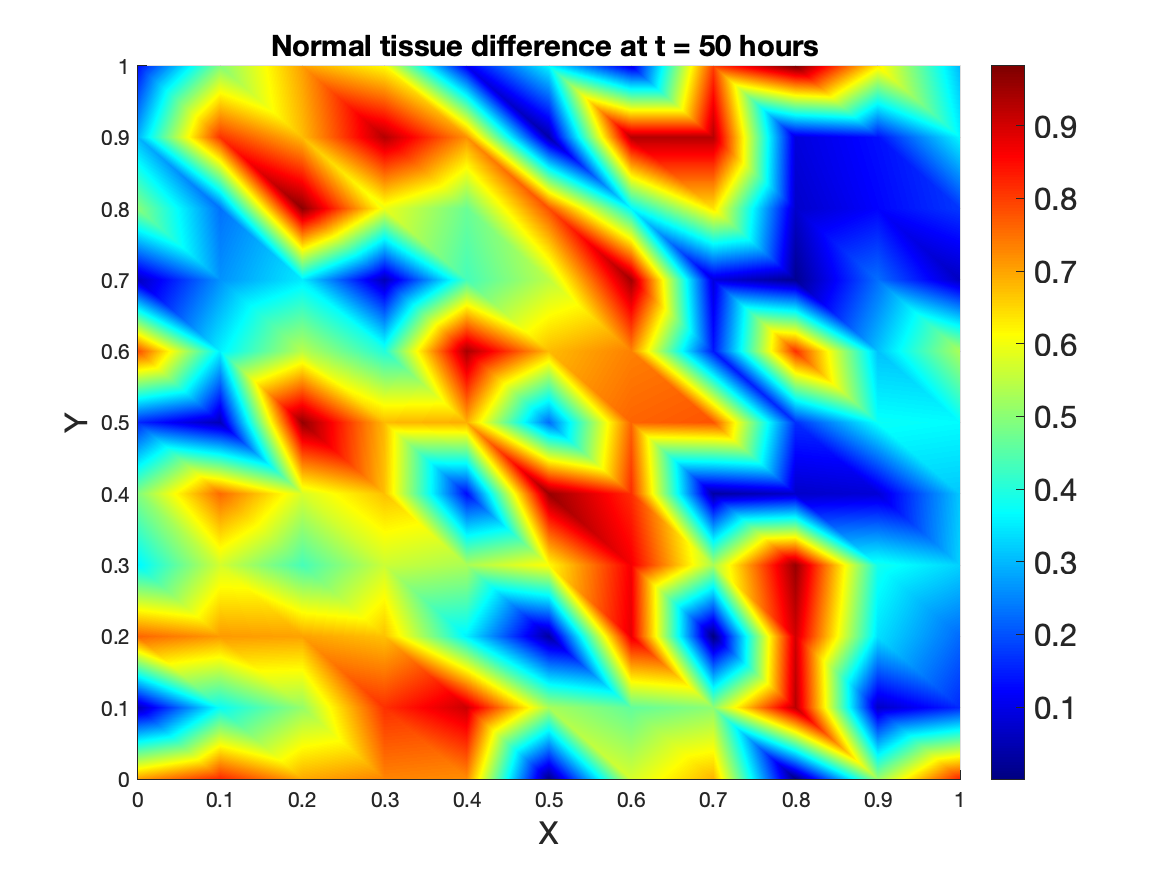}
		\caption{Normal tissue at \\ $t=50$}
	\end{subfigure}
	\begin{subfigure}{0.24\textwidth} 		\centering
		\includegraphics[width=\textwidth]{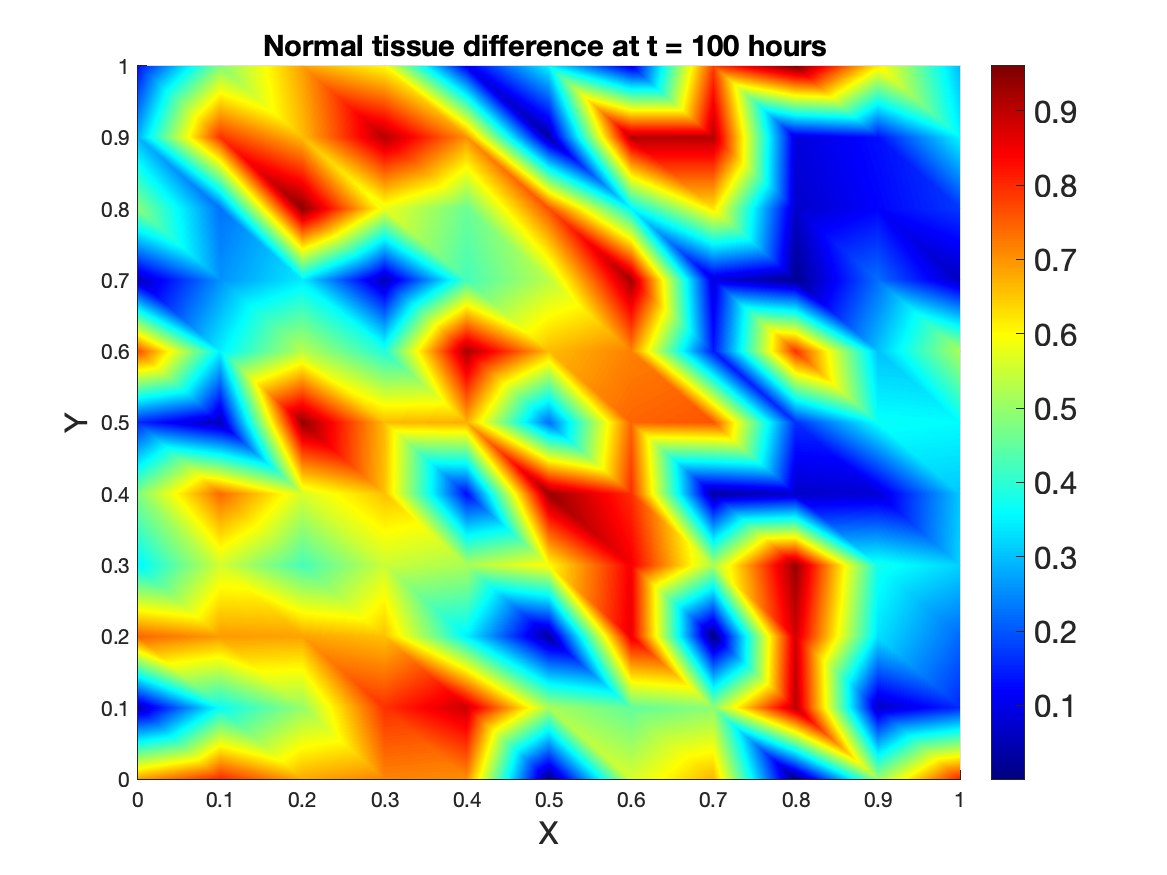}
		\caption{Normal tissue at \\ $t=100$}
	\end{subfigure}
	\begin{subfigure}{0.24\textwidth} 		\centering
		\includegraphics[width=\textwidth]{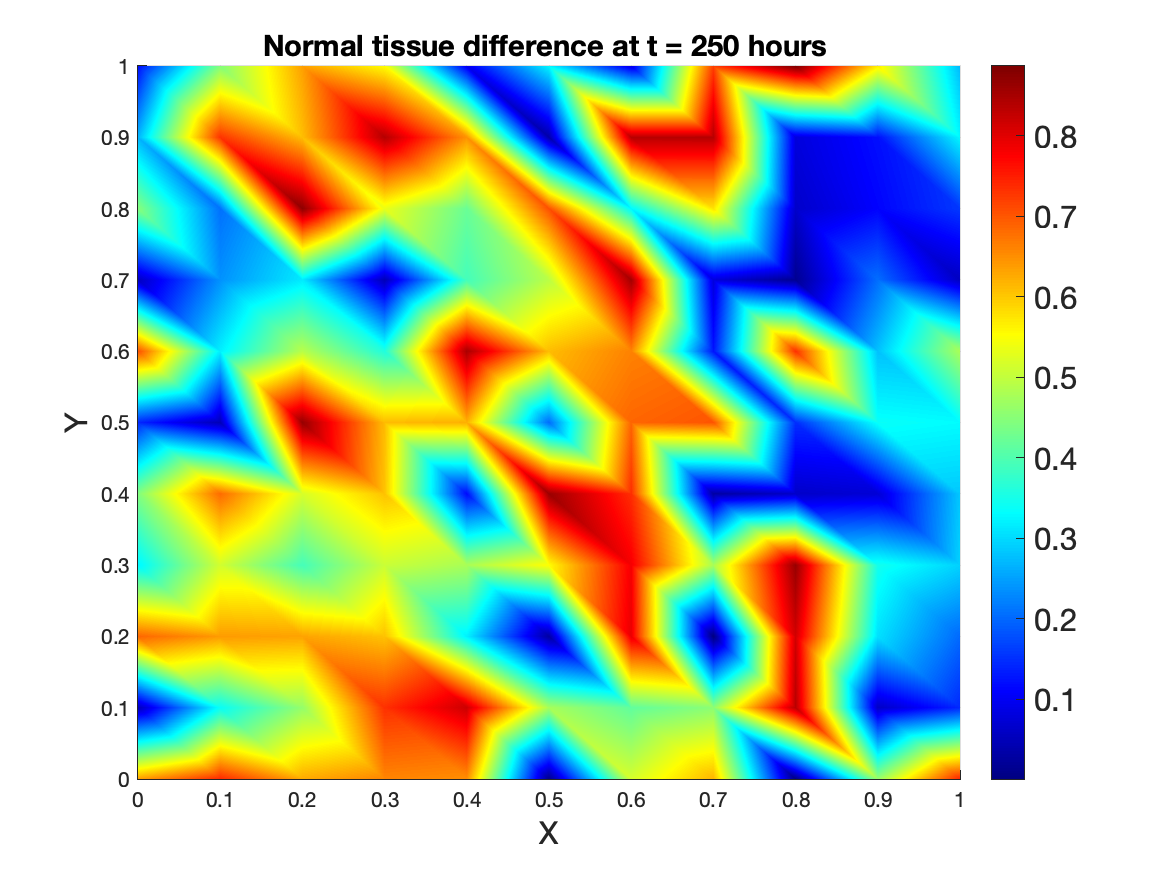}
		\caption{Normal tissue at \\ $t=250$}
	\end{subfigure}\\
	
	\begin{subfigure}{0.24\textwidth} 		\centering
		\includegraphics[width=\textwidth]{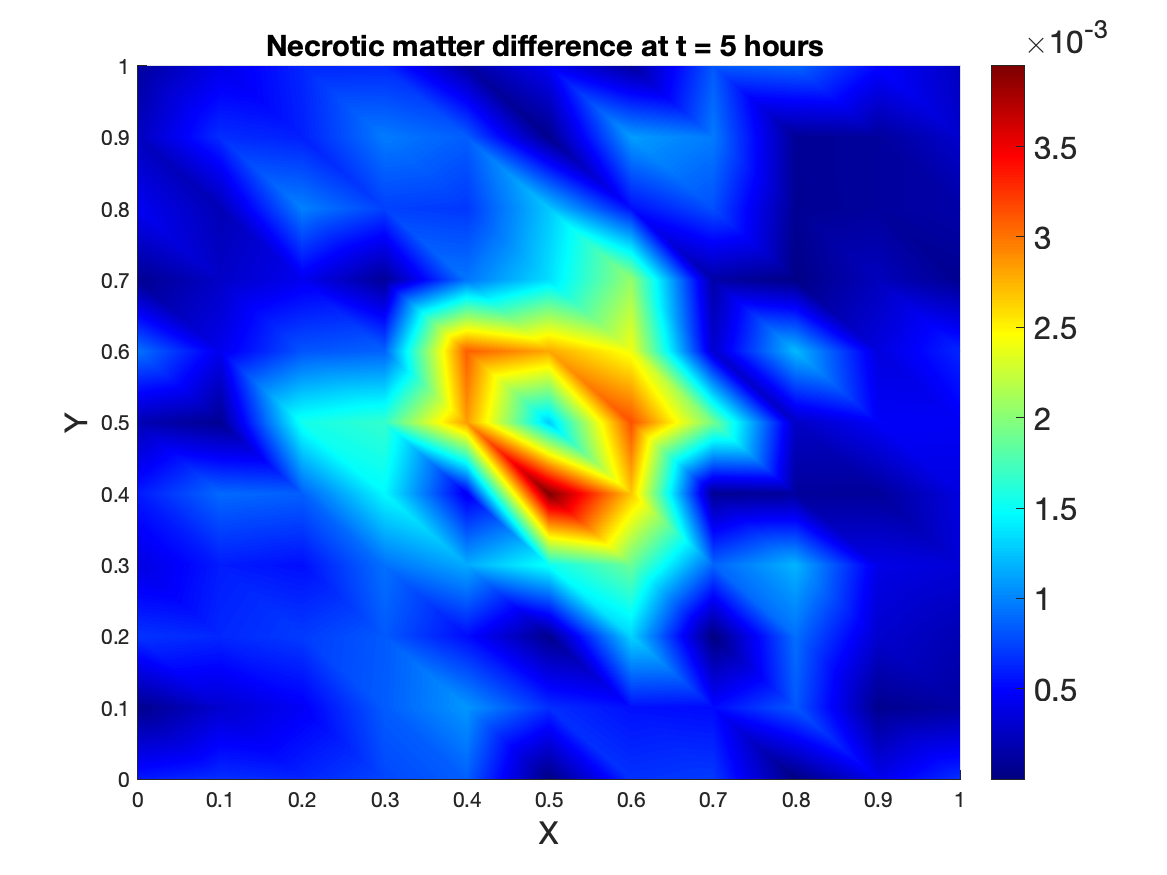}
		\caption{Necrotic matter \\ at $t=5$}
		\label{}
	\end{subfigure}  	
	\begin{subfigure}{0.24\textwidth} 		\centering
		\includegraphics[width=\textwidth]{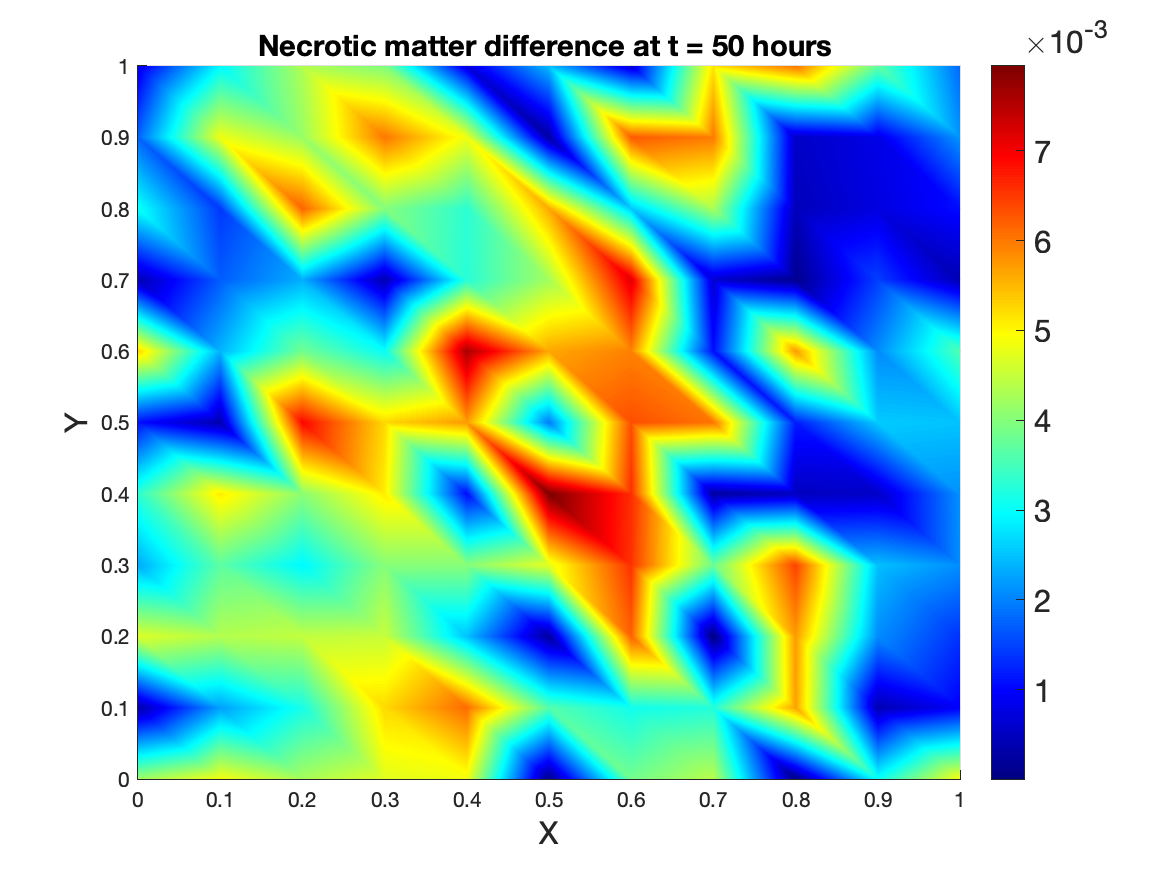}
		\caption{Necrotic matter at \\ $t=50$}
		
	\end{subfigure} 
	\begin{subfigure}{0.24\textwidth} 		\centering
		\includegraphics[width=\textwidth]{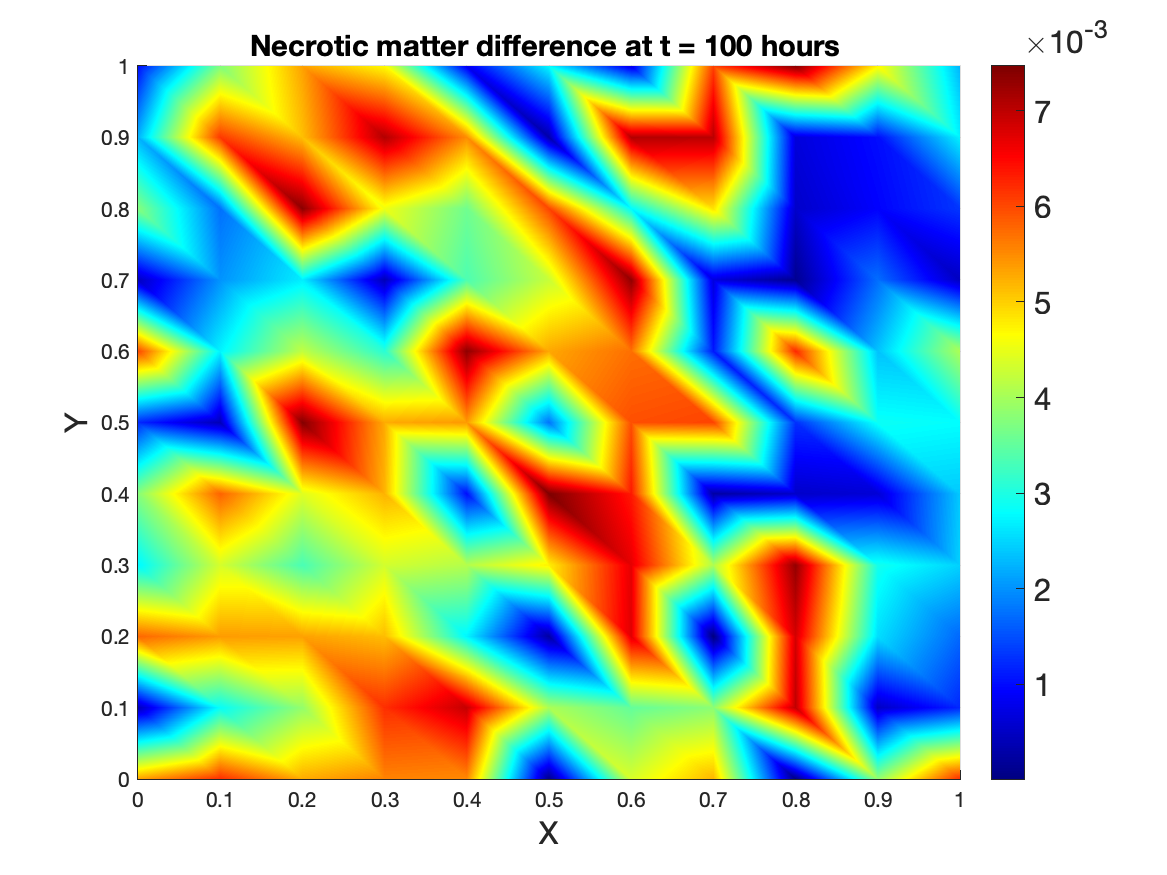}
		\caption{Necrotic matter at \\ $t=100$}
		
	\end{subfigure} 
	\begin{subfigure}{0.24\textwidth} 		\centering
		\includegraphics[width=\textwidth]{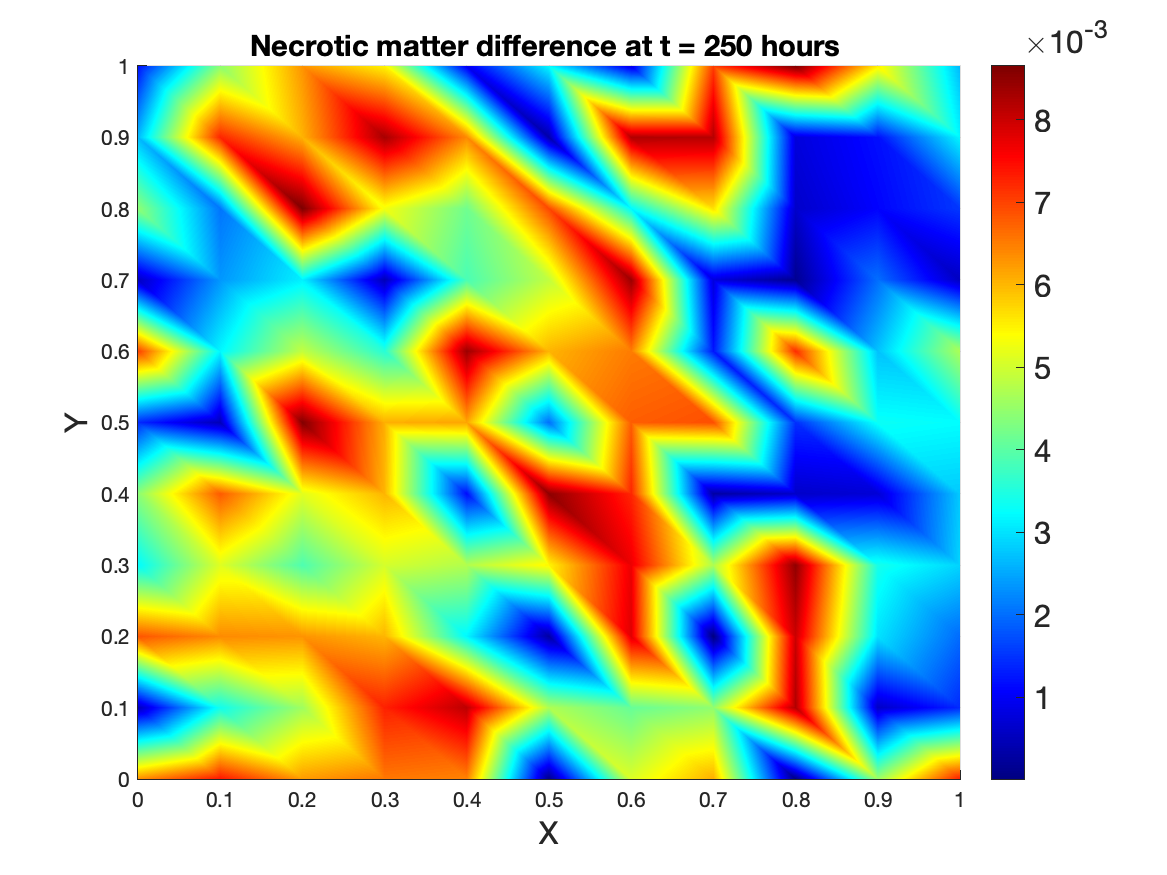}
		\caption{Necrotic matter at \\ $t=250$}
	\end{subfigure}     \vspace{0.5cm}     \\
	\caption{ Difference between densities of bacteria, mycolactone, normal tissue, and necrotic matter at different times for Scenario 1 and Scenario 5 (Scenario 1 - Scenario 5).}
	\label{fig:comparison3}
\end{figure}

\noindent
We also provide the comparison of Scenario 1 with Scenario 5 in Figure \ref{fig:comparison3}. The difference in the initial density of normal tissue is given in Figure \ref{fig:IC-comparison5and1}.\\[-2ex]

\noindent
In the comparison between Scenario 5 and Scenario 1, all entities (bacteria, mycolactone, normal tissue, and necrotic tissue) exhibit positive differences in regions of higher concentrations. This results from starting with a lower initial amount of normal tissue in Scenario 5, leading to reduced growth of all other entities compared to Scenario 1, which begins with an adequate amount of normal tissue. This highlights the significant impact of initial conditions on the progression and dynamics of bacteria infection and tissue response.\\[-2ex]

\noindent
From the scenarios discussed above, it becomes evident that the tactic sensitivity of bacteria towards normal tissue, along with the initial density of normal tissue, significantly influence the progression of bacterial colonisation and the development of Buruli ulcer. 
For reduced tactic sensitivity (as in the case of $\gamma_1>\gamma_2$), the production and spread of bacteria, as well as other associated entities, are diminished compared to scenarios where the tactic sensitivity $\gamma_2$ is higher. This difference underscores the pivotal role of bacteria motility in dictating the extent of infection spread and subsequent tissue damage.
\\[-2ex]

\noindent
Likewise, lower amounts of normal tissue lead to a distinct outcome: bacterial degradation over time rather than proliferation. This observation highlights the critical importance of reduced normal tissue in slowing bacterial growth and tissue degradation. 
Considering the body's natural ability to produce normal tissue, this suggests a potential avenue for intervention, wherein strategies aimed at making the normal tissue unfavourable e.g., by intensely heating the affected areas could mitigate the progression of Buruli ulcer.

\subsection{Simulations of system \eqref{burulinonlinear}}

\noindent
In this subsection we use simulations to investigate the effect of the modified motility coefficients in \eqref{burulinonlinear}, which were obtained by the position jump approach from Subsection \ref{subsec:PJA}. \\[-2ex]

\noindent
Before proceeding we nondimenionalize \eqref{burulinonlinear} by taking 
\begin{equation}\label{40}\bar{u}=\frac{u}{K_u}, \ \tilde{v}=\frac{v}{K_v}, \ \tilde{n}=\frac{n}{K_n}, \ \tilde{m}=\frac{m}{K_m}, \  \tilde{x}={x}\sqrt{\frac{\alpha_u}{D_m}}, \ \tilde{t}=t \alpha_u.\end{equation}
This leads to 
\begin{equation}
	\begin{aligned}\label{nondimsysburulinonlinear}
		& \bar{u}_{\tilde{t}}=\nabla \cdot \left(\frac{\tilde{D}}{(1+\bar{u}\tilde{v})^2}\nabla \bar{u}\right)-\nabla \cdot \left(\frac{\tilde{D}\bar{u}^2}{(1+\bar{u}\tilde{v})^2}\nabla \tilde{v}\right)-\nabla \cdot \left(\frac{\tilde{\chi}_n\bar{u}}{(1+\tilde{n})^2 (1+un)}\nabla \tilde{n}\right) \\
		&\quad \quad \quad 
		+ \frac{\tilde{n}}{1+\tilde{n}} \bar{u} \left(1-\bar{u}-\tilde{v}-\tilde{n}\right),
		\\[5pt]
		&\tilde{m}_{\tilde{t}}=\Delta \tilde{m}+\frac{\tilde{\delta} \bar{u}}{1+\bar{u}}-\tilde{\lambda} \tilde{m},
		\\[5pt]
		& \tilde{v}_{\tilde{t}}=-\tilde{\beta}_1\tilde{m}\tilde{v},
		\\[5pt]
		& \tilde{n}_{\tilde{t}}=\tilde{\beta}_2\tilde{m}\tilde{v}-\tilde{\gamma}\tilde{n},
	\end{aligned}
\end{equation}
where 
\begin{equation*}
	\begin{aligned} 
		\tilde{D}&=\frac{D}{2*D_m}, \quad \tilde{\chi}_n=\frac{K_1}{D_m}, \quad \tilde{\beta}_1=\frac{{\beta}_1}{\alpha_u}, 
		\\[5pt]
		\tilde{\beta_2} &=\frac{\beta_2 }{\alpha_u K_N},\quad \tilde{\gamma} =\frac{\gamma}{\alpha_u}, \quad \tilde{\delta} =\frac{\delta}{K_m \alpha_u}, \quad  \tilde{\lambda} =\frac{\lambda} {\alpha_u}
	\end{aligned}
\end{equation*}
\noindent
To simplify notations, we omit all tildes and bars in system \eqref{nondimsysburulinonlinear} and continue with these equations.\\[-2ex]

\noindent
For the discretization we use the same method as above, and the same parameters are considered in the implementation. We show in this subsection only the differences between the simulations obtained with system \eqref{14} and those resulting from system 
\eqref{nondimsysburulinonlinear} above, with no-flux boundary conditions and initail conditions as for \eqref{14}.\\[-2ex]

\begin{figure}[htbp!]
	\centering
	\begin{subfigure}{0.24 \textwidth}
		\centering
		\includegraphics[width=\textwidth]{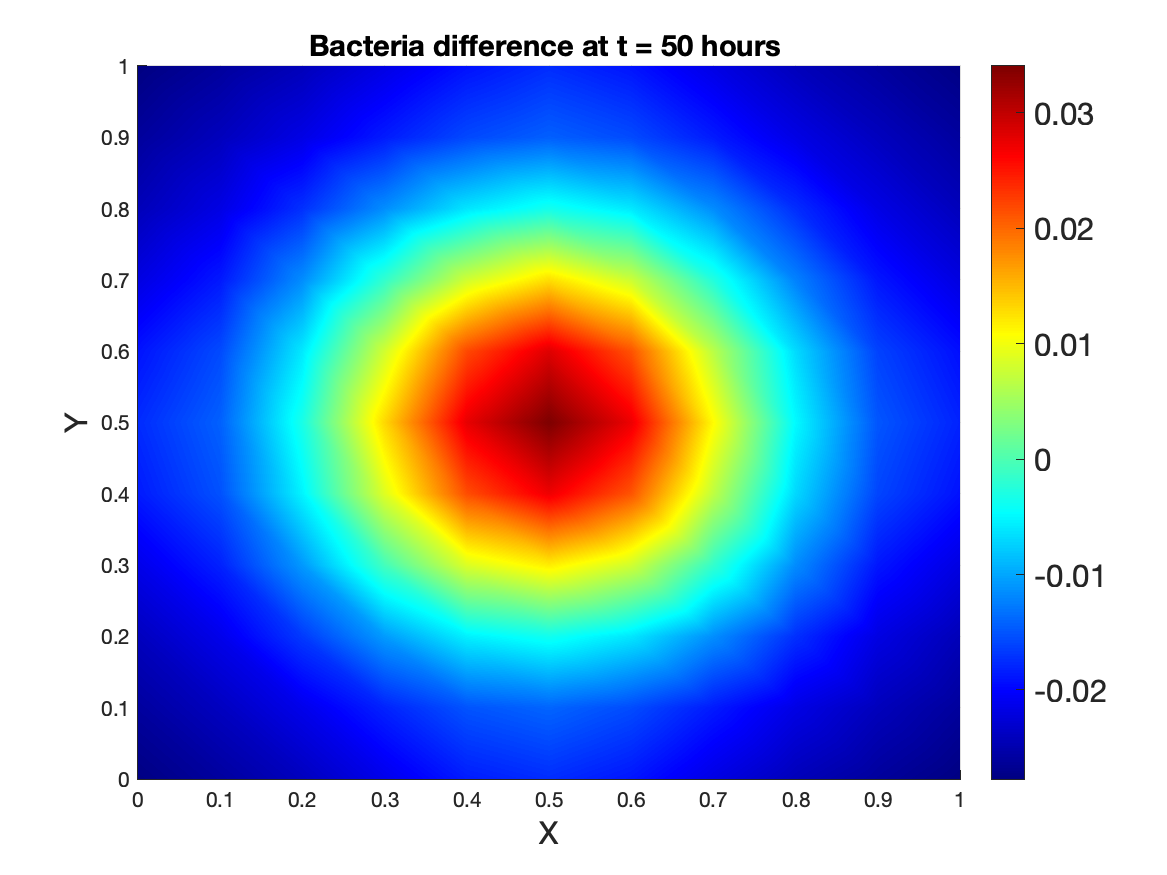}
		\caption*{\footnotesize{Bacteria \\at t=50}}
		\label{}
	\end{subfigure} 
	\begin{subfigure}{0.24 \textwidth}
		\centering
		\includegraphics[width=\textwidth]{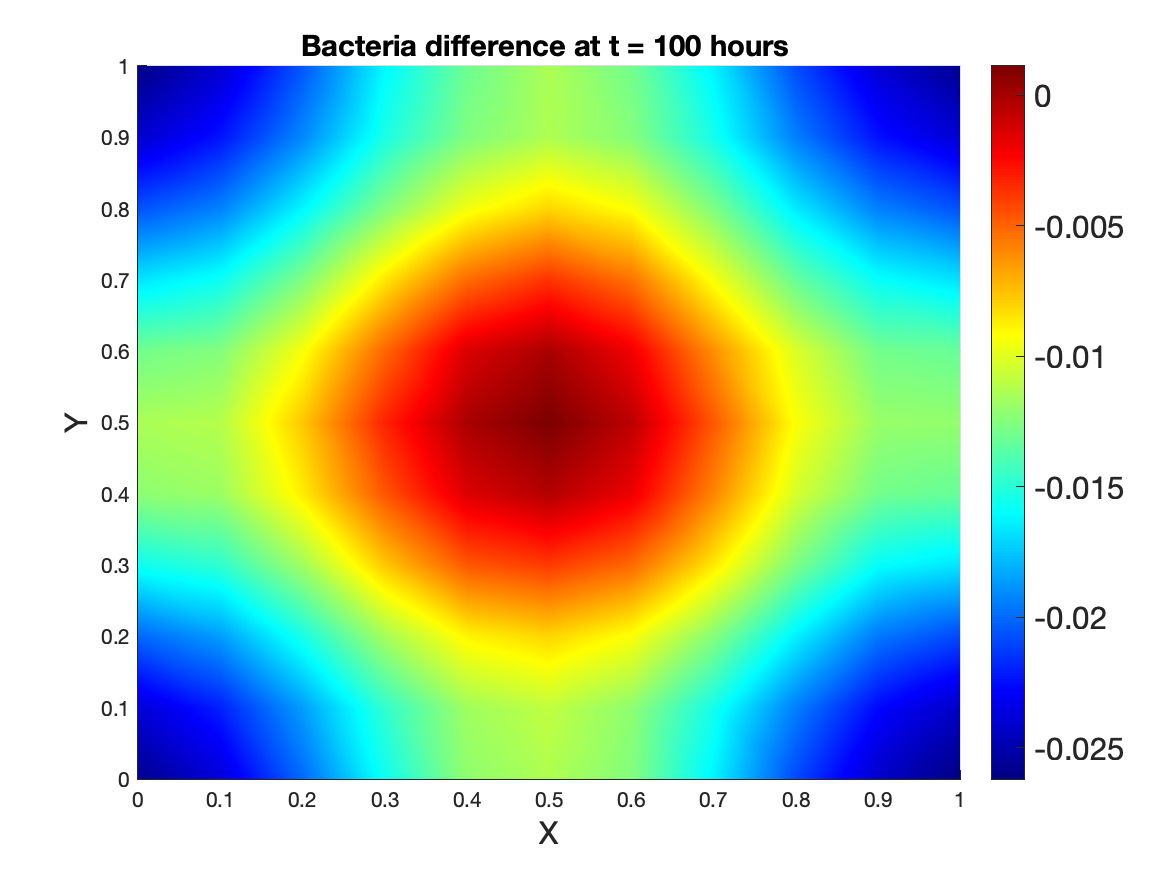}
		\caption*{\footnotesize{Bacteria \\at t=100}}
	\end{subfigure}
	\begin{subfigure}{0.24 \textwidth}
		\centering
		\includegraphics[width=\textwidth]{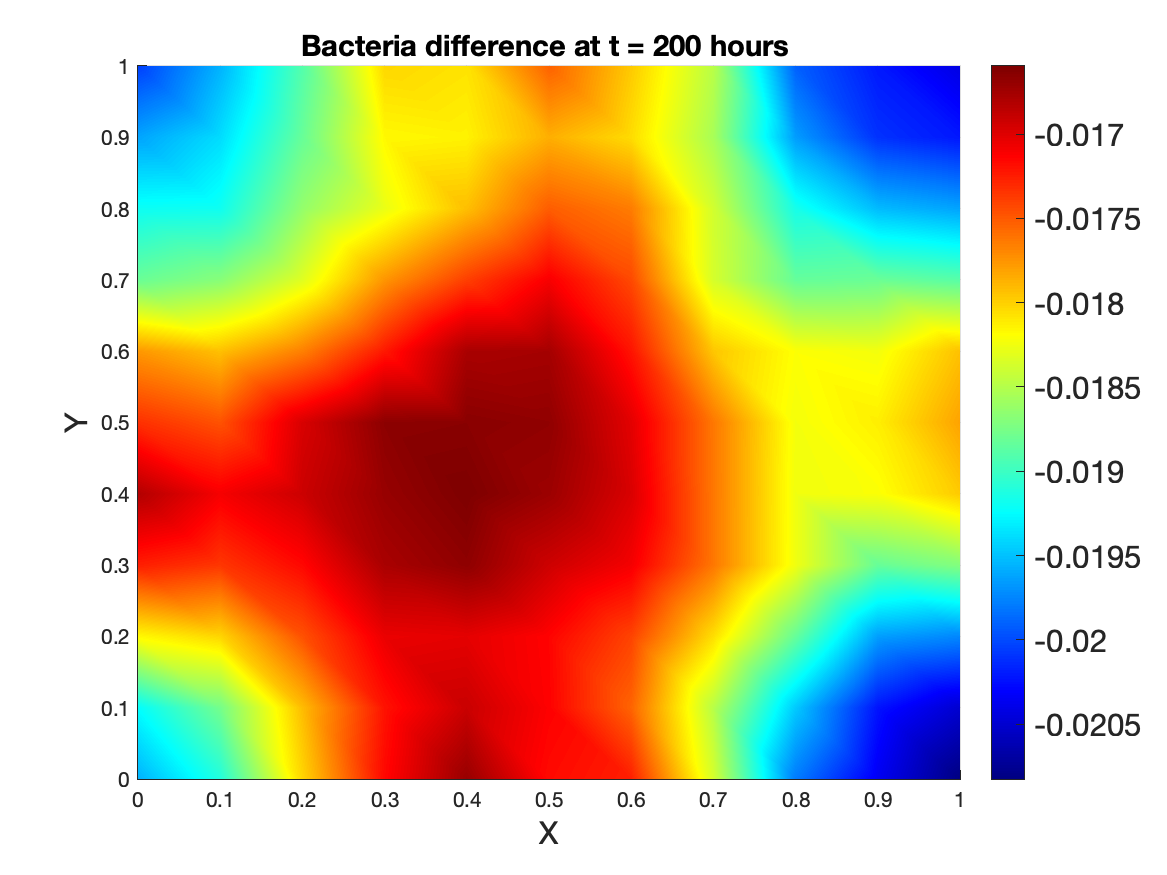}
		\caption*{\footnotesize{Bacteria \\at t=200}}
	\end{subfigure} 
	\begin{subfigure}{0.24 \textwidth}
		\centering
		\includegraphics[width=\textwidth]{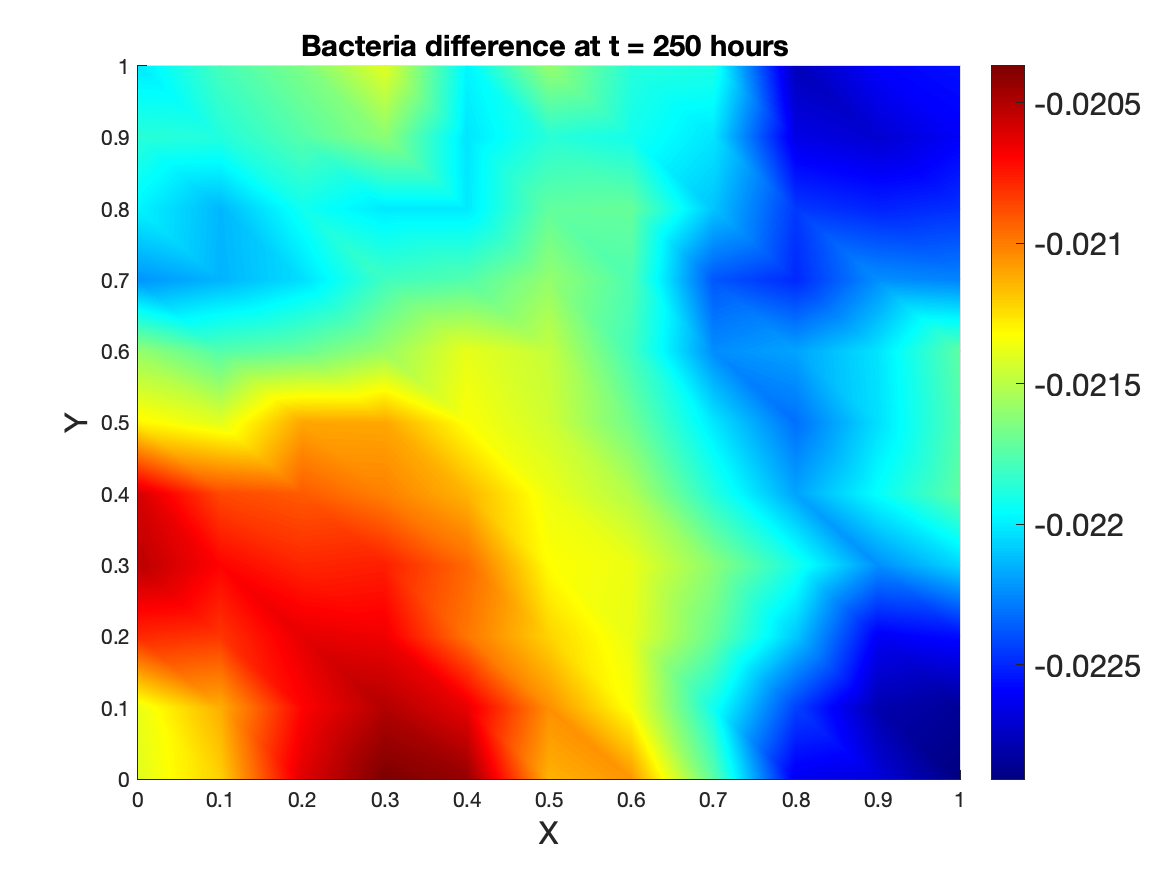}
		\caption*{\footnotesize{Bacteria \\at t=250}}
	\end{subfigure}     \vspace{0.5cm}     \\
	
	\begin{subfigure}{0.24 \textwidth}
		\centering
		\includegraphics[width=\textwidth]{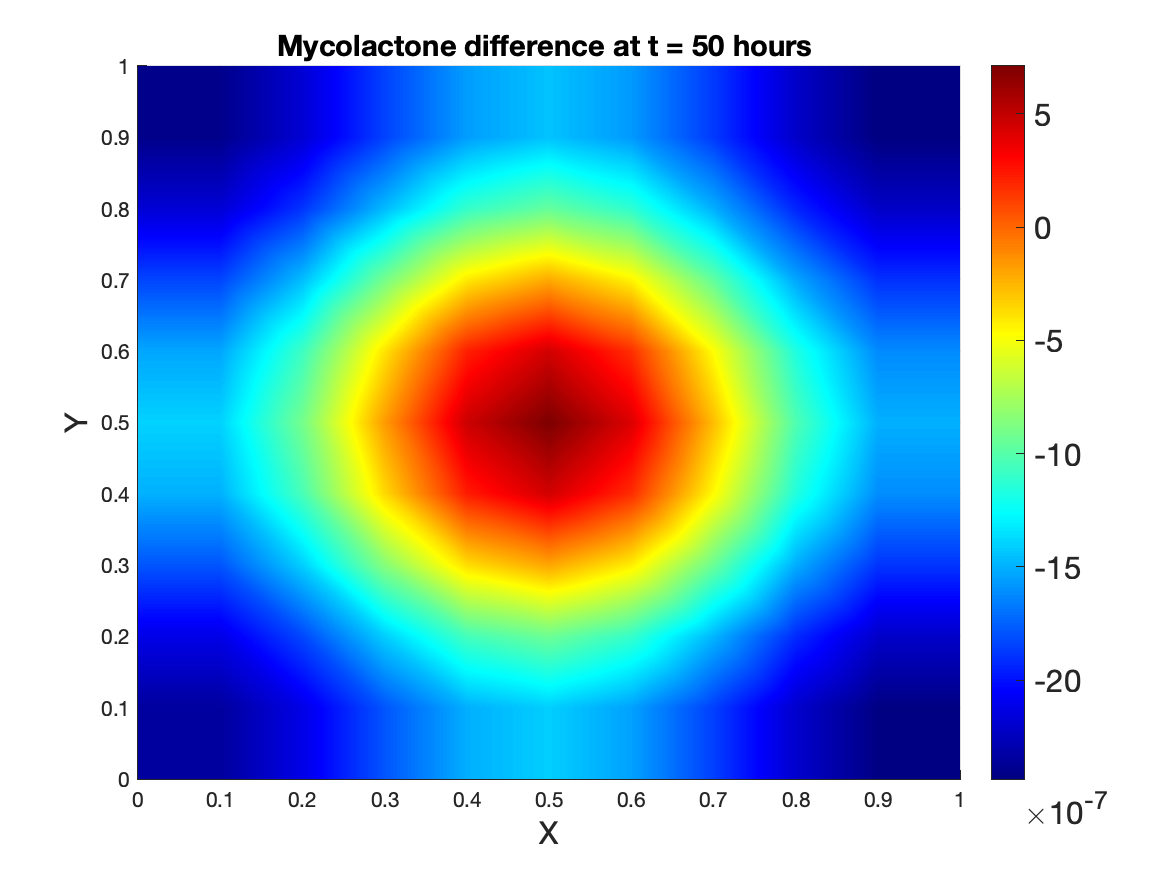}
		\caption*{\footnotesize{Mycolactone \\ at t=50}}
		\label{}
	\end{subfigure}
	\begin{subfigure}{0.24 \textwidth}
		\centering
		\includegraphics[width=\textwidth]{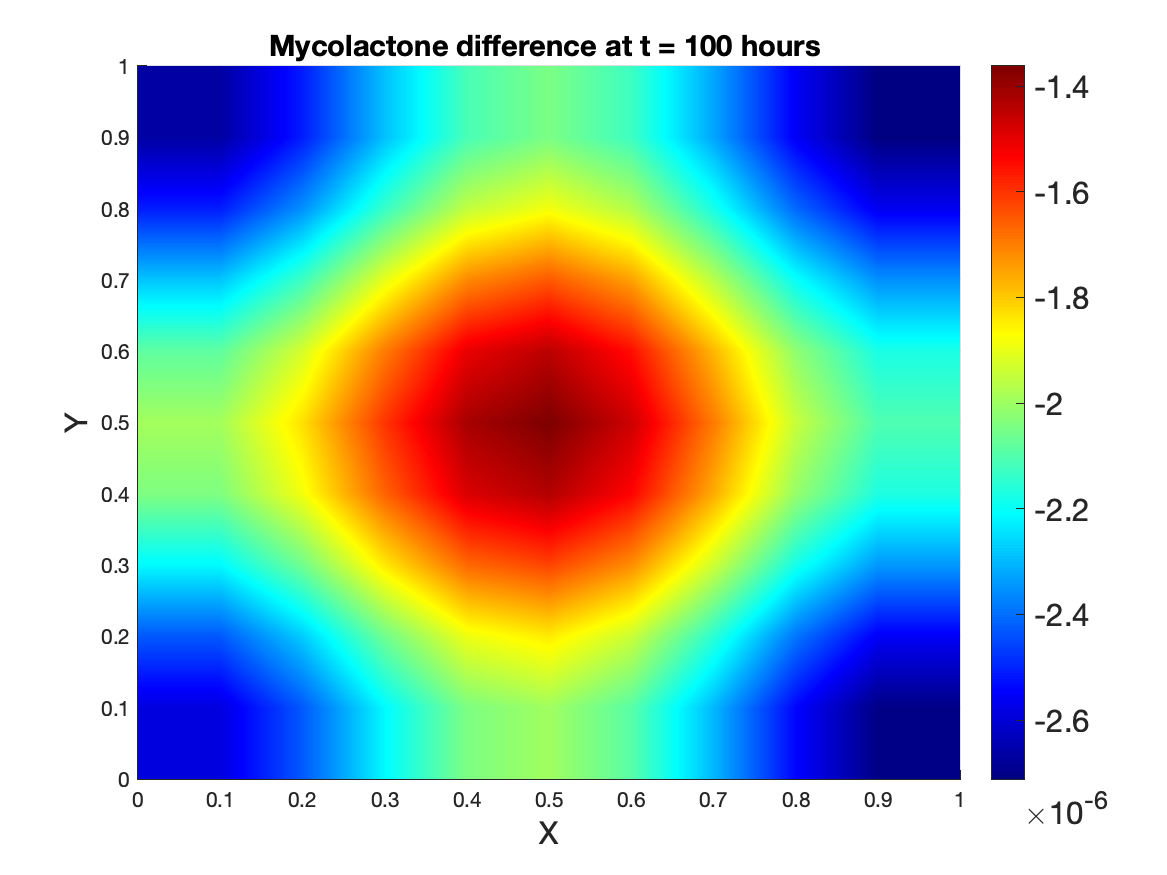}
		\caption*{\footnotesize{Mycolactone \\ at t=100}}
	\end{subfigure} 
	\begin{subfigure}{0.24 \textwidth}
		\centering
		\includegraphics[width=\textwidth]{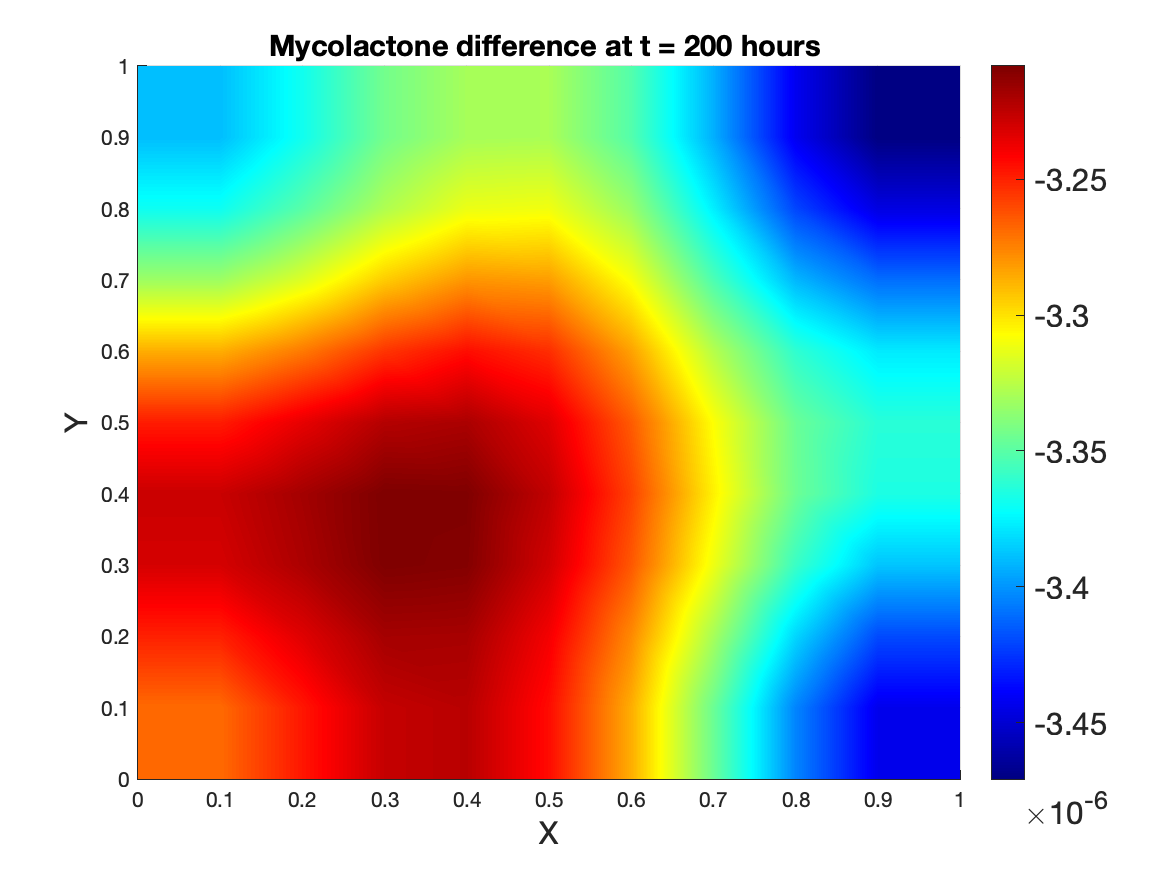}
		\caption*{\footnotesize{Mycolactone \\ at t=200}}
	\end{subfigure} 
	\begin{subfigure}{0.24 \textwidth}
		\centering
		\includegraphics[width=\textwidth]{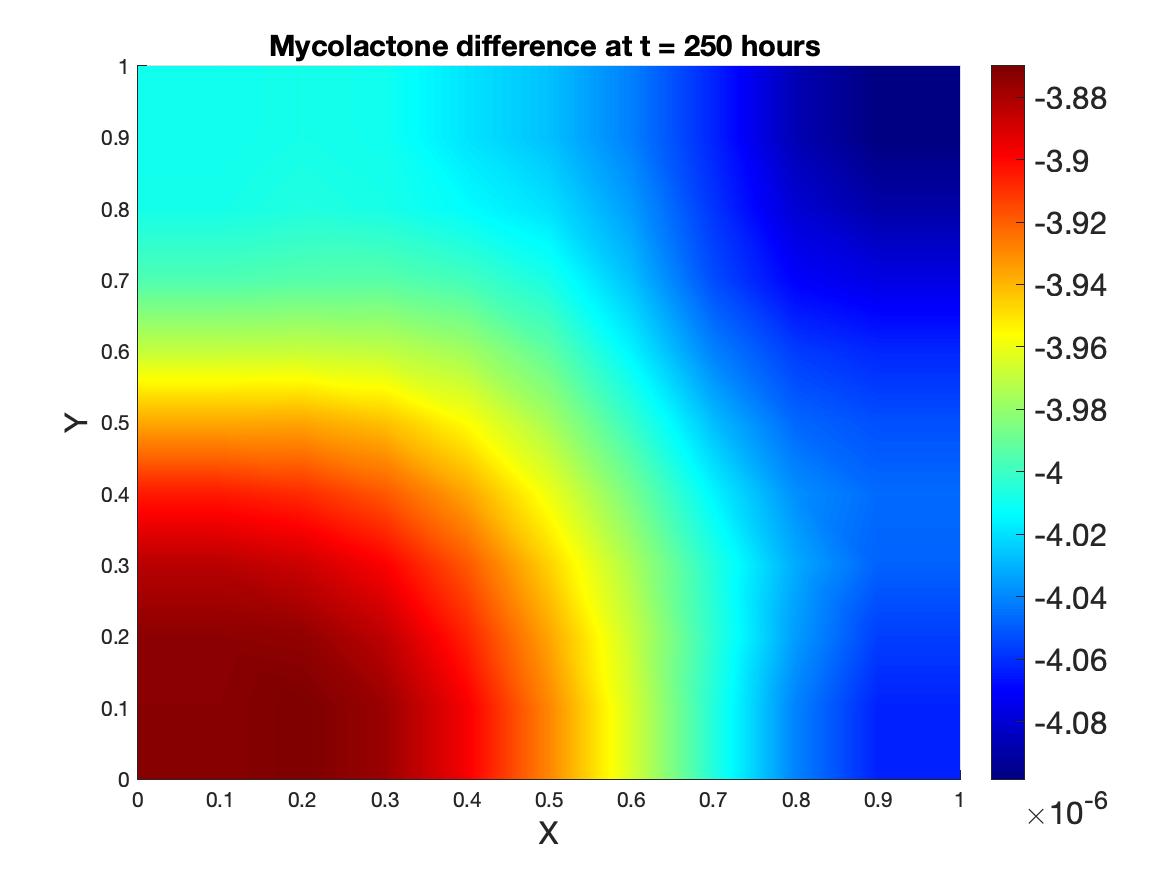}
		\caption*{\footnotesize{Mycolactone \\ at t=250}}
	\end{subfigure}     \vspace{0.5cm}     \\
	
	\begin{subfigure}{0.24 \textwidth}
		\centering
		\includegraphics[width=\textwidth]{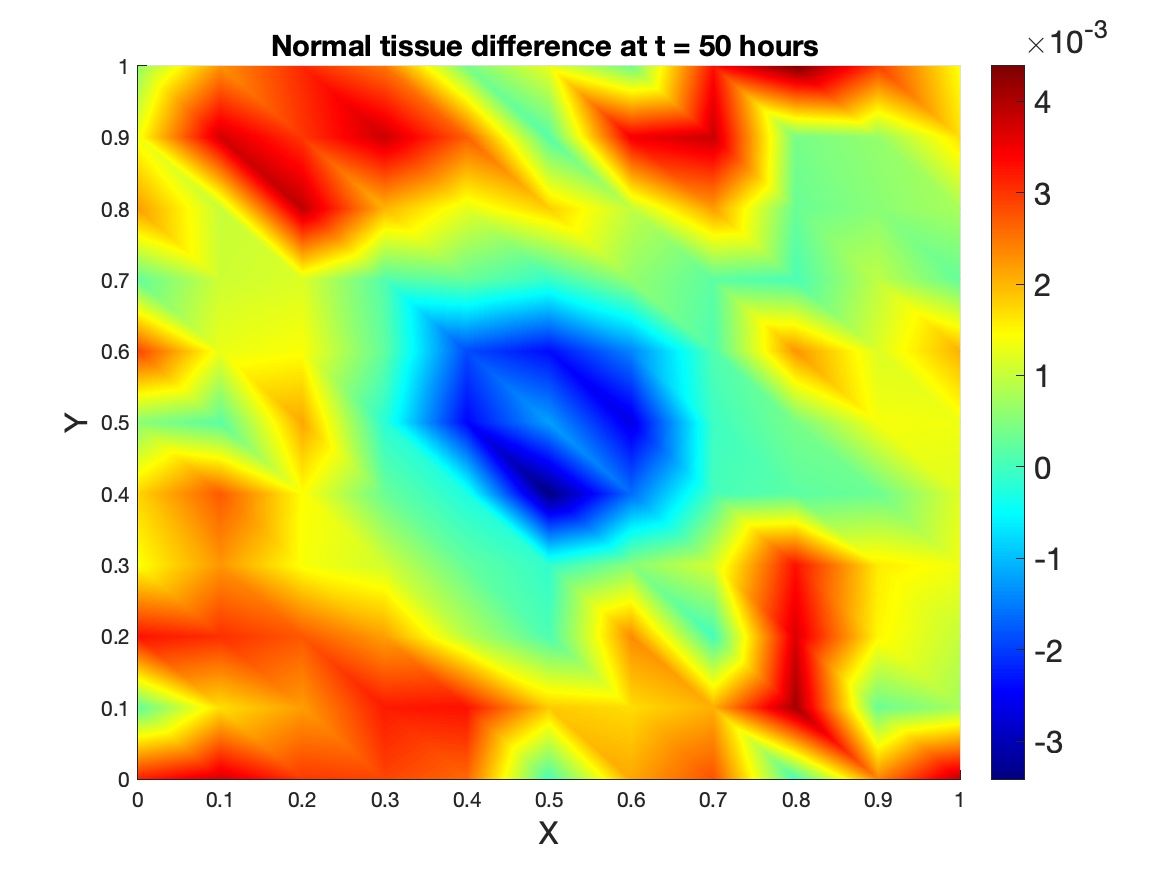}
		\caption*{\footnotesize{Normal tissue \\ at t=50}}
		\label{}
	\end{subfigure} 
	\begin{subfigure}{0.24 \textwidth}
		\centering
		\includegraphics[width=\textwidth]{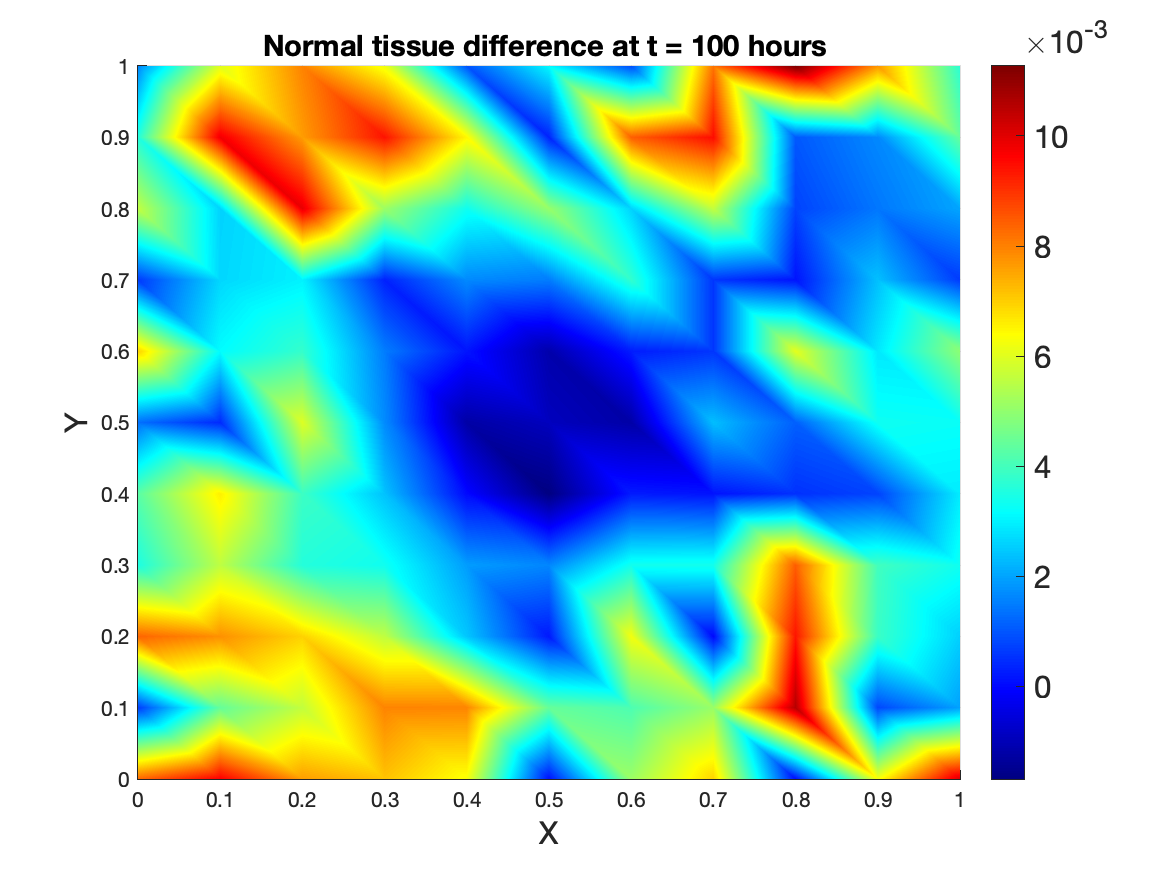}
		\caption*{\footnotesize{Normal tissue \\ at t=100}}
	\end{subfigure} 
	\begin{subfigure}{0.24 \textwidth}
		\centering
		\includegraphics[width=\textwidth]{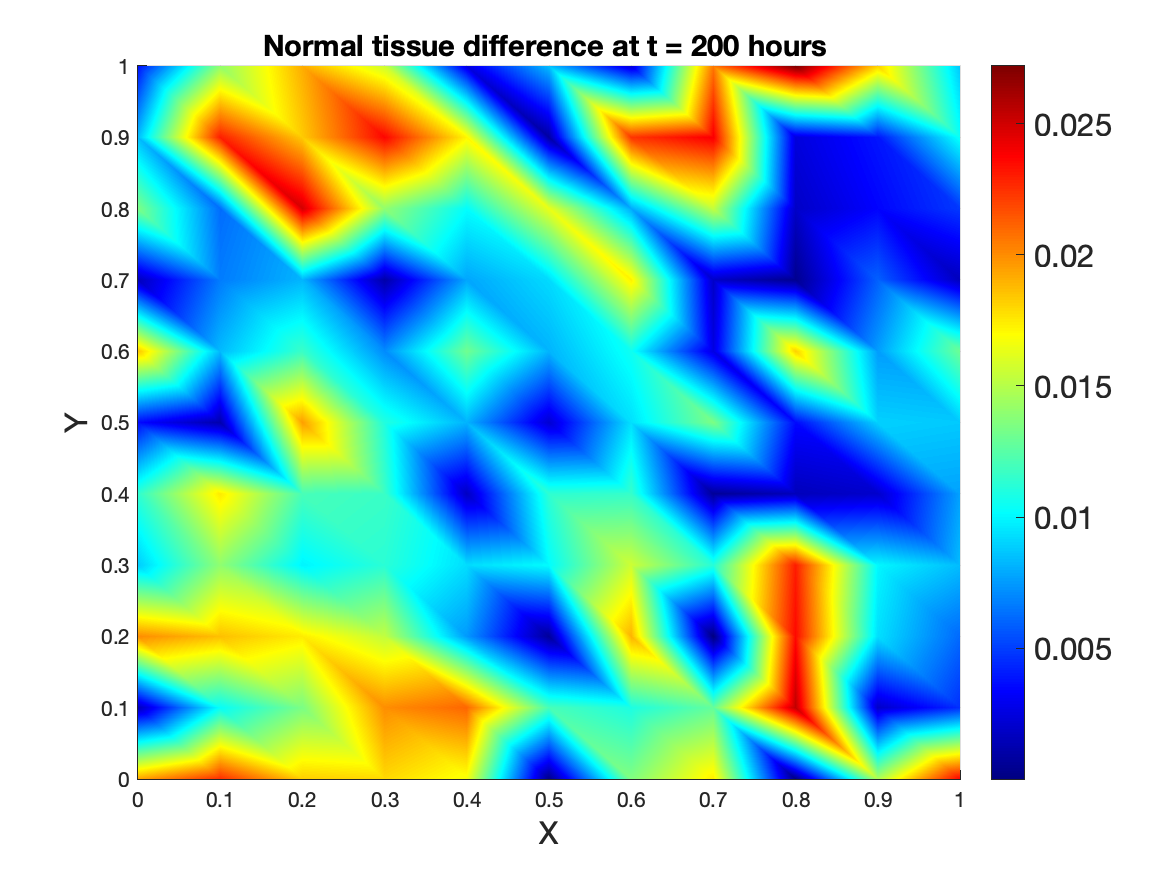}
		\caption*{\footnotesize{Normal tissue \\ at t=200}}
	\end{subfigure} 
	\begin{subfigure}{0.24 \textwidth}
		\centering
		\includegraphics[width=\textwidth]{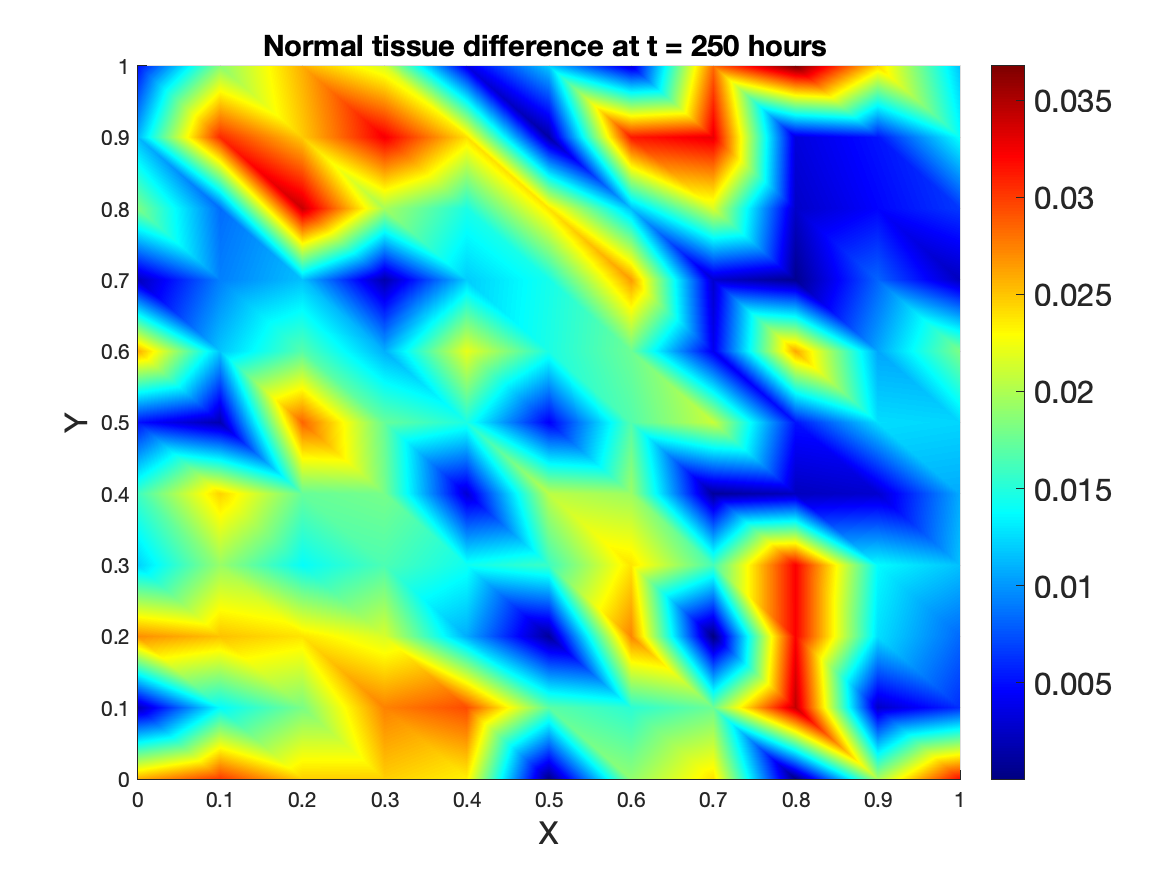}
		\caption*{\footnotesize{Normal tissue \\ at t=250}}
	\end{subfigure}
	\vspace{0.5cm} 
	\\
	
	\begin{subfigure}{0.24 \textwidth}
		\centering
		\includegraphics[width=\textwidth]{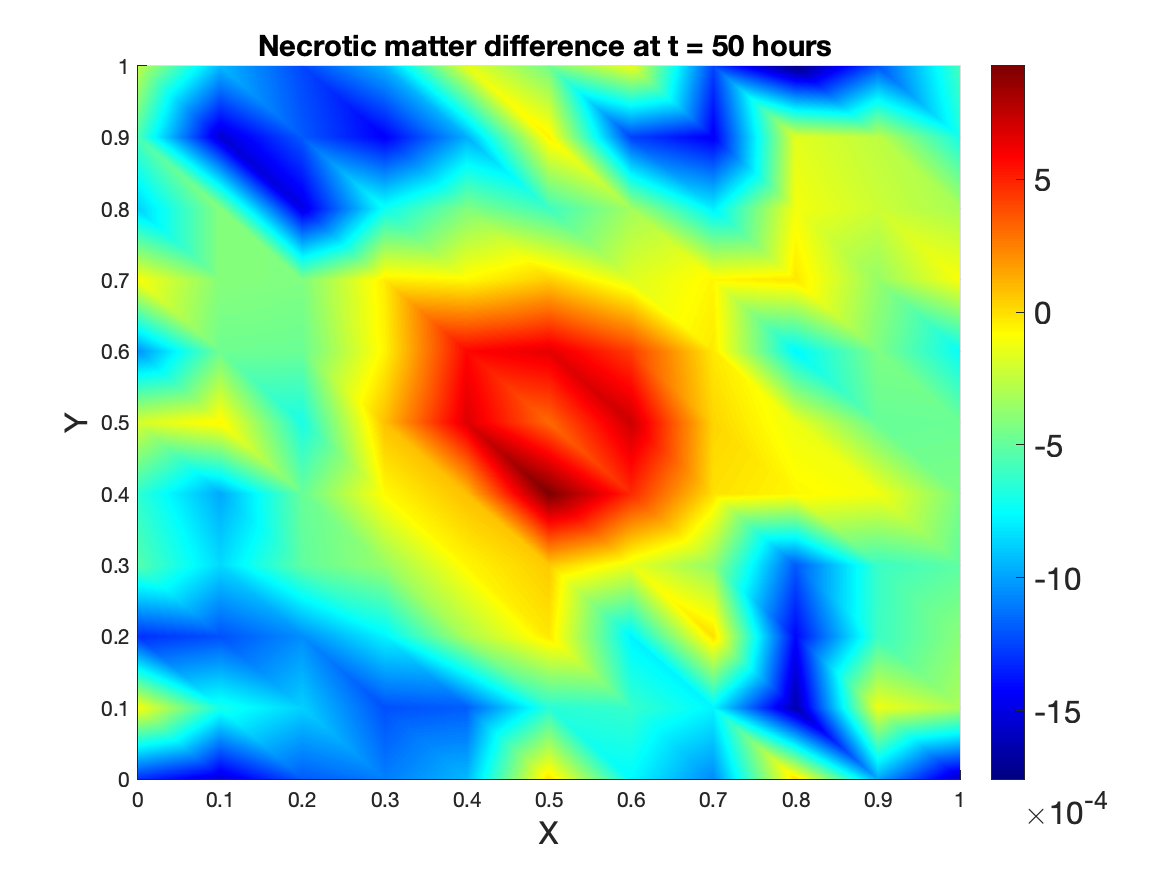}
		\caption*{\footnotesize{Necrotic matter \\ at t=50}}
		\label{}
	\end{subfigure} 
	\begin{subfigure}{0.24 \textwidth}
		\centering
		\includegraphics[width=\textwidth]{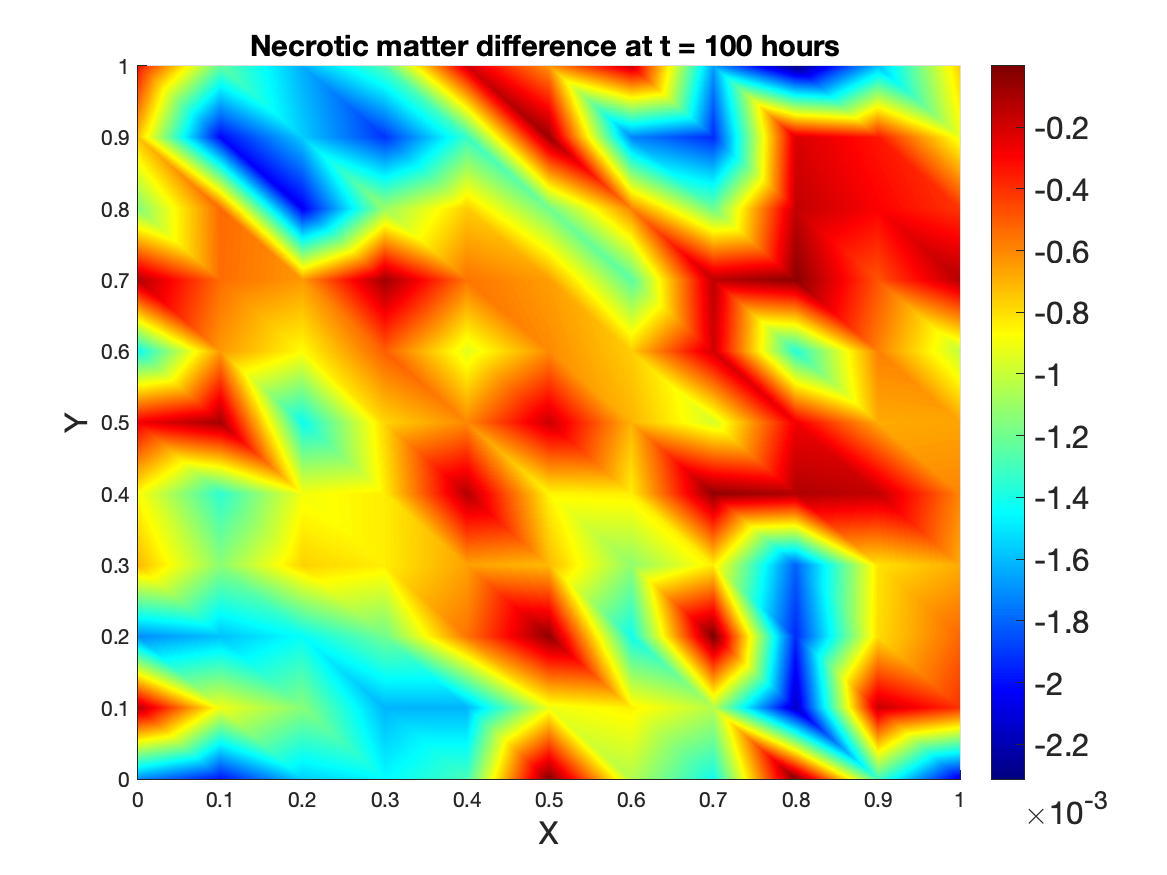}
		\caption*{\footnotesize{Necrotic matter \\ at t=100}}
		\label{}
	\end{subfigure} 
	\begin{subfigure}{0.24 \textwidth}
		\centering
		\includegraphics[width=\textwidth]{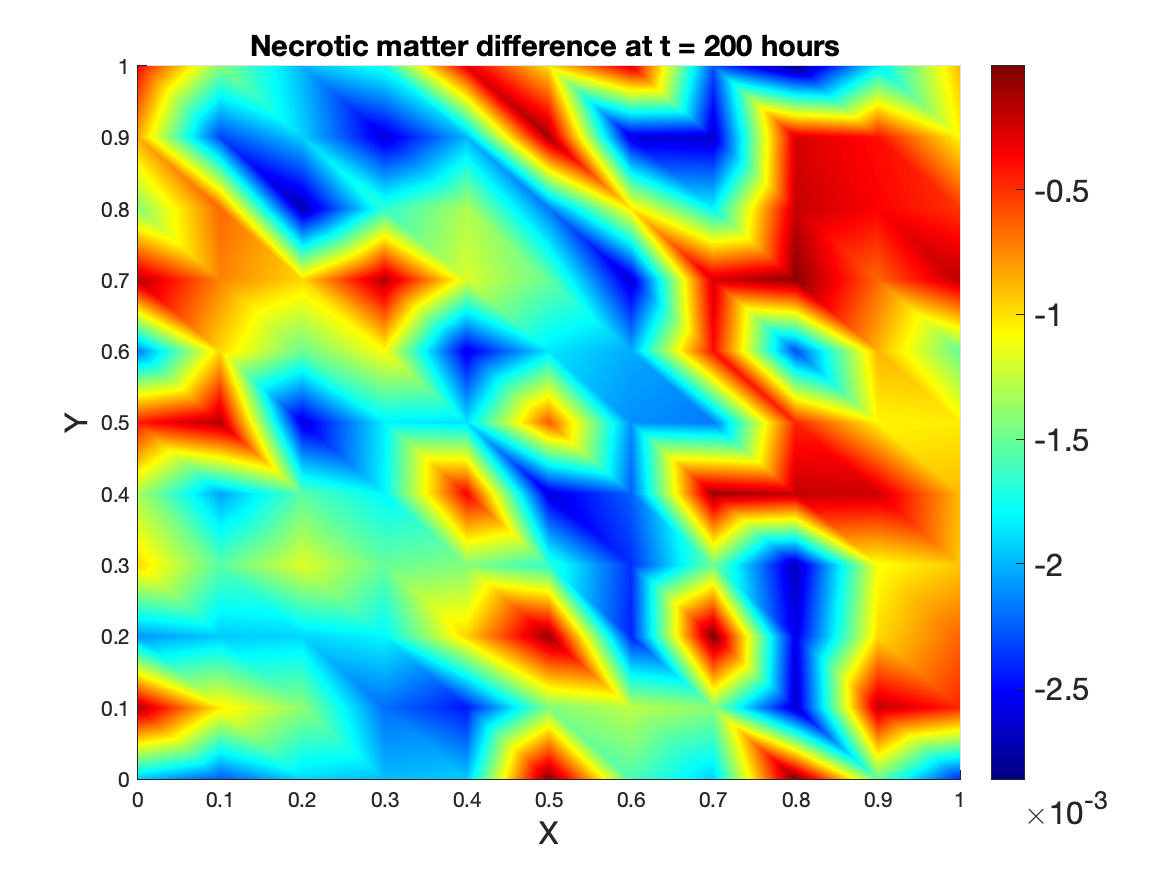}
		\caption*{\footnotesize{Necrotic matter \\ at t=200}}
		\label{}
	\end{subfigure} 
	\begin{subfigure}{0.24\textwidth}
		\centering
		\includegraphics[width=\textwidth]{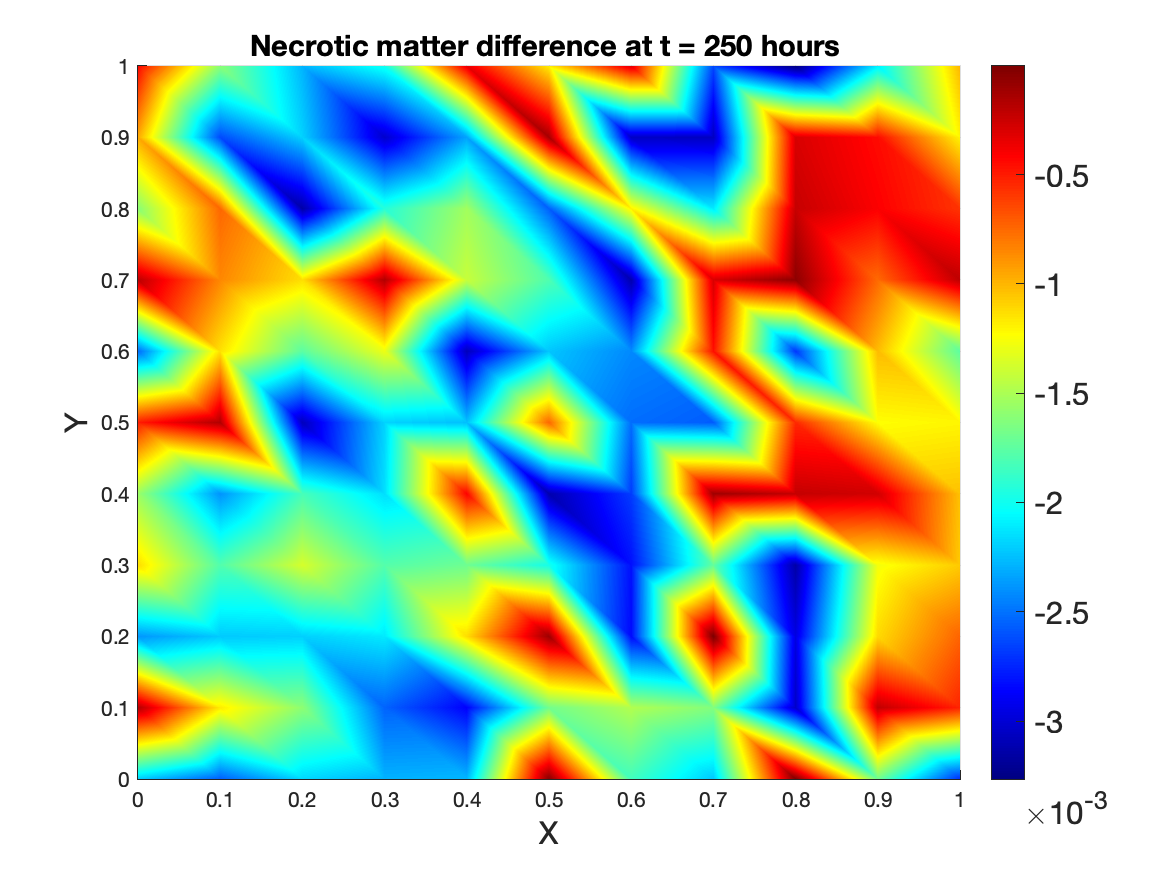}
		\caption*{\footnotesize{Necrotic matter \\ at t=250}}
		\label{}
	\end{subfigure} 
	\caption{Difference in density of bacteria, mycolactone concentration, normal tissue and necrotic matter at different times for systems \eqref{nondimsysburulinonlinear}  and \eqref{14} (\eqref{nondimsysburulinonlinear}-\eqref{14}: Scenario 1)}
	\label{fig:comparisonnonlinearburuli}
\end{figure}

\noindent
Fig \ref{fig:comparisonnonlinearburuli} compares the densities of bacteria, mycolactone concentration, normal tissue, and necrotic matter at various time points. Initially, differences between models show a peak of around 0.03 at early time points (t=50 to t=100), which reduces to 0.0225 by t=250, suggesting that bacterial growth and spread decelerates over time in Chapter 4. This phenomenon may be attributed to the diffusion coefficient in \eqref{nondimsysburulinonlinear} being dependent on both bacteria and normal tissue concentrations, potentially retarding the initial spread of bacteria at such high densities. Similarly, the taxis coefficients could alter bacterial movement patterns and reduce their accumulation in favourable areas. These more sophisticated representations of bacteria movement and spread in this model likely contribute to a more gradual disease progression by introducing additional limiting factors on bacteria proliferation and distribution. Mycolactone concentration differences mirror bacteria growth patterns, as mycolactone production is biologically linked to bacteria activity. \\[-2ex]

\noindent A key observation is the positive difference in normal tissue, which suggests that the model \eqref{nondimsysburulinonlinear} predicts either slower tissue damage or more preservation of healthy tissue compared to the results of \eqref{14}. By t=250, differences reach 0.035, indicating that tissue damage is either mitigated or delayed in certain regions. This might be due to the nonlinear diffusion and taxis coefficients. \\[-2ex]

\noindent In conclusion, the motility coefficients in model \eqref{nondimsysburulinonlinear} result in reduced bacteria spread, better preservation of normal tissue, and less necrosis. This outcome may be attributed to decreased tactic sensitivity towards necrotic matter and normal tissue, and reduced diffusion.\\[-2ex]

\noindent Our simulations suggest that system \eqref{14} should be preferred to \eqref{nondimsysburulinonlinear}, as it needs fewer ad-hoc assumptions regarding bacteria responses to signals in their environment. Furthermore, it has the potential to incorporate more detailed information about the spatial distribution of tissue, thus allowing for more realistic descriptions than the uniform distribution considered in this work. At the same time, the predictions of system \eqref{nondimsysburulinonlinear} might underestimate the wound spread, leading to insufficient resection of the affected area.\\[-2ex]

\noindent
Overall, our findings underscore the complex interplay between bacteria dynamics, tissue health, and the body's regenerative capacity in the pathogenesis of Buruli ulcer. Further research into these mechanisms could inform the development of novel therapeutic approaches for managing this debilitating disease.
\\

\section*{Acknowledgement} SCM acknowledges funding by the DAAD, the Research Initiative MathApp, and the TU-Nachwuchsring at the RPTU in Kaiserslautern. NNP thanks the DAAD for granting her a short-term visit at the TU Kaiserslautern during June-July 2019.

	
	\phantomsection
	\printbibliography

@book{bellomo2017quest,
title={A quest towards a mathematical theory of living systems},
author={Bellomo, Nicola and Bellouquid, Abdelghani and Gibelli, Livio and Outada, Nisrine},
year={2017},
doi={10.1007/978-3-319-57436-3},
publisher={Springer},
}

@article{conte2023mathematical,
title={Mathematical modeling of glioma invasion and therapy approaches via kinetic theory of active particles},
author={Conte, Martina and Dzierma, Yvonne and Knobe, Sven and Surulescu, Christina},
journal={Mathematical Models and Methods in Applied Sciences},
volume = {33},
number = {05},
pages = {1009-1051},
year = {2023},
publisher={World Scientific},
doi={10.1142/S0218202523500227}
}

@article{conte2021mathematical,
title={Mathematical modeling of glioma invasion: acid-and vasculature mediated go-or-grow dichotomy and the influence of tissue anisotropy},
author={Conte, Martina and Surulescu, Christina},
journal={Applied Mathematics and Computation},
volume={407},
pages={126305},
year={2021},
publisher={Elsevier},
doi={10.1016/j.amc.2021.126305}
}

@article{corbin2021modeling,
title={Modeling glioma invasion with anisotropy-and hypoxia-triggered motility enhancement: From subcellular dynamics to macroscopic PDEs with multiple taxis},
author={Corbin, Gregor and Klar, Axel and Surulescu, Christina and Engwer, Christian and Wenske, Michael and Nieto, Juanjo and Soler, Juan},
journal={Mathematical Models and Methods in Applied Sciences},
volume={31},
number={01},
pages={177--222},
year={2021},
publisher={World Scientific},
doi={10.1142/S0218202521500056}
}

@article{engwer2015glioma,
title={Glioma follow white matter tracts: a multiscale DTI-based model},
author={Engwer, Christian and Hillen, Thomas and Knappitsch, Markus and Surulescu, Christina},
journal={Journal of Mathematical Biology},
volume={71},
pages={551--582},
year={2015},
publisher={Springer},
doi={10.1007/s00285-014-0822-7}
}

@article{engwer2016effective,
title={Effective equations for anisotropic glioma spread with proliferation: a multiscale approach and comparisons with previous settings},
author={Engwer, Christian and Hunt, Alexander and Surulescu, Christina},
journal={Mathematical medicine and biology: a journal of the IMA},
volume={33},
number={4},
pages={435--459},
year={2016},
publisher={Oxford University Press},
doi={10.1093/imammb/dqv030}
}

@article{kumar2021multiscale,
title={Multiscale modeling of glioma pseudopalisades: contributions from the tumor microenvironment},
author={Kumar, Pawan and Li, Jing and Surulescu, Christina},
journal={Journal of Mathematical Biology},
volume={82},
pages={1--45},
year={2021},
publisher={Springer},
doi={10.1007/s00285-021-01599-x}
}

@article{painter2013mathematical,
title = {Mathematical modelling of glioma growth: The use of Diffusion Tensor Imaging (DTI) data to predict the anisotropic pathways of cancer invasion},
journal = {Journal of Theoretical Biology},
volume = {323},
pages = {25-39},
year = {2013},
issn = {0022-5193},
doi = {10.1016/j.jtbi.2013.01.014},
author = {K.J. Painter and T. Hillen},
publisher={Elsevier}
}

@article{plaza,
doi = {10.1007/s00285-018-1323-x},
year = {2019},
publisher = {Springer Science and Business Media {LLC}},
volume = {78},
number = {6},
pages = {1681--1711},
author = {Ram{\'{o}}n G. Plaza},
title = {Derivation of a bacterial nutrient-taxis system with doubly degenerate cross-diffusion as the parabolic limit of a velocity-jump process},
journal = {Journal of Mathematical Biology}
}

@article{ABT02,
title={Buruli ulcer in Ghana: results of a national case search},
  author={Amofah, George and Bonsu, Frank and Tetteh, Christopher and Okrah, Jane and Asamoa, Kwame and Asiedu, Kingsley and Addy, Jonathan},
  journal={Emerging infectious diseases},
  volume={8},
  number={2},
  pages={167},
  year={2002},
  publisher={Centers for Disease Control and Prevention},
doi = {10.3201/eid0802.010119}
}

@article{Asiedu,
author = {Asiedu, K. and Etuaful, S.},
  title = {Socioeconomic implications of \textit{Buruli ulcer} in Ghana: a three-year review},
  journal = {Transactions of the Royal Society of Tropical Medicine and Hygiene},
  volume = {92},
  number = {6},
  pages = {1015-1022},
  year = {1998},
  doi = {10.4269/ajtmh.1998.59.1015}
}

@article{JSSPMPHHA,
title={Buruli ulcer (M. ulcerans infection): new insights, new hope for disease control},
  author={Johnson, Paul D R and Stinear, Timothy and Small, Pamela L C and Pluschke, Gerd and Merritt, Richard W and Portaels, Francoise and Huygen, Kris and Hayman, John A and Asiedu, Kingsley},
  journal={PLoS medicine},
  volume={2},
  number={4},
  pages={e108},
  year={2005},
doi = {10.1371/journal.pmed.0020108}
}

@article{Muelder,
  author = {Muelder, K. and Nourou, A.},
  title = {Buruli ulcer in Benin},
  journal = {Lancet},
  volume = {336},
  pages = {1109-1111},
  year = {1990},
  doi={10.1016/0140-6736(90)92581-2}
}

@article{Oliviera,
  author={Oliveira, M. S and et al.},
  title = {Infection with Mycobacterium ulcerans induces persistent inflammatory responses in mice},
  journal = {Infection and Immunity},
  volume = {73},
  number = {10},
  pages = {6299-6310},
  year = {2005},
 doi={10.1128/IAI.73.10.6299-6310.2005}
}

@article{Walsh,
    author = {Walsh, Douglas S. and Portaels, Françoise and Meyers, Wayne M.},
    title = "{Buruli ulcer (Mycobacterium ulcerans infection)}",
    journal = {Transactions of The Royal Society of Tropical Medicine and Hygiene},
    volume = {102},
    number = {10},
    pages = {969-978},
    year = {2008},
    issn = {0035-9203},
    doi = {10.1016/j.trstmh.2008.06.006}
}

@article{Wansbrough,
title = {Buruli ulcer: emerging from obscurity},
journal = {The Lancet},
volume = {367},
number = {9525},
pages = {1849-1858},
year = {2006},
issn = {0140-6736},
doi ={10.1016/S0140-6736(06)68807-7},
author = {Mark Wansbrough-Jones and Richard Phillips}
}

@article{winkler,
title = {Aggregation vs. global diffusive behavior in the higher-dimensional Keller-Segel model},
journal = {Journal of Differential Equations},
volume = {248},
number = {12},
pages = {2889-2905},
year = {2010},
issn = {0022-0396},
doi = {10.1016/j.jde.2010.02.008},
author = {Michael Winkler}
}

@article{ab,
author = {Agbenorku, Pius and Saunderson, Paul},
year = {2013},
month = {12},
pages = {4 pages},
title = {Mycobacterium ulcerans Disease of the Face: The Fate of the Victims},
volume = {3},
journal = {Mycobacterial Diseases},
doi = {10.4172/2161-1068.1000133}
}

@article{bar,
  author = {Barker, D.J.P.},
  title = {Epidemiology of \textit{Mycobacterium ulcerans} infection},
  journal = {Transactions of the Royal Society of Tropical Medicine and Hygiene},
  volume = {67},
  pages = {43-7},
  year = {1973},
  doi = {10.1016/0035-9203(73)90317-9}
}

@book{LAD,
  author = {Ladyzenskaja, O. A. and Solonnikov, V. A. and Ural'ceva, N. N.},
  title = {Linear and Quasi-linear Equations of Parabolic Type},
  publisher = {American Mathematical Society Translations},
  volume = {23},
  address = {Providence, RI},
  year = {1968}
}

@article{Portaels,
title={First cultivation and characterization of Mycobacterium ulcerans from the environment},
  author={Portaels, Fran{\c{c}}oise and Meyers, Wayne M and Ablordey, Anthony and Castro, Ant{\'o}nio G and Chemlal, Karim and de Rijk, Pim and Elsen, Pierre and Fissette, Krista and Fraga, Alexandra G and Lee, Richard and others},
  journal={PLoS neglected tropical diseases},
  volume={2},
  number={3},
  pages={e178},
  year={2008},
  publisher={Public Library of Science San Francisco, USA},
doi = {10.1371/journal.pntd.0000178}
}

@article{TAOchemohapto,
author = {Tao, Youshan and Winkler, Michael},
title = {A Chemotaxis-Haptotaxis Model: The Roles of Nonlinear Diffusion and Logistic Source},
journal = {SIAM Journal on Mathematical Analysis},
volume = {43},
number = {2},
pages = {685-704},
year = {2011},
doi = {10.1137/100802943}
}

@mastersthesis{lenz,
          school = {Technische Universit{\"a}t Darmstadt},
            year = {2020},
         address = {Darmstadt},
           title = {Global Existence for a Tumor Invasion Model with Repellent Taxis and Therapy},
          author = {Jonas Lenz},
    doi = {10.25534/TUPRINTS-00011578}
}

@article{sym12111870,
AUTHOR = {Kumar, Pawan and Surulescu, Christina},
TITLE = {A Flux-Limited Model for Glioma Patterning with Hypoxia-Induced Angiogenesis},
JOURNAL = {Symmetry},
VOLUME = {12},
YEAR = {2020},
NUMBER = {11},
ARTICLE-NUMBER = {1870},
ISSN = {2073-8994},
DOI = {10.3390/sym12111870}
}

@article{randomjumpKPH,
author = {Painter, Kevin and Hillen, Thomas},
year = {2002},
pages = {501-543},
title = {Volume-filling and quorum-sensing in models for chemosensitive movement},
volume = {10},
journal = {Can. Appl. Math. Q.}
}

@article{othmer2002diffusion,
  title={The diffusion limit of transport equations II: Chemotaxis equations},
  author={Othmer, Hans G and Hillen, Thomas},
  journal={SIAM Journal on Applied Mathematics},
  volume={62},
  number={4},
  pages={1222--1250},
  year={2002},
  doi = {10.1137/S0036139900382772},
  publisher={SIAM}
}

@article{randomjumpSAO,
author = {Stevens, Angela and Othmer, Hans G.},
title = {Aggregation, Blowup, and Collapse: The ABC's of Taxis in Reinforced Random Walks},
journal = {SIAM Journal on Applied Mathematics},
volume = {57},
number = {4},
pages = {1044-1081},
year = {1997},
doi = {10.1137/S0036139995288976}
}

@article{doi:10.1098/rsif.2008.0014,
author = {Codling, Edward A  and Plank, Michael J  and Benhamou, Simon },
title = {Random walk models in biology},
journal = {Journal of The Royal Society Interface},
volume = {5},
number = {25},
pages = {813-834},
year = {2008},
doi = {10.1098/rsif.2008.0014}
}

@article{Nyabadza2015,
    author = {Nyabadza, Farai and Bonyah, Ebenezer} ,
    title = {On the transmission dynamics of Buruli ulcer in Ghana: Insights through a mathematical model},
    journal = {BMC Research Notes},
    volume={8},
  year={2015},
doi={10.1186/s13104-015-1619-5}
}

@article{neglectedburuli,
author = {Jeannette Guarner},
title = {Buruli Ulcer: Review of a Neglected Skin Mycobacterial Disease},
journal = {Journal of Clinical Microbiology},
volume = {56},
number = {4},
year = {2018},
doi = {10.1128/jcm.01507-17},
}

@article{Roberts2019,
  title={Mathematical model predicts anti-adhesion-antibiotic-debridement combination therapies can clear an antibiotic resistant infection},
  author={Roberts, Paul A and Huebinger, Ryan M and Keen, Emma and Krachler, Anne-Marie and Jabbari, Sara},
  journal={PLoS Comput Biol},
  volume={15},
  number={7},
  pages={e1007211},
  year={2019},
  publisher={Public Library of Science},
  doi={10.1371/journal.pcbi.1007211}
}

@article{Shi2016,
  title={Mathematical Model of Innate and Adaptive Immunity of Sepsis: A Modeling and Simulation Study of Infectious Disease},
  author={Shi, Zhenzhen and Wu, Chih-Hang J. and Ben-Arieh, David and Simpson, Steven Q.},
  journal={Computational and Mathematical Methods in Medicine},
  volume={2016},
  pages={1-10},
  year={2016},
  publisher={Hindawi},
  doi={10.1155/2015/504259}
}

@article{King2003,
  title={Modelling host tissue degradation by extracellular bacterial pathogens},
  author={King, J. R. and Koerber, A. J. and Croft, J. M. and Ward, J. P. and Williams, P. and Sockett, R. E.},
  journal={Mathematical Medicine and Biology: A Journal of the IMA},
  volume={20},
  number={3},
  pages={227--260},
  year={2003},
  publisher={Oxford University Press},
  doi={10.1093/imammb/20.3.227},
  pmid={14667046}
}

@misc{alghamdi2024combinedexperimentalmathematicalstudy,
      title={A Combined Experimental and Mathematical Study of The Evolution of Microbial Community Composed of Interacting Staphylococcus Strains}, 
      author={Nouf Alghamdi and Mal Horsburgh and Bakhtier Vasiev},
      year={2024},
      eprint={2402.04939},
      archivePrefix={arXiv},
      primaryClass={q-bio.PE},
      doi={10.48550/arXiv.2402.04939}, 
}

@article{McGann2009,
  title={Buruli ulcer in United Kingdom tourist returning from Latin America},
  author={McGann, Hugh and Stragier, Pieter and Portaels, Francoise and Gascoyne Binzi, Deborah and Collyns, Timothy and Lucas, Sebastian and Mawer, Damian},
  journal={Emerging Infectious Diseases},
  volume={15},
  number={11},
  pages={1827--1829},
  year={2009},
  publisher={Centers for Disease Control and Prevention (CDC)},
  doi={10.3201/eid1511.090460},
  pmid={19891876},
  pmcid={PMC2857232}
}

@article{Licata2016,
  title = {Diffusion of Bacterial Cells in Porous Media},
  author = {Licata, Nicholas A. and Mohari, Bitan and Fuqua, Clay and Setayeshgar, Sima},
  journal = {Biophysical Journal},
  volume = {110},
  number = {1},
  pages = {247--257},
  year = {2016},
  doi = {10.1016/j.bpj.2015.09.035},
  issn = {1542-0086},
 }

@misc{mohanan2024mathematicalmodeltissueregeneration,
      title={On a mathematical model for tissue regeneration}, 
      author={Chettiparambil Mohanan,Shimi and Mohan,Nishith  and Surulescu, Christina },
      year={2024},
      eprint={2403.04516},
      archivePrefix={arXiv},
      primaryClass={math.AP},
      doi={10.48550/arXiv.2403.04516}
}

@article{MARINO2019154,
title = {Moser iteration applied to elliptic equations with critical growth on the boundary},
journal = {Nonlinear Analysis},
volume = {180},
pages = {154-169},
year = {2019},
issn = {0362-546X},
doi = {10.1016/j.na.2018.10.002},
author = {Marino,Greta  and Winkert,Patrick },
}

@article{Giesselmann,
title = {Existence and uniqueness of global classical solutions to a two dimensional two species cancer invasion haptotaxis model},
journal = {Discrete and Continuous Dynamical Systems - B},
volume = {23},
number = {10},
pages = {4397-4431},
year = {2018},
issn = {1531-3492},
doi = {10.3934/dcdsb.2018169},
author = {Jan Giesselmann and Niklas Kolbe and Maria Lukacova-Medvidova  and Nikolaos Sfakianakis},
}

@article{zhao2021mathematical,
  title={A mathematical model for the coinfection of buruli ulcer and cholera},
  author={Zhao, Jin-Qiang and Bonyah, Ebenezer and Yan, Bing and Khan, Muhammad Altaf and Okosun, KO and Alshahrani, Mohammad Y and Muhammad, Taseer},
  journal={Results in Physics},
  volume={29},
  pages={104746},
  year={2021},
  publisher={Elsevier},
  doi = {10.1016/j.rinp.2021.104746}
}

@article{palomino1998effect,
  title={Effect of oxygen on growth of Mycobacterium ulcerans in the BACTEC system},
  author={Palomino, JC and Obiang, AM and Realini, L and Meyers, WM and Portaels, F},
  journal={Journal of clinical microbiology},
  volume={36},
  number={11},
  pages={3420--3422},
  year={1998},
  publisher={Am Soc Microbiol},
  doi={10.1128/jcm.36.11.3420-3422.1998}
}

@book{protter2012maximum,
  title={Maximum principles in differential equations},
  author={Protter, Murray H and Weinberger, Hans F},
  year={2012},
  publisher={Springer Science \& Business Media},
  doi={10.1007/978-1-4612-5282-5}
}

@article{nyarko,
title={Mathematical Modeling and Numerical Simulation for Mycobacterium Ulcerans Tissue Invasion: A Macroscopic Model for The Buruli Ulcer Disease},
author={Nyarko, P. R. and Adikorley, I. J. N. A. and Amoako-Yirenkyi, P. and Dontwi, I. K.},
journal={International Journal of Advances in Science, Engineering and Technology},
volume={5},
pages={41-48},
year={2017},
}

@phdthesis{diss-shimi,
	school = {RPTU Kaiserslautern-Landau},
	year = {2025},
	title = {Mathematical models for cell triggered tissue
	regeneration and wound formation},
	author = {Chettiparambil Mohanan, Shimi},
	doi = {10.26204/KLUEDO/8808}
}
	
\end{document}